\documentclass[reqno]{amsart}
\usepackage{setspace}
\usepackage{amssymb}
\usepackage{color}
\usepackage{comment}
\usepackage[margin=1in]{geometry}
\usepackage{graphicx}
\usepackage{float}
\usepackage{mathrsfs}
\usepackage{yfonts}
\usepackage{hyperref}
\usepackage{enumitem}
\usepackage{lscape}
\usepackage{subfigure}

\numberwithin{figure}{section}
\numberwithin{equation}{section}
\newtheorem{rhp}{Riemann-Hilbert Problem}[subsection]
\newtheorem{theorem}{Theorem}[section]

\theoremstyle{remark}
\newtheorem{remark}{Remark}[section]

\newcommand*{\dif}{\mathop{}\!\mathrm{d}}

\begin{document}

\title[Riemann problem of the defocusing NLS hydrodynamics]{\bf Rigorous asymptotic analysis for the Riemann problem of the defocusing nonlinear Schr\"{o}dinger hydrodynamics}

\author{Deng-Shan Wang}
\address{Deng-Shan Wang: Laboratory of Mathematics and Complex Systems (Ministry of Education), School of Mathematical Sciences, Beijing Normal University, Beijing 100875, China}
\email{dswang@bnu.edu.cn}
\author{Peng Yan}
\address{Peng Yan: Laboratory of Mathematics and Complex Systems (Ministry of Education), School of Mathematical Sciences, Beijing Normal University, Beijing 100875, China}
\email{yp2023@mail.bnu.edu.cn}

\date{\today}

\begin{abstract}

The rigorous asymptotic analysis for the Riemann problem of the defocusing nonlinear Schr\"{o}dinger hydrodynamics is a very interesting problem with many challenges. To date, the full analysis of this problem remains open. In this work, the long-time asymptotics for the defocusing nonlinear Schr\"{o}dinger equation with general step-like initial data is investigated by the Whitham modulation theory and Riemann-Hilbert formulation. The Whitham modulation theory shows that there are six cases for the initial discontinuity problem according to the orders of the Riemann invariants. The leading-order terms and the corresponding error estimates for each region of the six cases are formulated by the Deift-Zhou nonlinear steepest descent method for oscillatory Riemann-Hilbert problems. It is demonstrated that the long-time asymptotic solutions match very well with the results from Whitham modulation theory and the numerical simulations.

\end{abstract}

\maketitle

\setcounter{tocdepth}{1}
\tableofcontents

\section{Introduction}

The nonlinear Schr\"{o}dinger (NLS) equation can model a wide range of physical phenomena, such as unidirectional propagation of nonlinear waves in nonlinear optics \cite{Sulem1999}, superconductivity \cite{Sharma1976}, surface gravity waves \cite{Onorato2001}, Bose-Einstein condensates \cite{Yau2010}, and so on. The one-dimensional NLS equation with cubic nonlinearity is one of the integrable systems possessing infinite conservation laws and can be exactly solvable by the inverse scattering transform \cite{GGKM,Ablowitz1991}.
\par
The NLS equation can be divided into focusing type and defocusing type according to the sign of the nonlinearity term. The focusing NLS equation \cite{Tovbis2004}-\cite{Lenells2021} exhibits abundant nonlinear structures such as rogue waves \cite{Tovbis2013} because of Benjamin-Feir instability \cite{Benjamin-Feir} or modulation instability \cite{Zakharov2009}. The defocusing NLS equation \cite{Dodson2019,Merle2022} can also describe many physical systems with nonlinear interactions, and compared with the focusing NLS equation, certain properties of the defocusing NLS equation may be more amenable to mathematical analysis  due to the stability of the dark solitons \cite{Cuccagna2016}. It is worth noting out that the nonlinearity term in the defocusing NLS equation, while opposite in sign to that of the focusing type, does not necessarily mean that its physical meaning is completely opposite \cite{Cuccagna2016}. Instead, it leads to different physical and mathematical properties in some aspects \cite{El-Phys.D1995,Jenkins2015}. Moreover, although the defocusing nonlinearity can cause a reduction in the intensity of an optical beam as it propagates through a medium, it can also lead to the formation of stable soliton trains in certain conditions \cite{Kamchatnov2002}, which can be used to counteract the effect of linear dispersion, and thus enable long-distance transmission of optical signals without distortion.
\par
This work considers the one-dimensional defocusing NLS equation
\begin{equation}
 \mathrm{i} q_t+\frac{1}{2}q_{xx}-|q|^2q=0, ~~~ x\in \mathbb{R},~ t \geqslant 0,  \label{(NLS)}
\end{equation}
subject to the step-like initial data
\begin{equation} \label{(initial)}
q(x,0)=q_0(x):=
\left\{  \begin{array}{ll}
A_le^{-2\mathrm{i} \mu_l x}, ~~~~~x<0,\\
A_re^{-2\mathrm{i} \mu_r x}, ~~~~~x>0,\\
\end{array} \right.
\end{equation}
where $ A_l$, $A_r$, $\mu_l$ and $\mu_r$ are four real constants with $A_l>0$ and $A_r>0$. This family of step-like initial data naturally corresponds to the plane wave solutions $q^{\mathrm{p}}_j (x,t)=A_j e^{-2\mathrm{i} \mu_j x-\mathrm{i}(A_j^2+2\mu_j^2) t},~ j\in\{l, r\}$.
To ensure that the initial value problem (\ref{(NLS)}) with (\ref{(initial)}) makes sense, the following additional conditions should be satisfied:
\begin{equation}
  \int_{0}^{\pm \infty} |q(x,t) - q^{\mathrm{p}}_j (x,t) | \,\dif x < \infty, \qquad \forall t \geqslant 0, \qquad j\in\{l, r\}.
\end{equation}
\par

Introducing the Madelung transform
\begin{equation}
\label{Madelung} q(x,t)=\sqrt{\rho(x,t)}e^{\mathrm{i}\phi(x,t)},\quad \partial_x\phi(x,t)=v(x,t),
\end{equation}
the defocusing NLS equation (\ref{(NLS)}) is converted into the hydrodynamic form
\begin{align}\label{h1}
&\rho_t+(\rho v)_x =0,\\
& v_t+vv_x+\rho_x=
\frac{1}{4}\partial_x\left(\frac{\partial_{xx}^2\rho}{\rho}-\frac{(\partial_x\rho)^2}{2\rho^2}\right).
\label{h2}
\end{align}
The zero-dispersion limit of the defocusing NLS hydrodynamics yields the Euler equation for the ideal compressible fluid with velocity $v$ and density $\rho$ associated with pressure $P=\rho^2/2.$ With the Madelung transform (\ref{Madelung}) in mind, the step-like initial data (\ref{(initial)}) for $q(x,0)$ can be transformed into the step-like initial data for $\rho(x,0)$ and $v(x,0)$ as
\begin{equation} \label{initial-rho-v}
\rho(x,0)=
\left\{  \begin{array}{ll}
\rho_l, ~~~x<0,\\
\rho_r, ~~~x>0,\\
\end{array} \right.\quad \quad
v(x,0)=
\left\{  \begin{array}{ll}
v_l, ~~~x<0,\\
v_r, ~~~x>0,\\
\end{array} \right.
\end{equation}
where $\rho_l, \rho_r$ are positive constants and $v_l, v_r$ are arbitrary constants. This is usually called the Riemann problem of the defocusing NLS hydrodynamics.
\par
Dispersive systems with discontinuous initial conditions may experience rapidly oscillatory waves, i.e., dispersive shock waves (DSWs) \cite{Hoefer2014}. In fact, the DSW can be regarded as the counterpart of the classical shock wave in the dispersive case. In the long-time range, the DSW can be described by slowly modulated periodic waves, which is embodied in the Whitham modulation theory proposed by  Whitham \cite{Whitham1974}. Another important wave structure observed in dispersive systems is the rarefaction wave \cite{Kamchatnov2021}, which is continuous and different from the DSW and does not involve breaking. The study of the Riemann problem for dispersive hydrodynamics dates back to the work of Gurevich and Pitaevskii \cite{GP1974} for the KdV equation, in which a general approach was formulated to investigate the problem of initial discontinuity based on the Whitham modulation theory. Thus, the Riemann problem of dispersive systems is often named the Gurevich-Pitaevskii problem. Much work in this direction has been carried out for integrable systems and nonintegrable systems \cite{Congy2021}.
\par
The exploration of the Riemann problem for defocusing NLS hydrodynamics (\ref{(NLS)}) with (\ref{(initial)}) is a significant and challenging topic in both physics and mathematics. El et al. \cite{El-Phys.D1995} gave a complete classification of the solutions to the special Riemann problem of the defocusing NLS equation based on Whitham modulation theory, and demonstrated that there are six possible wave structure types. Jenkins \cite{Jenkins2015} reexamined this problem by the Riemann-Hilbert formulation and gave the zero-dispersion limit/long-time solution of one type of wave structure. Trillo et al. experimentally investigated the dam-break Riemann problem \cite{Trillo2017} and piston Riemann problem  \cite{Trillo2022} in quantum hydrodynamics modeled by the defocusing NLS equation and clearly observed the propagations of DSWs and rarefaction waves.
\par
The Riemann-Hilbert approach based on the Deift-Zhou nonlinear steepest descent technique \cite{Deift-Zhou1993} is an effective way to study the long-time asymptotic behaviors of the solution in integrable systems \cite{Bilman-Miller-2019}-\cite{MaWX2022}. Recently, a great deal of work has been done to investigate the step-like initial problems of the NLS equation. For example, Buckingham and Venakides \cite{Buckingham2007} considered the long-time asymptotics of the shock problem for the focusing NLS equation. Boutet de Monvel, Kotlyarov and Shepelsky \cite{Shepelsky2011} studied the long-time asymptotics of the solution to the pure step-like initial value problem for the focusing NLS equation and found three different solution regions. Jenkins and McLaughlin \cite{Jenkins2014} examined the zero-dispersion limit of the box problem for the focusing NLS equation. Biondini and Mantzavinos \cite{Biondini2017} analyzed the long-time asymptotic behaviors of the focusing NLS equation with a symmetric
nonvanishing boundary condition and then considered the same problem in the presence of
the discrete spectrum \cite{Biondini2021}. Boutet de Monvel, Lenells and Shepelsky \cite{Lenells2021} gave a full asymptotic analysis of the general step-like initial problem for the focusing NLS equation. There are only a few studies in the literatures on the step-like initial problems of the defocusing NLS equation and the most remarkable one belongs to the work of Jenkins \cite{Jenkins2015}. Until now, rigorous asymptotic analysis for the general step-like initial problem (\ref{(initial)}) (i.e., the Riemann problem) of the defocusing NLS equation (\ref{(NLS)}) is still open. Thus in this work, we make use of the nonlinear steepest descent method for Riemann-Hilbert problems (RHPs) toward addressing this problem.
\par
To make this work self-contained, some results of Whitham modulation theory for the Riemann problem of the defocusing NLS equation should be reviewed before presenting the rigorous asymptotic analysis. In what follows, the zero-phase modulated solution and one-phase modulated elliptic wave along with the corresponding Whitham equation for the defocusing NLS equation are listed \cite{Gong-Wang2023}. The complete classification of the initial value problem (\ref{(initial)}) of equation (\ref{(NLS)}) is described in two kinds of coordinates. Finally, some notations are given briefly in advance.

\subsection{The modulated (unmodulated) solutions and Whitham equations}\

This subsection lists the modulated (unmodulated) solutions and the corresponding Whitham equations for the defocusing NLS equation (\ref{(NLS)}) based on Whitham modulation theory \cite{Congy-AMK-2017,Gong-Wang2023}.

\subsubsection{The zero-phase modulated solution and Whitham equations}
\
\newline
\indent
The plane wave solution of equation (\ref{(NLS)}) can be obtained easily, so the zero-phase modulated solution is
\begin{equation} \label{g=0-solution}
q(x,t)=\sqrt{\rho_0}e^{\mathrm{i} \phi_0(x,t)},
\end{equation}
where the phase $\phi_0(x,t)=v_{0}x-(v_{0}^{2}/2+\rho_0)t$ is a fast variable, while the amplitude $\sqrt{\rho_0}$ and wavenumber $v_{0}$, which depend on the Whitham equations of genus zero, are slowly varying. Taking the Riemann invariants $\lambda_1$ and $\lambda_2$ to satisfy $\lambda_1+\lambda_2=-v_0$ and $\lambda_1\lambda_2=v_0^2/4-\rho_0$, the Whitham equations for $\lambda_1$ and $\lambda_2$ are given by
\begin{equation}
\label{Whitham-0-phase}
\begin{aligned}
&\frac{\partial \lambda_1}{\partial t}-(\frac{3}{2}\lambda_1 +\frac{1}{2}\lambda_2)
\frac{\partial \lambda_1}{\partial x}=0, \\
&\frac{\partial \lambda_2}{\partial t}-(\frac{1}{2}\lambda_1 +\frac{3}{2}\lambda_2)
\frac{\partial \lambda_2}{\partial x}=0,
\end{aligned}
\end{equation}
where $\lambda_1>\lambda_2$. Solving the Whitham equations (\ref{Whitham-0-phase}) for the Riemann invariants $\lambda_1$ and $\lambda_2$ and then for the modulated amplitude $\sqrt{\rho_0}$ and wavenumber $v_{0}$, the zero-phase modulated solution in (\ref{g=0-solution}) can be formulated explicitly.

\subsubsection{The one-phase modulated elliptic wave and Whitham equations}
\
\newline
\indent
When considering the one-phase modulated elliptic wave, the number of Riemann invariants increases to four and the variables of the defocusing NLS hydrodynamics (\ref{h1})-(\ref{h2}) are expressed by the Jacobi elliptic function as
\begin{align}\label{rho-equation-solution-new}
&\rho(x,t)=\rho_2-(\rho_2-\rho_3) {\rm cn}^2(\sqrt {\rho_1-\rho_3}(x-Vt),m),\\
& v(x,t)=V+\frac{\sqrt{\rho_1\rho_2\rho_3}}{\rho(x,t)},
\end{align}
with the parameters $\rho_1=\frac{1}{4}(\lambda_1+\lambda_2-\lambda_3-\lambda_4)^2, \rho_2=\frac{1}{4}(\lambda_1-\lambda_2+\lambda_3-\lambda_4)^2, \rho_3=\frac{1}{4}(\lambda_1-\lambda_2-\lambda_3+\lambda_4)^2,$ $V=-\frac{1}{2}\sum_{i=1}^{4}\lambda_i$ and the modulus of the Jacobi elliptic function $m=\sqrt{ {(\rho_2-\rho_3)}/{(\rho_1-\rho_3)}}$. When $m\rightarrow 1$, the one-phase modulated elliptic wave degenerates into a soliton solution, while $m\rightarrow 0$ produces the harmonic front. Here, the Riemann invariants $\lambda_1, \lambda_2, \lambda_3$ and $\lambda_4$ with $\lambda_1>\lambda_2>\lambda_3>\lambda_4$ solve the Whitham equations of genus one as
\begin{equation} \label{Whitham1234-1234}
\frac{\partial \lambda_i}{\partial t}+v_i(\lambda_i)\frac{\partial \lambda_i}{\partial x}=0,\quad i=1, 2, 3, 4,
\end{equation}
where the characteristic velocities $v_i$~$(i=1,2,3,4)$ are
\begin{equation} \label{Whitham1234-1234-charact-velocity}
 \begin{split}
&v_1=-\frac{1}{2}\sum_{i=1}^{4}\lambda_i
-\frac{(\lambda_1-\lambda_4)(\lambda_1-\lambda_2)~K(m)}{(\lambda_1-\lambda_4)~K(m)+(\lambda_4-\lambda_2)~E(m)},\\
&v_2=-\frac{1}{2}\sum_{i=1}^{4}\lambda_i
+\frac{(\lambda_2-\lambda_3)(\lambda_1-\lambda_2)~K(m)}{(\lambda_2-\lambda_3)~K(m)+(\lambda_3-\lambda_1)~E(m)},\\
&v_3=-\frac{1}{2}\sum_{i=1}^{4}\lambda_i
-\frac{(\lambda_2-\lambda_3)(\lambda_3-\lambda_4)~K(m)}{(\lambda_2-\lambda_3)~K(m)+(\lambda_4-\lambda_2)~E(m)},\\
&v_4=-\frac{1}{2}\sum_{i=1}^{4}\lambda_i
+\frac{(\lambda_1-\lambda_4)(\lambda_3-\lambda_4)~K(m)}{(\lambda_1-\lambda_4)~K(m)+(\lambda_3-\lambda_1)~E(m)},
\end{split}
\end{equation}
where $K(m)$ and $E(m)$ are the first kind and second kind complete elliptic integrals, respectively.

\subsubsection{The unmodulated elliptic wave}
\
\newline
\indent
Since the Whitham equations (\ref{Whitham1234-1234}) with (\ref{Whitham1234-1234-charact-velocity}) are strictly hyperbolic \cite{Kodama2006}, the number of nonconstant Riemann invariants (called soft edges) is at most one.
Thus, there may exist an unmodulated elliptic wave region when taking special values in the step-like initial data (\ref{(initial)}). In this case, the Riemann invariants $\lambda_1, \lambda_2, \lambda_3$ and $\lambda_4$ are all constants (called hard edges), which satisfy the Whitham equations (\ref{Whitham1234-1234}) automatically. This case appears in the middle region of Case A below.

\subsection{The complete classification of asymptotic solutions to the initial problem}
\
\newline
\indent
In fact, the Whitham equations of genus zero and the corresponding Riemann invariants can also be determined by the zero-dispersion limit of the defocusing NLS hydrodynamics (\ref{h1})-(\ref{h2}), i.e.,
the following zero-dispersion hydrodynamic system
\begin{equation}\label{non-dispersion}
\begin{pmatrix} \rho \\ v \end{pmatrix}_t + \begin{pmatrix} v & \rho \\ 1 & v
\end{pmatrix} \begin{pmatrix} \rho \\ v \end{pmatrix}_x=0,
\end{equation}
which can be used to describe surface wave motion in shallow water. Diagonalizing the system (\ref{non-dispersion}), yields the same Whitham equations of genus zero as (\ref{Whitham-0-phase}), where the Riemann invariants are $\lambda_1=-{v}/{2}+\sqrt{\rho}$ and $\lambda_2=-{v}/{2}-\sqrt{\rho}$.
For the step-like initial data (\ref{initial-rho-v}) (corresponding to (\ref{(initial)})), the initial Riemann invariants are given by
\begin{equation}
  \begin{aligned}
  \label{Riin}
  &\lambda_l^+=-\frac{v_l}{2}+\sqrt{\rho_l} = \mu_l+A_l, &\qquad \lambda_l^- =-\frac{v_l}{2}-\sqrt{\rho_l} = \mu_l-A_l, \\
  &\lambda_r^+=-\frac{v_r}{2}+\sqrt{\rho_r} = \mu_r+A_r, &\qquad \lambda_r^-=-\frac{v_r}{2}-\sqrt{\rho_r} = \mu_r-A_r.
\end{aligned}
\end{equation}
\par

\begin{figure}[htbp]
  \centering
  \subfigure[]{\includegraphics[width=7cm]{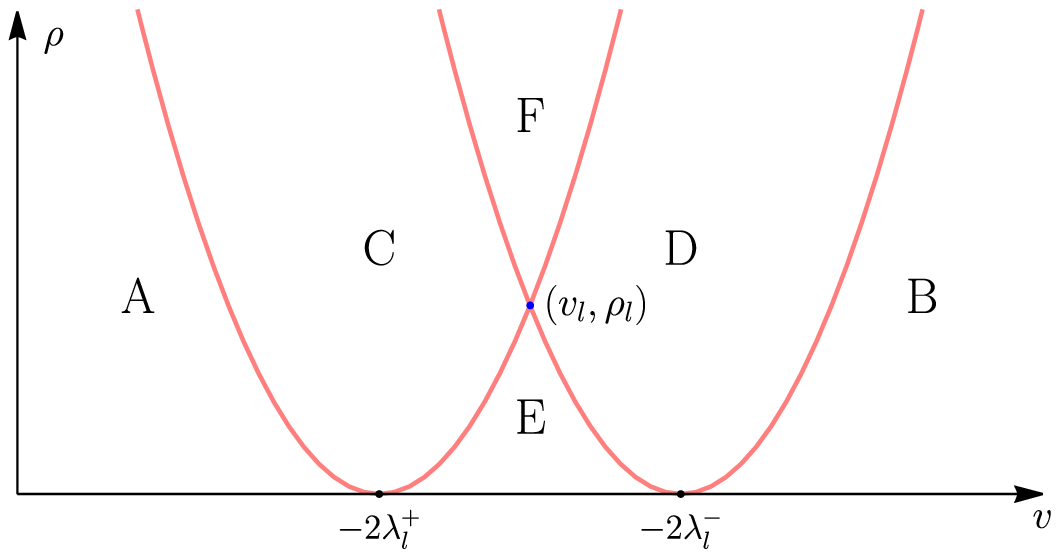}}\quad
  \subfigure[]{\includegraphics[width=7cm]{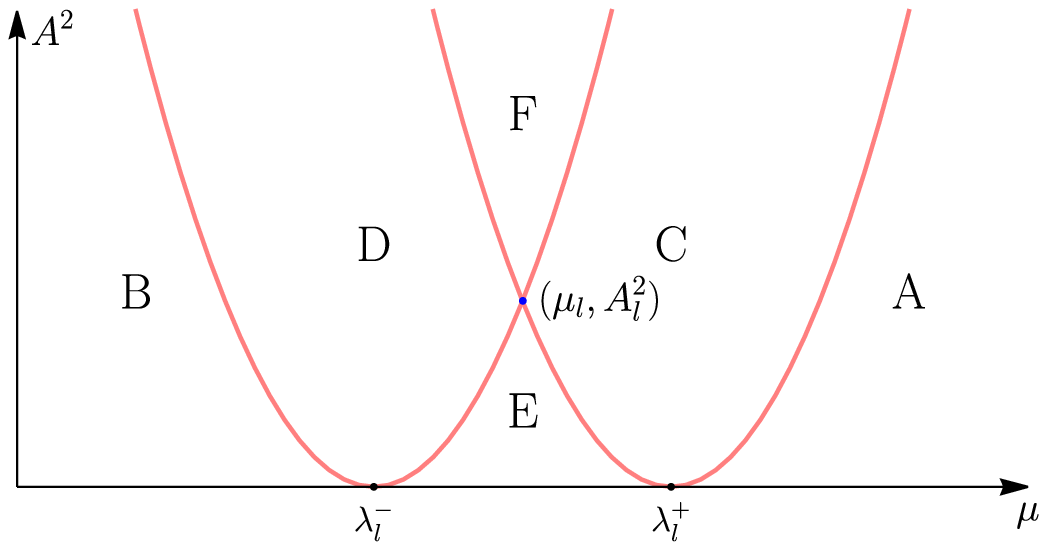}}
  \caption{{\protect\small  Classification diagrams of the asymptotic solutions in upper $(v, \rho)$-plane and upper $(\mu, A^2)$-plane, respectively.}}\label{Classifaction}
\end{figure}

Based on the above analysis, the general step-like initial conditions (\ref{(initial)}) of equation (\ref{(NLS)}) are fully classified \cite{Gong-Wang2023}. To give the complete classification of the asymptotic solutions, the left initial value $(v_l,\rho_l)$ is fixed first, and then it is naturally observed that the following two upward parabolas defined by
\begin{equation}\label{NLS-classification-1-1}
\rho=\frac{1}{4}({v}+2\lambda_l^+)^2\quad {\rm and} \quad \rho=\frac{1}{4}({v}+2\lambda_l^-)^2,
\end{equation}
along with the $v$-axis divide the upper $(v, \rho)$-plane into six regions, as shown in Figure \ref{Classifaction} (a). The two solid parabolas in Figure \ref{Classifaction} (a) correspond to the fixed Riemann invariants $\lambda_l^+$ and $\lambda_l^-$ with respect to the fixed left boundary. Similarly, the classification diagram for the initial data (\ref{(initial)}) in the upper $(\mu, A^2)$-plane is shown in Figure \ref{Classifaction} (b). Additionally, six regions are found by arranging the Riemann invariants in order, which are listed as follows:

\begin{equation}\label{Classification}
\begin{aligned}
&{\rm Case~ A:}~   \lambda_r^+>\lambda_r^->\lambda_l^+>\lambda_l^-, \quad
{\rm Case~ B:}~   \lambda_l^+>\lambda_l^->\lambda_r^+>\lambda_r^-, \\
&{\rm Case~ C:}~   \lambda_r^+>\lambda_l^+ >\lambda_r^->\lambda_l^-, \quad
{\rm Case~ D:}~  \lambda_l^+>\lambda_r^+>\lambda_l^->\lambda_r^-, \\
&{\rm Case~ E:} ~  \lambda_r^+>\lambda_l^+> \lambda_l^->\lambda_r^-, \quad
{\rm Case~ F:}~   \lambda_l^+>\lambda_r^+>\lambda_r^->\lambda_l^-. \\
\end{aligned}
\end{equation}
\par
The main purpose of this work is to formulate the long-time asymptotics of the solutions of the defocusing NLS equation (\ref{(NLS)}) in all six cases above. In particular, the theoretical asymptotic solutions obtained by the Riemann-Hilbert formulation will be compared with the results given by the Whitham modulation theory and direct numerical simulation. To the best of our knowledge, this is the first work to check the correctness of long-time asymptotic solutions in two different ways.

\subsection{Notation}
\
\newline
\indent
For brevity, denote the intervals by $\mathcal{I}_l=(\lambda_l^-, \lambda_l^+)$ and $\mathcal{I}_r=(\lambda_r^-, \lambda_r^+)$.
For a scalar function $f(z)$, let $f^*(z)=\overline{f(\overline{z})}$ denote the Schwarz reflection about the real axis.
Given a piecewise oriented smooth contour $\Sigma$ and a function $f(z)$ analytic in $z \in \mathbb{C} \backslash \Sigma$, denote by $f_{\pm}(z)$ for $z\in \Sigma$ the boundary values from the left
and right sides of $\Sigma$, repectively. We make use of the Pauli matrix $\sigma_3=\begin{pmatrix}  1 & 0 \\ 0 & -1 \end{pmatrix}$ and the matrix  notation $f(z)^{\sigma_3}=\begin{pmatrix}  f(z) & 0 \\ 0 & f^{-1}(z) \end{pmatrix}$.
Given two real numbers $\lambda_1$ and $\lambda_2$ such that $\lambda_1 > \lambda_2$,
we introduce the functions
\begin{equation}
  \begin{aligned}
  \beta(z;\lambda_1,  \lambda_2)= \left( \frac{z-\lambda_1}{z-\lambda_2} \right)^{1/4}, \qquad \mathcal{R}(z;\lambda_1,  \lambda_2)=\sqrt{(z-\lambda_1)(z-\lambda_2)}, \\
  \mathcal{E} (z;\lambda_1,  \lambda_2)=\begin{pmatrix} \frac{\beta(z;\lambda_1,  \lambda_2)+\beta^{-1}(z;\lambda_1,  \lambda_2)}{2} & -\frac{\beta(z;\lambda_1,  \lambda_2)-\beta^{-1}(z;\lambda_1,  \lambda_2)}{2i}
    \\ \frac{\beta(z;\lambda_1,  \lambda_2)-\beta^{-1}(z;\lambda_1,  \lambda_2)}{2i} & \frac{\beta(z;\lambda_1,  \lambda_2)+\beta^{-1}(z;\lambda_1,  \lambda_2)}{2}\end{pmatrix},
  \end{aligned}
  \end{equation}
where the branches are chosen such that these functions are analytic in $\mathbb{C} \backslash [\lambda_2, \lambda_1]$ and satisfy the large $z$ asymptotics of the forms
\begin{equation}
  \beta(z;\lambda_1,  \lambda_2)=1+ \mathcal{O}(z^{-1}), \quad \mathcal{R}(z;\lambda_1,  \lambda_2)=z -\frac{\lambda_1 +\lambda_2}{2} + \mathcal{O}(z^{-1}), \quad
  \mathcal{E} (z;\lambda_1,  \lambda_2)=I+\mathcal{O}(z^{-1}).
\end{equation}
For a real-valued vector $\boldsymbol{\lambda}=(\lambda_1, \lambda_2, \lambda_3, \lambda_4)$ such that $\lambda_1 > \lambda_2 > \lambda_3 > \lambda_4$, define
\begin{equation}
  \mathcal{R}(z;\boldsymbol{\lambda})=\sqrt{(z-\lambda_1)(z-\lambda_2)(z-\lambda_3)(z-\lambda_4)},
\end{equation}
which is cut along $[\lambda_4, \lambda_3] \cup [\lambda_2, \lambda_1]$ such that $\mathcal{R}(z;\boldsymbol{\lambda})>0$ for $z>\lambda_1$.

\section{The main theorems}

\begin{figure}[tbp]
  \centering
  \subfigure[${\rm Case~ A:}~   \lambda_r^+>\lambda_r^->\lambda_l^+>\lambda_l^-$]{\includegraphics[width=7.5cm]{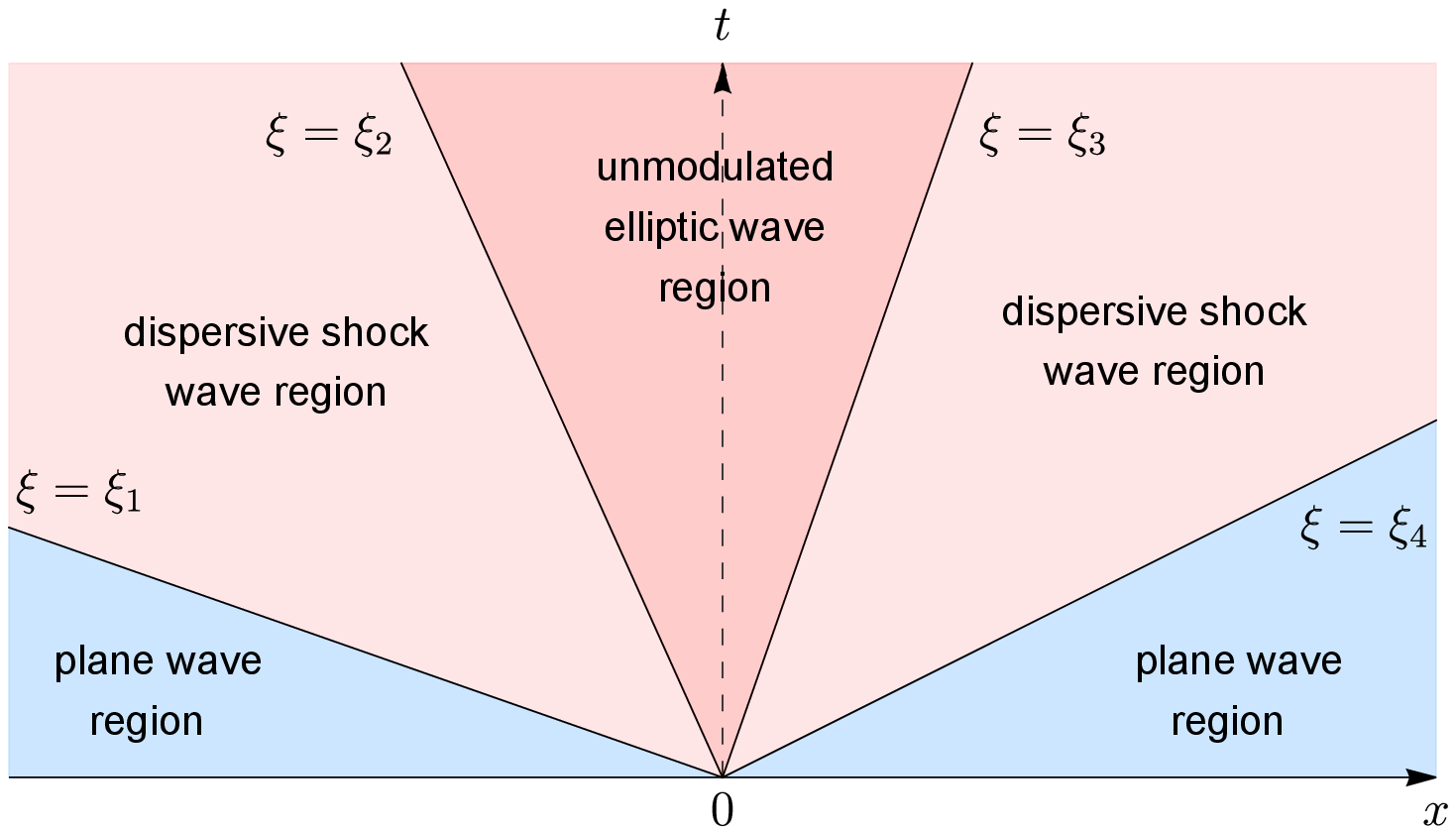}} \subfigure[${\rm Case~ B:}~   \lambda_l^+>\lambda_l^->\lambda_r^+>\lambda_r^-$]{\includegraphics[width=7.5cm]{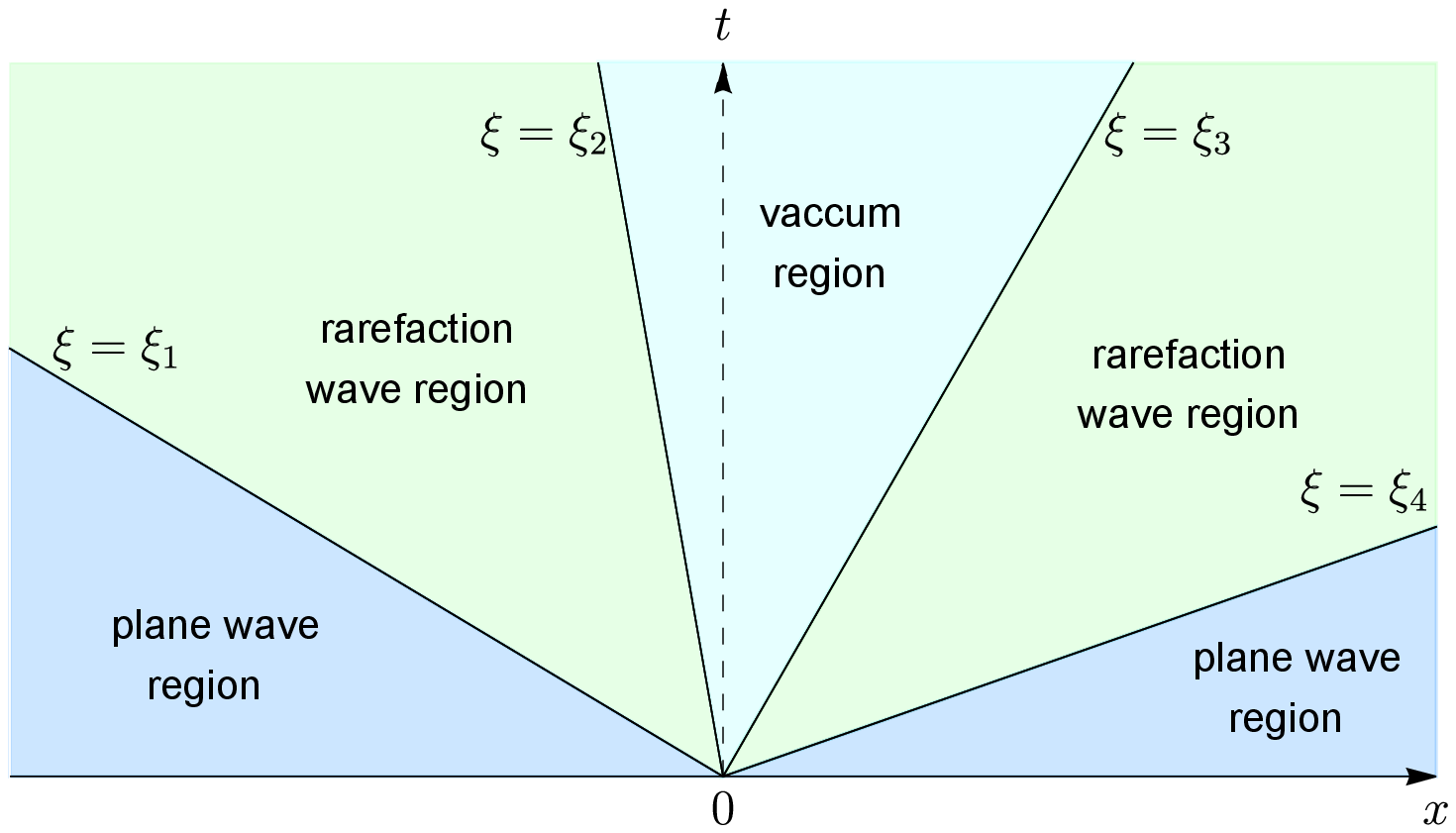}}

  \subfigure[${\rm Case~ C:}~   \lambda_r^+>\lambda_l^+ >\lambda_r^->\lambda_l^-$]{\includegraphics[width=7.5cm]{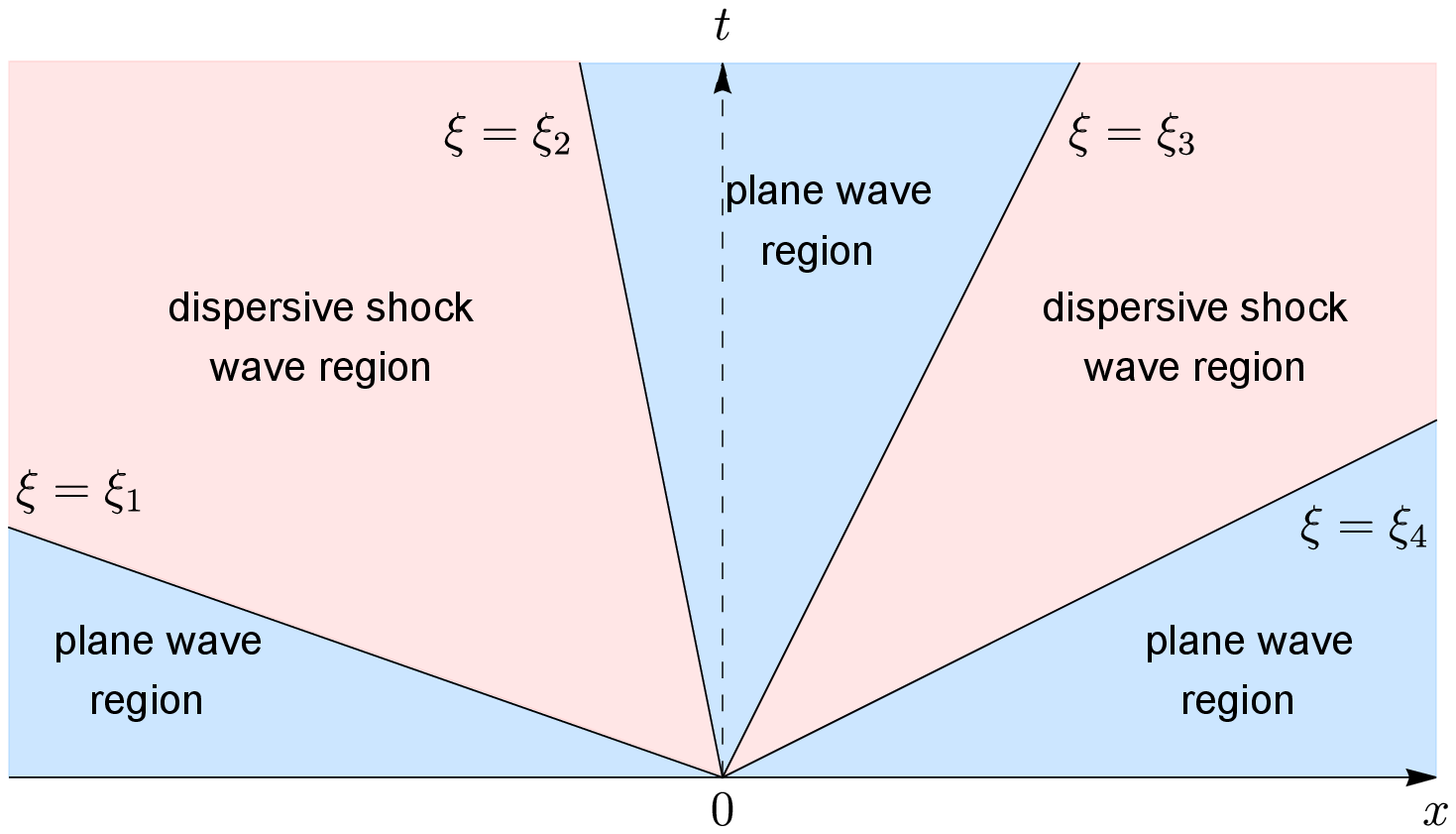}} \subfigure[${\rm Case~ D:}~  \lambda_l^+>\lambda_r^+>\lambda_l^->\lambda_r^-$]{\includegraphics[width=7.5cm]{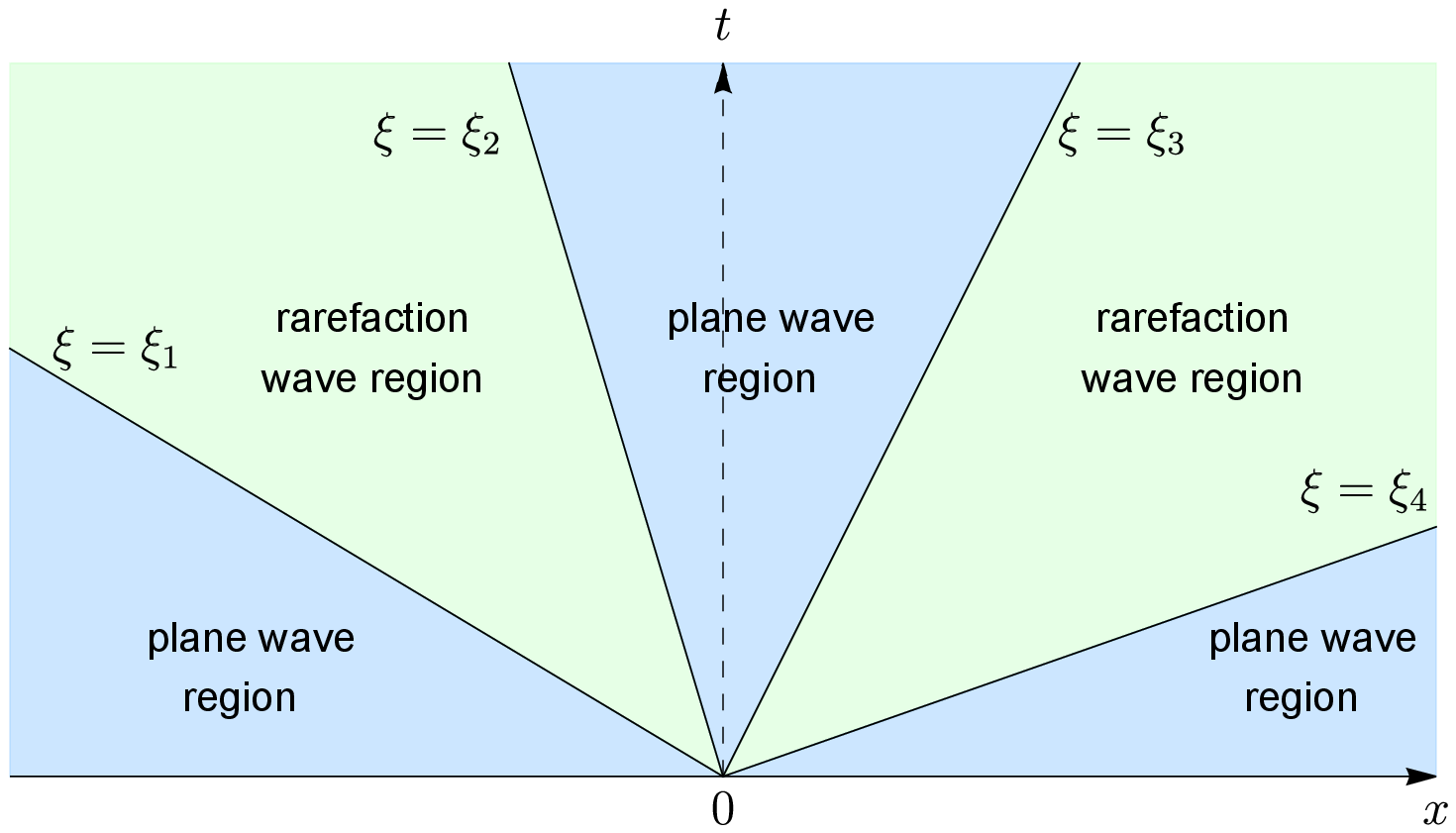}}

  \subfigure[${\rm Case~ E:} ~  \lambda_r^+>\lambda_l^+> \lambda_l^->\lambda_r^-$]{\includegraphics[width=7.5cm]{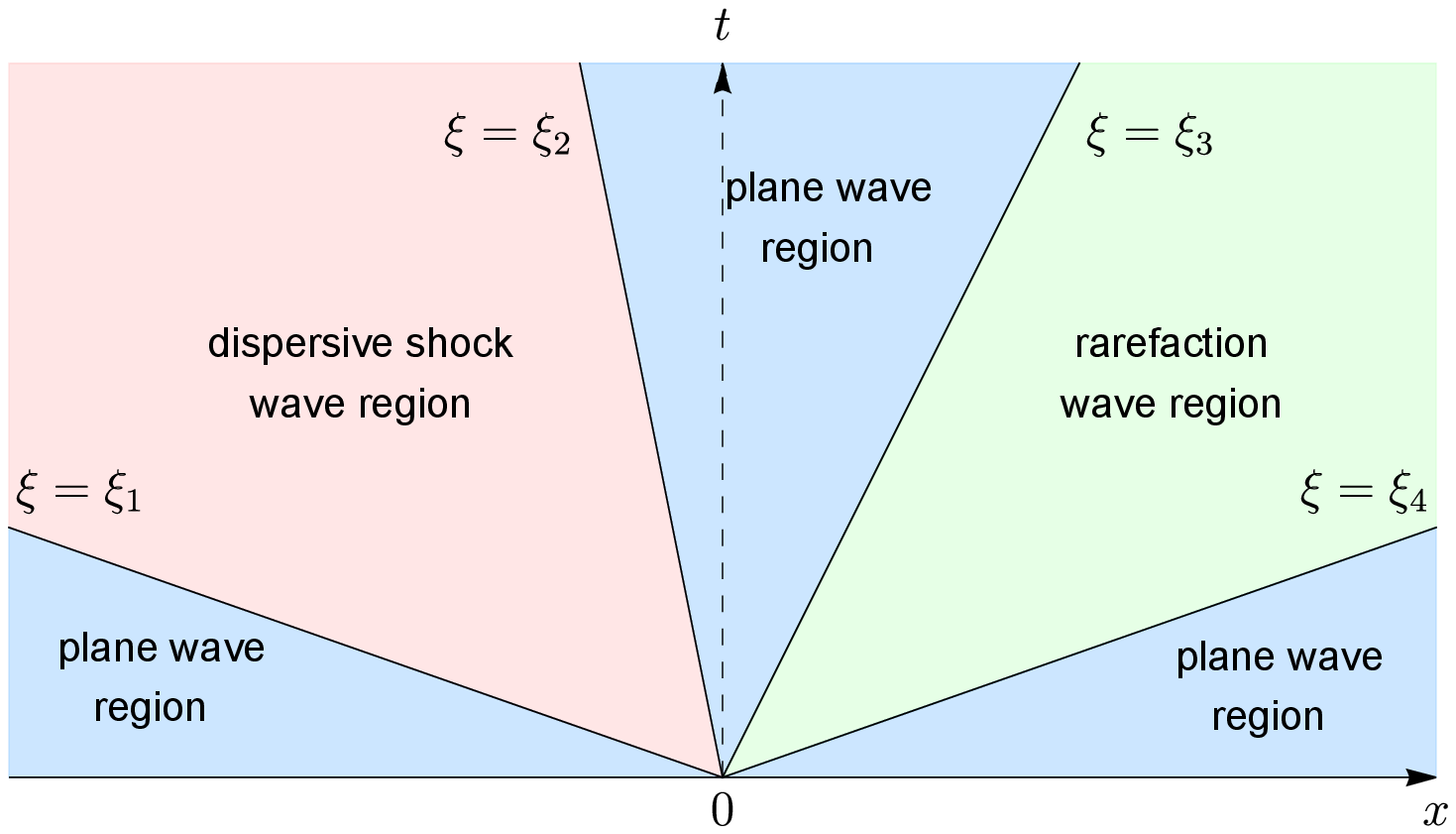}} \subfigure[${\rm Case~ F:}~   \lambda_l^+>\lambda_r^+>\lambda_r^->\lambda_l^-$]{\includegraphics[width=7.5cm]{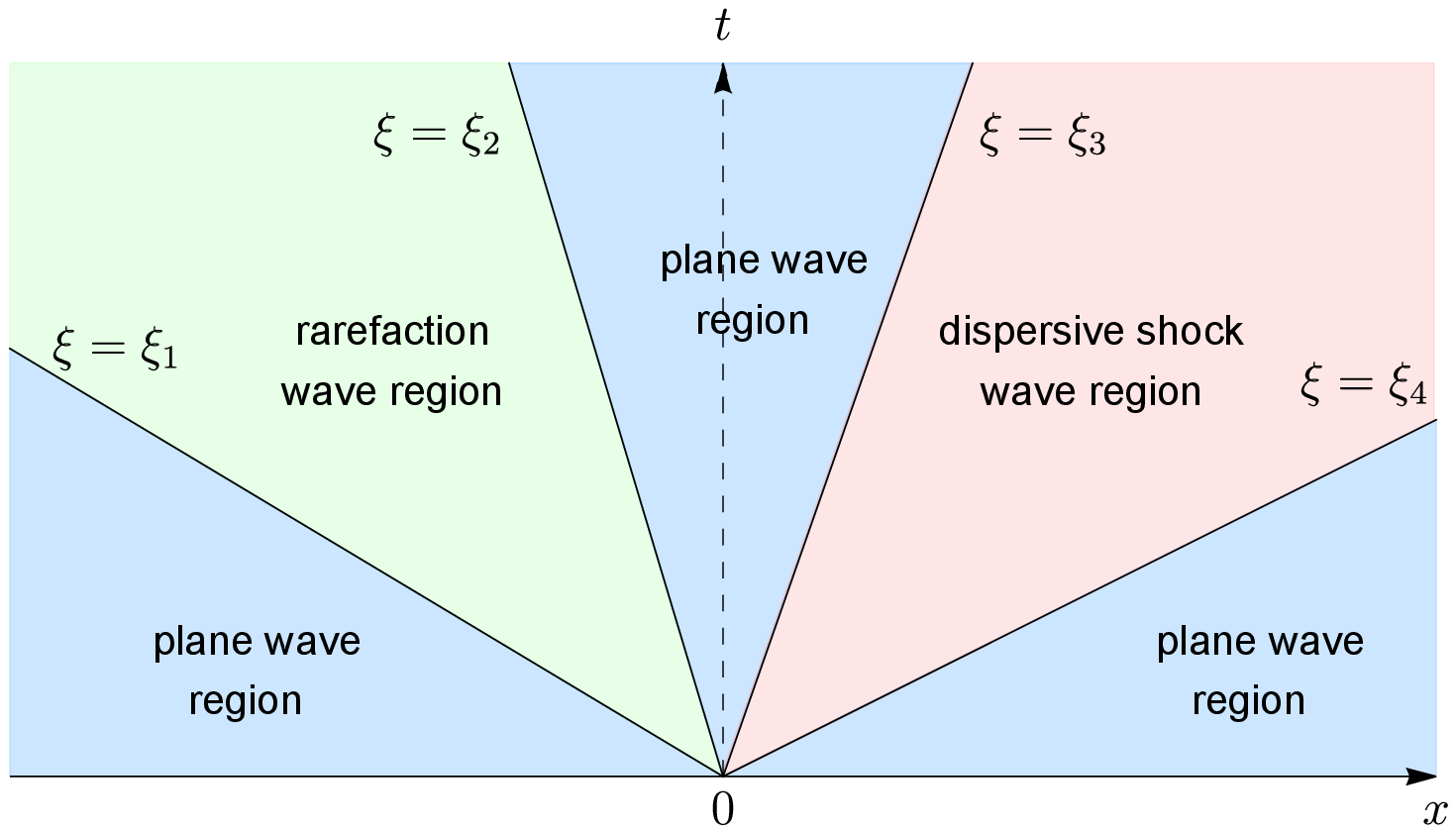}}
  \caption{{\protect\small The solution regions for each case of the classification in (\ref{Classification}) on the upper $(x,t)$-half-plane. }}\label{XT-1}
\end{figure}

\begin{figure}[htbp]
  \centering

  \subfigure[${\rm Case~ A:}~   \lambda_r^+>\lambda_r^->\lambda_l^+>\lambda_l^-$]{\includegraphics[width=8cm]{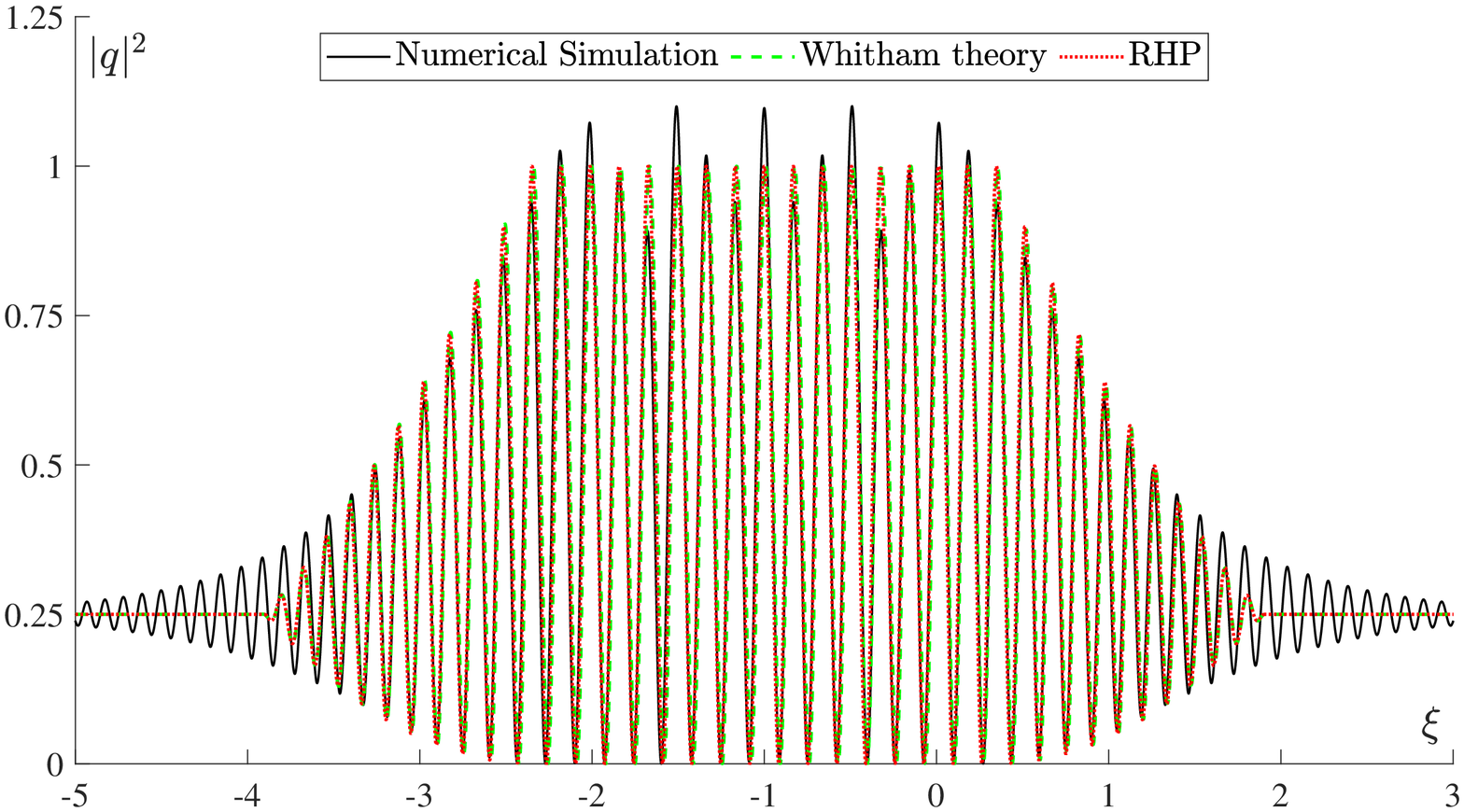}} \subfigure[${\rm Case~ B:}~   \lambda_l^+>\lambda_l^->\lambda_r^+>\lambda_r^-$]{\includegraphics[width=8cm]{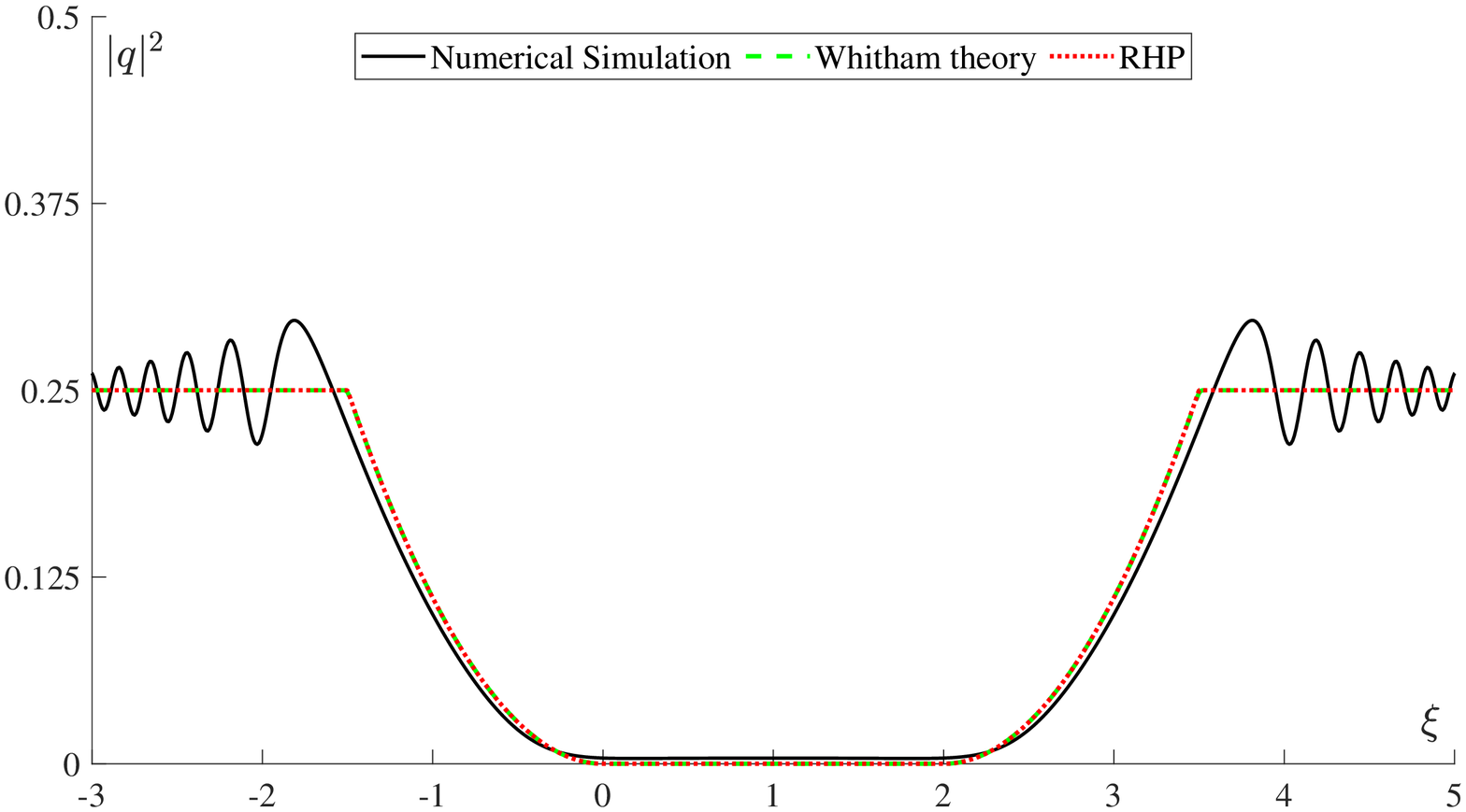}}

  \subfigure[${\rm Case~ C:}~   \lambda_r^+>\lambda_l^+ >\lambda_r^->\lambda_l^-$]{\includegraphics[width=8cm]{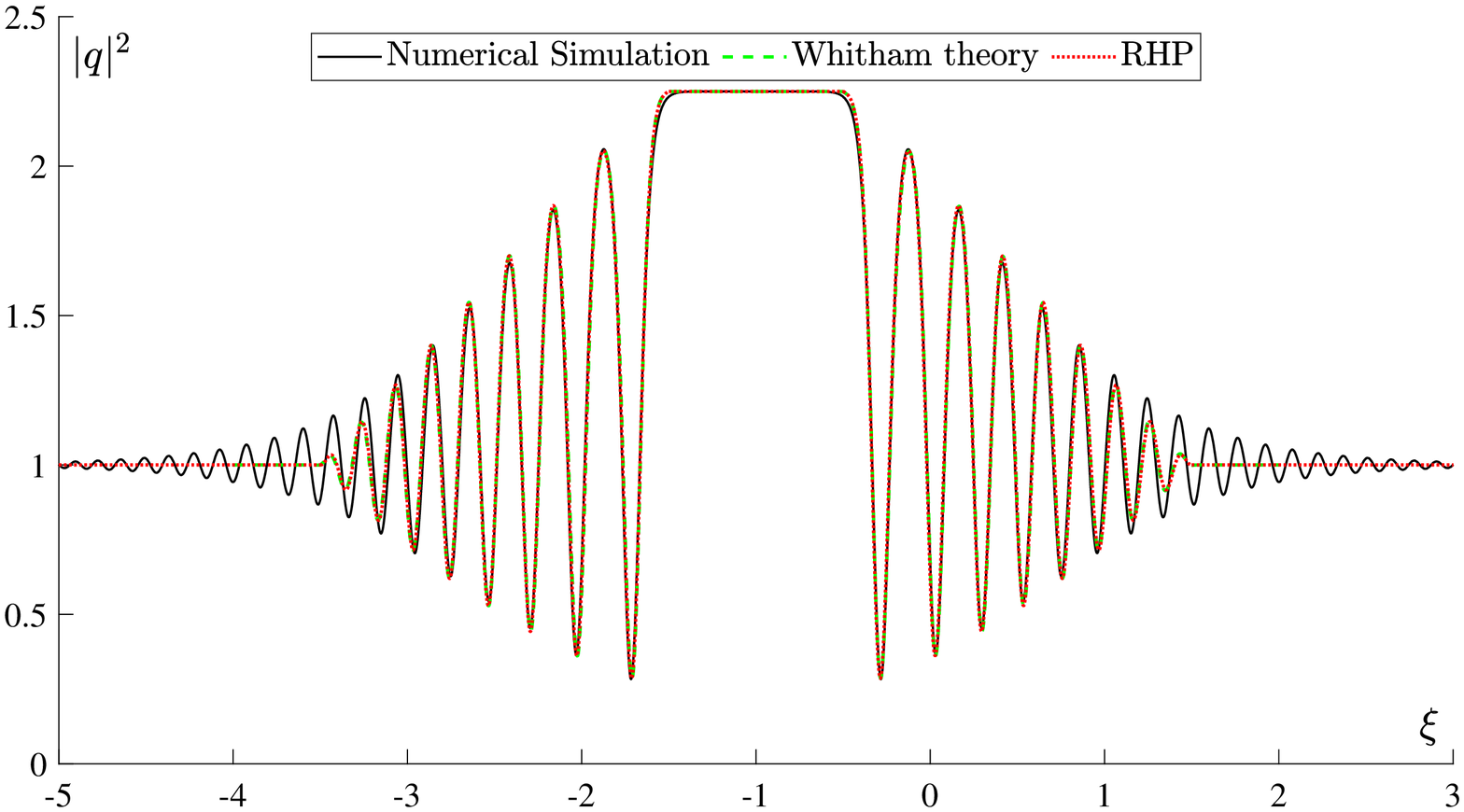}} \subfigure[${\rm Case~ D:}~  \lambda_l^+>\lambda_r^+>\lambda_l^->\lambda_r^-$]{\includegraphics[width=8cm]{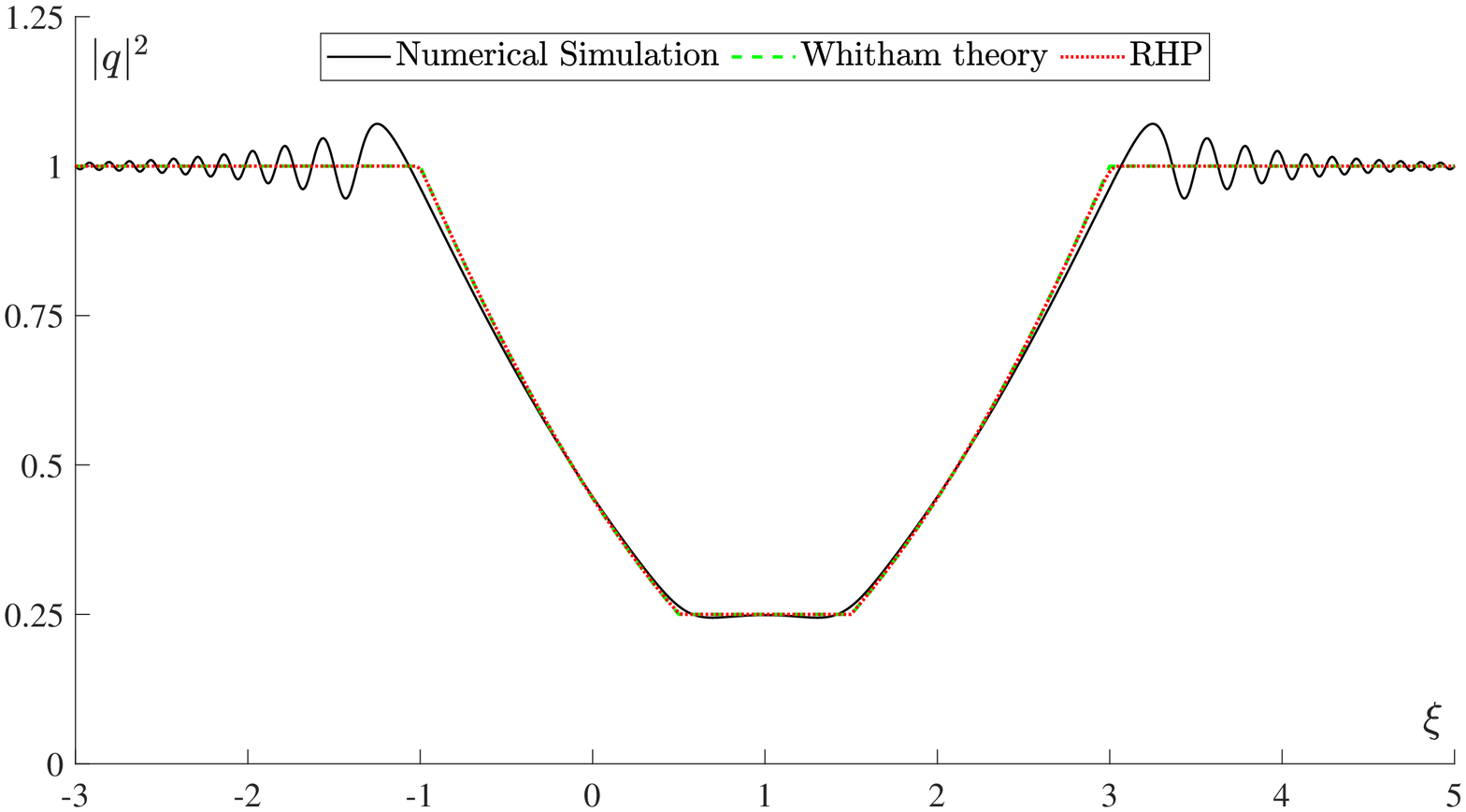}}

  \subfigure[${\rm Case~ E:} ~  \lambda_r^+>\lambda_l^+> \lambda_l^->\lambda_r^-$]{\includegraphics[width=8cm]{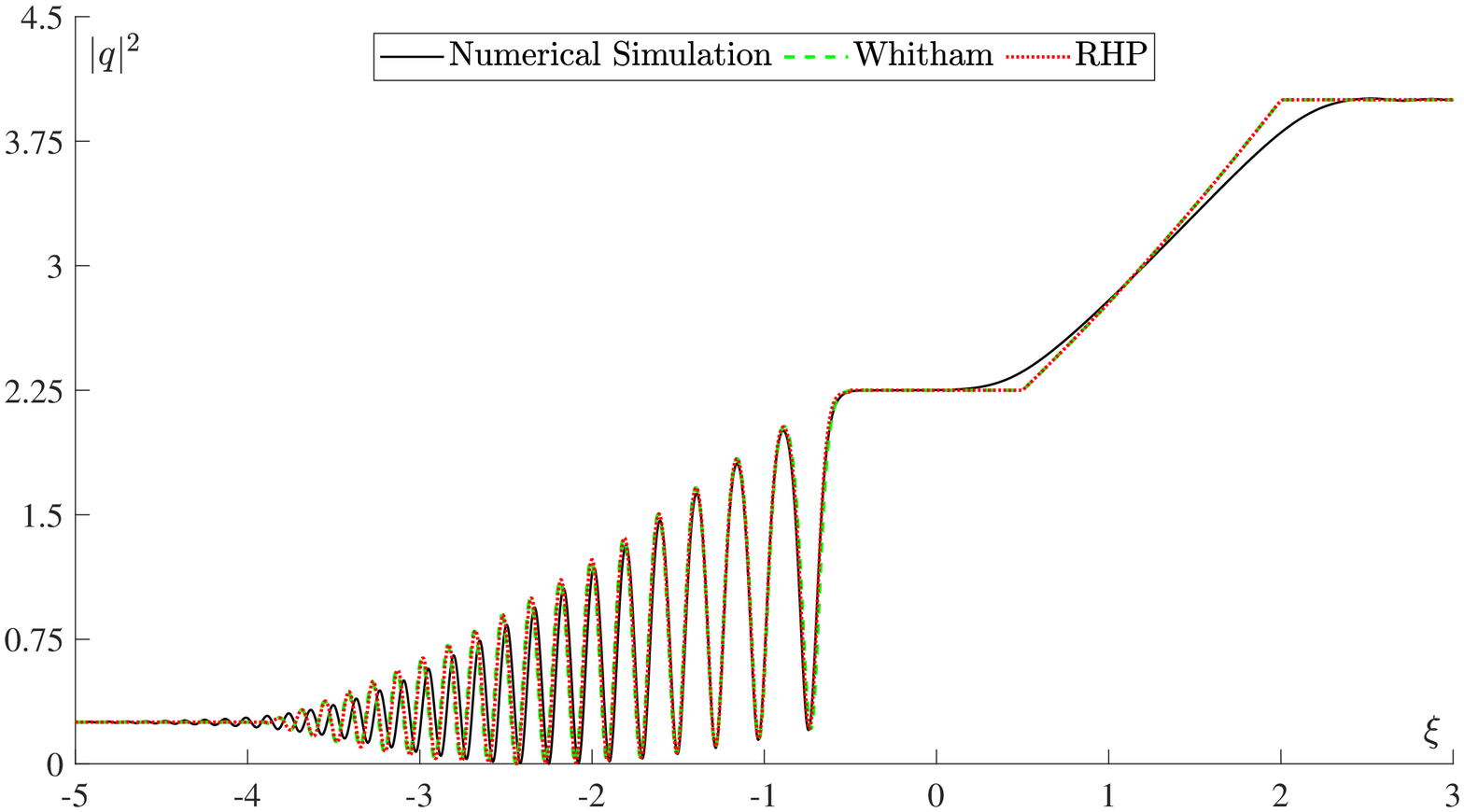}} \subfigure[${\rm Case~ F:}~   \lambda_l^+>\lambda_r^+>\lambda_r^->\lambda_l^-$]{\includegraphics[width=8cm]{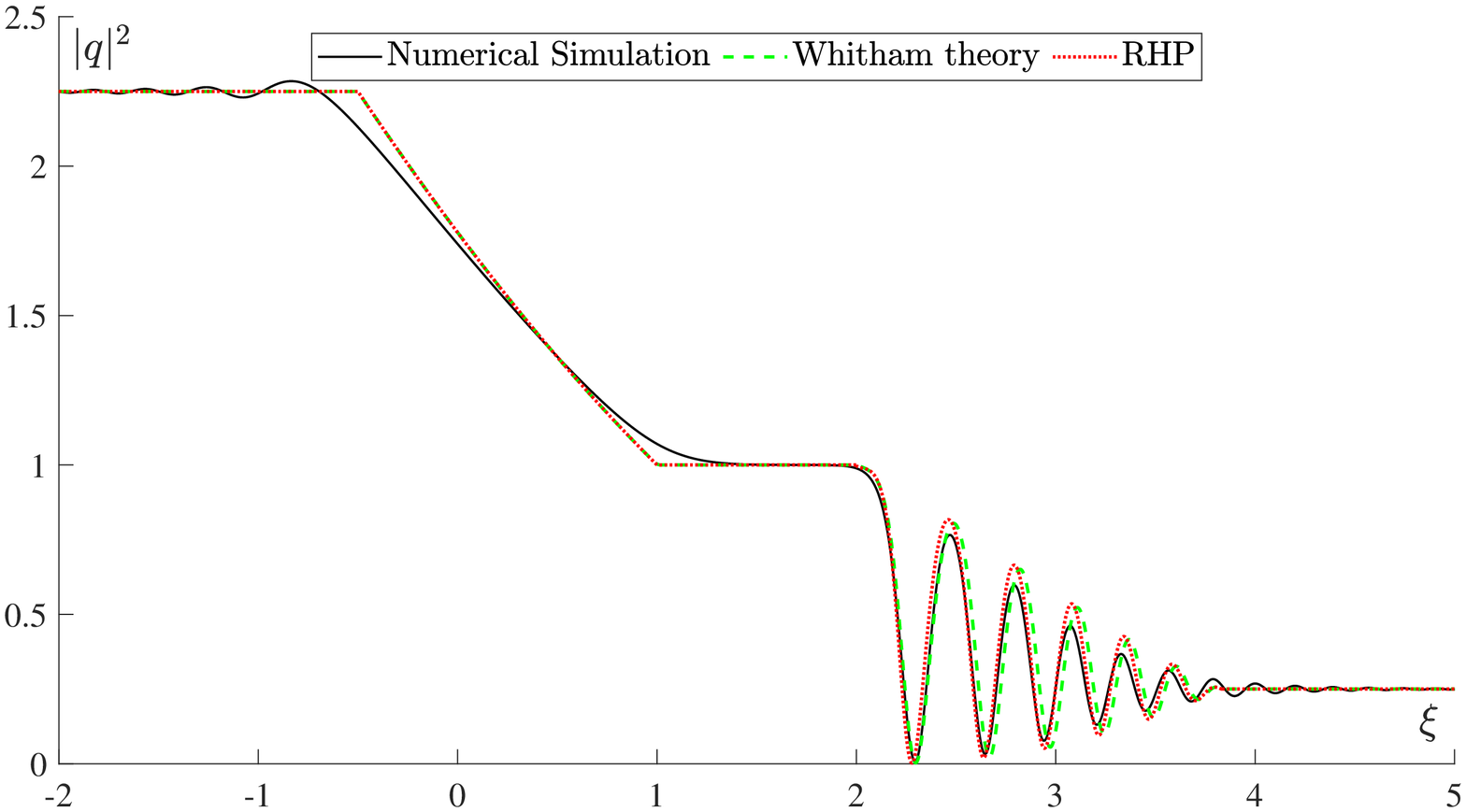}}
  \caption{{\protect\small The comparisons among the leading-order term from Riemann-Hilbert formulation in Theorems \ref{Shock-un}-\ref{Rare-Shock} (red dotted line), the results of Whitham modulation theory (green dashed line) and the numerical simulations (black solid line) for each case.
  The parameters in the initial data (\ref{(initial)}) are selected as follows:
  (a) $A_l=0.5, ~A_r=0.5, ~\mu_l=-0.5, ~\mu_r=1.5$; (b) $A_l=0.5, ~A_r=0.5, ~\mu_l=0.5, ~\mu_r=-1.5$;
  (c) $A_l=1, ~A_r=1, ~\mu_l=0, ~\mu_r=1$; (d) $A_l=1, ~A_r=1, ~\mu_l=0, ~\mu_r=-1$;
  (e)  $A_l=0.5, ~A_r=2, ~\mu_l=-0.5, ~\mu_r=0$; (f) $A_l=1.5, ~A_r=0.5, ~\mu_l=-0.5, ~\mu_r=-0.5$.
  }}
  \label{XT-2}
\end{figure}

This section proposes the main theorems that state the long-time asymptotic solutions of the defocusing NLS equation (\ref{(NLS)}) with the step-like initial data (\ref{(initial)}) for the six cases listed in (\ref{Classification}). It is shown that there are five types of regions in the classification, i.e., the plane wave region, rarefaction wave region, dispersive shock wave region, vacuum region and unmodulated elliptic wave region. Each case is regarded as the combination of the five types of regions that are stitched together in a particular way. Theorems \ref{Shock-un}-\ref{Rare-Shock} elaborate on the main results of this work, which are illustrated in Figure \ref{XT-1} and Figure \ref{XT-2}, respectively.

\begin{theorem}
  \label{Shock-un}
 {In Case {\rm A} $(\lambda_r^+>\lambda_r^->\lambda_l^+>\lambda_l^-)$,  the long-time asymptotic solution $q(x,t)$ of the defocusing NLS equation $(\ref{(NLS)})$ with initial data (\ref{(initial)}) behaves differently in five regions depending on the range of the self-similar variable $\xi=x/t$, as shown in Figure \ref{XT-1} (a) and Figure \ref{XT-2} (a). As $t\to \infty$, the leading-order term and error term in each region are formulated specifically}:
\par
{\rm (i)}. In the left plane wave region $(i.e., \xi<\xi_{A1}$ with $\xi_{A1}$ defined by equation $(\ref{A1B}))$, as $t \to \infty$,
 \begin{equation}\label{A1}
    q(x,t)=A_l e^{-2\mathrm{i} \mu_l x-\mathrm{i} (2 \mu_l^2+A_l^2)t}e^{-\mathrm{i}\phi_{lp}(\xi) }+\mathcal{O}(t^{-\frac{1}{2}}),
 \end{equation}
where the phase shift $\phi_{lp}(\xi)$ is given by equation $(\ref{Aphi1})$ and depends only on the parameters in the initial data $q(x,0)$ and the self-similar variable $\xi$.

{\rm (ii)}. In the dispersive shock wave region $(i.e., \xi_{A1}<\xi<\xi_{A2}$ with $\xi_{A1}$ and $\xi_{A2}$ defined by equation  $(\ref{A2B}))$, as $t \to \infty$,
 \begin{equation}\label{A2}
  q(x,t)=\frac{\lambda_r^+ - \lambda_{s}(\xi) + \lambda_l^+ -\lambda_l^-}{2}\frac{\Theta (0)~\Theta (2\mathcal{A}(\infty)+(\gamma + \hat{\gamma})/{2 \pi})}{\Theta ((\gamma + \hat{\gamma})/{2 \pi})~\Theta (2\mathcal{A}(\infty))} e^{2\mathrm{i} (t g_\infty+\hat{g}_\infty-\varphi^{(0)})}+\mathcal{O}(t^{-1}),
 \end{equation}
where $\Theta$ is the one-phase Riemann theta function defined by equation (\ref{Theta}). Here $\lambda_{s}$ is a movable Riemann invariant (or soft edge) given by the unique solution of equation (\ref{vv2}). The real quantities $\gamma$, $g_\infty$, $\varphi^{(0)}$, $\hat{\gamma}$ and $\hat{g}_\infty$ are given by equations (\ref{gamma}), (\ref{A2GINF}), (\ref{phi0}), (\ref{gammah}) and (\ref{AB2ghatINF}), respectively.
The pure imaginary quantity $\mathcal{A}(\infty)$ is given by equation (\ref{Abelinf}). These quantities depend only on the parameters in the initial data $q(x,0)$ and the self-similar variable $\xi$.

In addition, the square modulus of the leading-order term in (\ref{A2}) can be expressed in terms of the Jacobi elliptic $\rm{cn}$ function as:
\begin{equation}\label{A2-density}
  |q_{\mathrm{as}}(x,t)|^2= \rho_2 - (\rho_2- \rho_3) \mathrm{cn}^2 \left(\sqrt{\rho_1-\rho_3}\left(x- Vt + \varphi^{(1)}\left(\frac{x}{t}\right)\right) -K(m),m \right),
\end{equation}
where  $\rho_1, \rho_2$ and $\rho_3$ are given by equation (\ref{rho}) in Appendix \ref{CC}, and $\varphi^{(1)}$ and $V$ are given by equations (\ref{phi111}) and (\ref{VVV}), respectively. It is remarked that the leading-order term of the density in (\ref{A2-density}) is consistent with the one-phase modulated elliptic wave solution in (\ref{rho-equation-solution-new}) from the Whitham modulation theory.
\par
{\rm (iii)}. In the unmodulated elliptic wave region $(i.e., \xi_{A2}<\xi<\xi_{A3}$ with $\xi_{A2}$ and $\xi_{A3}$ defined by equation $(\ref{aB3b}))$, as $t \to \infty$,
 \begin{equation}\label{A3}
  q(x,t)=(A_l+A_r)\frac{\Theta (0)~\Theta (2\mathcal{A}(\infty)+(\gamma + \hat{\gamma})/{2 \pi})}{\Theta ((\gamma + \hat{\gamma})/{2 \pi})~\Theta (2\mathcal{A}(\infty))}e^{2\mathrm{i} (t g_\infty+\hat{g}_\infty-\varphi^{(0)})}+\mathcal{O}(t^{-1/2}),
\end{equation}
where the real quantities $\gamma$, $g_\infty$, $\hat{g}_\infty$, $\varphi^{(0)}$ and $\hat{\gamma}$ are given by equations (\ref{gamma1}), (\ref{A3GINF}), (\ref{aB3mdeltaINF}), (\ref{aB3mdeltaINF}) and (\ref{gammah1}), respectively, which depend only on the parameters in the initial data $q(x,0)$ and the self-similar variable $\xi$. The leading-order term in (\ref{A3}) is a one-phase periodic wave that can convert into a Jacobi elliptic wave.
\par
{\rm (iv)}. In the dispersive shock wave region $(i.e.,  \xi_{A3}<\xi<\xi_{A4}$ with $\xi_{A3}$ and $\xi_{A4}$ defined by equation $(\ref{AB4b}))$, as $t \to \infty$,
 \begin{equation}\label{A4}
  q(x,t)=\frac{\lambda_r^+ -\lambda_r^- + \lambda_s(\xi) -\lambda_l^-}{2}\frac{\Theta (0)~\Theta (2\mathcal{A}(\infty)+(\gamma + \hat{\gamma})/{2 \pi})}{\Theta ((\gamma + \hat{\gamma})/{2 \pi})~\Theta (2\mathcal{A}(\infty))} e^{2\mathrm{i} (t g_\infty+\hat{g}_\infty-\varphi^{(0)})}+\mathcal{O}(t^{-1}),
\end{equation}
where the real quantities $\gamma$, $g_\infty$, $\hat{g}_\infty$, $\varphi^{(0)}$, $\hat{\gamma}$ are given by equations (\ref{gamma2}), (\ref{A4GINF}), (\ref{AB4mdeltaINF}), (\ref{AB4mdeltaINF}) and (\ref{gammah2}), respectively.
Here $\lambda_{s}$ is a movable Riemann invariant (or soft edge) given by the unique solution of equation (\ref{vv3}). These quantities depend only on the parameters in the initial data $q(x,0)$ and the self-similar variable $\xi$.
\par
{\rm (iv)}. In the right plane wave region $(i.e., \xi > \xi_{A4}$ with $\xi_{A4}$ defined by equation $(\ref{A5B}))$, as $t \to \infty$,
\begin{equation}\label{A5}
    q(x,t)=A_r e^{-2\mathrm{i} \mu_r x-\mathrm{i}(2 \mu_r^2+A_r^2)t}e^{-\mathrm{i}\phi_{rp}(\xi) }+\mathcal{O}(t^{-\frac{1}{2}}),
\end{equation}
where the phase shift $\phi_{rp}(\xi)$ is given by equation $(\ref{A5phi1})$ and depends only on the parameters in the initial data $q(x,0)$ and the self-similar variable $\xi$.

\end{theorem}

\begin{theorem}
 \label{Rare-vacuum}
 {In Case {\rm B} $(\lambda_l^+>\lambda_l^->\lambda_r^+>\lambda_r^-)$, the long-time asymptotic solution $q(x,t)$ of the defocusing NLS equation $(\ref{(NLS)})$ with initial data (\ref{(initial)}) behaves differently in five regions depending on the range of the self-similar variable $\xi=x/t$, as shown in Figure \ref{XT-1} (b) and Figure \ref{XT-2} (b). As $t\to \infty$, the leading-order term and error term in each region are formulated specifically}:
\par
{\rm (i)}. In the left plane wave region $(i.e., \xi<\xi_{B1}$ with $\xi_{B1}$ defined by equation $(\ref{B1b}))$, the asymptotics of the solution $q(x,t)$ as $t \to \infty$ is the same as (\ref{A1}) up to the phase shift $\phi_{lp}(\xi)$ now given by $(\ref{philp})$.
\par
{\rm (ii)}. In the rarefaction wave region $(i.e., \xi_{B1} < \xi < \xi_{B2}$ with $\xi_{B1}$ and $\xi_{B2}$ defined by equation $(\ref{B2b}))$,
 as $t \to \infty$,
 \begin{equation}\label{B2}
    q(x,t)=-\frac{x+2\lambda_l^- t}{3t}e^{-\mathrm{i}(2{\lambda_l^-}^2t+2\lambda_l^- x-x^2/t)/3}e^{-\mathrm{i}\phi_{lr}(\xi) }+\mathcal{O}(t^{-1}),
 \end{equation}
 where the phase shift $\phi_{lr}(\xi)$ is given by equation $(\ref{Bphilr})$.
\par
{\rm (iii)}. In the vacuum region $(i.e., \xi_{B2} < \xi < \xi_{B3}$ with $\xi_{B2}$ and $\xi_{B3}$ defined by equation $(\ref{B3b}))$,
 as $t \to \infty$,
 \begin{equation}\label{B3}
   q(x,t)=\frac{\nu(-\xi)}{\sqrt{t}}e^{\mathrm{i}(t{\xi^2}/{2}- \nu(-\xi/2) \ln t + \phi_{va} (\xi) )} +o(t^{-\frac{1}{2}})=\mathcal{O}(t^{-\frac{1}{2}}),
 \end{equation}
 where $\nu(-\xi)$ and $\phi_{va} (\xi)$ are given by equations (\ref{B3nu}) and (\ref{phiva}), respectively.
\par
{\rm (iv)}. In the rarefaction wave region $(i.e., \xi_{B3} < \xi < \xi_{B4}$ with $\xi_{B3}$ and $\xi_{B4}$ defined by equation $(\ref{B4b}))$,
 as $t \to \infty$,
 \begin{equation}\label{B4}
   q(x,t)=\frac{x+2\lambda_r^+ t}{3t}e^{-\mathrm{i}(2{\lambda_r^+}^2t+2\lambda_r^+ x-x^2/t)/3}e^{-\mathrm{i}\phi_{rr}(\xi) }+\mathcal{O}(t^{-1}),
 \end{equation}
 where the phase shift $\phi_{rr}(\xi)$ is given by equation $(\ref{phirr})$.
\par
{\rm (v)}. In the left plane wave region $(i.e., \xi>\xi_{B4}$ with $\xi_{B4}$ defined by equation $(\ref{B5b}))$, the asymptotics of the solution $q(x,t)$ as $t \to \infty$ is the same as (\ref{A5}).

\end{theorem}

\begin{theorem}
  \label{Shock-two}
{In Case {\rm C} $(\lambda_r^+>\lambda_l^+ >\lambda_r^->\lambda_l^-)$, the long-time asymptotic solution $q(x,t)$ of the defocusing NLS equation $(\ref{(NLS)})$ with initial data (\ref{(initial)}) behaves differently in five regions depending on the range of the self-similar variable $\xi=x/t$, as shown in Figure \ref{XT-1} (c) and Figure \ref{XT-2} (c). As $t\to \infty$, the leading-order term and error term in each region are formulated specifically}:
\par
{\rm (i)}. In the left plane wave region $(i.e., \xi<\xi_{C1}$ with $\xi_{C1}$ defined by equation $(\ref{C1b}))$, the asymptotics of the solution $q(x,t)$ as $t \to \infty$ is the same as (\ref{A1}) up to the phase shift $\phi_{lp}(\xi)$ now given by $(\ref{philp})$.
\par
{\rm (ii)}. In the dispersive shock wave region $(i.e., \xi_{C1} < \xi<\xi_{C2}$ with $\xi_{C1}$ and $\xi_{C2}$ defined by equation $(\ref{C2b}))$, the asymptotics of the solution $q(x,t)$ as $t \to \infty$ is the same as (\ref{A2}) up to the phase shift $\varphi^{(0)}$ now given by $(\ref{C2phi0})$.
\par
{\rm (iii)}. In the middle plane wave region $(i.e., \xi_{C2} < \xi<\xi_{C3}$ with $\xi_{C2}$ and $\xi_{C3}$ defined by equation $(\ref{C3b}))$, as $t \to \infty$,
  \begin{equation}\label{C3}
    \begin{aligned}
      &q(x,t)=\frac{\lambda_r^+-\lambda_l^-}{2} e^{\mathrm{i}(kx- \omega t)}e^{-\mathrm{i}\phi_{mp} }+\mathcal{O}(e^{-ct}),\\
      &k=-(\lambda_r^++\lambda_l^-), \qquad \omega=\frac{(\lambda_r^++\lambda_l^-)^2}{2}+\frac{(\lambda_r^+ - \lambda_l^-)^2}{4},
    \end{aligned}
  \end{equation}
  for some constant $c>0$, where the phase shift $\phi_{mp}$ is given by equation $(\ref{phimp})$ and depends only on the parameters in the initial data $q(x,0)$.
\par
{\rm (iv)}. In the dispersive shock wave region $(i.e., \xi_{C3} < \xi<\xi_{C4}$ with $\xi_{C3}$ and $\xi_{C4}$ defined by equation $(\ref{C4b}))$, the asymptotics of the solution $q(x,t)$ as $t \to \infty$ is the same as (\ref{A4}).
\par
{\rm (v)}. In the right plane wave region $(i.e., \xi>\xi_{C4}$ with $\xi_{C4}$ defined by equation $(\ref{C5b}))$, the asymptotics of the solution $q(x,t)$ as $t \to \infty$ is the same as (\ref{A5}).

\end{theorem}

\begin{theorem}
  \label{Rare-two}
{In Case {\rm D} $(\lambda_l^+>\lambda_r^+>\lambda_l^->\lambda_r^-)$, the long-time asymptotic solution $q(x,t)$ of the defocusing NLS equation $(\ref{(NLS)})$ with initial data (\ref{(initial)}) is formulated according to the range of the self-similar variable $\xi=x/t$ as shown in Figure \ref{XT-1} (d) and Figure \ref{XT-2} (d).}
\par
{\rm (i)}. In the left plane wave region $(i.e., \xi<\xi_{D1}$ with $\xi_{D1}$ defined by equation $(\ref{D1b}))$, the asymptotics of the solution $q(x,t)$ as $t \to \infty$ is the same as (\ref{A1}) up to the phase shift $\phi_{lp}(\xi)$ now given by $(\ref{philp})$.
\par
{\rm (ii)}. In the rarefaction wave region $(i.e., \xi_{D1} < \xi<\xi_{D2}$ with $\xi_{D1}$ and $\xi_{D2}$ defined by equation $(\ref{D2b}))$, the asymptotics of the solution $q(x,t)$ as $t \to \infty$ is the same as (\ref{B2}) up to the phase shift $\phi_{lr}(\xi)$ now given by $(\ref{philrlr})$.
\par
{\rm (iii)}. In the middle plane wave region $(i.e., \xi_{D2} < \xi<\xi_{D3}$ with $\xi_{D2}$ and $\xi_{D3}$ defined by equation $(\ref{D3b}))$, the asymptotics of the solution $q(x,t)$ as $t \to \infty$ is the same as (\ref{C3}) up to the phase shift $\phi_{mp}$ now given by $(\ref{phimpmp})$.
\par
{\rm (iv)}. In the dispersive shock wave region $(i.e., \xi_{D3} < \xi<\xi_{D4}$ with $\xi_{D3}$ and $\xi_{D4}$ defined by equation $(\ref{D4b}))$, the asymptotics of the solution $q(x,t)$ as $t \to \infty$ is the same as (\ref{A4}).
\par
{\rm (v)}. In the right plane wave region $(i.e., \xi>\xi_{D4}$ with $\xi_{D4}$ defined by equation $(\ref{D5b}))$, the asymptotics of the solution $q(x,t)$ as $t \to \infty$ is the same as (\ref{A5}).

\end{theorem}

\begin{theorem}
  \label{Shock-Rare}
{In Case {\rm E} $(\lambda_r^+>\lambda_l^+> \lambda_l^->\lambda_r^-
)$, the long-time asymptotic solution $q(x,t)$ of the defocusing NLS equation $(\ref{(NLS)})$ with initial data (\ref{(initial)}) is formulated according to the range of the self-similar variable $\xi=x/t$ as shown in Figure \ref{XT-1} (e) and Figure \ref{XT-2} (e).}
\par
{\rm (i)}. In the left plane wave region $(i.e., \xi<\xi_{E1}$ with $\xi_{E1}$ defined by equation $(\ref{E1b}))$, the asymptotics of the solution $q(x,t)$ as $t \to \infty$ is the same as (\ref{A1}) up to the phase shift $\phi_{lp}(\xi)$ now given by $(\ref{philp})$.
\par
{\rm (ii)}. In the dispersive shock wave region $(i.e., \xi_{E1} < \xi<\xi_{E2}$ with $\xi_{E1}$ and $\xi_{E2}$ defined by equation $(\ref{E2b}))$, the asymptotics of the solution $q(x,t)$ as $t \to \infty$ is the same as (\ref{A2}) up to the phase shift $\varphi^{(0)}$ now given by $(\ref{C2phi0})$.
\par
{\rm (iii)}. In the middle plane wave region $(i.e., \xi_{E2} < \xi<\xi_{E3}$ with $\xi_{E2}$ and $\xi_{E3}$ defined by equation $(\ref{E3b}))$, the asymptotics of the solution $q(x,t)$ as $t \to \infty$ is the same as (\ref{C3}) up to the phase shift $\phi_{mp}$ now given by $(\ref{phimpmp})$.
\par
{\rm (iv)}. In the rarefaction wave region $(i.e., \xi_{E3} < \xi<\xi_{E4}$ with $\xi_{E3}$ and $\xi_{E4}$ defined by equation $(\ref{E4b}))$, the asymptotics of the solution $q(x,t)$ as $t \to \infty$ is the same as (\ref{B4}).
\par
{\rm (v)}. In the right plane wave region $(i.e., \xi>\xi_{E4}$ with $\xi_{E4}$ defined by equation $(\ref{E5b}))$, the asymptotics of the solution $q(x,t)$ as $t \to \infty$ is the same as (\ref{A5}).

\end{theorem}

\begin{theorem}
  \label{Rare-Shock}

{In Case {\rm F} $(\lambda_l^+>\lambda_r^+>\lambda_r^->\lambda_l^-)$, the long-time asymptotic solution $q(x,t)$ of the defocusing NLS equation $(\ref{(NLS)})$ with initial data (\ref{(initial)}) is formulated according to the range of the self-similar variable $\xi=x/t$ as shown in Figure \ref{XT-1} (f) and Figure \ref{XT-2} (f).}
\par
{\rm (i)}. In the left plane wave region $(i.e., \xi<\xi_{F1}$ with $\xi_{F1}$ defined by equation $(\ref{F1b}))$, the asymptotics of the solution $q(x,t)$ as $t \to \infty$ is the same as (\ref{A1}) up to the phase shift $\phi_{lp}(\xi)$ now given by $(\ref{philp})$.
\par
{\rm (ii)}. In the rarefaction wave region $(i.e., \xi_{F1} < \xi<\xi_{F2}$ with $\xi_{F1}$ and $\xi_{F2}$ defined by equation $(\ref{F2b}))$, the asymptotics of the solution $q(x,t)$ as $t \to \infty$ is the same as (\ref{B2}) up to the phase shift $\phi_{lr}(\xi)$ now given by $(\ref{philrlr})$.
\par
{\rm (iii)}. In the middle plane wave region $(i.e., \xi_{F2} < \xi<\xi_{F3}$ with $\xi_{F2}$ and $\xi_{F3}$ defined by equation $(\ref{F3b}))$, the asymptotics of the solution $q(x,t)$ as $t \to \infty$ is the same as (\ref{C3}) up to the phase shift $\phi_{mp}$ now given by $(\ref{phimpmp})$.
\par
{\rm (iv)}. In the dispersive shock wave region $(i.e., \xi_{F3} < \xi<\xi_{F4}$ with $\xi_{F3}$ and $\xi_{F4}$ defined by equation $(\ref{F4b}))$, the asymptotics of the solution $q(x,t)$ as $t \to \infty$ is the same as (\ref{A4}).
\par
{\rm (v)}. In the right plane wave region $(i.e., \xi>\xi_{F4}$ with $\xi_{F4}$ defined by equation $(\ref{F5b}))$, the asymptotics of the solution $q(x,t)$ as $t \to \infty$ is the same as (\ref{A5}).

\end{theorem}

\begin{remark}
  It is meaningful to consider the physically interesting symmetric initial data \cite{Buckingham2007}
  \begin{equation}\label{ini}
    q(x,0)=q_0(x):= A e^{-2\mathrm{i} \mu |x|}, \qquad x \in \mathbb{R},
  \end{equation}
  where $A$ and $\mu$ are nonzero real parameters with $A>0$, since it can be considered as the collision or separation (determined by the sign of $\mu$) of two plane waves corresponding to the initial data at the origin point for $t=0$.
  Case $\mu>0$ is called a shock problem, which corresponds to our Case A or Case C depending on the relationship between parameters $A$ and $\mu$.
  Case $\mu<0$ is called a rarefaction problem, which corresponds to our  Case B or Case D.
  In Figure \ref{XT-2} (a)-(d), it can be intuitively seen that the asymptotic behaviors of the solutions inherit the symmetry of the initial data (\ref{ini}).
\end{remark}

\section{The proofs of the main theorems}

This section provides the proof of Theorems \ref{Shock-un}-\ref{Rare-Shock} by adopting the Deift-Zhou nonlinear steepest descent technique \cite{Deift-Zhou1993} for the Riemann-Hilbert problem \ref{RHP-rhp-0} given in Appendix \ref{AAC}. Initially, we briefly describe how to analyze the inverse problems through a series of equivalent transformations. For the six cases listed in (\ref{Classification}), the procedures are similar up to the introduced $g$-functions and the deformations of the contours.
\par
First, take the so-called ``$g$-function mechanism'' \cite{Deift-Zhou-vanakides1994} to renormalize the oscillatory or exponentially large entries of the jump matrices. In
different regions, the appropriate $g$-functions are selected to carry out asymptotic analysis with common properties:
\begin{enumerate} [label=(\roman*)]
  \item  $g(z)$ is analytic in $z \in \mathbb{C} \backslash J $, where $J$  is the union of several bands (closed intervals),
  \item  the asymptotics: $g(z)=\theta(z)+g_\infty+ \mathcal{O}(z^{-1})$ as $z\to \infty$,
  \item  the Schwarz symmetry: $g(z)=g^*(z) $.
\end{enumerate}
Then the first transformation $M(z;x,t) \mapsto M^{(1)}(z;x,t)$ is defined by
\begin{equation}\label{g-fun}
  M^{(1)}(z;x,t)=e^{-\mathrm{i} tg_\infty \sigma_3 } M(z;x,t) e^{-\mathrm{i}t(\theta(z) -g(z))\sigma_3}.
\end{equation}
\par
Second, deform the contours to the steepest descent contours $L_i$ $(i=1, 2 , 3, 4),$ on which the jumps decay to identity for a large time, which is called ``opening lenses''. The contours and the real axis will divide the complex plane into six regions $D_i$ $(i=1,2, \cdots, 6)$ (for example, see Figure \ref{figA1}(c)), where the transformations to ``open lenses'' will be defined below. Here we need the following  upper/lower and lower/upper matrix factorizations on specific intervals.
These matrix factorizations can be directly verified by using (\ref{rraa1}) and (\ref{rraa2}).
For $z\in \mathbb{R}\backslash (\overline{\mathcal{I}_l \cup \mathcal{I}_r} ) $, it follows
\begin{subequations}\label{f1}
  \begin{align}
    \begin{pmatrix} 1-r(z)r^*(z) & -r^*(z) \\ r(z) & 1  \end{pmatrix} &= \begin{pmatrix} 1 & -r^*(z) \\ 0 & 1  \end{pmatrix}\begin{pmatrix} 1 & 0 \\ r(z) & 1  \end{pmatrix}  \\
    &=\begin{pmatrix} 1 & 0 \\ \frac{r(z)}{1-r(z)r^*(z)} & 1  \end{pmatrix} (1-r(z)r^*(z))^{\sigma_3} \begin{pmatrix} 1 & \frac{-r^*(z)}{1-r(z)r^*(z)} \\ 0 & 1  \end{pmatrix} .
  \end{align}
\end{subequations}
For $z\in \mathcal{I}_l \backslash \overline{\mathcal{I}_r}$, one has
\begin{subequations}\label{f2}
  \begin{align}
    \begin{pmatrix} 0  & -r^*_-(z) \\ r_+(z) & 1  \end{pmatrix} &= \begin{pmatrix} 1 & -r_-^*(z) \\ 0 & 1  \end{pmatrix}\begin{pmatrix} 1 & 0 \\ r_+(z) & 1  \end{pmatrix}  \\
    &=\begin{pmatrix} 1 & 0 \\ \frac{r_-(z)}{1-r_-(z)r_-^*(z)} & 1  \end{pmatrix}\begin{pmatrix} 0 & -r^*_-(z) \\ r_+(z) & 0  \end{pmatrix}\begin{pmatrix} 1 & \frac{-r_+^*(z)}{1-r_+(z)r_+^*(z)} \\ 0 & 1  \end{pmatrix}.
  \end{align}
\end{subequations}
For $z\in \mathcal{I}_r \backslash \overline{\mathcal{I}_l} $, one has
\begin{subequations}\label{f3}
  \begin{align}
    \begin{pmatrix} \frac{1}{a_+(z)a_-^*(z)}  & -1 \\ 1 & 0  \end{pmatrix} &= \begin{pmatrix} 1 & -r_-^*(z) \\ 0 & 1  \end{pmatrix} \begin{pmatrix} 0 & -1 \\ 1 & 0  \end{pmatrix} \begin{pmatrix} 1 & 0 \\ r_+(z) & 1  \end{pmatrix}  \\
    &=\begin{pmatrix} 1 & 0 \\ \frac{r_-(z)}{1-r_-(z)r_-^*(z)} & 1  \end{pmatrix} (a_+(z)a_-^*(z))^{-\sigma_3} \begin{pmatrix} 1 & \frac{-r_+^*(z)}{1-r_+(z)r_+^*(z)} \\ 0 & 1  \end{pmatrix}.
  \end{align}
\end{subequations}
For $z \in \mathcal{I}_l \cap \mathcal{I}_r$, one has
\begin{subequations}\label{f4}
  \begin{align}
    \begin{pmatrix} 0  & -1 \\ 1 & 0  \end{pmatrix} &= \begin{pmatrix} 1 & -r^*(z) \\ 0 & 1  \end{pmatrix}\begin{pmatrix} 0 & -1 \\ 1 & 0  \end{pmatrix}\begin{pmatrix} 1 & 0 \\ r(z) & 1  \end{pmatrix}  \\
    &=\begin{pmatrix} 1 & 0 \\ \frac{r(z)}{1-r(z)r^*(z)} & 1  \end{pmatrix} \begin{pmatrix} 0 & -1 \\ 1 & 0  \end{pmatrix} \begin{pmatrix} 1 & \frac{-r^*(z)}{1-r(z)r^*(z)} \\ 0 & 1  \end{pmatrix}.
  \end{align}
\end{subequations}
These matrix factorizations imply the following piecewise-analytic transformation $M^{(1)}(z;x,t) \mapsto M^{(2)}(z;x,t)$ defined by
\begin{equation}\label{open}
  M^{(2)}(z;x,t)=\left\{  \begin{array}{ll}
    M^{(1)}(z;x,t) \begin{pmatrix} 1 & 0 \\ -r(z)e^{2\mathrm{i} tg(z)} & 1\end{pmatrix}, & z \in D_1, \\ \\
    M^{(1)}(z;x,t), &z\in D_2 \cup D_5,\\
    M^{(1)}(z;x,t) \begin{pmatrix} 1 & \frac{r^*(z)}{1-r(z)r^*(z)}e^{-2\mathrm{i}tg(z)} \\ 0 & 1 \end{pmatrix}, & z \in D_3, \\ \\
    M^{(1)}(z;x,t) \begin{pmatrix} 1 & 0 \\ \frac{r(z)}{1-r(z)r^*(z)}e^{2\mathrm{i}tg(z)} & 1 \end{pmatrix}, & z \in D_4, \\ \\
    M^{(1)}(z;x,t) \begin{pmatrix} 1 & -r^*(z)e^{-2\mathrm{i}tg(z)} \\ 0 & 1 \end{pmatrix}, &z \in D_6. \end{array} \right.
\end{equation}
\par
Third, introduce a scalar undetermined function $\delta(z)$ and the third transformation $M^{(2)}(z;x,t) \mapsto M^{(3)}(z;x,t)$ defined by
\begin{equation} \label{2to3}
	M^{(3)}(z;x,t)=\delta _{\infty}^{\sigma_{3}} M^{(2)}(z;x,t) \delta(z) ^{-\sigma_{3}}
\end{equation}
to reduce the jump matrices across the real axis to matrices independent of $z$, which leads to a solvable limiting problem.
Then determine a global parametrix $M^{\mathrm{(par)}}(z;x,t)$ consisting of a local parametrix $M^{\mathrm{(loc)}}(z;x,t)$ near the stationary phase point
and an outer model $M^{\mathrm{(mod)}}(z;x,t)$ that solves the limiting problem.
\par
Finally, estimate the error term by considering the error matrix defined by
\begin{equation}
	M^{\mathrm{(err)}}(z;x,t)= M^{(3)}(z;x,t) {M^{\mathrm{(par)}}(z;x,t)}^{-1},
\end{equation}
which corresponds to the small-norm RHPs whose solutions can be expressed by Cauchy singular integrals. Inverting the above
transformations and using the reconstruction formula (\ref{recf}),  the long-time asymptotics of $q(x,t)$ with a leading-order term and an error estimate in each region is finally obtained as
\begin{equation}\label{rec}
   \begin{aligned}
      q(x,t) &= 2 \mathrm{i} \lim_{z \to \infty, z \in D_{2}} (z M(z;x,t))_{12}\\
      &=2 \mathrm{i} \lim_{z \to \infty, z \in D_{2}}(z e^{\mathrm{i} g_\infty \sigma_3 }\delta _{\infty}^{-\sigma_{3}}
       M^{(\mathrm{err})}(z;x,t) M^{(\mathrm{mod})}(z;x,t) \delta(z) ^{\sigma_{3}}e^{\mathrm{i}t(\theta(z) - g(z))\sigma_3})_{12}\\
      &=2 \mathrm{i} e^{2\mathrm{i} t g_\infty}\delta _{\infty}^{-2} (M^{(\mathrm{mod})}_1(x,t)+ M^{(\mathrm{err})}_1(x,t) )_{12}.
   \end{aligned}
\end{equation}
\par
In what follows, we repeat the above procedure to deform the RHPs in each case and eventually prove Theorems \ref{Shock-un}-\ref{Rare-Shock}.

\subsection{\rm Case A: $\lambda_r^+>\lambda_r^->\lambda_l^+>\lambda_l^-$}
\
\newline
\indent
In this case, the two intervals $\mathcal{I}_l$ and $\mathcal{I}_r$  do not intersect, so we only consider the factorizations (\ref{f1})-(\ref{f3}).
As the self-similar variable $\xi=x/t$ increases, the stationary phase points of the corresponding $g$-functions change continuously at different intervals on $\mathbb{R}$,
which implies that there are five different regions: the left plane wave region, dispersive shock wave region, unmodulated elliptic wave region, dispersive shock wave region and the right plane wave region.
The fact that two plane waves corresponding to the initial data collide at the origin point for $t=0$ is consistent with the appearance of elliptic wave regions.

\subsubsection{The left plane wave region: $\xi<-\frac{2\lambda_r^++\lambda_l^++\lambda_l^-}{2}+\frac{2(\lambda_r^+ - \lambda_l^+)
(\lambda_r^+ - \lambda_l^-)}{-2\lambda_r^+ + \lambda_l^+ + \lambda_l^-}$} \label{Alpw}
\
\newline

\begin{figure}[htbp]
  \centering
  \subfigure[]{\includegraphics[width=5.5cm]{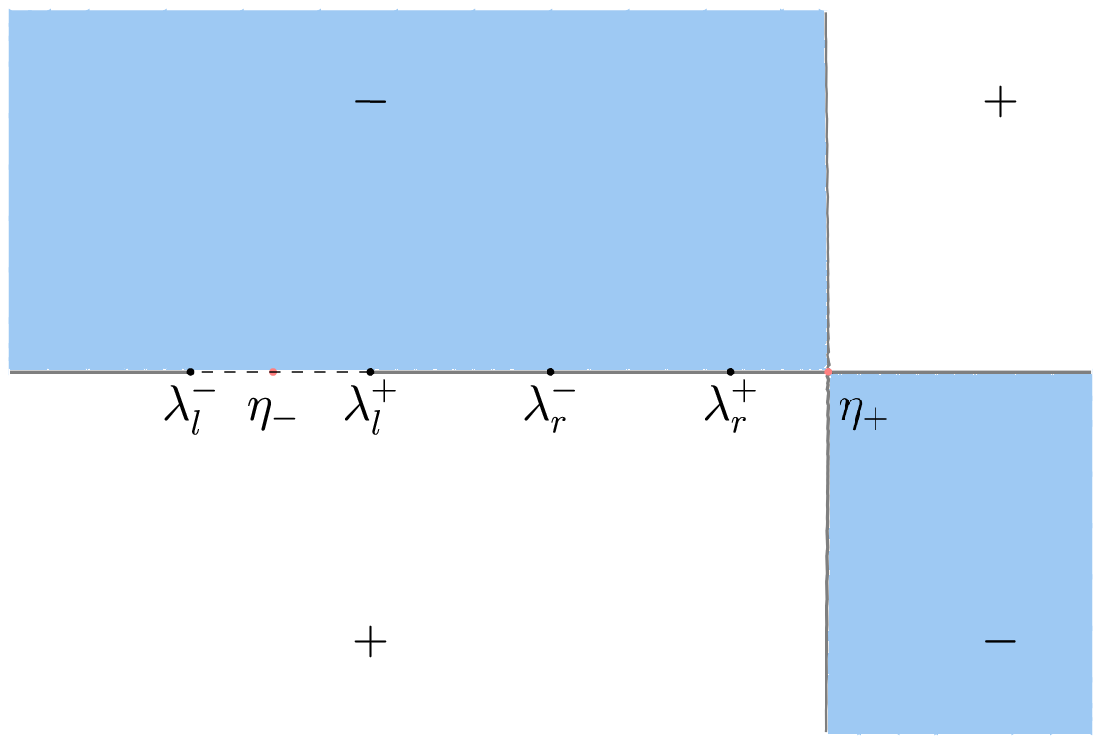}}
  \subfigure[]{\includegraphics[width=5cm]{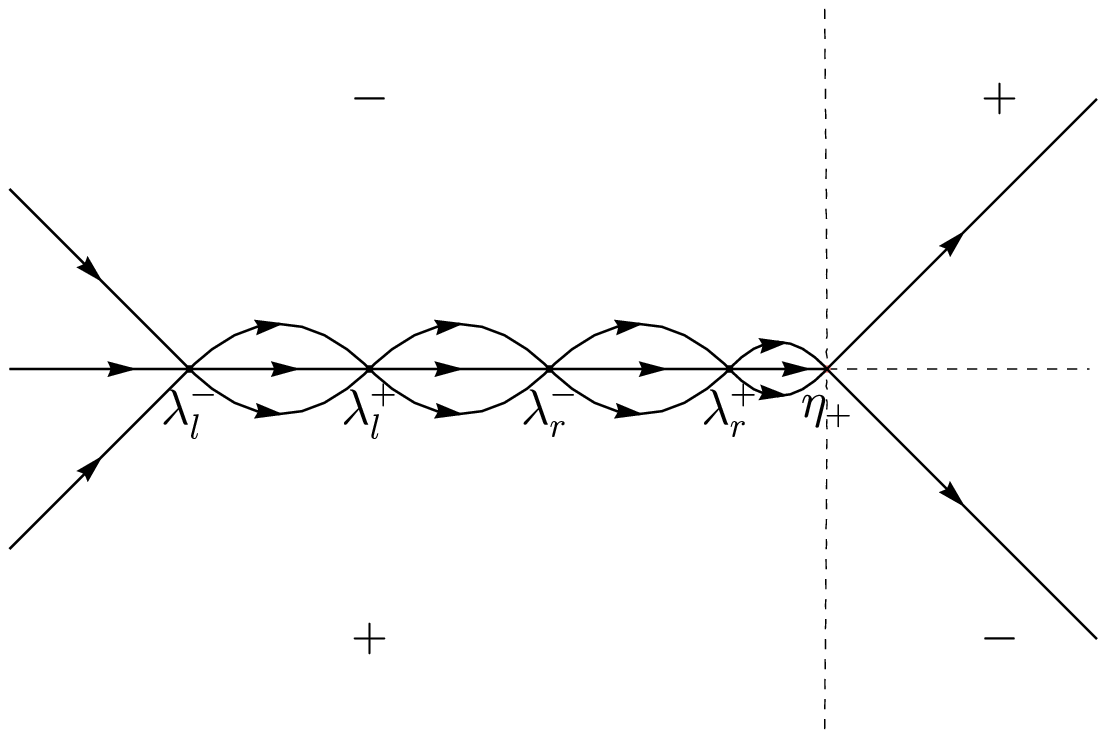}}
  \subfigure[]{\includegraphics[width=5.5cm]{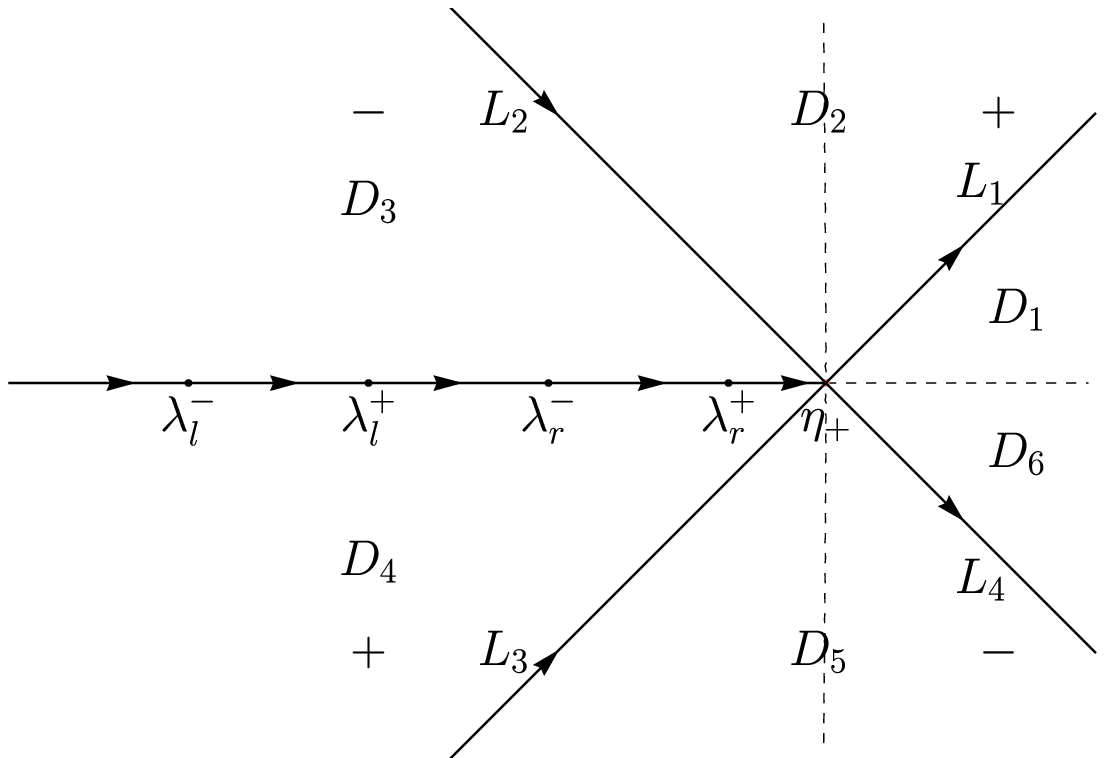}}
\caption{{\protect\small (a) Sign structure of Im$(g(z))$; (b) Opening lenses; (c) The jump contours of $M^{(2)}(z;x,t)$.}}
\label{figA1}
\end{figure}

In this region,  the phase function of the explicit eigenfunction $\Psi_l^{\mathrm{p}}(z;x,t)$ defined by (\ref{psip}) is the appropriate $g$-function
\begin{equation} \label{A1g}
  g(z)= \mathcal{R} (z;\lambda_l^+,\lambda_l^-)(z+\frac{1}{2}(\lambda_l^++\lambda_l^-)+\xi),
\end{equation}
which is consistent with the genus-zero case in which the Riemann invariants $\lambda_l^+$ and $\lambda_l^-$ are two hard edges in Whitham modulation theory.
The function $g(z)$ is analytic for $z \in \mathbb{C} \backslash \overline{\mathcal{I}_l}$ with
\begin{equation}\label{A1ginf}
   g_{\infty}=- \mu_l \xi-(2 \mu_l^2+A_l^2)/2.
\end{equation}
The signature table for $\mathrm{Im}(g(z))$ (see Figure \ref{figA1}(a)) is determined by the differential
\begin{equation} \label{A1dg}
  \dif g(z)= 2 \frac{(z-\eta_+(\xi))(z-\eta_-(\xi))}{\mathcal{R} (z;\lambda_l^+,\lambda_l^-)} \,\dif z,
\end{equation}
where $\eta_{\pm}(\xi)$ are the stationary phase points given by
\begin{equation}\label{lpweta}
  \eta_{\pm}(\xi)=\frac{\lambda_l^++\lambda_l^--\xi \pm \sqrt{(\lambda_l^++\lambda_l^-+\xi)^2+2(\lambda_l^+-\lambda_l^-)^2}}{4}.
\end{equation}
The boundary of the left plane wave region is characterized by
\begin{equation}\label{A1B}
  \eta_+ (\xi) > \lambda_r^+ \qquad \mathrm{iff}  \qquad  \xi< \xi_{A1} = -\frac{2\lambda_r^++\lambda_l^++\lambda_l^-}{2}+\frac{2(\lambda_r^+ - \lambda_l^+)
  (\lambda_r^+ - \lambda_l^-)}{-2\lambda_r^+ + \lambda_l^+ + \lambda_l^-}.
\end{equation}

Then we can open lenses from the real axis to the steepest descent contours through $\eta_+(\xi)$ (as shown in Figure \ref{figA1} (b)-(c)) by the transformations (\ref{g-fun}) and (\ref{open}),
which yields the following RHP.

\begin{rhp}
  Find a $2 \times 2$ matrix-valued function $M^{(2)}(z;x,t)$ with the following properties:
\begin{enumerate} [label=(\roman*)]
  \item  $M^{(2)}(z;x,t)$ is analytic in $z \in \mathbb{C} \backslash \Sigma^{(2)}$, where $\Sigma^{(2)}=(-\infty, \eta_{+}(\xi)] \cup L_1 \cup L_2 \cup  L_3 \cup L_4.$
  \item  $M^{(2)}(z;x,t)=I+\mathcal{O} (z^{-1}) $ as $z\to \infty$.
  \item  $M^{(2)}(z;x,t)$ achieves the continuous boundary values (CBVs) $M^{(2)}_{+}(z;x,t)$ and $M^{(2)}_{-}(z;x,t)$  on $\Sigma^{(2)}$ away from self-intersection points and
  branch points that satisfy the jump condition $M^{(2)}_{+}(z;x,t)=M^{(2)}_{-}(z;x,t)V^{(2)}(z;x,t)$, where
  \begin{equation}\label{A1V2}
    \begin{aligned}
    V^{(2)}(z;x,t)&=\left\{  \begin{array}{ll}
     (1-r(z)r^{*}(z))^{\sigma_{3}},  ~& z\in (-\infty, \lambda_l^-) \cup (\lambda_l^+, \lambda_r^-) \cup (\lambda_r^+, \eta_+(\xi)), \\ \\
     (a_+(z) a_-^*(z))^{-\sigma_{3}},  ~&z \in  (\lambda_r^-, \lambda_r^+), \\ \\
     \begin{pmatrix} 0 & -r_{-}^{*}(z) \\ r_{+}(z) & 0 \end{pmatrix}, ~&z \in (\lambda_l^-, \lambda_l^+), \\ \\
     \begin{pmatrix} 1 & 0 \\ r(z)e^{2\mathrm{i} tg(z)} & 1 \end{pmatrix}, ~ &z \in L_1, \\ \\
     \begin{pmatrix} 1 & \frac{-r^*(z)}{1-r(z)r^*(z)}e^{-2\mathrm{i} tg(z)} \\0 & 1 \end{pmatrix},  ~&z \in L_2,\\ \\
     \begin{pmatrix} 1 & 0 \\ \frac{r(z)}{1-r(z)r^*(z)}e^{2\mathrm{i} tg(z)} & 1 \end{pmatrix},  ~&z \in L_3,\\ \\
     \begin{pmatrix} 1 & -r^*(z)e^{-2\mathrm{i} tg(z)} \\ 0  & 1 \end{pmatrix}, ~ &z \in L_4.\end{array} \right.\\
    \end{aligned}
  \end{equation}
\end{enumerate}
\end{rhp}

To reduce jump matrices across the real axis to constant matrices (especially the identity matrix), a scalar function $\delta(z)$ should be introduced,
which is analytic in $z \in \mathbb{C} \backslash(-\infty, \eta_{+}(\xi)]$ and satisfies the following jump conditions:
\begin{equation*}
  \begin{aligned}
\left\{  \begin{array}{ll}
   \delta _{+}(z)=(1-r(z)r^{*}(z))\delta _{-}(z),  ~& z\in (-\infty, \lambda_l^-) \cup (\lambda_l^+, \lambda_r^-) \cup (\lambda_r^+, \eta_+(\xi)), \\ \\
   \delta _{+}(z)=(a_+(z) a_-^*(z))^{-1}\delta _{-}(z),  ~& z\in (\lambda_r^-, \lambda_r^+), \\ \\
   \delta _{+}(z)\delta _{-}(z)=r_{+}(z),  ~& z\in (\lambda_l^-, \lambda_l^+).
  \end{array} \right.\\
  \end{aligned}
\end{equation*}
By the Plemelj formulae \cite{Buckingham2007}, we have
\begin{equation} \label{A1delta}
  \begin{aligned}
  \delta(z) = \exp    &  \left\{ \frac{\mathcal{R} (z;\lambda_l^+,\lambda_l^-)}{2 \pi \mathrm{i}} \left[  \left(\int_{-\infty}^{\lambda_l^-} + \int_{\lambda_l^+}^{\lambda_r^-}+
  \int_{\lambda_r^+}^{\eta_+(\xi)}\right) \frac{\ln (1-|r(\zeta)|^2)}{\mathcal{R} (\zeta;\lambda_l^+,\lambda_l^-) } \,\frac{\dif\zeta}{\zeta - z} \right. \right.   \\
  & \left.  \left. + \int_{\lambda_r^-}^{\lambda_r^+} \frac{-\ln(a_+(\zeta)a_-^*(\zeta))}{\mathcal{R} (\zeta;\lambda_l^+,\lambda_l^-)} \,\frac{\dif\zeta}{\zeta - z}
  + \int_{\lambda_l^-}^{\lambda_l^+} \frac{\ln (r_+(\zeta))}{\mathcal{R}_+ (\zeta;\lambda_l^+,\lambda_l^-)} \, \frac{\dif\zeta}{\zeta - z} \right]  \right\}
  \end{aligned}
\end{equation}
with large $z$ asymptotic behavior $\delta(z)=\delta_{\infty}+\mathcal{O} (z^{-1}) $ as $z \to \infty$, where
\begin{equation}\label{A1deltaINF }
  \begin{aligned}
  \delta_{\infty }= \exp   & \left\{ \frac{ \mathrm{i}}{2\pi} \left[ \left(\int_{-\infty}^{\lambda_l^-} + \int_{\lambda_l^+}^{\lambda_r^-}+
  \int_{\lambda_r^+}^{\eta_+(\xi)}\right)  \frac{\ln (1-|r(\zeta)|^2)}{\mathcal{R} (\zeta;\lambda_l^+,\lambda_l^-) } \dif \zeta \right. \right. \\
  & \left.  \left. + \int_{\lambda_r^-}^{\lambda_r^+} \frac{-\ln(a_+(\zeta)a_-^*(\zeta))}{\mathcal{R} (\zeta;\lambda_l^+,\lambda_l^-)} \dif \zeta
  + \int_{\lambda_l^-}^{\lambda_l^+} \frac{\ln (r_+(\zeta))}{\mathcal{R}_+ (\zeta;\lambda_l^+,\lambda_l^-)} \dif \zeta \right]  \right\}.
  \end{aligned}
\end{equation}
This results in the following transformation
\begin{equation}
	M^{(3)}(z;x,t)=\delta _{\infty}^{\sigma_{3}} M^{(2)}(z;x,t) \delta(z) ^{-\sigma_{3}}
\end{equation}
and the following new RHP.
\begin{rhp} Find a $2 \times 2$ matrix-valued function $M^{(3)}(z;x,t)$ with the following properties:
\begin{enumerate} [label=(\roman*)]
  \item $M^{(3)}(z;x,t)$ is analytic in $z \in \mathbb{C} \backslash \Sigma^{(3)}$, where $\Sigma^{(3)}=[\lambda_l^-, \lambda_l^+] \cup L_1 \cup L_2 \cup  L_3 \cup L_4.$
  \item $M^{(3)}(z;x,t)=I+\mathcal{O} (z^{-1}) $ as $z\to \infty$.
  \item $M^{(3)}(z;x,t)$ achieves the CBVs $M^{(3)}_{+}(z;x,t)$ and $M^{(3)}_{-}(z;x,t)$ on $\Sigma^{(3)}$ away from self-intersection points and
branch points that satisfy the jump condition $M^{(3)}_{+}(z;x,t)=M^{(3)}_{-}(z;x,t)V^{(3)}(z;x,t)$, where
\begin{equation}\label{A1V3}
  \begin{aligned}
V^{(3)}(z;x,t)&=\left\{  \begin{array}{ll}
   \begin{pmatrix} 0 & -1 \\ 1 & 0 \end{pmatrix}, ~&z \in (\lambda_l^-, \lambda_l^+),\\ \\
   \begin{pmatrix} 1 & 0 \\ r(z) \delta ^{-2}(z) e^{2\mathrm{i} tg(z)} & 1 \end{pmatrix}, ~ &z \in L_1, \\ \\
   \begin{pmatrix} 1 & \frac{-r^*(z)}{1-r(z)r^*(z)} \delta ^2(z) e^{-2\mathrm{i}t g(z)} \\0 & 1 \end{pmatrix},  ~&z \in L_2,\\ \\
   \begin{pmatrix} 1 & 0 \\ \frac{r(z)}{1-r(z)r^*(z)} \delta ^{-2}(z) e^{2\mathrm{i} tg(z)} & 1 \end{pmatrix},  ~&z \in L_3,\\ \\
   \begin{pmatrix} 1 & -r^*(z) \delta ^2(z) e^{-2\mathrm{i} tg(z)} \\ 0  & 1 \end{pmatrix}, ~ &z \in L_4. \end{array} \right.\\
  \end{aligned}
\end{equation}
\end{enumerate}
\end{rhp}

Outside a small  neighborhood $\mathcal{U} $ of $\eta_+(\xi)$, the jump matrices $V^{(3)}(z;x,t)$ on the steepest descent contours uniformly converge
to the identity matrix. Therefore, consider the limiting problem
\begin{equation}
  \left\{  \begin{array}{ll}
  M^{(\mathrm{mod})}_+(z;x,t)=M^{(\mathrm{mod})}_-(z;x,t) V^{(\mathrm{mod})}, \quad  z \in (\lambda_l^-, \lambda_l^+), \\\\
  V^{(\mathrm{mod})}=\begin{pmatrix} 0 & -1 \\ 1 & 0 \end{pmatrix},
\end{array} \right.
\end{equation}
whose solution is exactly given by
\begin{equation} \label{A1mod}
  M^{(\mathrm{mod})}(z;x,t) = \mathcal{E}_l (z),
\end{equation}
where $\mathcal{E}_l (z)$ is given by (\ref{EPsilon}) in Appendix \ref{AAC}.
\par
\begin{figure}[htbp]
  \centering
  \subfigure[]{\includegraphics[width=7cm]{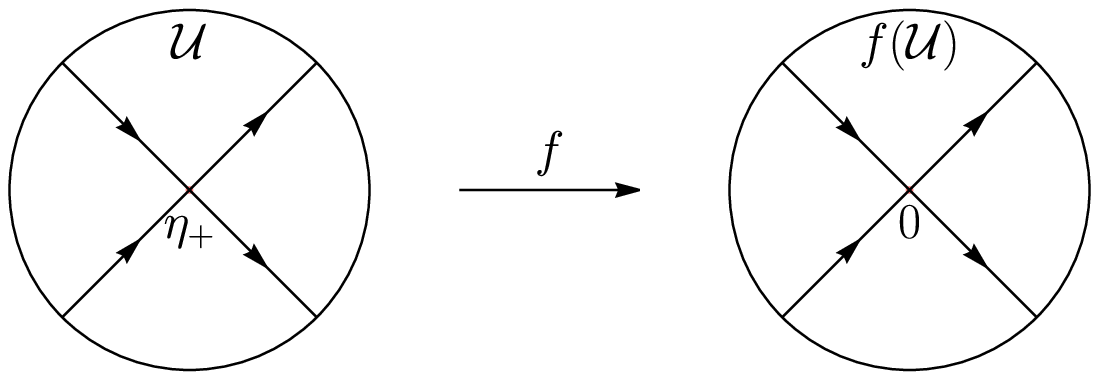}} \qquad
  \subfigure[]{\includegraphics[width=7cm]{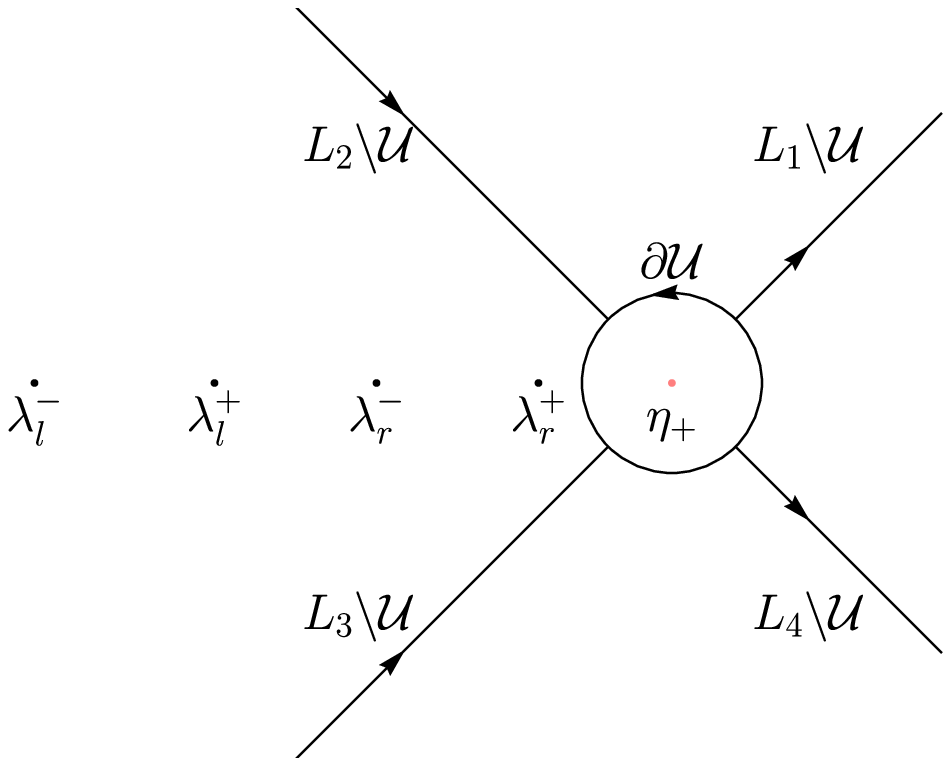}}
\caption{{\protect\small (a) The conformal mapping $f$ on $\mathcal{U}$; (b) The jump contours of error RHP.}}
\label{figA11}
\end{figure}

Inside  $\mathcal{U}$, the jump matrices $V^{(3)}(z;x,t)$ cannot  uniformly converge  to the identity matrix due to the quadratic vanishment
of the phase $g(z)$. Thus, we construct a local parametrix $M^{\mathrm{(loc)}}(z;x,t)$ that matches the same jumps as $M^{(3)}(z;x,t)$
inside $\mathcal{U}$ with the help of the parabolic cylinder model $M^{(\mathrm{PC})}(k;r)$ defined in Appendix \ref{pc} to determine a global parametrix $M^{(\mathrm{par})}(z;x,t)$. The construction is standard and we only give a brief overview.
Fix the small neighborhood $\mathcal{U}$ and define the mapping
\begin{equation}
  k=\sqrt{t} f(z)=2\sqrt{t} (g(z)-g(\eta_+(\xi)))^{1/2},
\end{equation}
where the branch is chosen such that $f'(\eta_+(\xi))>0$, which is conformal from $\mathcal{U}$ (in $z$) to a neighborhood of the origin (in $k$), as shown in Figure \ref{figA11} (a).
Define the local parametrix by
\begin{equation}
  M^{\mathrm{(loc)}}(z;x,t)=M^{(\mathrm{mod})}(z;x,t)h(z)^{-\sigma_3}M^{(\mathrm{PC})}(\sqrt{t} f(z);r(z))h(z)^{\sigma_3},
\end{equation}
where
\begin{equation}
  h(z)=e^{\mathrm{i}g(\eta_+(\xi))}(\sqrt{t} f(z))^{\mathrm{i}\nu(z)}\delta^{-1}(z)
\end{equation}
is principally branched. It is easy to verify that this local parametrix satisfies our requirements.

To consider the error estimate, define the error matrix $M^{\mathrm{(err)}}(z;x,t)$ by
\begin{equation}\label{A1para}
  \begin{aligned}
  M^{(3)}(z;x,t) =M^{(\mathrm{err})}(z;x,t)M^{(\mathrm{par})}(z;x,t)
  &=\left\{  \begin{array}{ll}M^{(\mathrm{err})}(z;x,t) M^{\mathrm{(loc)}}(z;x,t), ~ & \mathrm{inside} ~ \mathcal{U}, \\  \\
    M^{(\mathrm{err})}(z;x,t) M^{(\mathrm{mod})}(z;x,t), ~& \mathrm{outside} ~ \mathcal{U},
   \end{array} \right.\\
  \end{aligned}
\end{equation}
which leads to the following RHP.
\begin{rhp}\label{eRHP} Find a $2 \times 2$ matrix-valued function $M^{(\mathrm{err})}(z;x,t)$ with the following properties:
  \begin{enumerate} [label=(\roman*)]
    \item $M^{(\mathrm{err})}(z;x,t)$ is analytic in $z \in \mathbb{C} \backslash \Sigma^{(\mathrm{err})}$, where $\Sigma^{(\mathrm{err})}=((L_1 \cup L_2 \cup  L_3 \cup L_4) \backslash\mathcal{U})\cup \partial \mathcal{U}$.
    \item $M^{(\mathrm{err})}(z;x,t)=I+\mathcal{O} (z^{-1}) $ as $z\to \infty$.
    \item $M^{(\mathrm{err})}(z;x,t)$ achieves the CBVs $M^{(\mathrm{err})}_{+}(z;x,t)$ and $M^{(\mathrm{err})}_{-}(z;x,t)$ on $\Sigma^{(\mathrm{err})}$, which satisfy the jump condition $M^{(\mathrm{err})}_{+}(z;x,t)=M^{(\mathrm{err})}_{-}(z;x,t)V^{(\mathrm{err})}(z;x,t)$, where
  \begin{equation}\label{A1Verr}
    \begin{aligned}
  V^{(\mathrm{err})}(z;x,t)&=\left\{  \begin{array}{ll}
    M^{(\mathrm{mod})}(z;x,t) \begin{pmatrix} 1 & 0 \\ r(z) \delta ^{-2}(z) e^{2\mathrm{i} tg(z)} & 1 \end{pmatrix} {M^{(\mathrm{mod})}(z;x,t)}^{-1}, ~ &z \in L_1\backslash\mathcal{U}, \\ \\
    M^{(\mathrm{mod})}(z;x,t) \begin{pmatrix} 1 & \frac{-r^*(z)}{1-r(z)r^*(z)} \delta ^2(z) e^{-2\mathrm{i}t g(z)} \\0 & 1 \end{pmatrix} {M^{(\mathrm{mod})}(z;x,t)}^{-1},  ~&z \in L_2\backslash\mathcal{U},\\ \\
    M^{(\mathrm{mod})}(z;x,t) \begin{pmatrix} 1 & 0 \\ \frac{r(z)}{1-r(z)r^*(z)} \delta ^{-2}(z) e^{2\mathrm{i} tg(z)} & 1 \end{pmatrix} {M^{(\mathrm{mod})}(z;x,t)}^{-1},  ~&z \in L_3\backslash\mathcal{U},\\ \\
    M^{(\mathrm{mod})}(z;x,t) \begin{pmatrix} 1 & -r^*(z) \delta ^2(z) e^{-2\mathrm{i} tg(z)} \\ 0  & 1 \end{pmatrix} {M^{(\mathrm{mod})}(z;x,t)}^{-1}, ~ &z \in L_4\backslash\mathcal{U},\\ \\
    M^{(\mathrm{mod})}(z;x,t) {M^{(\mathrm{loc})}(z;x,t)}^{-1}, ~ &z \in \partial \mathcal{U}.\end{array} \right.
    \end{aligned}
  \end{equation}
  \end{enumerate}
\end{rhp}
In Appendix \ref{appb} it is shown that the error estimate is $M^{(\mathrm{err})}_1(x,t)=\mathcal{O}(t^{-1/2})$ as $t\to \infty$.
\par
Using equations (\ref{rec}), (\ref{A1ginf}), (\ref{A1deltaINF }) and (\ref{A1mod}), it is concluded that the long-time asymptotics of the defocusing NLS equation (\ref{(NLS)})
in the left plane wave region is given by
\begin{equation}
  \begin{aligned}
   q(x,t)&=A_l e^{-2\mathrm{i} \mu_l x-\mathrm{i}(2 \mu_l^2+A_l^2)t}e^{-\mathrm{i}\phi_{lp}(\xi) }+\mathcal{O}(t^{-\frac{1}{2}}),
  \end{aligned}
\end{equation}
where
\begin{equation}\label{Aphi1}
  \begin{aligned}
  \phi_{lp}(\xi)=\frac{1}{\pi} & \left[ \left(\int_{-\infty}^{\lambda_l^-} +\int_{\lambda_l^+}^{\lambda_r^-}+
   \int_{\lambda_r^+}^{\eta_+(\xi)}\right) \frac{\ln (1-|r(\zeta)|^2)}{\sqrt{(\zeta-\lambda_l^+)(\zeta-\lambda_l^-)} } \dif \zeta \right.\\
   & \left. +\int_{\lambda_r^-}^{\lambda_r^+} \frac{-\ln(a_+(\zeta)a_-^*(\zeta))}{\sqrt{(\zeta-\lambda_l^+)(\zeta-\lambda_l^-)}} \dif \zeta
   + \int_{\lambda_l^-}^{\lambda_l^+} \frac{\arg (r_+(\zeta))}{\sqrt{(\lambda_l^+-\zeta)(\zeta-\lambda_l^-)}} \dif \zeta \right].\\
  \end{aligned}
\end{equation}

\subsubsection{Dispersive shock wave region: $-\frac{2\lambda_r^++\lambda_l^++\lambda_l^-}{2}+\frac{2(\lambda_r^+ - \lambda_l^+)
(\lambda_r^+ - \lambda_l^-)}{-2\lambda_r^+ + \lambda_l^+ + \lambda_l^-} < \xi < v_2(\lambda_r^+, \lambda_r^-, \lambda_l^+, \lambda_l^-)$} \label{ALSH}
\
\newline
\indent
When the stationary phase point $\eta_+ (\xi)$ moves inside $\mathcal{I}_r$, the exponential oscillation $e^{\pm 2 \mathrm{i} t \theta(z)}$ appears in the jump matrix $V^{(2)}(z;x,t)$
on the interval $(\eta_+(\xi), \lambda_r^+) \subset \mathcal{I}_r$. To remove the oscillation, the $g$-function should be modified by including two bands
$[\lambda_l^-, \lambda_l^+]$ and $[\lambda_s(\xi), \lambda_r^+]$. For $\boldsymbol{\lambda}=(\lambda_r^+, \lambda_s(\xi), \lambda_l^+, \lambda_l^-)$, the soft edge $\lambda_s(\xi)$ is uniquely determined by the following implicit function (solvable for $\lambda_s(\xi) \in ( \lambda_l^+, \lambda_r^+ )$) \cite{Kodama2006}:
\begin{equation}\label{vv2}
  \xi = v_2(\boldsymbol{\lambda}) := -\frac{\lambda_r^+ + \lambda_s(\xi) + \lambda_l^+ + \lambda_l^-}{2}+ \frac{ \lambda_r^+- \lambda_s(\xi)  }{1- \frac{\lambda_r^+-\lambda_l^+}{\lambda_s(\xi)-\lambda_l^+} \frac{E(m(\xi))}{K(m(\xi))}}
\end{equation}
with the elliptic modulus
\begin{equation}
  m(\xi)= \sqrt{\frac{(\lambda_r^+ - \lambda_s(\xi))(\lambda_l^+ -\lambda_l^-)}{(\lambda_r^+ -\lambda_l^+ )(\lambda_s(\xi) -\lambda_l^-)}}.
\end{equation}
The boundaries of this region are characterized by degeneration below
\begin{equation}\label{A2B}
  \begin{aligned}
  &\xi \to \xi_{A1} = v_2(\lambda_r^+, \lambda_r^+, \lambda_l^+, \lambda_l^-)= -\frac{2\lambda_r^++\lambda_l^++\lambda_l^-}{2}+\frac{2(\lambda_r^+ - \lambda_l^+)
  (\lambda_r^+ - \lambda_l^-)}{-2\lambda_r^+ + \lambda_l^+ + \lambda_l^-}, \quad &\mathrm{as} \quad \lambda_s(\xi) \to \lambda_r^+ ,\\
  &\xi \to \xi_{A2} =  v_2(\lambda_r^+, \lambda_r^-, \lambda_l^+, \lambda_l^-), \quad &\mathrm{as} \quad \lambda_s(\xi) \to \lambda_r^-.
  \end{aligned}
\end{equation}
\begin{figure}[htbp]
  \centering
 \includegraphics[width=10cm]{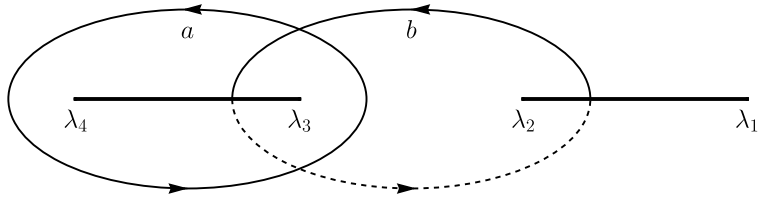}
\caption{{\protect\small The basis $\{a, b\}$  of the genus-one Riemann surface  $X$.}}
\label{figRS}
\end{figure}

To determine the two-band $g$-function, consider the Riemann surface $X$ of genus one with the branch points $ \boldsymbol{\lambda}=(\lambda_1, \lambda_2, \lambda_3, \lambda_4)$ defined by
\begin{equation}
  X:=\{ P=(z,w) : w^2=\prod_{k=1}^4 (z-\lambda_k )=:\mathcal{R} (z ;\boldsymbol{\lambda})^2 \},
\end{equation}
whose upper and lower sheets are denoted by $X_+$ and $X_-$ using the relations:
\begin{equation}
  \mathcal{R} (z ;\boldsymbol{\lambda})= \pm z^2 (1+ \mathcal{O}(z^{-1})), \qquad \mathrm{as} \quad z=\pi (P) \to \infty, \quad P \in X_{\pm}.
\end{equation}
Here the standard projection $\pi (P)=z$ defines $X$ as a two-sheet cover of the Riemann sphere $\mathbb{C} \mathbb{P} ^1$. Denote the preimages of $\infty$ on $X_{\pm}$ by $\infty_\pm$, respectively. Choose the base
$\{a, b\}$ of the homology group $H_1(X)$ so that the $a$-cycle becomes a counterclockwise oval entirely on the upper sheet around the interval $(\lambda_4, \lambda_3)$,
while the $b$-cycle starts from $(\lambda_2, \lambda_1)$, goes counterclockwise on the upper sheet to $(\lambda_4, \lambda_3)$,
and returns to the starting point on the lower sheet, as shown in Figure \ref{figRS}. Now, define two polynomials
\begin{equation}\label{polynomial}
  \begin{aligned}
    &P_{0} (z ;\boldsymbol{\lambda}):= z^2- \frac{1}{2} \sigma_1(\boldsymbol{\lambda})z + \gamma_0,\\
    &P_{1} (z ;\boldsymbol{\lambda}):= z^3- \frac{1}{2} \sigma_1(\boldsymbol{\lambda})z^2 + (\frac{1}{2} \sigma_2(\boldsymbol{\lambda})-\frac{1}{8} \sigma_1(\boldsymbol{\lambda})^2 )z+ \gamma_1,
  \end{aligned}
\end{equation}
where $  \sigma_1(\boldsymbol{\lambda})=\sum _{k=1} ^{4}\lambda_k$ and $\sigma_2(\boldsymbol{\lambda})= \sum _{1 \leq k < l \leq 4 } \lambda_k\lambda_l$ are elementary symmetric polynomials
and
\begin{equation}
  \begin{aligned}
    &\gamma_{0}=\frac{1}{2}(\lambda_1\lambda_2+\lambda_3\lambda_4)-\frac{1}{2}(\lambda_1-\lambda_3)(\lambda_2-\lambda_4)\frac{E(m)}{K(m)},\\
    &\gamma_{1}=\frac{1}{8}(\lambda_1\lambda_2-\lambda_3\lambda_4)(\lambda_1+\lambda_2-\lambda_3-\lambda_4)-\frac{1}{8}\sigma_1(\boldsymbol{\lambda})(\lambda_1-\lambda_3)(\lambda_2-\lambda_4)\frac{E(m)}{K(m)}
  \end{aligned}
\end{equation}
are uniquely determined by the conditions
\begin{equation} \label{acv}
  \begin{aligned}
  &\int_{\lambda_4}^{\lambda_3} \frac{P_0 (z ;\boldsymbol{\lambda})}{\mathcal{R}_+(z ;\boldsymbol{\lambda})}\,\dif z =0, \\
  &\int_{\lambda_4}^{\lambda_3} \frac{P_1 (z ;\boldsymbol{\lambda})}{\mathcal{R}_+(z ;\boldsymbol{\lambda})}\,\dif z =0.
  \end{aligned}
\end{equation}
Then define the differential $\dif g(z)$ as the Abelian differential of the second kind on $X$ as
\begin{equation}
  \dif g(z)=\frac{2 P_1 (z ;\boldsymbol{\lambda})+ \xi P_0 (z ;\boldsymbol{\lambda})}{\mathcal{R}(z ;\boldsymbol{\lambda})}\,\dif z,
\end{equation}
and we obtain the two-band $g$-function
\begin{equation}
  g(z)= \int_{\lambda_1}^{z} \frac{2 P_1 (\zeta ;\boldsymbol{\lambda})+ \xi P_0 (\zeta ;\boldsymbol{\lambda})}{\mathcal{R}(\zeta ;\boldsymbol{\lambda})}\dif \zeta,
\end{equation}
which is analytic in $z\in \mathbb{C} \backslash ([\lambda_4,\lambda_3] \cup [\lambda_2,\lambda_1])$ and satisfies the jump conditions
\begin{equation}
  \begin{aligned}
\left\{  \begin{array}{ll}
   g _{+}(z)+g _{-}(z)=0,  ~& z\in (\lambda_2,\lambda_1), \\ \\
  g _{+}(z)+g _{-}(z)= \oint_b \dif g ,  ~& z\in (\lambda_4,\lambda_3).
  \end{array} \right.
  \end{aligned}
\end{equation}

Now, go back to this region where the branch points $\boldsymbol{\lambda}=(\lambda_r^+, \lambda_s(\xi), \lambda_l^+, \lambda_l^-)$ have one soft edge $\lambda_s(\xi)$ and three hard edges.
Then $\lambda_s(\xi)$ is the simple zero of the cubic polynomial $2 P_1 (z ;\boldsymbol{\lambda})+ \xi P_0 (z ;\boldsymbol{\lambda})$ since $\lambda_k$ is
a soft edge if and only if the Whitham velocity $v_k(\boldsymbol{\lambda})$  is equal to the self-similar variable $\xi=x/t$ for some $k$ \cite{grava2002generation}:
\begin{equation}
  \xi=v_k(\boldsymbol{\lambda})= -2 \frac{P_1 (\lambda_k ;\boldsymbol{\lambda})}{P_0 (\lambda_k ;\boldsymbol{\lambda})}.
\end{equation}
Due to the $a$-cycle vanishment condition (\ref{acv}) of $\dif g(z)$, the other two zeros of $2 P_1 (z ;\boldsymbol{\lambda})+ \xi P_0 (z ;\boldsymbol{\lambda})$ are located on two bands, respectively.
Further, label the zeros $\eta_-(\xi) \in (\lambda_4,\lambda_3)$ and $\eta_+(\xi) \in (\lambda_2,\lambda_1)$ and  rewrite
$2 P_1 (z ;\boldsymbol{\lambda})+ \xi P_0 (z ;\boldsymbol{\lambda})=2(z-\eta_+(\xi))(z-\lambda_s(\xi))(z-\eta_-(\xi))$, where $\eta_+(\xi)$ and $\eta_-(\xi)$ can be determined directly by comparing the coefficients.
\par
\begin{figure}[tbp]
  \centering
  \subfigure[]{\includegraphics[width=5.5cm]{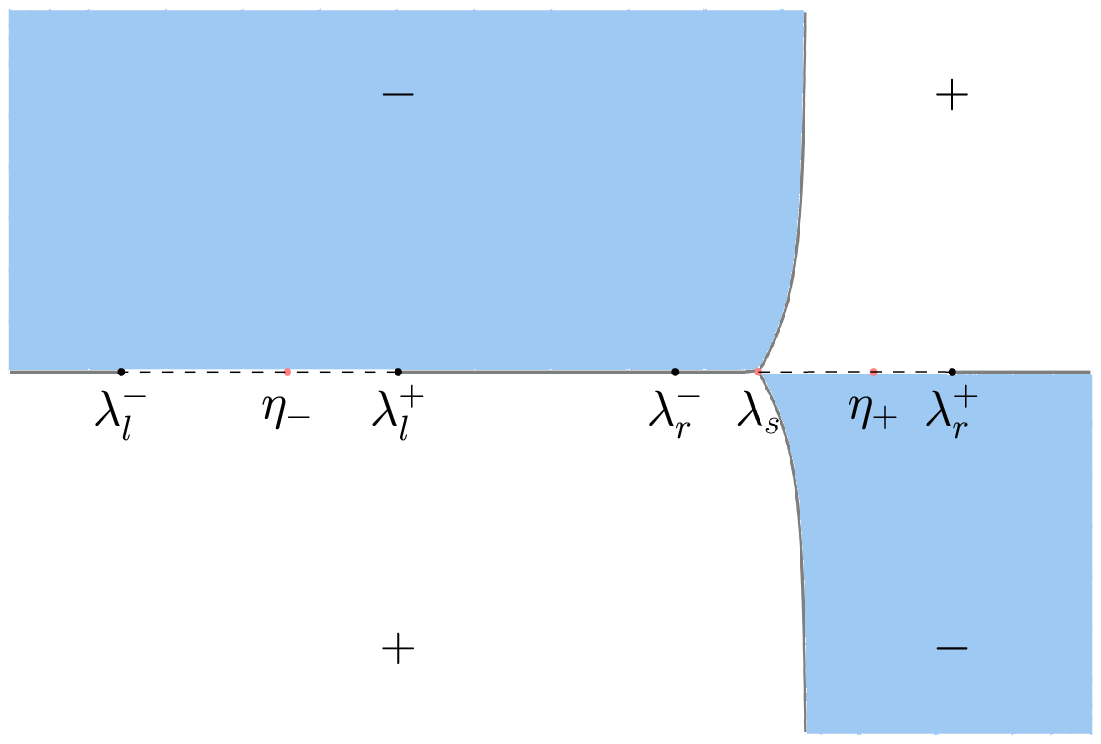}}
  \subfigure[]{\includegraphics[width=5cm]{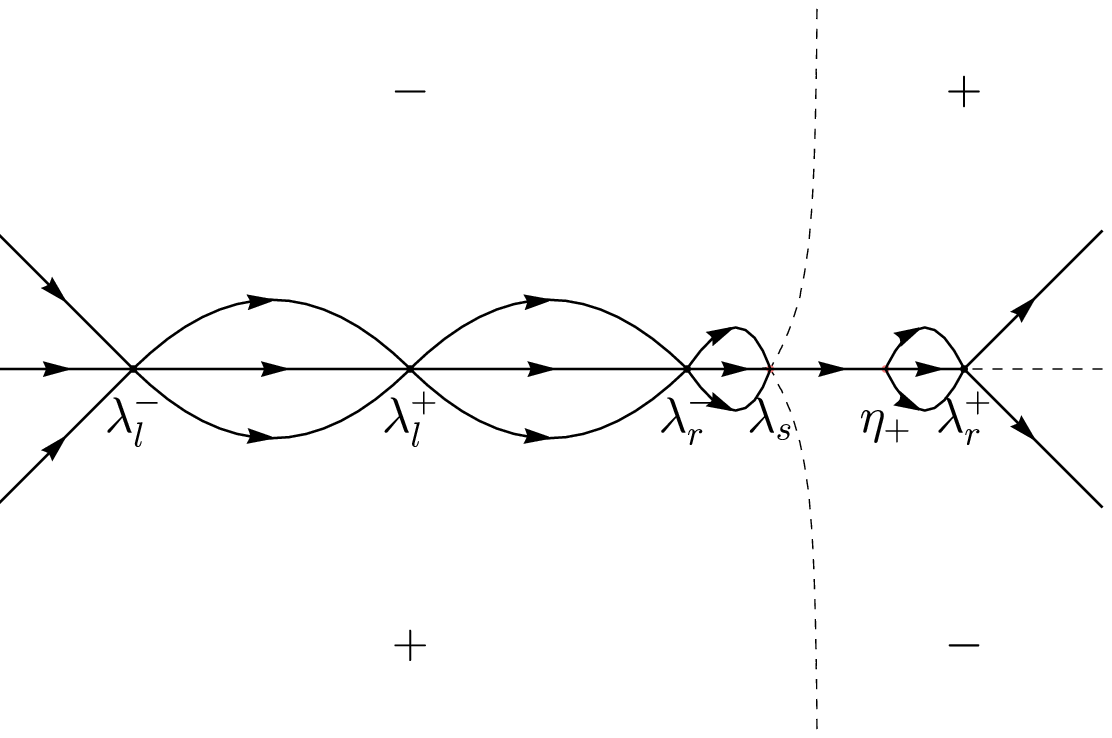}}
  \subfigure[]{\includegraphics[width=5.5cm]{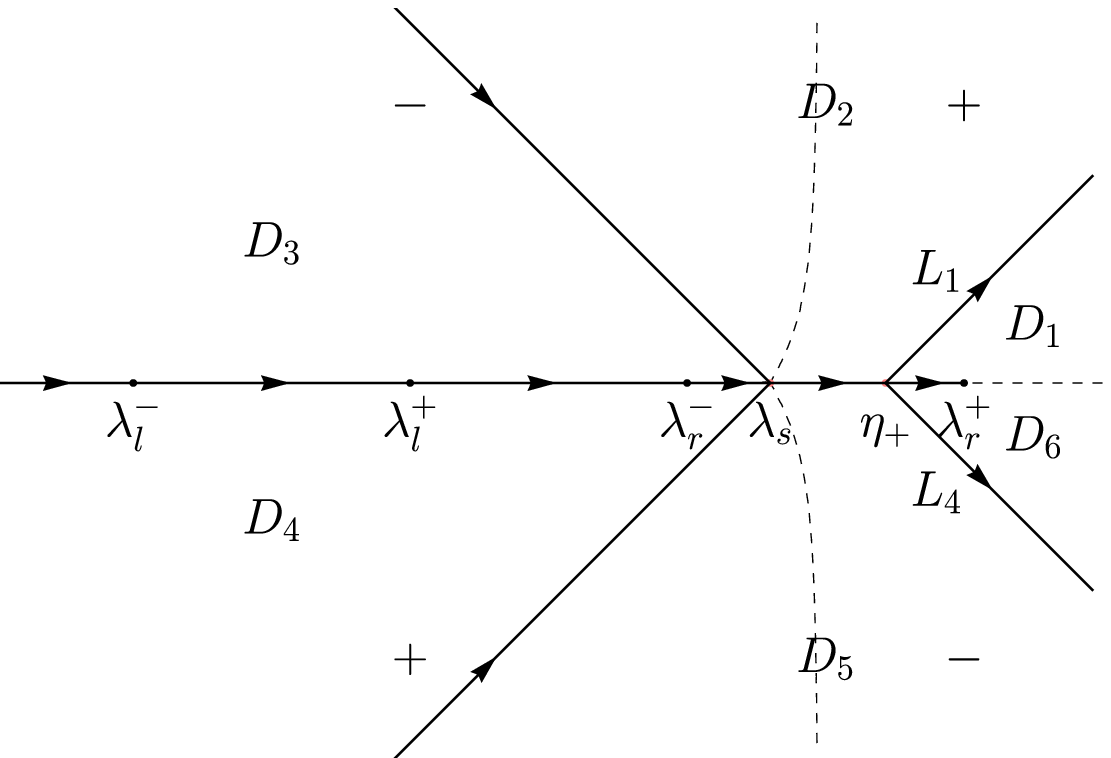}}
\caption{{\protect\small (a) The sign structure of Im$(g(z))$; (b) Opening lenses; (c) The jump contours of $M^{(2)}(z;x,t)$.}}
\label{figA2}
\end{figure}
The two-band $g$-function with one soft edge (or shock $g$-function) is given by
\begin{equation}\label{A2G}
  g(z)= \int_{\lambda_r^+}^{z} \frac{2(\zeta-\eta_+(\xi))(\zeta-\lambda_s(\xi))(\zeta-\eta_-(\xi))}{\mathcal{R}(\zeta ;\lambda_r^+, \lambda_s(\xi), \lambda_l^+, \lambda_l^-)}\dif \zeta,
\end{equation}
which is analytic in $z\in \mathbb{C} \backslash ([\lambda_l^-,\lambda_l^+] \cup [\lambda_s(\xi),\lambda_r^+])$ and satisfies the jump conditions
\begin{equation} \label{gamma}
  \begin{aligned}
\left\{  \begin{array}{ll}
   g _{+}(z)+g _{-}(z)=0,  ~& z\in (\lambda_s(\xi),\lambda_r^+), \\ \\
  g _{+}(z)+g _{-}(z)= \gamma :=\oint_b \dif g ,  ~& z\in (\lambda_l^-,\lambda_l^+).
  \end{array} \right.
  \end{aligned}
\end{equation}
Furthermore, it is seen that
\begin{equation} \label{A2GINF}
    g_\infty= \int_{\lambda_r^+}^{\infty} \left(\frac{2(\zeta-\eta_+(\xi))(\zeta-\lambda_s(\xi))(\zeta-\eta_-(\xi))}{\mathcal{R}(\zeta ;\lambda_r^+, \lambda_s(\xi), \lambda_l^+, \lambda_l^-)}-2\zeta-\xi\, \right)\,\dif \zeta - \theta(\lambda_r^+)
\end{equation}
is well-defined since the relation
\begin{equation}
  \frac{2(\zeta-\eta_+(\xi))(\zeta-\lambda_s(\xi))(\zeta-\eta_-(\xi))}{\mathcal{R}(\zeta ;\lambda_r^+, \lambda_s(\xi), \lambda_l^+, \lambda_l^-)}-2\zeta-\xi=\mathcal{O}(\frac{1}{\zeta^2}), \qquad \zeta\to \infty,
\end{equation}
ensures that the integral converges.

Then we open lenses from intervals to the steepest descent contours through $\lambda_s(\xi)$ and $\eta_+(\xi)$  (as shown in Figure \ref{figA2})  by the transformations
$M(z;x,t) \mapsto M^{(1)}(z;x,t) \mapsto M^{(2)}(z;x,t)$, which yields the following RHP.

\begin{rhp}
  Find a $2 \times 2$ matrix-valued function $M^{(2)}(z;x,t)$ with the following properties:
\begin{enumerate} [label=(\roman*)]
  \item  $M^{(2)}(z;x,t)$ is analytic in $z \in \mathbb{C} \backslash \Sigma^{(2)}$, where $\Sigma^{(2)}=(-\infty, \lambda_r^+] \cup L_1 \cup L_2 \cup  L_3 \cup L_4.$
  \item  $M^{(2)}(z;x,t)=I+\mathcal{O} (z^{-1}) $ as $z\to \infty$.
  \item  $M^{(2)}(z;x,t)$ achieves the CBVs  $M^{(2)}_{+}(z;x,t)$ and $M^{(2)}_{-}(z;x,t)$  on $\Sigma^{(2)}$ away from self-intersection points and
  branch points that satisfy the jump condition $M^{(2)}_{+}(z;x,t)=M^{(2)}_{-}(z;x,t)V^{(2)}(z;x,t)$, where
  \begin{equation}\label{AB2V2}
    \begin{aligned}
    V^{(2)}(z;x,t)&=\left\{  \begin{array}{ll}
     (1-r(z)r^{*}(z))^{\sigma_{3}},  ~& z\in (-\infty, \lambda_l^-) \cup (\lambda_l^+, \lambda_r^-), \\ \\
     (a_+(z) a_-^*(z))^{-\sigma_{3}},  ~&z \in  (\lambda_r^-, \lambda_s(\xi)), \\ \\
     \begin{pmatrix} 0 & -r_{-}^{*}(z) e^{-\mathrm{i} \gamma}  \\ r_{+}(z) e^{\mathrm{i} \gamma}  & 0 \end{pmatrix}, ~&z \in (\lambda_l^-, \lambda_l^+), \\ \\
     \begin{pmatrix} \frac{e^{ 2 \mathrm{i} t g_+(z)}}{a_+(z) a_-^*(z)} & -1 \\1 & 0 \end{pmatrix}, ~&z \in (\lambda_s(\xi), \eta_+(\xi)), \\ \\
     \begin{pmatrix} 0 & -1 \\1 & 0 \end{pmatrix}, ~&z \in (\eta_+(\xi), \lambda_r^+), \\ \\
     \begin{pmatrix} 1 & 0 \\ r(z)e^{2\mathrm{i} tg(z)} & 1 \end{pmatrix}, ~ &z \in L_1, \\ \\
     \begin{pmatrix} 1 & \frac{-r^*(z)}{1-r(z)r^*(z)}e^{-2\mathrm{i} tg(z)} \\0 & 1 \end{pmatrix},  ~&z \in L_2,\\ \\
     \begin{pmatrix} 1 & 0 \\ \frac{r(z)}{1-r(z)r^*(z)}e^{2\mathrm{i} tg(z)} & 1 \end{pmatrix},  ~&z \in L_3,\\ \\
     \begin{pmatrix} 1 & -r^*(z)e^{-2\mathrm{i} tg(z)} \\ 0  & 1 \end{pmatrix}, ~ &z \in L_4.\end{array} \right.\\
    \end{aligned}
  \end{equation}
\end{enumerate}
\end{rhp}

In comparison with the plane wave region, the above jump matrices $V^{(2)}(z;x,t)$ on two disjoint intervals $(\lambda_l^-, \lambda_l^+)$ and $(\lambda_s(\xi),\lambda_r^+)$ converge to two different constant matrices for large $t$.
To arrive at the solvable limiting RHP with piecewise constant jumps across two bands, introduce a scalar function
\begin{equation} \label{AB2delta}
  \begin{aligned}
     \delta(z)= \exp & \left\{  \frac{\mathcal{R} (z;\boldsymbol{\lambda})}{2 \pi \mathrm{i}}  \left[ \left(\int_{-\infty}^{\lambda_l^-} + \int_{\lambda_l^+}^{\lambda_r^-} \right) \frac{\ln (1-|r(\zeta)|^2)}{\mathcal{R} (\zeta;\boldsymbol{\lambda}) } \,\frac{\dif\zeta}{\zeta - z} \right. \right. \\
      &\left. \left. + \int_{\lambda_r^-}^{\lambda_s(\xi)} \frac{-\ln(a_+(\zeta)a_-^*(\zeta))}{\mathcal{R} (\zeta;\boldsymbol{\lambda})} \,\frac{\dif\zeta}{\zeta - z} + \int_{\lambda_l^-}^{\lambda_l^+} \frac{\ln (r_+(\zeta))}{\mathcal{R}_+ (\zeta;\boldsymbol{\lambda})} \, \frac{\dif\zeta}{\zeta - z}
      \right]  \right\},
    \end{aligned}
\end{equation}
where $\boldsymbol{\lambda}=(\lambda_r^+, \lambda_s(\xi), \lambda_l^+, \lambda_l^-)$. The function $\delta(z)$ satisfies the following properties:
\begin{enumerate} [label=(\roman*)]
  \item $\delta(z)$ is analytic in $z \in \mathbb{C} \backslash (-\infty, \lambda_r^+]$;
  \item $\delta(z)=\delta_\infty(z)(1+\mathcal{O} (z^{-1}))$ as $z\to \infty$, where
    \begin{equation}\label{AB2deltaINF}
     \begin{aligned}
     \delta_{\infty }(z)= \exp \{\mathrm{i}(\varphi^{(1)}z +\varphi^{(0)}) \}
     \end{aligned}
    \end{equation}
  and $\varphi^{(1)} $ and $ \varphi^{(0)}$ are real coefficients
  \begin{equation} \label{phi111}
    \begin{aligned}
     \varphi^{(1)}=\frac{1}{2 \pi}\left[ \left(\int_{-\infty}^{\lambda_l^-} + \int_{\lambda_l^+}^{\lambda_r^-} \right) \frac{\ln (1-|r(\zeta)|^2) }{\mathcal{R} (\zeta;\boldsymbol{\lambda}) } \dif \zeta
      +\int_{\lambda_r^-}^{\lambda_s(\xi)} \frac{-\ln(a_+(\zeta)a_-^*(\zeta))}{\mathcal{R} (\zeta;\boldsymbol{\lambda})} \dif \zeta + \int_{\lambda_l^-}^{\lambda_l^+} \frac{\ln (r_+(\zeta))}{\mathcal{R}_+ (\zeta;\boldsymbol{\lambda})} \,{\dif\zeta}\right],\\
    \end{aligned}
  \end{equation}
  \begin{equation} \label{phi0}
    \begin{aligned}
      \varphi^{(0)}=\frac{1}{2 \pi} & \left[\left(\int_{-\infty}^{\lambda_l^-} + \int_{\lambda_l^+}^{\lambda_r^-} \right) \frac{\ln (1-|r(\zeta)|^2)  (\zeta+ V) }{\mathcal{R} (\zeta;\boldsymbol{\lambda}) } \dif \zeta +  \int_{\lambda_r^-}^{\lambda_s(\xi)} \frac{-\ln(a_+(\zeta)a_-^*(\zeta))(\zeta+ V )}{\mathcal{R} (\zeta;\boldsymbol{\lambda})} \dif \zeta \right.\\
      & \left. +  \int_{\lambda_l^-}^{\lambda_l^+} \frac{\ln (r_+(\zeta)) (\zeta+ V)}{\mathcal{R}_+ (\zeta;\boldsymbol{\lambda})} \,{\dif\zeta} \right],
    \end{aligned}
  \end{equation}
  where
  \begin{equation} \label{VVV}
    V=-\frac{1}{2}\sigma_1(\boldsymbol{\lambda})=-\frac{1}{2}(\lambda_r^+ + \lambda_s(\xi)+ \lambda_l^++\lambda_l^-).
  \end{equation}
  \item $\delta(z)$ achieves the CBVs on  $(-\infty, \lambda_r^+ )$ satisfying the jump conditions
  \begin{equation*}
    \begin{aligned}
  \left\{  \begin{array}{ll}
     \delta _{+}(z)=(1-r(z)r^{*}(z))\delta _{-}(z),  ~& z\in (-\infty, \lambda_l^-) \cup (\lambda_l^+, \lambda_r^-) , \\ \\
     \delta _{+}(z)=(a_+(z) a_-^*(z))^{-1}\delta _{-}(z),  ~& z\in (\lambda_r^-, \lambda_s(\xi)), \\ \\
     \delta _{+}(z)\delta _{-}(z)=r_{+}(z) ,  ~& z\in (\lambda_l^-, \lambda_l^+),\\ \\
     \delta _{+}(z)\delta _{-}(z)=1,  ~& z\in (\lambda_s(\xi), \lambda_r^+).
    \end{array} \right.\\
    \end{aligned}
  \end{equation*}
\end{enumerate}

However, $\delta(z)$ has an essential singularity at infinity as shown in (\ref{AB2deltaINF}), which leads to a new oscillation  in the jump matrix by the previous
transformation (\ref{2to3}). We  use the $g$-function mechanism again to modify $\delta(z)$ to eliminate the singularity. Introduce  the Abelian integral
$\hat{g} (z) $ of the second kind as the new $g$-function:
\begin{equation}
  \hat{g} (z) = \varphi^{(1)} \int_{\lambda_r^+}^{z} \frac{P_0 (\zeta ;\lambda_r^+, \lambda_s(\xi), \lambda_l^+, \lambda_l^-)}{\mathcal{R}(\zeta ;\lambda_r^+, \lambda_s(\xi), \lambda_l^+, \lambda_l^-)}\,\dif \zeta,
\end{equation}
which is also analytic in $z\in \mathbb{C} \backslash ([\lambda_l^-,\lambda_l^+] \cup [\lambda_s(\xi),\lambda_r^+])$ and satisfies the jump conditions
\begin{equation} \label{gammah}
  \begin{aligned}
\left\{  \begin{array}{ll}
  \hat{g} _{+}(z)+\hat{g}_{-}(z)=0,  ~& z\in (\lambda_s(\xi),\lambda_r^+), \\ \\
  \hat{g}_{+}(z)+\hat{g} _{-}(z)= \hat{\gamma} :=\oint_b \, \dif\hat{g} ,  ~& z\in (\lambda_l^-,\lambda_l^+).
  \end{array} \right.
  \end{aligned}
\end{equation}
The large-$z$ expansion of $\hat{g} (z)$ is of the form
\begin{equation}
  \hat{g} (z) = \varphi^{(1)} z +\hat{g}_\infty +\mathcal{O}(z^{-1}),
\end{equation}
where
\begin{equation} \label{AB2ghatINF}
  \hat{g}_\infty = \varphi^{(1)} \int_{\lambda_r^+}^{\infty} \left(\frac{P_0 (z ;\lambda_r^+, \lambda_s(\xi), \lambda_l^+, \lambda_l^-)}{\mathcal{R}(z ;\lambda_r^+, \lambda_s(\xi), \lambda_l^+, \lambda_l^-)} -1 \right)\,\dif z - \varphi^{(1)} \lambda_r^+.
\end{equation}

Now, introduce the modified function
\begin{equation}\label{AB2mdelta}
  \hat{\delta} (z) = \delta  (z) e^{- \mathrm{i} \hat{g} (z) }
\end{equation}
with the following properties:
\begin{enumerate} [label=(\roman*)]
  \item $\hat{\delta}(z)$ is analytic in $z \in \mathbb{C} \backslash (-\infty, \lambda_r^+]$;
  \item $\hat{\delta}(z)=\hat{\delta}_\infty+\mathcal{O} (z^{-1})$ as $z\to \infty$, where
    \begin{equation}\label{AB2mdeltaINF}
     \begin{aligned}
     \hat{\delta}_{\infty }= \exp \{\mathrm{i}(\varphi^{(0)} -\hat{g}_\infty ) \};
     \end{aligned}
    \end{equation}
  \item $\hat{\delta}(z)$ achieves the CBVs on  $(-\infty, \lambda_r^+ )$ satisfying the jump conditions
  \begin{equation*}
    \begin{aligned}
  \left\{  \begin{array}{ll}
     \hat{\delta} _{+}(z)=(1-r(z)r^{*}(z))\hat{\delta} _{-}(z),  ~& z\in (-\infty, \lambda_l^-) \cup (\lambda_l^+, \lambda_r^-) , \\ \\
     \hat{\delta} _{+}(z)=(a_+(z) a_-^*(z))^{-1}\hat{\delta} _{-}(z),  ~& z\in (\lambda_r^-, \lambda_s(\xi)), \\ \\
     \hat{\delta} _{+}(z)\hat{\delta} _{-}(z)=r_{+}(z) e^{-\mathrm{i} \hat{\gamma}} ,  ~& z\in (\lambda_l^-, \lambda_l^+),\\ \\
     \hat{\delta} _{+}(z)\hat{\delta} _{-}(z)=1,  ~& z\in (\lambda_s(\xi), \lambda_r^+).
    \end{array} \right.\\
    \end{aligned}
  \end{equation*}
\end{enumerate}
Then define the modified transformation
\begin{equation}
	M^{(3)}(z;x,t)=\hat{\delta} _{\infty}^{\sigma_{3}} M^{(2)}(z;x,t) \hat{\delta}(z) ^{-\sigma_{3}},
\end{equation}
which leads to the following RHP.
\begin{rhp}\label{lsrhp} Find a $2 \times 2$ matrix-valued function $M^{(3)}(z;x,t)$ with the following properties:
  \begin{enumerate} [label=(\roman*)]
    \item $M^{(3)}(z;x,t)$ is analytic in $z \in \mathbb{C} \backslash \Sigma^{(3)}$, where $\Sigma^{(3)}= [\lambda_l^-,\lambda_l^+] \cup [\lambda_s(\xi),\lambda_r^+]\cup L_1 \cup L_2 \cup  L_3 \cup L_4.$
    \item $M^{(3)}(z;x,t)=I+\mathcal{O} (z^{-1}) $ as $z\to \infty$.
    \item $M^{(3)}(z;x,t)$ achieves the CBVs $M^{(3)}_{+}(z;x,t)$ and $M^{(3)}_{-}(z;x,t)$ on $\Sigma^{(3)}$ away from self-intersection points and
  branch points that satisfy the jump condition $M^{(3)}_{+}(z;x,t)=M^{(3)}_{-}(z;x,t)V^{(3)}(z;x,t)$, where
  \begin{equation}\label{AB2V3}
    \begin{aligned}
  V^{(3)}(z;x,t)&=\left\{  \begin{array}{ll}
     \begin{pmatrix} 0 & - e^{-\mathrm{i} (\gamma+\hat{\gamma})}  \\  e^{\mathrm{i} (\gamma+\hat{\gamma})}  & 0 \end{pmatrix}, ~&z \in (\lambda_l^-, \lambda_l^+), \\ \\
     \begin{pmatrix}  \frac{\hat{\delta}^{-1}_{+}(z)\hat{\delta} _{-}(z)}{a_+(z) a_-^*(z)}  e^{2 \mathrm{i} t g_+(z)} & -1 \\1 & 0 \end{pmatrix}, ~&z \in (\lambda_s(\xi), \eta_+(\xi)), \\ \\
     \begin{pmatrix} 0 & -1 \\1 & 0 \end{pmatrix}, ~&z \in (\eta_+(\xi), \lambda_r^+), \\ \\
     \begin{pmatrix} 1 & 0 \\ r(z) \hat{\delta} ^{-2}(z) e^{2\mathrm{i} tg(z)} & 1 \end{pmatrix}, ~ &z \in L_1, \\ \\
     \begin{pmatrix} 1 & \frac{-r^*(z)}{1-r(z)r^*(z)} \hat{\delta} ^2(z) e^{-2\mathrm{i}t g(z)} \\0 & 1 \end{pmatrix},  ~&z \in L_2,\\ \\
     \begin{pmatrix} 1 & 0 \\ \frac{r(z)}{1-r(z)r^*(z)} \hat{\delta} ^{-2}(z) e^{2\mathrm{i} tg(z)} & 1 \end{pmatrix},  ~&z \in L_3,\\ \\
     \begin{pmatrix} 1 & -r^*(z) \hat{\delta} ^2(z) e^{-2\mathrm{i} tg(z)} \\ 0  & 1 \end{pmatrix}, ~ &z \in L_4. \end{array} \right.\\
    \end{aligned}
  \end{equation}
  \end{enumerate}
\end{rhp}

The jump matrices for $M^{(3)}(z;x,t)$ uniformly converge to the identity matrix or constant matrices outside a fixed neighborhood $\mathcal{U}$ of $\lambda_s(\xi)$,
while the convergence is not uniform inside $\mathcal{U}$ due to the exponential phase function $g(z)=\mathcal{O}((z-\lambda_s(\xi))^{3/2})$.
We construct an outer model parametrix that is the solution of the two-band limiting problem
\begin{equation}
  M^{(\mathrm{mod})}_+(z;x,t)=M^{(\mathrm{mod})}_-(z;x,t) V^{(\mathrm{mod})},
\end{equation}
where
\begin{equation}\label{AB2vmod}
  \begin{aligned}
    V^{(\mathrm{mod})}&=\left\{  \begin{array}{ll}
   \begin{pmatrix} 0 & - e^{-\mathrm{i} (\gamma+\hat{\gamma})}  \\  e^{\mathrm{i} (\gamma+\hat{\gamma})}  & 0 \end{pmatrix}, ~&z \in (\lambda_l^-, \lambda_l^+), \\ \\
   \begin{pmatrix} 0 & -1 \\1 & 0 \end{pmatrix}, ~&z \in (\lambda_s(\xi), \lambda_r^+). \end{array} \right.\\
   \end{aligned}
\end{equation}
The model problem can be solved by using Riemann theta functions \cite{farkas1992riemann,Deift-Its-Zhou-1997} on the aforementioned Riemann surface $X$  with branch points $\boldsymbol{\lambda}=(\lambda_r^+, \lambda_s(\xi), \lambda_l^+, \lambda_l^-)$ and
the canonical homology basis $\{a, b\}$, as shown in Figure \ref{figRS} .

Define the normalized Abelian differential of the first kind
\begin{equation}
  \omega_1= \frac{d}{\mathcal{R} (z;\boldsymbol{\lambda})} \, \dif z,
\end{equation}
where the constant $d$ is determined by the normalization condition
\begin{equation}
  \oint _{a} \, \omega_1 = 1,
\end{equation}
and can be expressed in terms of the complete elliptic integral of the first kind $K(m)$ via the formula
\begin{equation} \label{d}
  d= \frac{\mathrm{i}\sqrt{(\lambda_r^+-\lambda_l^+)(\lambda_s(\xi)-\lambda_l^-)}}{4 K(m)}.
\end{equation}
Define the Riemann period
\begin{equation} \label{tau}
  \tau := \oint _{b} \, \omega_1 = \frac{\mathrm{i}K(\sqrt{1-m^2})}{4 K(m)}
\end{equation}
with $\tau$ being purely imaginary and $\mathrm{Im}~\tau>0$. Then, define the Riemann theta function \cite{farkas1992riemann,Deift-Its-Zhou-1997}
\begin{equation}
  \label{Theta}
  \Theta (z)=\Theta (z;\tau)=\sum _{n \in \mathbb{Z} }  e^{2\pi \mathrm{i} n z + \pi \mathrm{i} n^2 \tau}
\end{equation}
with the following properties:
\begin{enumerate}[label=(\roman*)]
  \item $\Theta (z)$ is entire in $z \in \mathbb{C}$;
  \item $\Theta (z+1)=\Theta (z)$;
  \item $\Theta (z+\tau)=\Theta (z) e^{-2\pi \mathrm{i} z - \pi \mathrm{i}  \tau}$;
  \item $\Theta (z) =0$ iff $z=1/2 + \tau /2 + n + m \tau $ and $n, m \in \mathbb{Z}$.
\end{enumerate}
Let $\mathcal{A}(z)$ define the constraint of the Abel map to the complex plane
\begin{equation}
  \mathcal{A}(z):= \int_{\lambda_1}^{z} \omega_1 = \int_{\lambda_r^+}^{z} \omega_1,
\end{equation}
where the integration path lies away from $[\lambda_4, \lambda_3]\cup[\lambda_2, \lambda_1]=[\lambda_l^-, \lambda_l^+] \cup [\lambda_s(\xi), \lambda_r^+]$. Furthermore, $\mathcal{A}(z)$ satisfies the following relations:
\begin{equation}\label{Abelmap}
  \begin{aligned}
\left\{  \begin{array}{ll}
  \mathcal{A} _{+}(z)+\mathcal{A} _{-}(z)=\tau ,  ~& z\in (\lambda_4, \lambda_3)=(\lambda_l^-, \lambda_l^+),\\ \\
  \mathcal{A} _{+}(z)+\mathcal{A} _{-}(z)=0,  ~& z\in (\lambda_2, \lambda_1)=(\lambda_s(\xi), \lambda_r^+).
  \end{array} \right.\\
  \end{aligned}
\end{equation}
Define
\begin{equation}
  \alpha(z; \boldsymbol{\lambda}):=\left({\frac{(z-\lambda_1)(z-\lambda_3)}{(z-\lambda_2)(z-\lambda_4)}}\right)^{1/4}=\left({\frac{(z-\lambda_r^+)(z-\lambda_l^+)}{(z-\lambda_s(\xi))(z-\lambda_l^-)}}\right)^{1/4}
\end{equation}
with the branch cuts  $[\lambda_4, \lambda_3]=[\lambda_l^-, \lambda_l^+] $ and $ [\lambda_2, \lambda_1]=[\lambda_s(\xi), \lambda_r^+]$ and the behavior $\alpha(z; \boldsymbol{\lambda}) \thicksim 1 $ as $z \to \infty$. It is observed that $\alpha_+(z)= \mathrm{i} \alpha_-(z)$ along the cuts.

Now, define the $2 \times 2$ matrix-valued function $N(z; \boldsymbol{\lambda})$ with entries
\begin{equation}\label{N}
   \begin{aligned}
    N_{11}(z; \boldsymbol{\lambda})=\frac{(\alpha(z; \boldsymbol{\lambda})+\alpha^{-1}(z; \boldsymbol{\lambda}))}{2}  \frac{\Theta (\mathcal{A}(z)-\mathcal{A}(\infty)-(\gamma + \hat{\gamma})/{2 \pi})}{\Theta (\mathcal{A}(z)-\mathcal{A}(\infty))}, \\
    N_{12}(z; \boldsymbol{\lambda})=\frac{(\alpha(z; \boldsymbol{\lambda})-\alpha^{-1}(z; \boldsymbol{\lambda}))}{-2\mathrm{i}}  \frac{\Theta (\mathcal{A}(z)+\mathcal{A}(\infty)+(\gamma + \hat{\gamma})/{2 \pi})}{\Theta (\mathcal{A}(z)+\mathcal{A}(\infty))}, \\
    N_{21}(z; \boldsymbol{\lambda})=\frac{(\alpha(z; \boldsymbol{\lambda})-\alpha^{-1}(z; \boldsymbol{\lambda}))}{2\mathrm{i}}  \frac{\Theta (\mathcal{A}(z)+\mathcal{A}(\infty)-(\gamma + \hat{\gamma})/{2 \pi})}{\Theta (\mathcal{A}(z)+\mathcal{A}(\infty))}, \\
    N_{22}(z; \boldsymbol{\lambda})=\frac{(\alpha(z; \boldsymbol{\lambda})+\alpha^{-1}(z; \boldsymbol{\lambda}))}{2}  \frac{\Theta (\mathcal{A}(z)-\mathcal{A}(\infty)+(\gamma + \hat{\gamma})/{2 \pi})}{\Theta (\mathcal{A}(z)-\mathcal{A}(\infty))}, \\
 \end{aligned}
\end{equation}
where
\begin{equation} \label{Abelinf}
  \mathcal{A} (\infty)=  \int_{\lambda_r^+}^{\infty} \omega_1.
\end{equation}
It is claimed that $N(z; \boldsymbol{\lambda})$ can be viewed as a function analytic in $z\in \mathbb{C} \backslash ([\lambda_l^-,\lambda_l^+] \cup [\lambda_s(\xi),\lambda_r^+])$.
Indeed, the only possible singularities of $N(z; \boldsymbol{\lambda})$, except for branch points, originate from the functions $\Theta (\mathcal{A}(z)\pm\mathcal{A}(\infty))$ on the denominator, but these singularities can be canceled
by the possible zeros of $\alpha(z; \boldsymbol{\lambda})\pm\alpha^{-1}(z; \boldsymbol{\lambda})$. We will show that the function $\Theta (\mathcal{A}(z)+\mathcal{A}(\infty))$ has a unique zero on the upper sheet $X_+$. As a consequence,
the function $\Theta (\mathcal{A}(z)-\mathcal{A}(\infty))$ is nonzero on $X_+$ (see \cite{farkas1992riemann,Deift-Its-Zhou-1997}). Specifically, consider
\begin{equation}
  h(z)=\alpha(z; \boldsymbol{\lambda})(\alpha^{-1}(z; \boldsymbol{\lambda})-\alpha(z; \boldsymbol{\lambda}))
\end{equation}
as a function on the Riemann surface $X$. The function $h(z)$ has simple poles at $\lambda_s(\xi)$ and $\lambda_l^-$, and simple zeros at $\infty_+$ and at the same unique finite zero $\lambda^*$ of the function
$\alpha(z; \boldsymbol{\lambda})-\alpha^{-1}(z; \boldsymbol{\lambda})$ given by
\begin{equation}
  \lambda^*=\frac{\lambda_r^+\lambda_s(\xi)-\lambda_l^+\lambda_l^-}{\lambda_r^++\lambda_s(\xi)-\lambda_l^+-\lambda_l^-}.
\end{equation}
Then, $h$ is a meromorphic function on $X$ with principal divisor
\begin{equation}
  (h)=\lambda^*+\infty_+-\lambda_s(\xi)-\lambda_l^-.
\end{equation}
By Abel's theorem \cite{farkas1992riemann}, we have $\mathcal{A}((h))=0$, and thus
\begin{equation}
  \mathcal{A}(\lambda^*)+\mathcal{A}(\infty)=\mathcal{A}(\lambda_s(\xi))+\mathcal{A}(\lambda_l^-)=\frac{1}{2}+\frac{\tau}{2}+\mathbb{Z}+\tau\mathbb{Z},
\end{equation}
which is exactly the zero of the function $\Theta(z)$. Hence,  the function has a unique zero at the preimage of $\lambda^*$ on the upper sheet $X_+$.
The properties of $\Theta(z)$ and the relations (\ref{Abelmap}) imply that $N(z; \boldsymbol{\lambda})$ satisfies the jump conditions (\ref{AB2vmod}) of the model problem. Thus the outer model parametrix is given by
\begin{equation} \label{Ab2mod}
  M^{(\mathrm{mod})}(z;x,t)=N^{-1}(\infty; \boldsymbol{\lambda}) N(z; \boldsymbol{\lambda}).
\end{equation}

\begin{figure}[htbp]
  \centering
  \subfigure[]{\includegraphics[width=7cm]{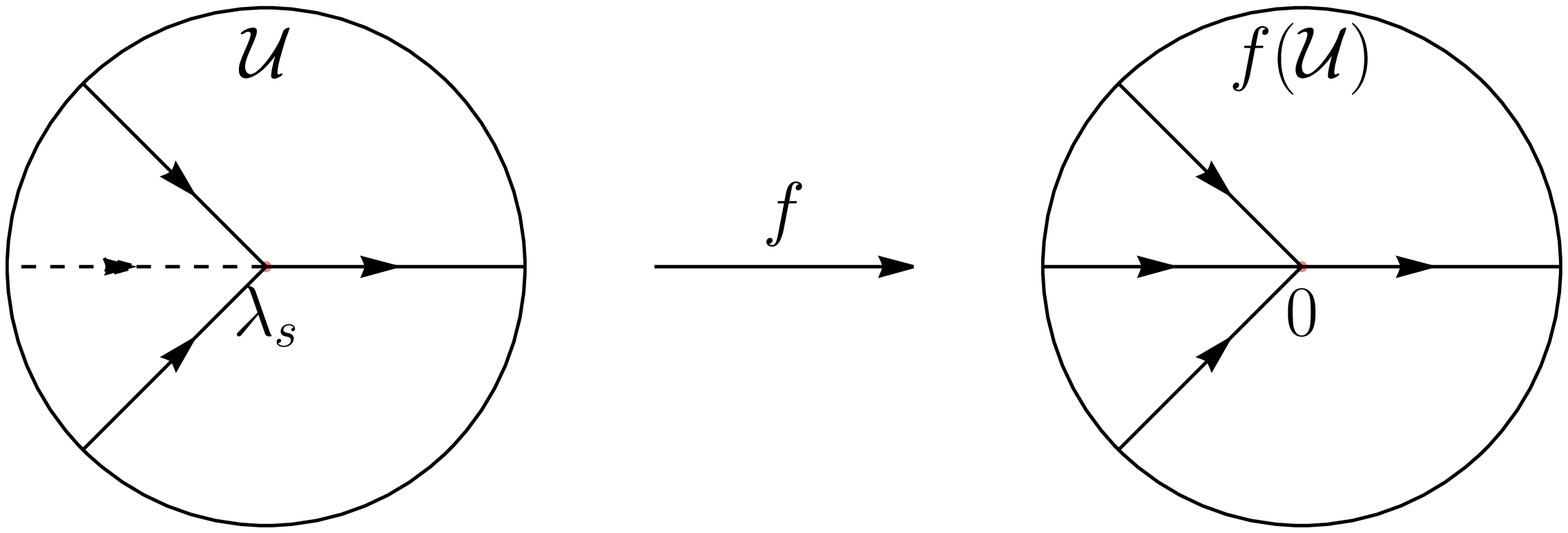}} \qquad
  \subfigure[]{\includegraphics[width=7cm]{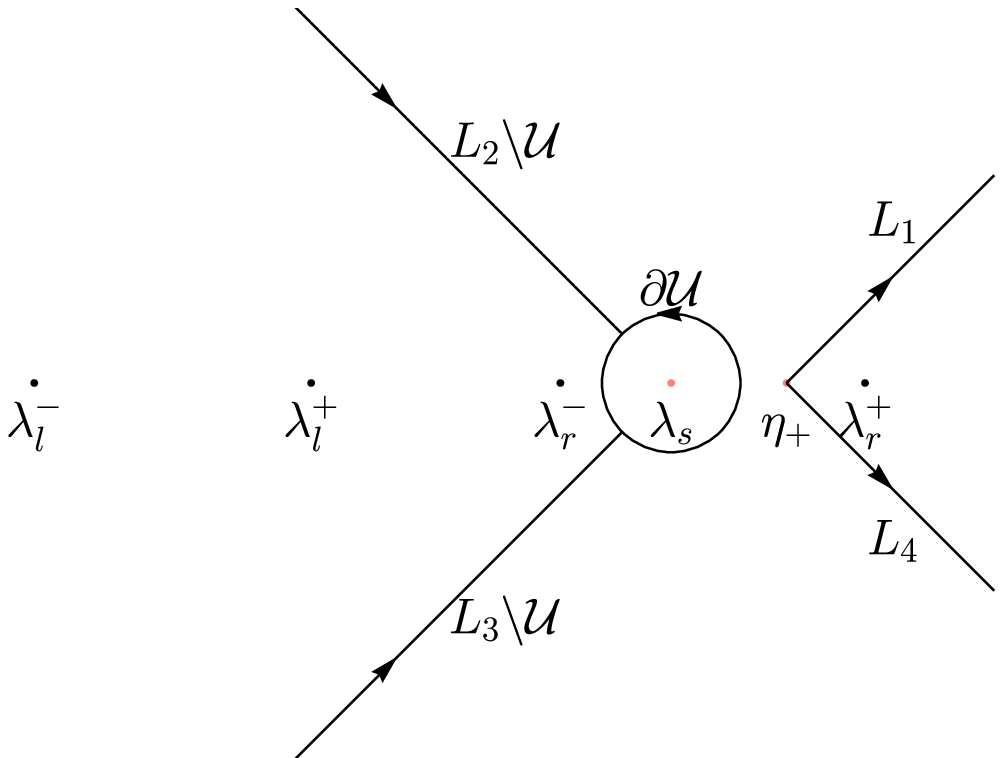}}
\caption{{\protect\small (a) The conformal mapping $f$ on $\mathcal{U}$; (b) The jump contours of error RHP.}}
\label{figA22}
\end{figure}

Inside $\mathcal{U}$, construct a local parametrix $M^{\mathrm{(loc)}}(z;x,t)$ with the same jump matrices as $M^{(3)}(z;x,t)$ with the help of the Airy local model $M^{(\mathrm{Ai})}(k)$ defined in Appendix \ref{airy}
and the change of variables $k=t^{2/3} f(z)$, where $f$ is a conformal mapping from $\mathcal{U}$ (in $z$) to a neighborhood of the origin (in $k$), as shown in Figure \ref{figA22}, defined by
\begin{equation}\label{Ab2f}
  \frac{2}{3}f^{3/2}(z)=-\mathrm{i\, sgn}(\mathrm{Im}\,z)\,(g(z)-\gamma/2)
\end{equation}
and the branch cut is chosen so that $f'(\lambda_s(\xi))>0$. Define the local parametrix
\begin{equation}
  \begin{aligned}
  M^{\mathrm{(loc)}}(z;x,t)=& \left\{
    \begin{array}{ll}
    T(z) M^{(\mathrm{Ai})}(t^{2/3} f(z)) h_1(z)^{\sigma_3} \begin{pmatrix} 0 & -1 \\ 1 & 0 \end{pmatrix}, ~ & z \in \mathcal{U} \cap \mathbb{C}^+, \\
    T(z) M^{(\mathrm{Ai})}(t^{2/3} f(z)) h_2(z)^{\sigma_3} , ~ & z \in \mathcal{U} \cap \mathbb{C}^-,
  \end{array} \right.
\end{aligned}
\end{equation}
where the matching factors are analytic in $\mathcal{U}$ given by
\begin{equation}
  \begin{aligned}
    T(z)=& \left\{
      \begin{array}{ll}
      e^{-\mathrm{i}\pi/12} \sqrt{\pi} M^{(\mathrm{mod})}(z;x,t) \begin{pmatrix} 0 & 1 \\ -1 & 0 \end{pmatrix} h_1(z)^{-\sigma_3} e^{\mathrm{i}\pi \sigma_3/4} \begin{pmatrix} 1 & -1 \\ 1 & 1 \end{pmatrix} (t^{2/3} f(z))^{\sigma_3/4}, ~ & z \in \mathcal{U} \cap \mathbb{C}^+, \\
      e^{-\mathrm{i}\pi/12} \sqrt{\pi} M^{(\mathrm{mod})}(z;x,t)  h_2(z)^{-\sigma_3} e^{\mathrm{i}\pi \sigma_3/4} \begin{pmatrix} 1 & -1 \\ 1 & 1 \end{pmatrix} (t^{2/3} f(z))^{\sigma_3/4}, ~ & z \in \mathcal{U} \cap \mathbb{C}^-,
    \end{array} \right.
  \end{aligned}
\end{equation}
and
\begin{equation}
  h_1(z)=\left(\frac{r^*(z)}{1-r(z)r^*(z)}\right)^{1/2} \hat{\delta}(z), \qquad h_2(z)=\left(\frac{r(z)}{1-r(z)r^*(z)}\right)^{1/2} {\hat{\delta}}^{-1}(z).
\end{equation}
By direct calculation, it is easy to verify that this local parametrix satisfies our requirements.
Then construct a global parametrix $M^{(\mathrm{par})}(z;x,t)$ by
\begin{equation}\label{Ab2para}
  \begin{aligned}
  M^{(\mathrm{par})}(z;x,t)
  &=\left\{  \begin{array}{ll} M^{\mathrm{(loc)}}(z;x,t), ~ & \mathrm{inside} ~ \mathcal{U}, \\
     M^{(\mathrm{mod})}(z;x,t), ~& \mathrm{outside} ~ \mathcal{U},
   \end{array} \right.
  \end{aligned}
\end{equation}
and define the error matrix as
$M^{(\mathrm{err})}(z;x,t)= M^{(3)}(z;x,t) {M^{(\mathrm{par})}(z;x,t)}^{-1}$.
This results in the following error RHP below.

\begin{rhp} \label{a2err}Find a $2 \times 2$ matrix-valued function $M^{(\mathrm{err})}(z;x,t)$ with the following properties:
  \begin{enumerate} [label=(\roman*)]
    \item $M^{(\mathrm{err})}(z;x,t)$ is analytic in $z \in \mathbb{C} \backslash \Sigma^{(\mathrm{err})}$, where $\Sigma^{(\mathrm{err})}=((L_1 \cup L_2 \cup  L_3 \cup L_4) \backslash\mathcal{U})\cup \partial \mathcal{U}$.
    \item $M^{(\mathrm{err})}(z;x,t)=I+\mathcal{O} (z^{-1}) $ as $z\to \infty$.
    \item $M^{(\mathrm{err})}(z;x,t)$ achieves the CBVs $M^{(\mathrm{err})}_{+}(z;x,t)$ and $M^{(\mathrm{err})}_{-}(z;x,t)$ on $\Sigma^{(\mathrm{err})}$, which satisfy the jump condition $M^{(\mathrm{err})}_{+}(z;x,t)=M^{(\mathrm{err})}_{-}(z;x,t)V^{(\mathrm{err})}(z;x,t)$, where
  \begin{equation}\label{Ab2Verr}
    \begin{aligned}
  V^{(\mathrm{err})}(z;x,t)&=\left\{  \begin{array}{ll}
    M^{(\mathrm{mod})}(z;x,t) \begin{pmatrix} 1 & 0 \\ r(z) \hat{\delta}  ^{-2}(z) e^{2\mathrm{i} tg(z)} & 1 \end{pmatrix} {M^{(\mathrm{mod})}(z;x,t)}^{-1}, ~ &z \in L_1, \\ \\
    M^{(\mathrm{mod})}(z;x,t) \begin{pmatrix} 1 & \frac{-r^*(z)}{1-r(z)r^*(z)} \hat{\delta}  ^2(z) e^{-2\mathrm{i}t g(z)} \\0 & 1 \end{pmatrix} {M^{(\mathrm{mod})}(z;x,t)}^{-1},  ~&z \in L_2 \backslash\mathcal{U},\\ \\
    M^{(\mathrm{mod})}(z;x,t) \begin{pmatrix} 1 & 0 \\ \frac{r(z)}{1-r(z)r^*(z)} \hat{\delta}  ^{-2}(z) e^{2\mathrm{i} tg(z)} & 1 \end{pmatrix} {M^{(\mathrm{mod})}(z;x,t)}^{-1},  ~&z \in L_3 \backslash\mathcal{U},\\ \\
    M^{(\mathrm{mod})}(z;x,t) \begin{pmatrix} 1 & -r^*(z) \hat{\delta}  ^2(z) e^{-2\mathrm{i} tg(z)} \\ 0  & 1 \end{pmatrix} {M^{(\mathrm{mod})}(z;x,t)}^{-1}, ~ &z \in L_4,\\ \\
    M^{(\mathrm{mod})}(z;x,t) {M^{(\mathrm{loc})}(z;x,t)}^{-1}, ~ &z \in \partial \mathcal{U}.\end{array} \right.\\
    \end{aligned}
  \end{equation}
  \end{enumerate}
\end{rhp}
Appendix \ref{appb} shows that the error estimate is $M^{(\mathrm{err})}_1(x,t)=\mathcal{O}(t^{-1})$ as $t\to \infty$.
Using equations (\ref{rec}) (with $\delta$ replaced by $\hat{\delta}$), (\ref{AB2mdeltaINF}) and (\ref{Ab2mod}), we conclude that the long-time asymptotics of the defocusing NLS equation (\ref{(NLS)})
in the dispersive shock wave region is given by
\begin{equation} \label{onephase}
  q(x,t)=\frac{\lambda_r^+ - \lambda_s(\xi) + \lambda_l^+ -\lambda_l^-}{2}\frac{\Theta (0)~\Theta (2\mathcal{A}(\infty)+(\gamma + \hat{\gamma})/{2 \pi})}{\Theta ((\gamma + \hat{\gamma})/{2 \pi})~\Theta (2\mathcal{A}(\infty))} e^{2\mathrm{i} (t g_\infty+\hat{g}_\infty-\varphi^{(0)})}+\mathcal{O}(t^{-1}),
\end{equation}
where the real quantities $\gamma$, $g_\infty$, $\varphi^{(0)}$, $\hat{\gamma}$, $\hat{g}_\infty$,  and the pure imaginary quantity $\mathcal{A}(\infty)$ are given by (\ref{gamma}), (\ref{A2GINF}), (\ref{phi0}), (\ref{gammah}), (\ref{AB2ghatINF}) and (\ref{Abelinf}), respectively.

In addition, the square modulus of the leading-order term can be expressed in terms of the Jacobi elliptic function as
\begin{equation}
  |q_{\mathrm{as}}(x,t)|^2= \rho_2 - (\rho_2- \rho_3) \mathrm{cn}^2 \left(\sqrt{\rho_1-\rho_3}\left(x- Vt + \varphi^{(1)}\left(\frac{x}{t}\right)\right) -K(m),m \right),
\end{equation}
which is in agreement with the genus-one self-similar solution (\ref{rho-equation-solution-new}) in Whitham modulation theory. See more details in Appendix \ref{CC}.

\subsubsection{The unmodulated elliptic wave region: $ v_2(\lambda_r^+, \lambda_r^-, \lambda_l^+, \lambda_l^-) < \xi < v_3(\lambda_r^+, \lambda_r^-, \lambda_l^+, \lambda_l^-) $} \label{aBmue}
\
\newline
\indent

\begin{figure}[htbp]
  \centering
  \subfigure[]{\includegraphics[width=5.5cm]{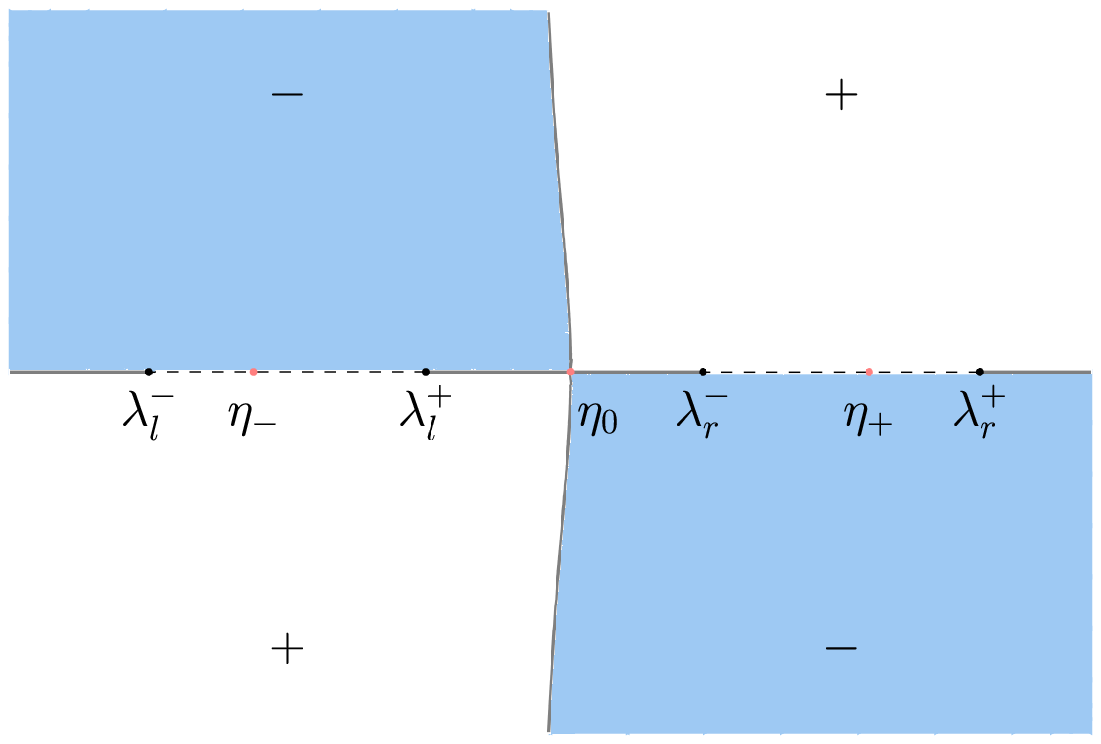}}
  \subfigure[]{\includegraphics[width=5cm]{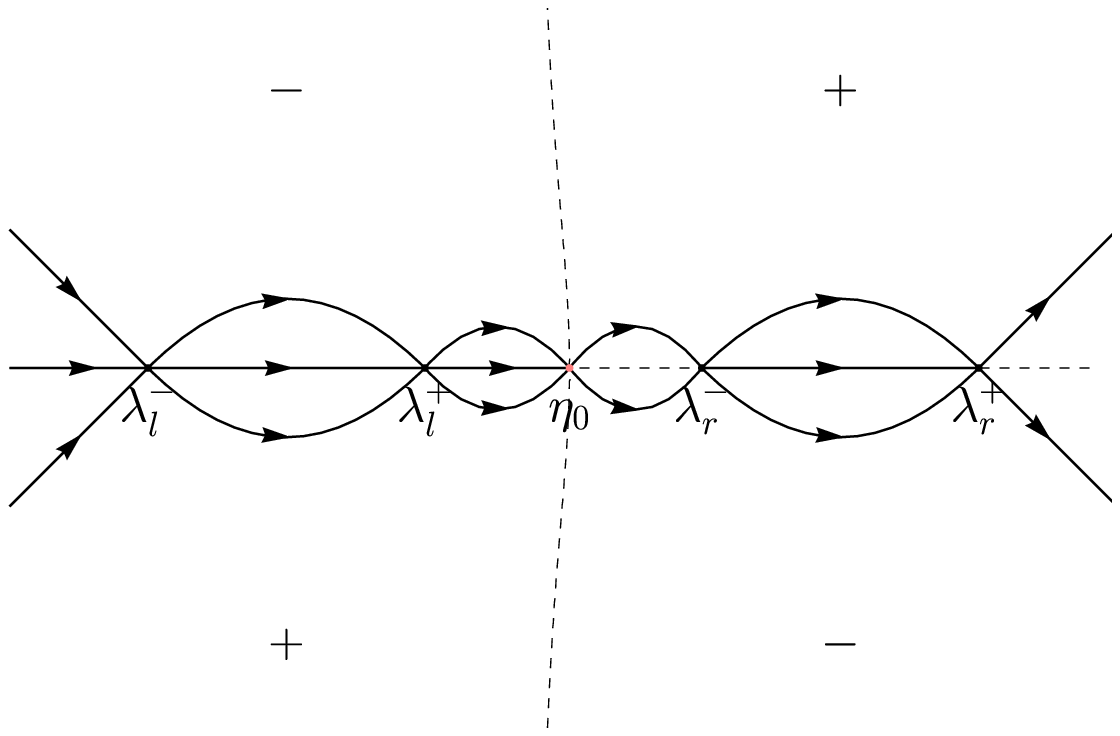}}
  \subfigure[]{\includegraphics[width=5.5cm]{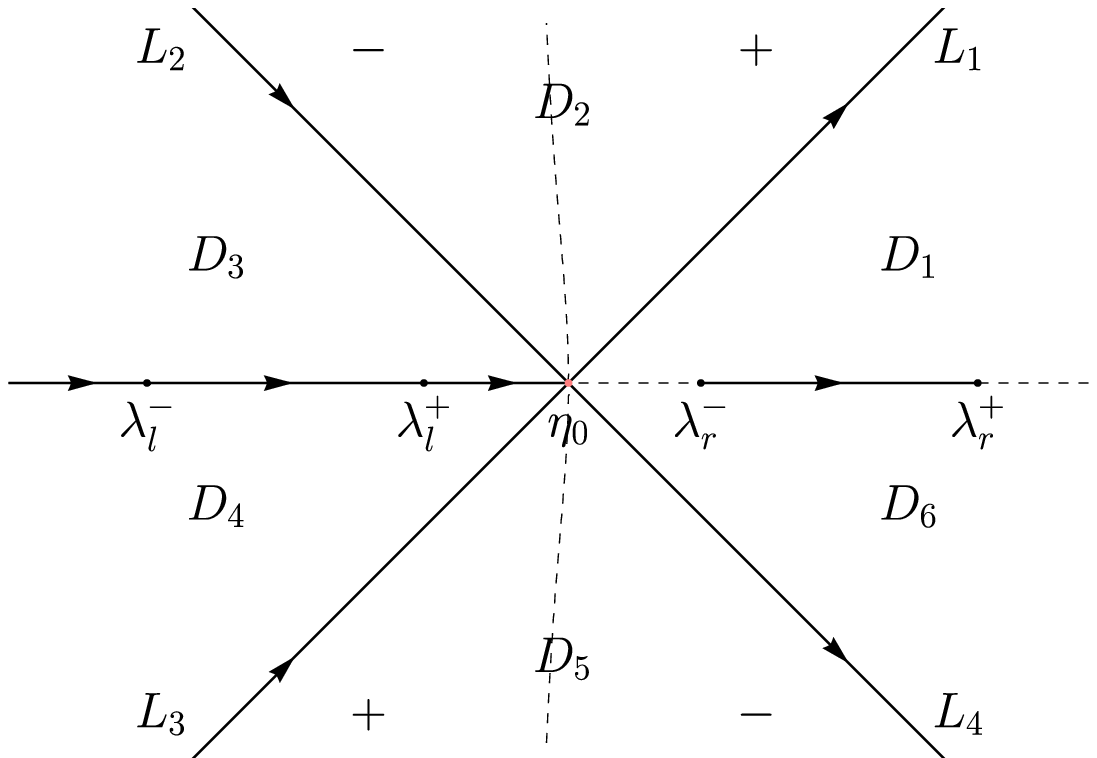}}
\caption{{\protect\small (a) The sign structure of Im$(g(z))$; (b) Opening lenses; (c) The jump contours of $M^{(2)}(z;x,t)$.}}
\label{figA3}
\end{figure}

If $\xi > v_2(\lambda_r^+, \lambda_r^-, \lambda_l^+, \lambda_l^-)$, the stationary phase point $\lambda_s(\xi)$ of the shock $g$-function (\ref{A2G}) is less than $\lambda_r^-$,
which makes exponentially large entries $e^{\pm 2 \mathrm{i}t g_+(z)}$ appear in the jump matrix (\ref{AB2V2}) on the interval $(\lambda_s(\xi), \lambda_r^-)$. Therefore, we need to modify the $g$-function so that $g(z)$ is analytic in $z\in \mathbb{C} \backslash ([\lambda_l^-, \lambda_l^+] \cup [\lambda_r^-, \lambda_r^+])$.
The construction of the $g(z)$ function is similar to that in Section \ref{ALSH}, except that the branch points are replaced by $\boldsymbol{\lambda}=(\lambda_r^+, \lambda_r^-, \lambda_l^+, \lambda_l^-)$ and the two bands are fixed. The new two-band $g$-function is given by
\begin{equation}\label{aB3g}
  g(z)= \int_{\lambda_r^+}^{z} \frac{2(\zeta-\eta_+(\xi))(\zeta-\eta_0(\xi))(\zeta-\eta_-(\xi))}{\mathcal{R}(\zeta ;\lambda_r^+, \lambda_r^-, \lambda_l^+, \lambda_l^-)}\dif \zeta,
\end{equation}
where $\eta_0(\xi)$, $\eta_+(\xi)$ and $\eta_-(\xi)$ are zeros of the polynomial $2 P_1 (z ;\boldsymbol{\lambda})+ \xi P_0 (z ;\boldsymbol{\lambda})$ defined by (\ref{polynomial}).
Along the bands, $g(z)$ satisfies the jump conditions
\begin{equation} \label{gamma1}
  \begin{aligned}
\left\{  \begin{array}{ll}
   g _{+}(z)+g _{-}(z)=0,  ~& z\in (\lambda_r^-,\lambda_r^+), \\ \\
  g _{+}(z)+g _{-}(z)= \gamma :=\oint_b \dif g ,  ~& z\in (\lambda_l^-,\lambda_l^+).
  \end{array} \right.
  \end{aligned}
\end{equation}
Note that $\eta_-(\xi) \in (\lambda_l^-,\lambda_l^+)$ and $\eta_+(\xi) \in (\lambda_r^-,\lambda_r^+)$.
Furthermore,
\begin{equation} \label{A3GINF}
    g_\infty= \int_{\lambda_r^+}^{\infty} \left(\frac{2(\zeta-\eta_+(\xi))(\zeta-\eta_0(\xi))(\zeta-\eta_-(\xi))}{\mathcal{R}(\zeta ;\lambda_r^+, \lambda_r^-, \lambda_l^+, \lambda_l^-)}-2\zeta-\xi\, \right)\dif \zeta - \theta(\lambda_r^+).
\end{equation}
With the help of Whitham modulation theory, the boundaries of this region are determined by Whitham velocities $v_2(\boldsymbol{\lambda})$ and $v_3(\boldsymbol{\lambda})$ as follows:
\begin{equation}\label{aB3b}
  \lambda_l^+ < \eta_0(\xi) < \lambda_r^- \qquad  \mathrm{iff}  \qquad  v_2(\lambda_r^+, \lambda_r^-, \lambda_l^+, \lambda_l^-)= \xi_{A2}  < \xi < \xi_{A3} = v_3(\lambda_r^+, \lambda_r^-, \lambda_l^+, \lambda_l^-).
\end{equation}
Then, open lenses from intervals to the steepest descent contours through $\eta_0(\xi)$  (as shown in Figure \ref{figA3})  by the transformations
$M(z;x,t) \mapsto M^{(1)}(z;x,t) \mapsto M^{(2)}(z;x,t)$, which yields the following RHP.
\begin{rhp}
  Find a $2 \times 2$ matrix-valued function $M^{(2)}(z;x,t)$ with the following properties:
\begin{enumerate} [label=(\roman*)]
  \item  $M^{(2)}(z;x,t)$ is analytic in $z \in \mathbb{C} \backslash \Sigma^{(2)}$, where $\Sigma^{(2)}=(-\infty, \eta_0(\xi)] \cup [\lambda_r^-, \lambda_r^+] \cup L_1 \cup L_2 \cup  L_3 \cup L_4.$
  \item  $M^{(2)}(z;x,t)=I+\mathcal{O} (z^{-1}) $ as $z\to \infty$.
  \item  $M^{(2)}(z;x,t)$ achieves the CBVs  $M^{(2)}_{+}(z;x,t)$ and $M^{(2)}_{-}(z;x,t)$  on $\Sigma^{(2)}$ away from self-intersection point and
  branch points that satisfy the jump condition $M^{(2)}_{+}(z;x,t)=M^{(2)}_{-}(z;x,t)V^{(2)}(z;x,t)$, where
  \begin{equation}\label{aB3V2}
    \begin{aligned}
    V^{(2)}(z;x,t)&=\left\{  \begin{array}{ll}
     (1-r(z)r^{*}(z))^{\sigma_{3}},  ~& z\in (-\infty, \lambda_l^-) \cup (\lambda_l^+, \eta_0(\xi)), \\ \\
     \begin{pmatrix} 0 & -r_{-}^{*}(z) e^{-\mathrm{i} \gamma}  \\ r_{+}(z) e^{\mathrm{i} \gamma}  & 0 \end{pmatrix}, ~&z \in (\lambda_l^-, \lambda_l^+), \\ \\
     \begin{pmatrix} 0 & -1 \\1 & 0 \end{pmatrix}, ~&z \in (\lambda_r^-, \lambda_r^+), \\ \\
     \begin{pmatrix} 1 & 0 \\ r(z)e^{2\mathrm{i} tg(z)} & 1 \end{pmatrix}, ~ &z \in L_1, \\ \\
     \begin{pmatrix} 1 & \frac{-r^*(z)}{1-r(z)r^*(z)}e^{-2\mathrm{i} tg(z)} \\0 & 1 \end{pmatrix},  ~&z \in L_2,\\ \\
     \begin{pmatrix} 1 & 0 \\ \frac{r(z)}{1-r(z)r^*(z)}e^{2\mathrm{i} tg(z)} & 1 \end{pmatrix},  ~&z \in L_3,\\ \\
     \begin{pmatrix} 1 & -r^*(z)e^{-2\mathrm{i} tg(z)} \\ 0  & 1 \end{pmatrix}, ~ &z \in L_4.\end{array} \right.\\
    \end{aligned}
  \end{equation}
\end{enumerate}
\end{rhp}

To arrive at the solvable limiting RHP with piecewise constant jumps across two fixed bands,  introduce a modified scalar function
\begin{equation}\label{aB3mdelta}
  \hat{\delta} (z) = \delta  (z) e^{- \mathrm{i} \hat{g} (z) },
\end{equation}
where
\begin{equation} \label{aB3delta}
  \delta(z)= \exp \left\{  \frac{\mathcal{R} (z;\boldsymbol{\lambda})}{2 \pi \mathrm{i}}  \left[ \left(\int_{-\infty}^{\lambda_l^-} + \int_{\lambda_l^+}^{\eta_0(\xi)} \right) \frac{\ln (1-|r(\zeta)|^2)}{\mathcal{R} (\zeta;\boldsymbol{\lambda}) } \,\frac{\dif\zeta}{\zeta - z}
   + \int_{\lambda_l^-}^{\lambda_l^+} \frac{\ln(r_+(\zeta))}{\mathcal{R}_+ (\zeta;\boldsymbol{\lambda})} \,\frac{\dif\zeta}{\zeta - z}
   \right]  \right\},
\end{equation}
and
\begin{equation}
  \begin{aligned}
    &\hat{g} (z)  = \varphi^{(1)} \int_{\lambda_r^+}^{z} \frac{P_0 (\zeta ;\boldsymbol{\lambda})}{\mathcal{R}(\zeta ;\boldsymbol{\lambda})}\,\dif \zeta,\\
    &\varphi^{(1)}=\frac{1}{2 \pi} \left[ \left(\int_{-\infty}^{\lambda_l^-} + \int_{\lambda_l^+}^{\eta_0(\xi)} \right) \frac{\ln (1-|r(\zeta)|^2) }{\mathcal{R} (\zeta;\boldsymbol{\lambda}) } \dif \zeta
    + \int_{\lambda_l^-}^{\lambda_l^+} \frac{\ln(r_+(\zeta))}{\mathcal{R}_+ (\zeta;\boldsymbol{\lambda})} \dif \zeta \right]. \\  \end{aligned}
\end{equation}
The function $\hat{\delta} (z)$ has the following properties:
\begin{enumerate} [label=(\roman*)]
  \item $\hat{\delta}(z)$ is analytic in $z \in \mathbb{C} \backslash ((-\infty, \eta_0(\xi)] \cup [\lambda_r^-, \lambda_r^+])$;
  \item $\hat{\delta}(z)=\hat{\delta}_\infty+\mathcal{O} (z^{-1})$ as $z\to \infty$, where
    \begin{equation}\label{aB3mdeltaINF}
     \begin{aligned}
     &\hat{\delta}_{\infty }= \exp \{\mathrm{i}(\varphi^{(0)} -\hat{g}_\infty ) \},\\
     &\varphi^{(0)}=\frac{1}{2 \pi} \left[ \left(\int_{-\infty}^{\lambda_l^-} + \int_{\lambda_l^+}^{\eta_0(\xi)} \right) \frac{\ln (1-|r(\zeta)|^2) (\zeta+ V) }{\mathcal{R} (\zeta;\boldsymbol{\lambda}) } \dif \zeta
     +  \int_{\lambda_l^-}^{\lambda_l^+} \frac{\ln(r_+(\zeta)) (\zeta+ V) }{\mathcal{R}_+(\zeta;\boldsymbol{\lambda})} \dif \zeta \right],\\
     &\hat{g}_\infty = \varphi^{(1)} \int_{\lambda_r^+}^{\infty} \left(\frac{P_0 (\zeta ;\boldsymbol{\lambda})}{\mathcal{R}(\zeta ;\boldsymbol{\lambda})} -1 \right)\,\dif \zeta - \varphi^{(1)} \lambda_r^+.
    \end{aligned}
    \end{equation}
  \item $\hat{\delta}(z)$ achieves CBVs on  $(-\infty, \eta_0(\xi)) \cup (\lambda_r^-, \lambda_r^+)$ satisfying the jump conditions
  \begin{equation*}
    \begin{aligned}
  \left\{  \begin{array}{ll}
     \hat{\delta} _{+}(z)=(1-r(z)r^{*}(z))\hat{\delta} _{-}(z),  ~& z\in (-\infty, \lambda_l^-) \cup (\lambda_l^+, \eta_0(\xi)) , \\ \\
     \hat{\delta} _{+}(z)\hat{\delta} _{-}(z)=r_{+}(z) e^{-\mathrm{i} \hat{\gamma}} ,  ~& z\in (\lambda_l^-, \lambda_l^+),\\ \\
     \hat{\delta} _{+}(z)\hat{\delta} _{-}(z)=1,  ~& z\in (\lambda_r^-, \lambda_r^+),
    \end{array} \right.\\
    \end{aligned}
  \end{equation*}
  where
  \begin{equation} \label{gammah1}
    \hat{\gamma} :=\oint_b \, \dif\hat{g}.
  \end{equation}

\end{enumerate}
Then, we can define the modified transformation
\begin{equation}
	M^{(3)}(z;x,t)=\hat{\delta} _{\infty}^{\sigma_{3}} M^{(2)}(z;x,t) \hat{\delta}(z) ^{-\sigma_{3}}
\end{equation}
which leads to the following RHP.
\begin{rhp} Find a $2 \times 2$ matrix-valued function $M^{(3)}(z;x,t)$ with the following properties:
  \begin{enumerate} [label=(\roman*)]
    \item $M^{(3)}(z;x,t)$ is analytic in $z \in \mathbb{C} \backslash \Sigma^{(3)}$, where $\Sigma^{(3)}= [\lambda_l^-,\lambda_l^+] \cup [\lambda_r^-,\lambda_r^+]\cup L_1 \cup L_2 \cup  L_3 \cup L_4.$
    \item $M^{(3)}(z;x,t)=I+\mathcal{O} (z^{-1}) $ as $z\to \infty$.
    \item $M^{(3)}(z;x,t)$ achieves CBVs $M^{(3)}_{+}(z;x,t)$ and $M^{(3)}_{-}(z;x,t)$ on $\Sigma^{(3)}$ away from the self-intersection point and
  branch points that satisfy the jump condition $M^{(3)}_{+}(z;x,t)=M^{(3)}_{-}(z;x,t)V^{(3)}(z;x,t)$, where
  \begin{equation}\label{aB3V3}
    \begin{aligned}
  V^{(3)}(z;x,t)&=\left\{  \begin{array}{ll}
     \begin{pmatrix} 0 & - e^{-\mathrm{i} (\gamma+\hat{\gamma})}  \\  e^{\mathrm{i} (\gamma+\hat{\gamma})}  & 0 \end{pmatrix}, ~&z \in (\lambda_l^-, \lambda_l^+), \\ \\
     \begin{pmatrix} 0 & -1 \\1 & 0 \end{pmatrix}, ~&z \in (\lambda_r^-, \lambda_r^+), \\ \\
     \begin{pmatrix} 1 & 0 \\ r(z) \hat{\delta} ^{-2}(z) e^{2\mathrm{i} tg(z)} & 1 \end{pmatrix}, ~ &z \in L_1, \\ \\
     \begin{pmatrix} 1 & \frac{-r^*(z)}{1-r(z)r^*(z)} \hat{\delta} ^2(z) e^{-2\mathrm{i}t g(z)} \\0 & 1 \end{pmatrix},  ~&z \in L_2,\\ \\
     \begin{pmatrix} 1 & 0 \\ \frac{r(z)}{1-r(z)r^*(z)} \hat{\delta} ^{-2}(z) e^{2\mathrm{i} tg(z)} & 1 \end{pmatrix},  ~&z \in L_3,\\ \\
     \begin{pmatrix} 1 & -r^*(z) \hat{\delta} ^2(z) e^{-2\mathrm{i} tg(z)} \\ 0  & 1 \end{pmatrix}, ~ &z \in L_4. \end{array} \right.\\
    \end{aligned}
  \end{equation}
  \end{enumerate}
\end{rhp}

Outside a small  neighborhood $\mathcal{U} $ of $\eta_0(\xi)$, the jump matrices $V^{(3)}(z;x,t)$ on the steepest descent contours converge
uniformly to the identity matrix. So we consider the limiting problem
\begin{equation}
  M^{(\mathrm{mod})}_+(z;x,t)=M^{(\mathrm{mod})}_-(z;x,t) V^{(\mathrm{mod})},
\end{equation}
where
\begin{equation}\label{aB3vmod}
  \begin{aligned}
    V^{(\mathrm{mod})}&=\left\{  \begin{array}{ll}
   \begin{pmatrix} 0 & - e^{-\mathrm{i} (\gamma+\hat{\gamma})}  \\  e^{\mathrm{i} (\gamma+\hat{\gamma})}  & 0 \end{pmatrix}, ~&z \in (\lambda_l^-, \lambda_l^+), \\ \\
   \begin{pmatrix} 0 & -1 \\1 & 0 \end{pmatrix}, ~&z \in (\lambda_r^-, \lambda_r^+). \end{array} \right.\\
   \end{aligned}
\end{equation}
The solution is  given by
\begin{equation} \label{AB3mod}
  M^{(\mathrm{mod})}(z;x,t)=N^{-1}(\infty; \boldsymbol{\lambda}) N(z; \boldsymbol{\lambda}),
\end{equation}
where $N(z; \boldsymbol{\lambda})$ is defined by (\ref{N}) with $\boldsymbol{\lambda}=(\lambda_r^+, \lambda_r^-, \lambda_l^+, \lambda_l^-)$.

Inside  $\mathcal{U}$, the jump matrices $V^{(3)}(z;x,t)$ cannot  uniformly converge  to the identity matrix due to the quadratic vanishment
of the phase $g(z)$. The construction of the parametrix proceeds as in the left plane wave region with $\eta_+(\xi)$  and $\delta(z)$ replaced by $\eta_0(\xi)$ and $\hat{\delta}(z)$, respectively.
Then we obtain the same error RHP as RHP \ref{eRHP} and thus the same error estimate $M^{(\mathrm{err})}_1(x,t)= \mathcal{O}(t^{-1/2})$.
Using the equations (\ref{rec}) (with $\delta$ replaced by $\hat{\delta}$), (\ref{aB3mdeltaINF}) and (\ref{AB3mod}), it is concluded that the long-time asymptotic behavior
of the defocusing NLS equation (\ref{(NLS)}) in the unmodulated elliptic wave region is given by
\begin{equation}
  q(x,t)=(A_l+A_r)\frac{\Theta (0)~\Theta (2\mathcal{A}(\infty)+(\gamma + \hat{\gamma})/{2 \pi})}{\Theta ((\gamma + \hat{\gamma})/{2 \pi})~\Theta (2\mathcal{A}(\infty))}e^{2\mathrm{i} (t g_\infty+\hat{g}_\infty-\varphi^{(0)})}+\mathcal{O}(t^{-1/2}),
\end{equation}
where the real quantities $\gamma$, $g_\infty$, $\hat{g}_\infty$, $\varphi^{(0)}$, $\hat{\gamma}$ and the pure imaginary quantity $\mathcal{A}(\infty)$ are given by (\ref{gamma1}), (\ref{A3GINF}), (\ref{aB3mdeltaINF}), (\ref{aB3mdeltaINF}), (\ref{gammah1}) and (\ref{Abelinf}), respectively.

\subsubsection{Dispersive shock wave region: $ v_3(\lambda_r^+, \lambda_r^-, \lambda_l^+, \lambda_l^-) < \xi < -\frac{\lambda_r^+ + \lambda_r^- + 2\lambda_l^-}{2}+\frac{2(\lambda_r^+ - \lambda_l^-)
(\lambda_r^- - \lambda_l^-)}{\lambda_r^+ + \lambda_r^- -2\lambda_l^-} $} \label{ABrsh}
\
\newline
\indent

\begin{figure}[htbp]
  \centering
  \subfigure[]{\includegraphics[width=5.5cm]{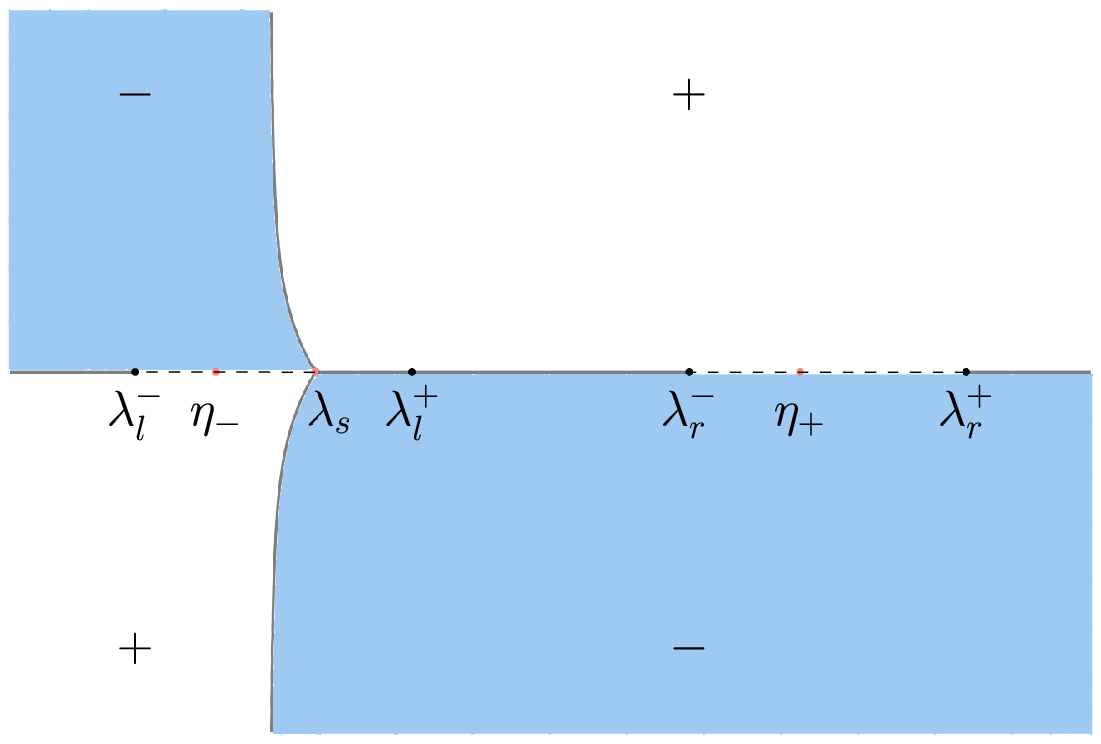}}
  \subfigure[]{\includegraphics[width=5cm]{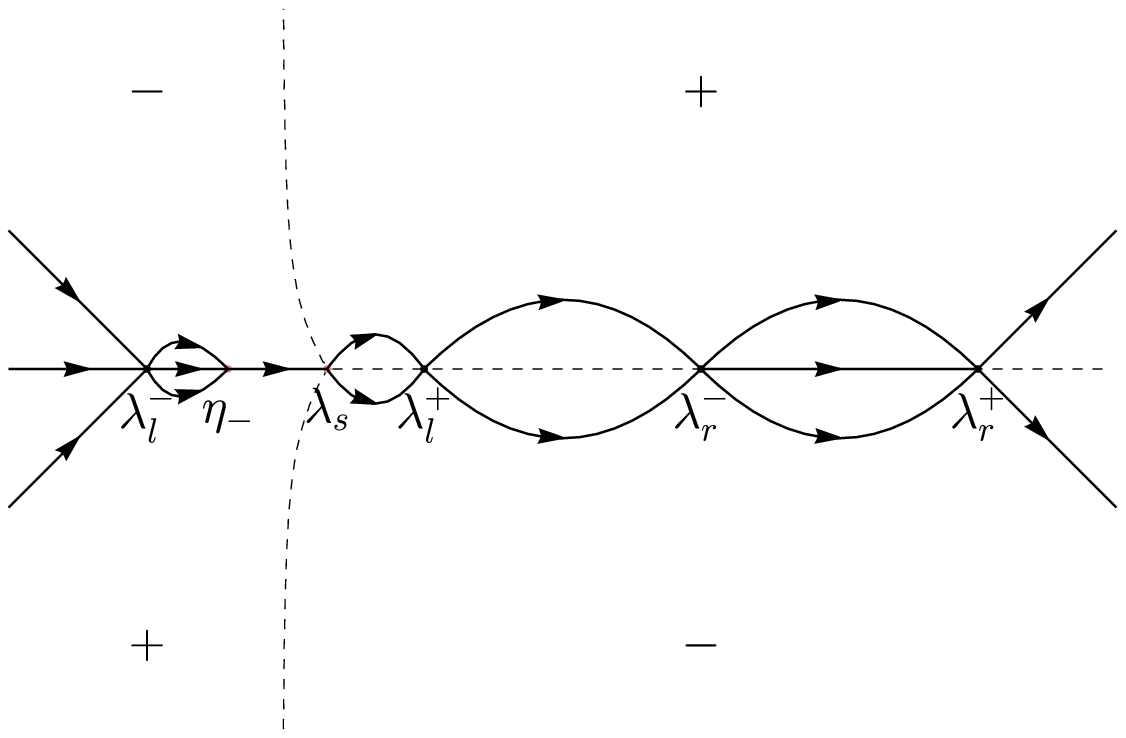}}
  \subfigure[]{\includegraphics[width=5.5cm]{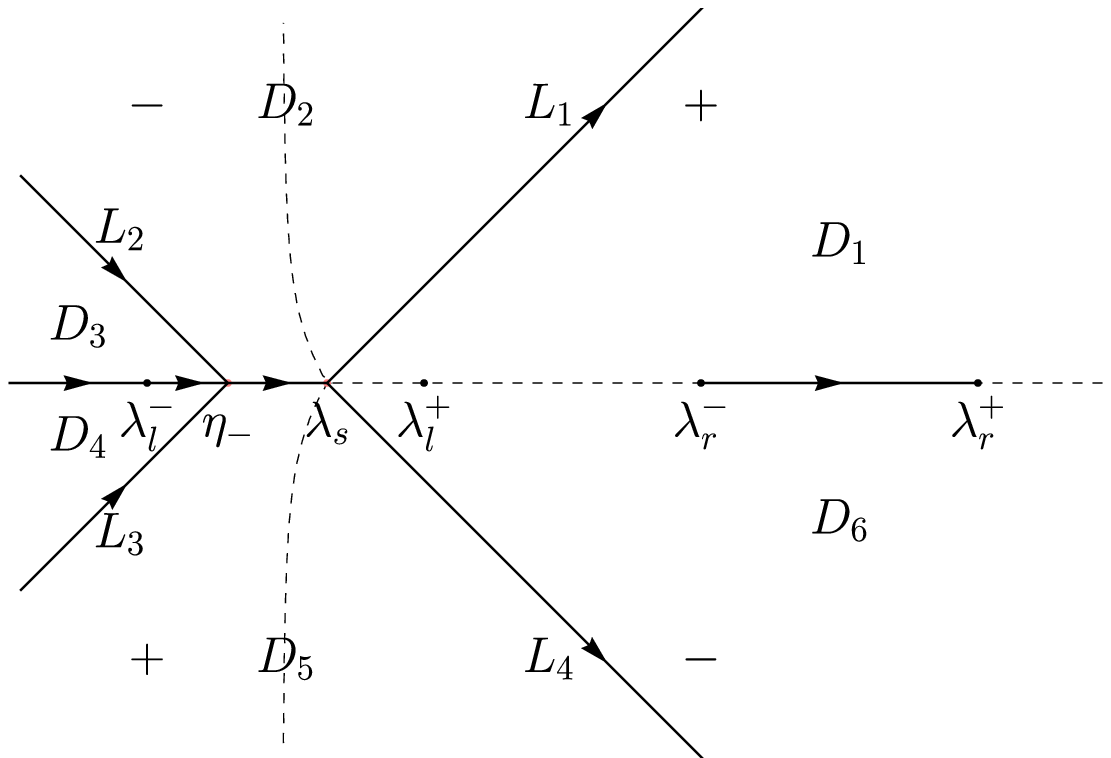}}
\caption{{\protect\small (a) The sign structure of Im$(g(z))$; (b) Opening lenses; (c) The jump contours of $M^{(2)}(z;x,t)$.}}
\label{figA4}
\end{figure}

If $\xi > v_3(\lambda_r^+, \lambda_r^-, \lambda_l^+, \lambda_l^-)$, the stationary phase point $\eta_0(\xi)$ of the unmodulated $g$-function (\ref{aB3g}) is less than $\lambda_l^+$. This leads to the emergence of the interval $(\eta_0(\xi), \lambda_l^+)$, on which the jump matrix (\ref{aB3V2}) grows exponentially with respect to  $t$.
To remove the exponentially large entries, we modify the $g$-function.
The construction of the shock $g$-function is similar to that in Section \ref{ALSH}, except that the branch points are replaced by $\boldsymbol{\lambda}=(\lambda_r^+, \lambda_r^-, \lambda_s(\xi), \lambda_l^-)$.
Here, the soft edge $\lambda_s(\xi)$ is uniquely determined by the following implicit function (solvable for $\lambda_s(\xi) \in ( \lambda_l^-, \lambda_r^- )$):
\begin{equation}\label{vv3}
  \xi = v_3(\boldsymbol{\lambda}) := -\frac{\lambda_r^+ +  \lambda_r^- + \lambda_s(\xi) + \lambda_l^-}{2}- \frac{ \lambda_s(\xi) -\lambda_l^-   }{1- \frac{\lambda_r^- -\lambda_l^-}{\lambda_r^- - \lambda_s(\xi)} \frac{E(m(\xi))}{K(m(\xi))}}
\end{equation}
with the elliptic modulus
\begin{equation}
  m(\xi)= \sqrt{\frac{(\lambda_r^+ -\lambda_r^- )(\lambda_s(\xi) -\lambda_l^-)}{(\lambda_r^+ - \lambda_s(\xi))(\lambda_r^- -\lambda_l^-)}}.
\end{equation}
The new two-band $g$-function is given by
\begin{equation}\label{AB4g}
  g(z)= \int_{\lambda_r^+}^{z} \frac{2(\zeta-\eta_+(\xi))(\zeta-\lambda_s(\xi))(\zeta-\eta_-(\xi))}{\mathcal{R}(\zeta ;\lambda_r^+, \lambda_r^-, \lambda_s(\xi), \lambda_l^-)}\dif \zeta,
\end{equation}
where $\eta_+(\xi) \in (\lambda_r^-, \lambda_r^+)$ and $\eta_-(\xi) \in  (\lambda_l^-, \lambda_s(\xi))$ are zeros of the polynomial $2 P_1 (z ;\boldsymbol{\lambda})+ \xi P_0 (z ;\boldsymbol{\lambda})$ defined by (\ref{polynomial}).
Along the bands, $g(z)$ satisfies the jump conditions
\begin{equation} \label{gamma2}
  \begin{aligned}
\left\{  \begin{array}{ll}
   g _{+}(z)+g _{-}(z)=0,  ~& z\in (\lambda_r^-,\lambda_r^+), \\ \\
  g _{+}(z)+g _{-}(z)= \gamma :=\oint_b \dif g ,  ~& z\in (\lambda_l^-,\lambda_s(\xi)).
  \end{array} \right.
  \end{aligned}
\end{equation}
Furthermore, we have
\begin{equation} \label{A4GINF}
    g_\infty= \int_{\lambda_r^+}^{\infty} \left(\frac{2(\zeta-\eta_+(\xi))(\zeta-\lambda_s(\xi))(\zeta-\eta_-(\xi))}{\mathcal{R}(\zeta ;\lambda_r^+, \lambda_r^-, \lambda_s(\xi), \lambda_l^-)}-2\zeta-\xi\, \right)\,\dif \zeta - \theta(\lambda_r^+).
\end{equation}
The boundaries of this region are characterized by the degeneration
\begin{equation}\label{AB4b}
  \begin{aligned}
  &\xi \to \xi_{A3} = v_3(\lambda_r^+, \lambda_r^-, \lambda_l^+, \lambda_l^-), \quad &\mathrm{as} \quad \lambda_s(\xi) \to \lambda_l^+,\\
  &\xi \to \xi_{A4} = v_3(\lambda_r^+, \lambda_r^-, \lambda_l^-, \lambda_l^-)= -\frac{\lambda_r^+ + \lambda_r^- + 2\lambda_l^-}{2}+\frac{2(\lambda_r^+ - \lambda_l^-)
  (\lambda_r^- - \lambda_l^-)}{\lambda_r^+ + \lambda_r^- -2\lambda_l^-}, \quad &\mathrm{as} \quad \lambda_s(\xi) \to \lambda_l^-.
  \end{aligned}
\end{equation}
Then open the steepest descent contours through $\lambda_s(\xi)$ and $\eta_-(\xi)$  (as shown in Figure \ref{figA4}) by the transformations $M(z;x,t) \mapsto M^{(1)}(z;x,t) \mapsto M^{(2)}(z;x,t)$, which yields the following RHP.

\begin{rhp} \label{rsrhp}
  Find a $2 \times 2$ matrix-valued function $M^{(2)}(z;x,t)$ with the following properties:
\begin{enumerate} [label=(\roman*)]
  \item  $M^{(2)}(z;x,t)$ is analytic in $z \in \mathbb{C} \backslash \Sigma^{(2)}$, where $\Sigma^{(2)}=(-\infty, \lambda_s(\xi)]  \cup [\lambda_r^-, \lambda_r^+] \cup L_1 \cup L_2 \cup  L_3 \cup L_4.$
  \item  $M^{(2)}(z;x,t)=I+\mathcal{O} (z^{-1}) $ as $z\to \infty$.
  \item  $M^{(2)}(z;x,t)$ achieves the CBVs  $M^{(2)}_{+}(z;x,t)$ and $M^{(2)}_{-}(z;x,t)$  on $\Sigma^{(2)}$ away from self-intersection points and
  branch points that satisfy the jump condition $M^{(2)}_{+}(z;x,t)=M^{(2)}_{-}(z;x,t)V^{(2)}(z;x,t)$, where
  \begin{equation}\label{AB4V2}
    \begin{aligned}
    V^{(2)}(z;x,t)&=\left\{  \begin{array}{ll}
     (1-r(z)r^{*}(z))^{\sigma_{3}},  ~& z\in (-\infty, \lambda_l^-), \\ \\
     \begin{pmatrix} 0 & -r_{-}^{*}(z) e^{-\mathrm{i} \gamma}  \\ r_{+}(z) e^{\mathrm{i} \gamma}  & 0 \end{pmatrix}, ~&z \in (\lambda_l^-, \eta_-(\xi)), \\ \\
     \begin{pmatrix} 0 & -r_{-}^{*}(z) e^{-\mathrm{i} \gamma}  \\ r_{+}(z) e^{\mathrm{i} \gamma} & e^{ -2 \mathrm{i} t g_+(z)} e^{\mathrm{i} \gamma} \end{pmatrix}, ~&z \in (\eta_-(\xi), \lambda_s(\xi)), \\ \\
     \begin{pmatrix} 0 & -1 \\1 & 0 \end{pmatrix}, ~&z \in (\lambda_r^-, \lambda_r^+), \\ \\
     \begin{pmatrix} 1 & 0 \\ r(z)e^{2\mathrm{i} tg(z)} & 1 \end{pmatrix}, ~ &z \in L_1, \\ \\
     \begin{pmatrix} 1 & \frac{-r^*(z)}{1-r(z)r^*(z)}e^{-2\mathrm{i} tg(z)} \\0 & 1 \end{pmatrix},  ~&z \in L_2,\\ \\
     \begin{pmatrix} 1 & 0 \\ \frac{r(z)}{1-r(z)r^*(z)}e^{2\mathrm{i} tg(z)} & 1 \end{pmatrix},  ~&z \in L_3,\\ \\
     \begin{pmatrix} 1 & -r^*(z)e^{-2\mathrm{i} tg(z)} \\ 0  & 1 \end{pmatrix}, ~ &z \in L_4.\end{array} \right.\\
    \end{aligned}
  \end{equation}
\end{enumerate}
\end{rhp}

To arrive at the solvable limiting RHP with piecewise constant jumps across one fixed band and one soft band,  introduce a modified scalar function
\begin{equation}\label{AB4mdelta}
  \hat{\delta} (z) = \delta  (z) e^{- \mathrm{i} \hat{g} (z) },
\end{equation}
where
\begin{equation} \label{AB4delta}
  \delta(z)= \exp \left\{  \frac{\mathcal{R} (z;\boldsymbol{\lambda})}{2 \pi \mathrm{i}}  \left[ \int_{-\infty}^{\lambda_l^-}  \frac{\ln (1-|r(\zeta)|^2)}{\mathcal{R} (\zeta;\boldsymbol{\lambda}) } \,\frac{\dif\zeta}{\zeta - z}
   + \int_{\lambda_l^-}^{\lambda_s(\xi)} \frac{\ln(r_+(\zeta))}{\mathcal{R}_+ (\zeta;\boldsymbol{\lambda})} \,\frac{\dif\zeta}{\zeta - z}
   \right]  \right\},
\end{equation}
and
\begin{equation}
  \begin{aligned}
    &\hat{g} (z)  = \varphi^{(1)} \int_{\lambda_r^+}^{z} \frac{P_0 (\zeta ;\boldsymbol{\lambda})}{\mathcal{R}(\zeta ;\boldsymbol{\lambda})}\,\dif \zeta,\\
    &\varphi^{(1)}=\frac{1}{2 \pi} \left[\int_{-\infty}^{\lambda_l^-}  \frac{\ln (1-|r(\zeta)|^2) }{\mathcal{R} (\zeta;\boldsymbol{\lambda}) } \dif \zeta
    +\int_{\lambda_l^-}^{\lambda_s(\xi)} \frac{\ln(r_+(\zeta))}{\mathcal{R}_+ (\zeta;\boldsymbol{\lambda})} \dif \zeta \right]. \\
   \end{aligned}
\end{equation}
The function $\hat{\delta} (z)$ has the following properties:
\begin{enumerate} [label=(\roman*)]
  \item $\hat{\delta}(z)$ is analytic in $z \in \mathbb{C} \backslash ((-\infty, \lambda_s(\xi)] \cup [\lambda_r^-, \lambda_r^+])$;
  \item $\hat{\delta}(z)=\hat{\delta}_\infty+\mathcal{O} (z^{-1})$ as $z\to \infty$, where
    \begin{equation}\label{AB4mdeltaINF}
     \begin{aligned}
     &\hat{\delta}_{\infty }= \exp \{\mathrm{i}(\varphi^{(0)} -\hat{g}_\infty ) \},\\
     &\varphi^{(0)}=\frac{1}{2 \pi} \left[\int_{-\infty}^{\lambda_l^-}  \frac{\ln (1-|r(\zeta)|^2) (\zeta+ V) }{\mathcal{R} (\zeta;\boldsymbol{\lambda}) } \dif \zeta
     +\int_{\lambda_l^-}^{\lambda_s(\xi)} \frac{\ln(r_+(\zeta)) (\zeta+ V) }{\mathcal{R}_+(\zeta;\boldsymbol{\lambda})} \dif \zeta \right],\\
     &\hat{g}_\infty = \varphi^{(1)}\left[ \int_{\lambda_r^+}^{\infty} \left(\frac{P_0 (\zeta ;\boldsymbol{\lambda})}{\mathcal{R}(\zeta ;\boldsymbol{\lambda})} -1 \right)\,\dif \zeta -  \lambda_r^+ \right].
    \end{aligned}
    \end{equation}
  \item $\hat{\delta}(z)$ achieves the CBVs on  $(-\infty, \lambda_s(\xi)) \cup (\lambda_r^-, \lambda_r^+)$ satisfying the jump conditions
  \begin{equation*}
    \begin{aligned}
  \left\{  \begin{array}{ll}
     \hat{\delta} _{+}(z)=(1-r(z)r^{*}(z))\hat{\delta} _{-}(z),  ~& z\in (-\infty, \lambda_l^-) , \\ \\
     \hat{\delta} _{+}(z)\hat{\delta} _{-}(z)=r_{+}(z) e^{-\mathrm{i} \hat{\gamma}} ,  ~& z\in (\lambda_l^-, \lambda_s(\xi)),\\ \\
     \hat{\delta} _{+}(z)\hat{\delta} _{-}(z)=1,  ~& z\in (\lambda_r^-, \lambda_r^+),
    \end{array} \right.\\
    \end{aligned}
  \end{equation*}
  where
  \begin{equation} \label{gammah2}
    \hat{\gamma} :=\oint_b \, \dif\hat{g}.
  \end{equation}
\end{enumerate}
Then define the modified transformation
\begin{equation}
	M^{(3)}(z;x,t)=\hat{\delta} _{\infty}^{\sigma_{3}} M^{(2)}(z;x,t) \hat{\delta}(z) ^{-\sigma_{3}}
\end{equation}
which leads to the following RHP.
\begin{rhp} Find a $2 \times 2$ matrix-valued function $M^{(3)}(z;x,t)$ with the following properties:
  \begin{enumerate} [label=(\roman*)]
    \item $M^{(3)}(z;x,t)$ is analytic in $z \in \mathbb{C} \backslash \Sigma^{(3)}$, where $\Sigma^{(3)}= [\lambda_l^-,\lambda_s(\xi)] \cup [\lambda_r^-,\lambda_r^+]\cup L_1 \cup L_2 \cup  L_3 \cup L_4.$
    \item $M^{(3)}(z;x,t)=I+\mathcal{O} (z^{-1}) $ as $z\to \infty$.
    \item $M^{(3)}(z;x,t)$ achieves the CBVs $M^{(3)}_{+}(z;x,t)$ and $M^{(3)}_{-}(z;x,t)$ on $\Sigma^{(3)}$ away from the self-intersection point and
  branch points that satisfy the jump condition $M^{(3)}_{+}(z;x,t)=M^{(3)}_{-}(z;x,t)V^{(3)}(z;x,t)$, where
  \begin{equation}\label{AB4V3}
    \begin{aligned}
  V^{(3)}(z;x,t)&=\left\{  \begin{array}{ll}
     \begin{pmatrix} 0 & - e^{-\mathrm{i} (\gamma+\hat{\gamma})}  \\  e^{\mathrm{i} (\gamma+\hat{\gamma})}  & 0 \end{pmatrix}, ~&z \in (\lambda_l^-, \eta_-(\xi)), \\ \\
     \begin{pmatrix} 0 & - e^{-\mathrm{i}  (\gamma+\hat{\gamma})}  \\ e^{\mathrm{i}  (\gamma+\hat{\gamma})} &  \hat{\delta}_{+}(z)\hat{\delta}^{-1} _{-}(z)e^{ -2 \mathrm{i} t g_+(z)} e^{\mathrm{i} \gamma}  \end{pmatrix}, ~&z \in (\eta_-(\xi), \lambda_s(\xi)), \\ \\
     \begin{pmatrix} 0 & -1 \\1 & 0 \end{pmatrix}, ~&z \in (\lambda_r^-, \lambda_r^+), \\ \\
     \begin{pmatrix} 1 & 0 \\ r(z) \hat{\delta} ^{-2}(z) e^{2\mathrm{i} tg(z)} & 1 \end{pmatrix}, ~ &z \in L_1, \\ \\
     \begin{pmatrix} 1 & \frac{-r^*(z)}{1-r(z)r^*(z)} \hat{\delta} ^2(z) e^{-2\mathrm{i}t g(z)} \\0 & 1 \end{pmatrix},  ~&z \in L_2,\\ \\
     \begin{pmatrix} 1 & 0 \\ \frac{r(z)}{1-r(z)r^*(z)} \hat{\delta} ^{-2}(z) e^{2\mathrm{i} tg(z)} & 1 \end{pmatrix},  ~&z \in L_3,\\ \\
     \begin{pmatrix} 1 & -r^*(z) \hat{\delta} ^2(z) e^{-2\mathrm{i} tg(z)} \\ 0  & 1 \end{pmatrix}, ~ &z \in L_4. \end{array} \right.\\
    \end{aligned}
  \end{equation}
  \end{enumerate}
\end{rhp}

The jump matrices for $M^{(3)}(z;x,t)$ uniformly converge to the identity matrix or constant matrices outside a fixed neighborhood $\mathcal{U}$ of $\lambda_s(\xi)$,
while the convergence is not uniform inside $\mathcal{U}$ due to the exponential phase function $g(z)=\mathcal{O}((z-\lambda_s(\xi))^{3/2})$.
We construct an outer model parametrix that is the solution of the two-band limiting problem
\begin{equation}
  M^{(\mathrm{mod})}_+(z;x,t)=M^{(\mathrm{mod})}_-(z;x,t) V^{(\mathrm{mod})},
\end{equation}
where
\begin{equation}\label{AB4vmod}
  \begin{aligned}
    V^{(\mathrm{mod})}&=\left\{  \begin{array}{ll}
   \begin{pmatrix} 0 & - e^{-\mathrm{i} (\gamma+\hat{\gamma})}  \\  e^{\mathrm{i} (\gamma+\hat{\gamma})}  & 0 \end{pmatrix}, ~&z \in (\lambda_l^-, \lambda_s(\xi)), \\ \\
   \begin{pmatrix} 0 & -1 \\1 & 0 \end{pmatrix}, ~&z \in ( \lambda_r^-, \lambda_r^+). \end{array} \right.\\
   \end{aligned}
\end{equation}
The solution is  given by
\begin{equation} \label{AB4mod}
  M^{(\mathrm{mod})}(z;x,t)=N^{-1}(\infty; \boldsymbol{\lambda}) N(z; \boldsymbol{\lambda}),
\end{equation}
where $N(z; \boldsymbol{\lambda})$ is defined by (\ref{N}) with $\boldsymbol{\lambda}=(\lambda_r^+, \lambda_r^-, \lambda_s(\xi), \lambda_l^-)$.

\begin{figure}[t]
  \centering
  \subfigure[]{\includegraphics[width=7cm]{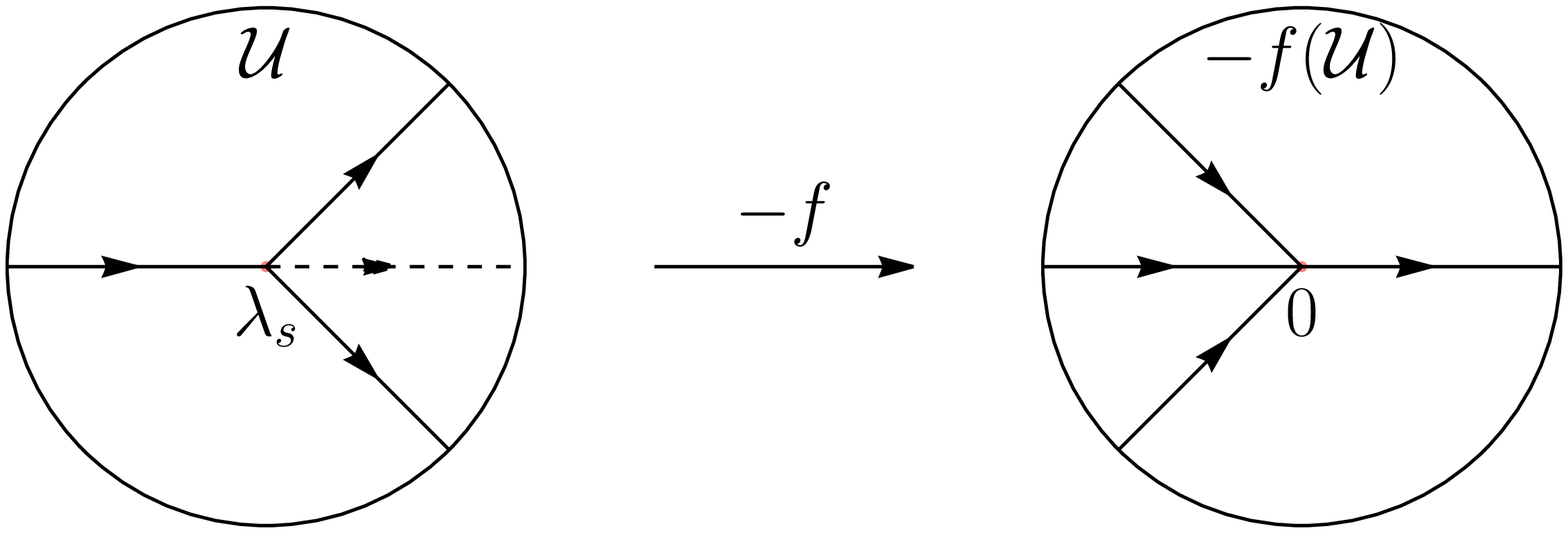}} \qquad
  \subfigure[]{\includegraphics[width=7cm]{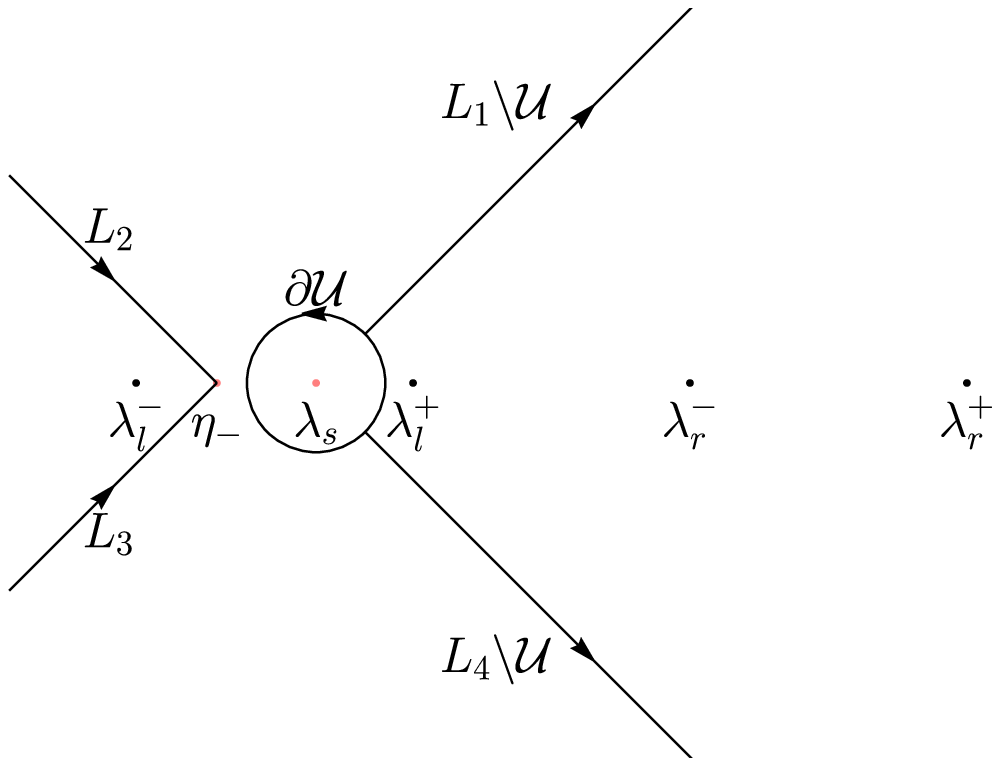}}
\caption{{\protect\small (a) The conformal mapping $-f$ on $\mathcal{U}$; (b) The jump contours of the error RHP.}}
\label{figA42}
\end{figure}

Inside $\mathcal{U}$, the convergence is not uniform since $g(z)=\mathcal{O}((z-\lambda_s(\xi))^{3/2})$ near $\lambda_s(\xi)$.
As before, construct a local parametrix $M^{\mathrm{(loc)}}(z;x,t)$ with the same jumps as $M^{(3)}(z;x,t)$
inside $\mathcal{U}$ by using the change in variables $k=-t^{2/3} f(z)$, where $-f$ is a conformal mapping from $\mathcal{U}$ (in $z$) to a neighborhood of the origin (in $k$), as shown in Figure \ref{figA42}, defined by
\begin{equation}\label{AB4f}
  \frac{2}{3}(-f(z))^{3/2}=\mathrm{i\, sgn}(\mathrm{Im}\,z)\,(g(z) - \gamma/2),
\end{equation}
and the branch is chosen such that $f'(\lambda_s(\xi))>0$. Define the local parametrix
\begin{equation}
  \begin{aligned}
  M^{\mathrm{(loc)}}(z;x,t)=& \left\{
    \begin{array}{ll}
    T(z) M^{(\mathrm{Ai})}(-t^{2/3} f(z)) \sigma_1 h_1(z)^{\sigma_3} \begin{pmatrix} 0 & -1 \\ 1 & 0 \end{pmatrix}, ~ & z \in \mathcal{U} \cap \mathbb{C}^+, \\
    T(z) M^{(\mathrm{Ai})}(-t^{2/3} f(z)) \sigma_1 h_2(z)^{\sigma_3} , ~ & z \in \mathcal{U} \cap \mathbb{C}^-,
  \end{array} \right.
\end{aligned}
\end{equation}
where the matching factors are analytic in $\mathcal{U}$ given by
\begin{equation}
  \begin{aligned}
    T(z)=& \left\{
      \begin{array}{ll}
      e^{-\mathrm{i}\pi/12} \sqrt{\pi} M^{(\mathrm{mod})}(z;x,t) \begin{pmatrix} 0 & 1 \\ -1 & 0 \end{pmatrix} h_1(z)^{-\sigma_3} \sigma_1 e^{\mathrm{i}\pi \sigma_3/4} \begin{pmatrix} 1 & -1 \\ 1 & 1 \end{pmatrix} (-t^{2/3} f(z))^{\sigma_3/4}, ~ & z \in \mathcal{U} \cap \mathbb{C}^+, \\
      e^{-\mathrm{i}\pi/12} \sqrt{\pi} M^{(\mathrm{mod})}(z;x,t)  h_2(z)^{-\sigma_3} \sigma_1 e^{\mathrm{i}\pi \sigma_3/4} \begin{pmatrix} 1 & -1 \\ 1 & 1 \end{pmatrix} (-t^{2/3} f(z))^{\sigma_3/4}, ~ & z \in \mathcal{U} \cap \mathbb{C}^-,
    \end{array} \right.
  \end{aligned}
\end{equation}
and
\begin{equation}
  h_1(z)=({r(z)})^{-1/2} \hat{\delta}(z), \qquad h_2(z)=({r^*(z)})^{-1/2} {\hat{\delta}}^{-1}(z).
\end{equation}
By direct calculation, it is easy to verify that this local parametrix satisfies our requirements.
We construct a global parametrix $M^{(\mathrm{par})}(z;x,t)$ by
\begin{equation}\label{A4para}
  \begin{aligned}
  M^{(\mathrm{par})}(z;x,t)
  &=\left\{  \begin{array}{ll} M^{\mathrm{(loc)}}(z;x,t), ~ & \mathrm{inside} ~ \mathcal{U}, \\  \\
     M^{(\mathrm{mod})}(z;x,t), ~& \mathrm{outside} ~ \mathcal{U},
   \end{array} \right.
  \end{aligned}
\end{equation}
and define the error matrix as
$M^{(\mathrm{err})}(z;x,t)= M^{(3)}(z;x,t) {M^{(\mathrm{par})}(z;x,t)}^{-1}$.
This leads to the following error RHP.

\begin{rhp} \label{rhp-rare} Find a $2 \times 2$ matrix-valued function $M^{(\mathrm{err})}(z;x,t)$ with the following properties:
  \begin{enumerate} [label=(\roman*)]
    \item $M^{(\mathrm{err})}(z;x,t)$ is analytic in $z \in \mathbb{C} \backslash \Sigma^{(\mathrm{err})}$, where $\Sigma^{(\mathrm{err})}=((L_1 \cup L_2 \cup  L_3 \cup L_4) \backslash\mathcal{U})\cup \partial \mathcal{U}$.
    \item $M^{(\mathrm{err})}(z;x,t)=I+\mathcal{O} (z^{-1}) $ as $z\to \infty$.
    \item $M^{(\mathrm{err})}(z;x,t)$ achieves the CBVs $M^{(\mathrm{err})}_{+}(z;x,t)$ and $M^{(\mathrm{err})}_{-}(z;x,t)$ on $\Sigma^{(\mathrm{err})}$, which satisfy the jump condition $M^{(\mathrm{err})}_{+}(z;x,t)=M^{(\mathrm{err})}_{-}(z;x,t)V^{(\mathrm{err})}(z;x,t)$, where
  \begin{equation}\label{A4Verr}
    \begin{aligned}
  V^{(\mathrm{err})}(z;x,t)&=\left\{  \begin{array}{ll}
    M^{(\mathrm{mod})}(z;x,t) \begin{pmatrix} 1 & 0 \\ r(z) \hat{\delta} ^{-2}(z) e^{2\mathrm{i} tg(z)} & 1 \end{pmatrix} {M^{(\mathrm{mod})}(z;x,t)}^{-1}, ~ &z \in L_1\backslash\mathcal{U}, \\ \\
    M^{(\mathrm{mod})}(z;x,t) \begin{pmatrix} 1 & \frac{-r^*(z)}{1-r(z)r^*(z)} \hat{\delta} ^2(z) e^{-2\mathrm{i}t g(z)} \\0 & 1 \end{pmatrix} {M^{(\mathrm{mod})}(z;x,t)}^{-1},  ~&z \in L_2,\\ \\
    M^{(\mathrm{mod})}(z;x,t) \begin{pmatrix} 1 & 0 \\ \frac{r(z)}{1-r(z)r^*(z)} \hat{\delta} ^{-2}(z) e^{2\mathrm{i} tg(z)} & 1 \end{pmatrix} {M^{(\mathrm{mod})}(z;x,t)}^{-1},  ~&z \in L_3,\\ \\
    M^{(\mathrm{mod})}(z;x,t) \begin{pmatrix} 1 & -r^*(z) \hat{\delta} ^2(z) e^{-2\mathrm{i} tg(z)} \\ 0  & 1 \end{pmatrix} {M^{(\mathrm{mod})}(z;x,t)}^{-1}, ~ &z \in L_4\backslash\mathcal{U},\\ \\
    M^{(\mathrm{mod})}(z;x,t) {M^{(\mathrm{loc})}(z;x,t)}^{-1}, ~ &z \in \partial \mathcal{U}.\end{array} \right.\\
    \end{aligned}
  \end{equation}
  \end{enumerate}
\end{rhp}
In the Appendix \ref{appb} we show that the error estimate is $M^{(\mathrm{err})}_1(x,t)=\mathcal{O}(t^{-1})$ as $t\to \infty$.
Using equations (\ref{rec}) (with $\delta$ replaced by $\hat{\delta}$), (\ref{AB4mdeltaINF}) and (\ref{AB4mod}), it is concluded that the long-time asymptotic behavior of the defocusing NLS equation (\ref{(NLS)}) in the dispersive shock wave region is given by
\begin{equation}\label{rsab}
  q(x,t)=\frac{\lambda_r^+ -\lambda_r^- + \lambda_s(\xi) -\lambda_l^-}{2}\frac{\Theta (0)~\Theta (2\mathcal{A}(\infty)+(\gamma + \hat{\gamma})/{2 \pi})}{\Theta ((\gamma + \hat{\gamma})/{2 \pi})~\Theta (2\mathcal{A}(\infty))} e^{2\mathrm{i} (t g_\infty+\hat{g}_\infty-\varphi^{(0)})}+\mathcal{O}(t^{-1}),
\end{equation}
where the real quantities $\gamma$, $g_\infty$, $\hat{g}_\infty$, $\varphi^{(0)}$, $\hat{\gamma}$ and the pure imaginary quantity $\mathcal{A}(\infty)$ are given by (\ref{gamma2}), (\ref{A4GINF}), (\ref{AB4mdeltaINF}), (\ref{AB4mdeltaINF}), (\ref{gammah2}) and (\ref{Abelinf}), respectively.

\subsubsection{The right plane wave region: $\xi > -\frac{\lambda_r^+ + \lambda_r^- + 2\lambda_l^-}{2}+\frac{2(\lambda_r^+ - \lambda_l^-)
(\lambda_r^- - \lambda_l^-)}{\lambda_r^+ + \lambda_r^- -2\lambda_l^-}$} \label{ABrpw}
\
\newline
\indent

\begin{figure}[htbp]
  \centering
  \subfigure[]{\includegraphics[width=5.5cm]{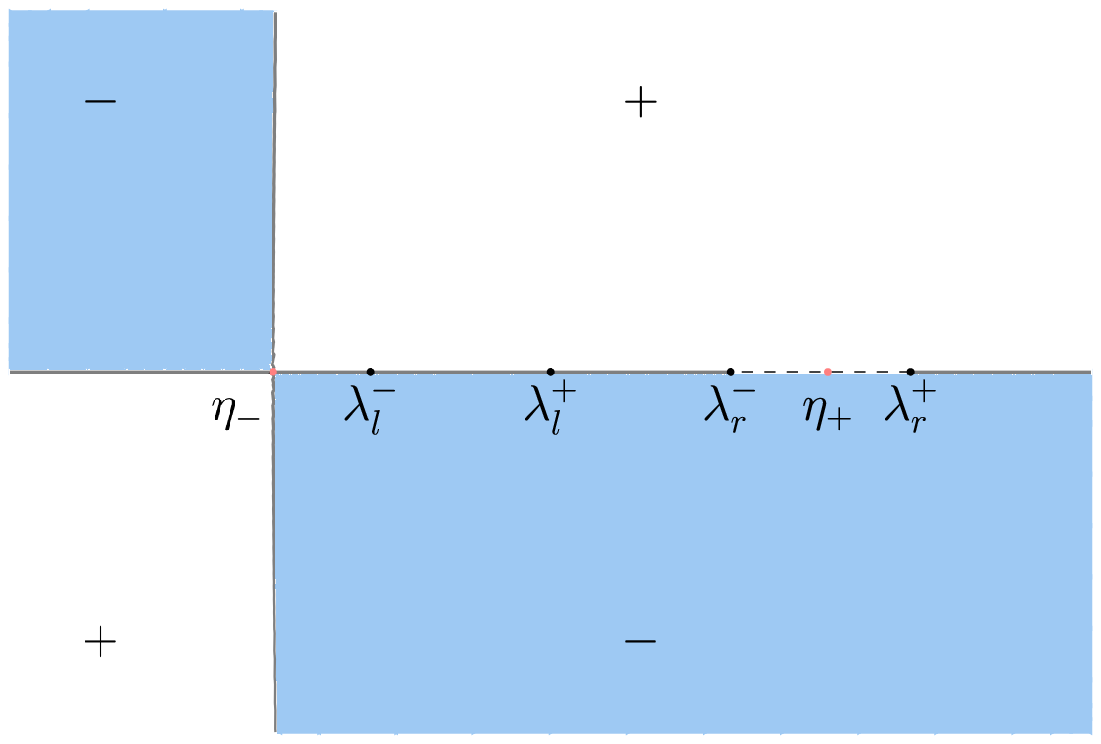}}
  \subfigure[]{\includegraphics[width=5cm]{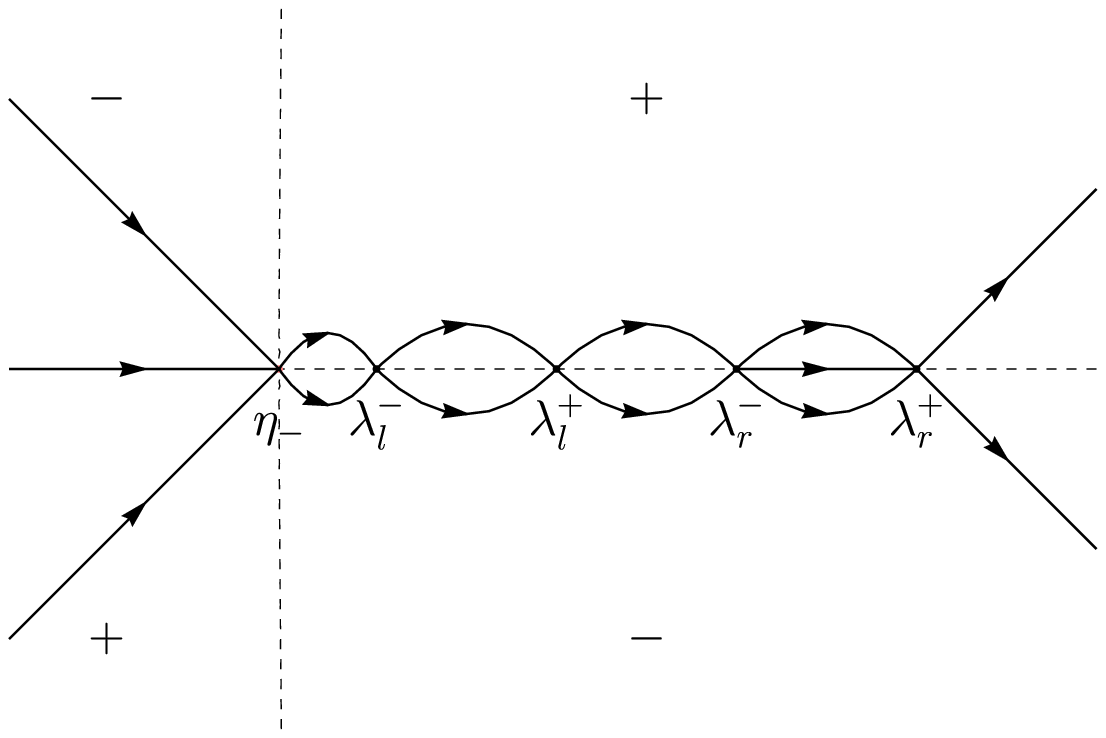}}
  \subfigure[]{\includegraphics[width=5.5cm]{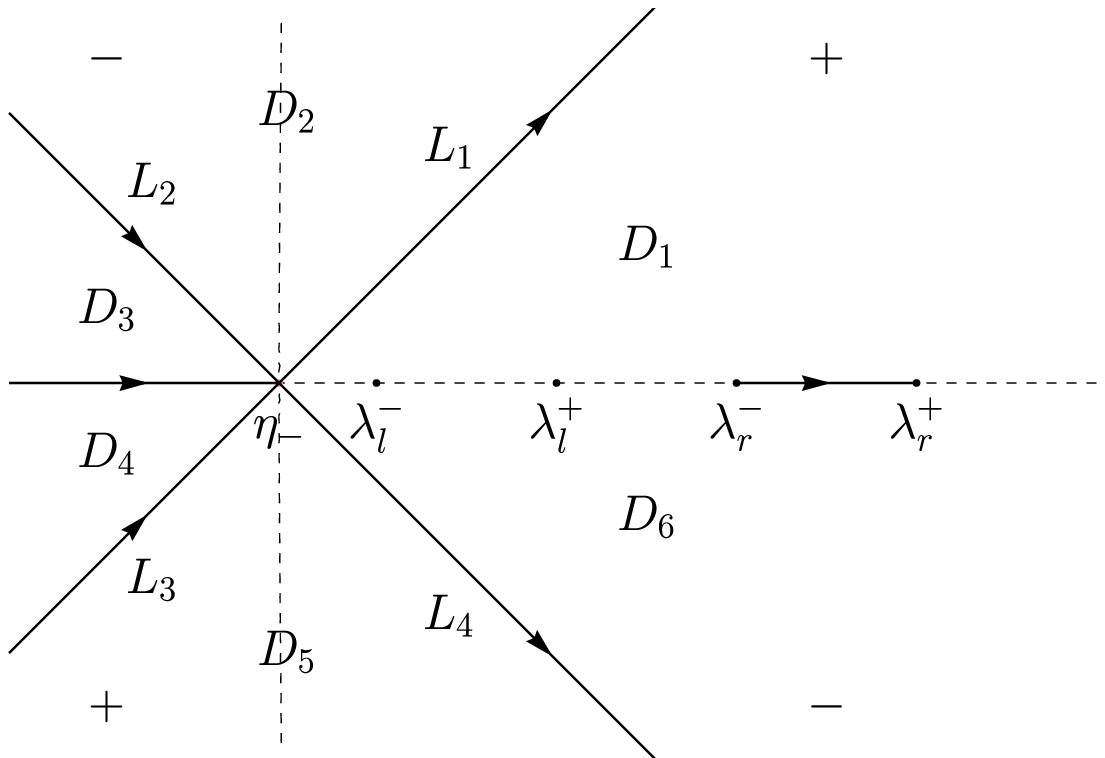}}
\caption{{\protect\small (a) The sign structure of Im$(g(z))$; (b) Opening lenses; (c) The jump contours of $M^{(2)}(z;x,t)$.}}
\label{figA5}
\end{figure}

When the stationary phase point $\lambda_s(\xi)$  coincides with $\lambda_l^-$, the two-band $g$-function degenerates to the one-band one.
Corresponding to the phase function of the explicit eigenfunction $\Psi_r^{\mathrm{p}}(z;x,t)$ defined by (\ref{psip}), the appropriate $g$-function is
\begin{equation} \label{A5g}
  g(z)= \mathcal{R} (z;\lambda_r^+,\lambda_r^-)(z+\frac{1}{2}(\lambda_r^++\lambda_r^-)+\xi),
\end{equation}
which is consistent with the genus-zero case in which the Riemann invariants $\lambda_r^+$ and $\lambda_r^-$
are two hard edges in Whitham modulation theory.
The function $g(z)$ is analytic in $z \in \mathbb{C} \backslash \overline{\mathcal{I}_r}$ with
\begin{equation}\label{A5ginf}
   g_{\infty}=- \mu_r \xi-(2 \mu_r^2+A_r^2)/2,
\end{equation}
and the differential of $g(z)$ is
\begin{equation} \label{A5dg}
  \dif g(z)= 2 \frac{(z-\eta_+(\xi))(z-\eta_-(\xi))}{\mathcal{R} (z;\lambda_r^+,\lambda_r^-)} \,\dif z,
\end{equation}
where $\eta_{\pm}(\xi)$ are the stationary phase points expressed by
\begin{equation} \label{rpweta}
  \eta_{\pm}(\xi)=\frac{\lambda_r^++\lambda_r^--\xi \pm \sqrt{(\lambda_r^++\lambda_r^-+\xi)^2+2(\lambda_r^+-\lambda_r^-)^2}}{4}.
\end{equation}
Furthermore, $\eta_- (\xi)$ determines the boundary of the right plane wave region as
\begin{equation}\label{A5B}
  \eta_- (\xi) < \lambda_l^- \qquad \mathrm{iff}  \qquad  \xi > \xi_{A4}=-\frac{\lambda_r^+ + \lambda_r^- + 2\lambda_l^-}{2}+\frac{2(\lambda_r^+ - \lambda_l^-)
  (\lambda_r^- - \lambda_l^-)}{\lambda_r^+ + \lambda_r^- -2\lambda_l^-},
\end{equation}
and this boundary is consistent with (\ref{AB4b}).

Opening contours  through $\eta_-(\xi)$  (as shown in Figure \ref{figA5})  by the transformations $M(z;x,t) \mapsto M^{(1)}(z;x,t) \mapsto M^{(2)}(z;x,t)$, we obtain the following RHP.

\begin{rhp}
  Find a $2 \times 2$ matrix-valued function $M^{(2)}(z;x,t)$ with the following properties:
\begin{enumerate} [label=(\roman*)]
  \item  $M^{(2)}(z;x,t)$ is analytic in $z \in \mathbb{C} \backslash \Sigma^{(2)}$, where $\Sigma^{(2)}=(-\infty, \eta_{-}(\xi)] \cup [\lambda_r^-, \lambda_r^+] \cup L_1 \cup L_2 \cup  L_3 \cup L_4.$
  \item  $M^{(2)}(z;x,t)=I+\mathcal{O} (z^{-1}) $ as $z\to \infty$.
  \item  $M^{(2)}(z;x,t)$ achieves the CBVs  $M^{(2)}_{+}(z;x,t)$ and $M^{(2)}_{-}(z;x,t)$  on $\Sigma^{(2)}$ away from self-intersection points and
  branch points that satisfy the jump condition $M^{(2)}_{+}(z;x,t)=M^{(2)}_{-}(z;x,t)V^{(2)}(z;x,t)$, where
  \begin{equation}\label{A5V2}
    \begin{aligned}
    V^{(2)}(z;x,t)&=\left\{  \begin{array}{ll}
     (1-r(z)r^{*}(z))^{\sigma_{3}},  ~& z\in (-\infty, \eta_{-}(\xi)), \\ \\
     \begin{pmatrix} 0 & -1 \\ 1 & 0 \end{pmatrix},  ~&z \in  (\lambda_r^-, \lambda_r^+), \\ \\
     \begin{pmatrix} 1 & 0 \\ r(z)e^{2\mathrm{i} tg(z)} & 1 \end{pmatrix}, ~ &z \in L_1, \\ \\
     \begin{pmatrix} 1 & \frac{-r^*(z)}{1-r(z)r^*(z)}e^{-2\mathrm{i} tg(z)} \\0 & 1 \end{pmatrix},  ~&z \in L_2,\\ \\
     \begin{pmatrix} 1 & 0 \\ \frac{r(z)}{1-r(z)r^*(z)}e^{2\mathrm{i} tg(z)} & 1 \end{pmatrix},  ~&z \in L_3,\\ \\
     \begin{pmatrix} 1 & -r^*(z)e^{-2\mathrm{i} tg(z)} \\ 0  & 1 \end{pmatrix}, ~ &z \in L_4.\end{array} \right.\\
    \end{aligned}
  \end{equation}
\end{enumerate}
\end{rhp}
The above problem is completely similar to that in the left plane wave region, so we put these transformations together:
$M^{(2)}(z;x,t) \to M^{(3)}(z;x,t) \to M^{(\mathrm{err})}(z;x,t)$, that is,
\begin{equation}\label{A5para}
  \begin{aligned}
  M^{(2)}(z;x,t)
  &=\left\{  \begin{array}{ll} \delta _{\infty}^{-\sigma_{3}} M^{(\mathrm{err})}(z;x,t) M^{\mathrm{(loc)}}(z;x,t)  \delta(z) ^{\sigma_{3}}, ~ & \mathrm{inside} ~ \mathcal{U}, \\  \\
    \delta _{\infty}^{-\sigma_{3}} M^{(\mathrm{err})}(z;x,t) M^{(\mathrm{mod})}(z;x,t)  \delta(z) ^{\sigma_{3}}, ~& \mathrm{outside} ~ \mathcal{U},
   \end{array} \right.\\
  \end{aligned}
\end{equation}
where
\begin{equation} \label{A5delta}
  \begin{aligned}
  \delta(z) = \exp & \left\{ \frac{\mathcal{R} (z;\lambda_r^+,\lambda_r^-)}{2 \pi \mathrm{i}}  \left(\int_{-\infty}^ {\eta_{-}(\xi) }  \frac{\ln (1-|r(\zeta)|^2)}{\mathcal{R} (\zeta;\lambda_r^+,\lambda_r^-) } \,\frac{\dif\zeta}{\zeta - z}   \right)  \right\}
  \end{aligned}
\end{equation}
and
\begin{equation}\label{A5deltaINF }
  \delta_{\infty }= \exp \left\{  \frac{ \mathrm{i}}{2\pi} \left(\int_{-\infty}^ {\eta_{-}(\xi) }  \frac{\ln (1-|r(\zeta)|^2)}{\mathcal{R} (\zeta;\lambda_r^+,\lambda_r^-) } \dif \zeta \right) \right\}.
\end{equation}
Here, fix a small neighborhood $\mathcal{U}$ of $\eta_-(\xi)$ and define
\begin{equation} \label{A5mod}
  M^{(\mathrm{mod})}(z;x,t)= \mathcal{E}_r(z).
\end{equation}
For the local parametrix $M^{(\mathrm{loc})}(z;x,t)$ inside $\mathcal{U}$, the construction is similar to that in the left plane wave region
\begin{equation}
  M^{\mathrm{(loc)}}(z;x,t)=M^{(\mathrm{mod})}(z;x,t)h(z)^{-\sigma_3}M^{(\mathrm{PC})}(\sqrt{t} f(z);r(z))h(z)^{\sigma_3},
\end{equation}
where
\begin{equation}
  k=\sqrt{t} f(z)=2\sqrt{t} (g(z)-g(\eta_-(\xi)))^{1/2},
\end{equation}
where the branch is chosen such that $f'(\eta_+(\xi))>0$, is a conformal mapping from $\mathcal{U}$ (in $z$) to a neighborhood of the origin (in $k$), and
\begin{equation}
  h(z)=e^{\mathrm{i}g(\eta_-(\xi))}(\sqrt{t} f(z))^{\mathrm{i}\nu(z)}\delta^{-1}(z)
\end{equation}
is principally branched. Hence, we again have the error estimate $M^{(\mathrm{err})}_1(x,t)= \mathcal{O}(t^{-1/2})$. The details of the above process can be found in Section \ref{Alpw}.
Using equations (\ref{rec}), (\ref{A5ginf}), (\ref{A5deltaINF }) and (\ref{A5mod}), it is concluded that the long-time asymptotic behavior of the defocusing NLS equation (\ref{(NLS)})
in the right plane wave region is given by
\begin{equation}
    q(x,t)=A_r e^{-2\mathrm{i} \mu_r x-\mathrm{i}(2 \mu_r^2+A_r^2)t}e^{-\mathrm{i}\phi_{rp}(\xi) }+\mathcal{O}(t^{-\frac{1}{2}}),
\end{equation}
where
\begin{equation}\label{A5phi1}
  \phi_{rp}(\xi)=\frac{1}{\pi}\int_{-\infty}^{\eta_{-}(\xi)} \frac{\ln (1-|r(\zeta)|^2)}{\sqrt{(\zeta-\lambda_r^+)(\zeta-\lambda_r^-)} } \dif \zeta.
\end{equation}

\subsection{\rm Case B: $\lambda_l^+>\lambda_l^->\lambda_r^+>\lambda_r^-$}
\
\newline
\indent
In this case, the two intervals $\mathcal{I}_l$ and $\mathcal{I}_r$ still do not intersect, so we only need the factorizations
(\ref{f1})-(\ref{f3}). In comparison with Case A, the positions of the two intervals are exchanged. This implies the absence of the elliptic (genus-one) wave region, as we will see later.
As the self-similar variable $\xi=x/t$ increases, the stationary phase points of the corresponding $g$-functions  change at different intervals on $\mathbb{R}$, which
implies that there are five different regions: the left plane wave region, rarefaction region, vacuum region, rarefaction region and the right plane wave region.
The fact that two plane waves corresponding to the initial data separate at the origin point for $t=0$ is consistent with the appearance of the middle vacuum region \cite{El-Phys.D1995}.

\subsubsection{The left plane wave region: $\xi<-\frac{3\lambda_l^++\lambda_l^-}{2}$}\label{Blpw}
\
\newline
\indent
As mentioned above, the leading order asymptotics in the left plane wave region  are consistent with the initial condition for $t = 0, x < 0$, up to the phase shift  $e^{-\mathrm{i}\phi_{lp}(\xi)}$.
The reason for this phenomenon is that the $g$-function is the same as (\ref{A1g}) for left plane wave regions in all six cases in (\ref{Classification}) and  the stationary phase point $\eta_+(\xi)$ of the left plane wave $g$-function  lies on the right side of both intervals $\mathcal{I}_l $ and $\mathcal{I}_r$.
The factorizations (\ref{f1})-(\ref{f4}) imply that only jump matrices on  $\mathcal{I}_l $  (sometimes consisting of $\mathcal{I}_l \backslash \overline{\mathcal{I}_r}  $  and $\mathcal{I}_l \cap \mathcal{I}_r $)   contribute to the asymptotic behaviors.
The boundary of the left plane wave region is characterized by
\begin{equation} \label{B1b}
  \eta_+(\xi) > \lambda_l^+    \qquad \mathrm{iff}  \qquad \xi< \xi_{B1}= -\frac{3\lambda_l^++\lambda_l^-}{2},
\end{equation}
where $\eta_+ (\xi)$ is given by (\ref{lpweta}).

Henceforth, a unified form of RHP for $M^{(2)}(z;x,t)$ is given by opening lenses from the real axis to the steepest descent contours through $\eta_+(\xi)$. The jumps depend on the specific forms of the intervals $(-\infty, \eta_{+}(\xi))\backslash (\overline{\mathcal{I}_l \cup \mathcal{I}_r} )$, $\mathcal{I}_l \backslash \overline{\mathcal{I}_r} $, $\mathcal{I}_r \backslash \overline{\mathcal{I}_l} $ and $\mathcal{I}_l \cap \mathcal{I}_r $ in the six cases. Therefore, the RHP below follows.

\begin{rhp}
  Find a $2 \times 2$ matrix-valued function $M^{(2)}(z;x,t)$ with the following properties:
\begin{enumerate} [label=(\roman*)]
  \item  $M^{(2)}(z;x,t)$ is analytic in $z \in \mathbb{C} \backslash \Sigma^{(2)}$, where $\Sigma^{(2)}=(-\infty, \eta_{+}(\xi)] \cup L_1 \cup L_2 \cup  L_3 \cup L_4.$
  \item  $M^{(2)}(z;x,t)=I+\mathcal{O} (z^{-1}) $ as $z\to \infty$.
  \item  $M^{(2)}(z;x,t)$ achieves the CBVs  $M^{(2)}_{+}(z;x,t)$ and $M^{(2)}_{-}(z;x,t)$  on $\Sigma^{(2)}$ away from self-intersection points and
  branch points that satisfy the jump condition $M^{(2)}_{+}(z;x,t)=M^{(2)}_{-}(z;x,t)V^{(2)}(z;x,t)$, where
  \begin{equation}\label{B1V2}
    \begin{aligned}
    V^{(2)}(z;x,t)&=\left\{  \begin{array}{ll}
     (1-r(z)r^{*}(z))^{\sigma_{3}},  ~& z\in (-\infty, \eta_{+}(\xi))\backslash (\overline{\mathcal{I}_l \cup \mathcal{I}_r} ), \\ \\
     (a_+(z) a_-^*(z))^{-\sigma_{3}},  ~&z \in  \mathcal{I}_r \backslash \overline{\mathcal{I}_l}, \\ \\
     \begin{pmatrix} 0 & -r_{-}^{*}(z) \\ r_{+}(z) & 0 \end{pmatrix}, ~&z \in \mathcal{I}_l \backslash \overline{\mathcal{I}_r}, \\ \\
     \begin{pmatrix} 0 & -1 \\ 1 & 0 \end{pmatrix}, ~&z \in \mathcal{I}_l \cap \mathcal{I}_r, \\ \\
     \begin{pmatrix} 1 & 0 \\ r(z)e^{2\mathrm{i} tg(z)} & 1 \end{pmatrix}, ~ &z \in L_1, \\ \\
     \begin{pmatrix} 1 & \frac{-r^*(z)}{1-r(z)r^*(z)}e^{-2\mathrm{i} tg(z)} \\0 & 1 \end{pmatrix},  ~&z \in L_2,\\ \\
     \begin{pmatrix} 1 & 0 \\ \frac{r(z)}{1-r(z)r^*(z)}e^{2\mathrm{i} tg(z)} & 1 \end{pmatrix},  ~&z \in L_3,\\ \\
     \begin{pmatrix} 1 & -r^*(z)e^{-2\mathrm{i} tg(z)} \\ 0  & 1 \end{pmatrix}, ~ &z \in L_4.\end{array} \right.\\
    \end{aligned}
  \end{equation}
\end{enumerate}
\end{rhp}
In this case,  $\mathcal{I}_l $ and $\mathcal{I}_r $ are separated, i.e., $\mathcal{I}_l \cap \mathcal{I}_r = \varnothing$, but we still keep this term  for the sake of uniformity.
As before, define a scalar function $\delta(z)$ to reduce jump matrices across the real axis to constant matrices independent of $z$,
which is analytic in $z \in \mathbb{C} \backslash(-\infty, \eta_{+}(\xi)]$ and satisfies the jump conditions
\begin{equation*}
  \begin{aligned}
\left\{  \begin{array}{ll}
   \delta _{+}(z)=(1-r(z)r^{*}(z))\delta _{-}(z),  ~& z\in (-\infty, \eta_{+}(\xi)) \backslash (\overline{\mathcal{I}_l \cup \mathcal{I}_r} ), \\ \\
   \delta _{+}(z)=(a_+(z) a_-^*(z))^{-1}\delta _{-}(z),  ~& z\in \mathcal{I}_r \backslash \overline{\mathcal{I}_l}, \\ \\
   \delta _{+}(z)\delta _{-}(z)=r_{+}(z),  ~& z\in \mathcal{I}_l \backslash \overline{\mathcal{I}_r},\\ \\
   \delta _{+}(z)\delta _{-}(z)=1,  ~& z\in \mathcal{I}_l \cap \mathcal{I}_r.
  \end{array} \right.\\
  \end{aligned}
\end{equation*}
As previously calculated, $\delta(z)$ is given by
\begin{equation} \label{B1delta}
  \begin{aligned}
  \delta(z) = \exp & \left\{ \frac{\mathcal{R} (z;\lambda_l^+,\lambda_l^-)}{2 \pi \mathrm{i}}  \left(\int_{(-\infty, \eta_{+}(\xi)) \backslash (\overline{\mathcal{I}_l \cup \mathcal{I}_r} )}  \frac{\ln (1-|r(\zeta)|^2)}{\mathcal{R} (\zeta;\lambda_l^+,\lambda_l^-) } \,\frac{\dif\zeta}{\zeta - z} \right. \right. \\
  & \left. \left.  + \int_{\mathcal{I}_r \backslash \overline{\mathcal{I}_l}} \frac{-\ln(a_+(\zeta)a_-^*(\zeta))}{\mathcal{R} (\zeta;\lambda_l^+,\lambda_l^-)} \,\frac{\dif\zeta}{\zeta - z}
  + \int_{\mathcal{I}_l \backslash \overline{\mathcal{I}_r}} \frac{\ln (r_+(\zeta))}{\mathcal{R}_+ (\zeta;\lambda_l^+,\lambda_l^-)} \, \frac{\dif\zeta}{\zeta - z} \right)  \right\}
  \end{aligned}
\end{equation}
with large $z$ asymptotic behavior $\delta(z)=\delta_{\infty}+\mathcal{O} (z^{-1}) $ as $z \to \infty$, where
\begin{equation}\label{B1deltaINF}
  \begin{aligned}
  \delta_{\infty }= \exp & \left\{ \frac{ \mathrm{i}}{2\pi} \left(\int_{(-\infty, \eta_{+}(\xi)) \backslash (\overline{\mathcal{I}_l \cup \mathcal{I}_r} )}  \frac{\ln (1-|r(\zeta)|^2)}{\mathcal{R} (\zeta;\lambda_l^+,\lambda_l^-) } \dif \zeta \right. \right. \\
  & \left. \left.  + \int_{\mathcal{I}_r \backslash \overline{\mathcal{I}_l}} \frac{-\ln(a_+(\zeta)a_-^*(\zeta))}{\mathcal{R} (\zeta;\lambda_l^+,\lambda_l^-)} \dif \zeta
  + \int_{\mathcal{I}_l \backslash \overline{\mathcal{I}_r}} \frac{\ln (r_+(\zeta))}{\mathcal{R}_+ (\zeta;\lambda_l^+,\lambda_l^-)} \dif \zeta \right)  \right\}.
  \end{aligned}
\end{equation}

The remaining steps are exactly the same as those in Section \ref{Alpw}, so we omit the details. Finally, a unified form of the long-time asymptotics of the defocusing NLS equation (\ref{(NLS)})
in the left plane wave region is given by
\begin{equation}
  q(x,t)=A_l e^{-2\mathrm{i} \mu_l x-\mathrm{i}(2 \mu_l^2+A_l^2)t}e^{-\mathrm{i}\phi_{lp}(\xi) }+\mathcal{O}(t^{-\frac{1}{2}}),
\end{equation}
where
\begin{equation}\label{philp}
  \begin{aligned}
   \phi_{lp}(\xi)=\dfrac{1}{\pi} & \left(\int_{(-\infty, \eta_{+}(\xi)) \backslash (\overline{\mathcal{I}_l \cup \mathcal{I}_r} )} \frac{\ln (1-|r(\zeta)|^2)}{\sqrt{(\zeta-\lambda_l^+)(\zeta-\lambda_l^-)} } \dif \zeta \right.\\
   &\left. + \int_{\mathcal{I}_r \backslash \overline{\mathcal{I}_l}} \frac{-\ln(a_+(\zeta)a_-^*(\zeta))}{\sqrt{(\zeta-\lambda_l^+)(\zeta-\lambda_l^-)}} \dif \zeta
   +  \int_{\mathcal{I}_l \backslash \overline{\mathcal{I}_r}} \frac{\arg (r_+(\zeta))}{\sqrt{(\lambda_l^+-\zeta)(\zeta-\lambda_l^-)}} \dif \zeta \right).\\
  \end{aligned}
\end{equation}

\begin{remark}
  In some cases the interval $\mathcal{I}_l \backslash \overline{\mathcal{I}_r}$ or ${\mathcal{I}_r \backslash \overline{\mathcal{I}_l}}$ may be  empty. When this happens, one only needs to omit the integral on the corresponding interval of $\phi_{lp}(\xi)$ in (\ref{philp}).
\end{remark}

\subsubsection{Rarefaction wave region: $-\frac{3\lambda_l^++\lambda_l^-}{2} <\xi <-2\lambda_l^-$} \label{Alrare}
\
\newline
\indent

\begin{figure}[h]
  \centering
  \subfigure[]{\includegraphics[width=5.5cm]{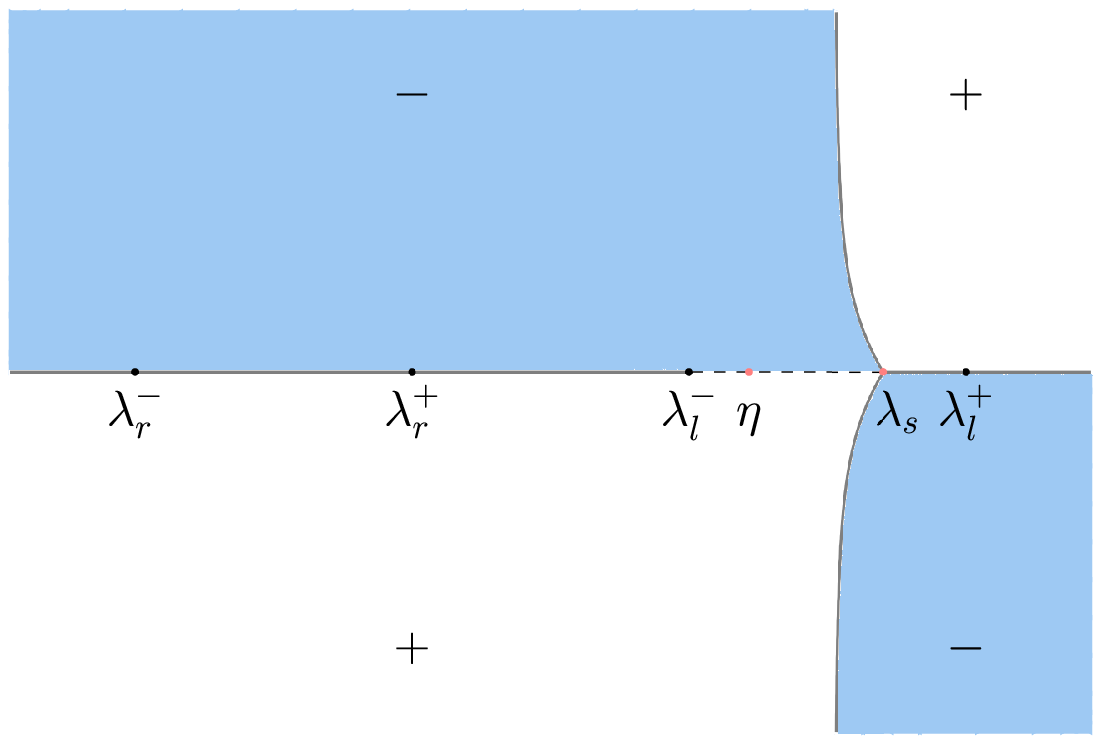}}
  \subfigure[]{\includegraphics[width=5cm]{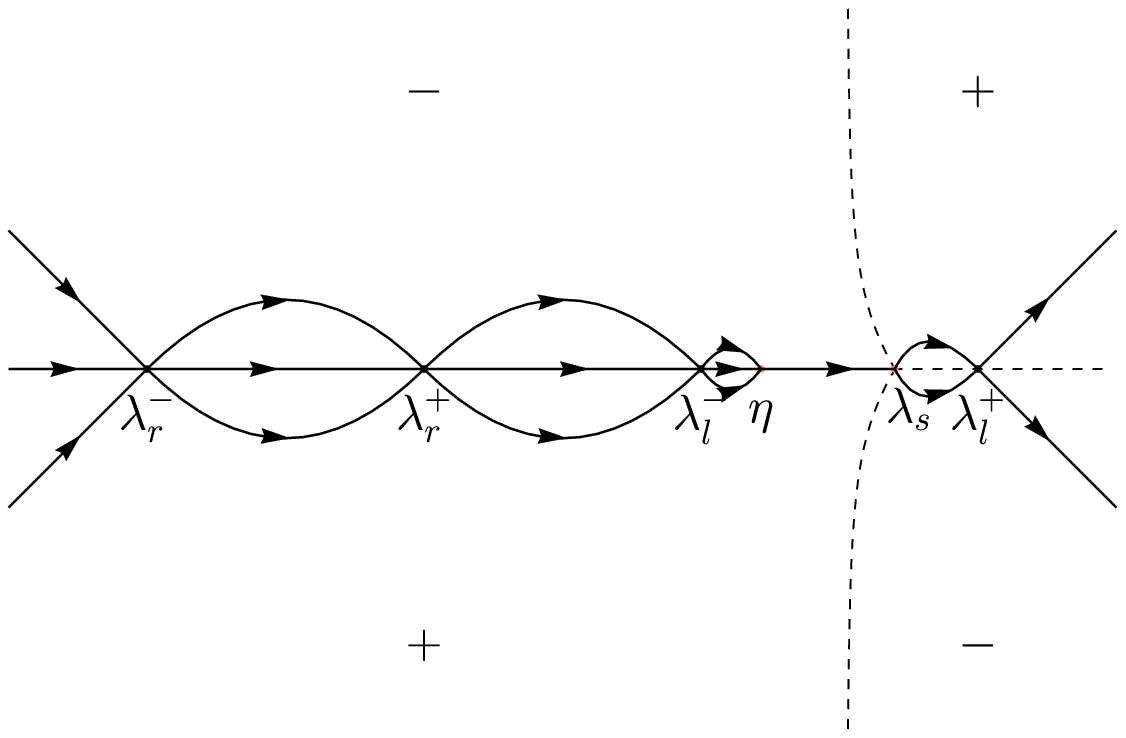}}
  \subfigure[]{\includegraphics[width=5.5cm]{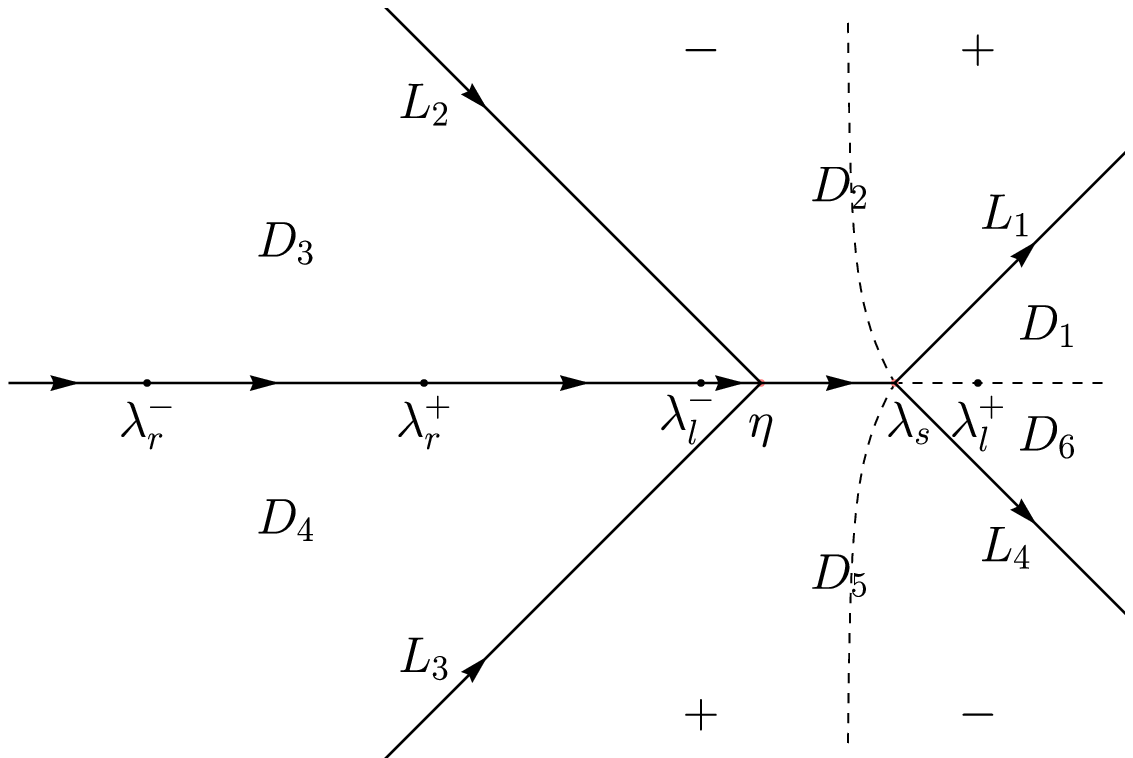}}
\caption{{\protect\small (a) The sign structure of Im$(g(z))$; (b) Opening lenses; (c) The jump contours of $M^{(2)}(z;x,t)$.}}
\label{figB2}
\end{figure}
As $\xi$ increases such that $\xi > -\frac{3\lambda_l^++\lambda_l^-}{2} $, the stationary phase point $\eta_+(\xi)$ of
the previous $g$-function moves inside $\mathcal{I}_l$, which contributes exponentially large diagonal entries $e^{\pm 2 \mathrm{i} t g_+(z)}$ of the jump matrix $V^{(2)}(z;x,t)$ on the interval $(\eta_+(\xi), \lambda_l^+) \subset
\mathcal{I}_l$ in (\ref{B1V2}). So a new $g$-function should be introduced, which is analytic in $z \in \mathbb{C} \backslash [\lambda_l^-, \lambda_s(\xi)]$ with the soft
edge $\lambda_s(\xi)=-(\lambda_l^- + 2 \xi)/3$ defined in Whitham modulation theory. The monotonically decreasing  function $\lambda_s(\xi)$ determines
the boundaries of this rarefaction wave region:
\begin{equation} \label{B2b}
  \lambda_l^- < \lambda_s(\xi) < \lambda_l^+ \qquad \mathrm{iff}  \qquad -\frac{3\lambda_l^++\lambda_l^-}{2} =\xi_{B1} < \xi < \xi_{B2} = -2\lambda_l^- .
\end{equation}

To determine $g(z)$, consider the differential
\begin{equation} \label{A2dg}
  \dif g(z)= 2 \frac{(z-\lambda_s(\xi))(z-\eta(\xi))}{\mathcal{R} (z;\lambda_s(\xi),\lambda_l^-)} \,\dif z,
\end{equation}
where $\eta(\xi)=(\lambda_s(\xi)+3\lambda_l^-)/4$ is determined by $\int_{\lambda_l^-}^{\lambda_s(\xi)} \,\dif g(z) =0$.
Hence, the rarefaction $g$-function is given by
\begin{equation} \label{A2g}
  g(z)= \mathcal{R} (z;\lambda_s(\xi),\lambda_l^-)(z-\lambda_s(\xi)),
\end{equation}
with
\begin{equation}\label{A2ginf}
  g_{\infty}=-(2{\lambda_l^-}^2+2\lambda_l^- \xi-\xi^2)/6.
\end{equation}

Then opening lenses to the steepest descent contours through $\lambda_s(\xi)$ and $\eta(\xi)$  (as shown in Figure \ref{figB2})
by the transformations (\ref{g-fun}) and (\ref{open}) yields the following RHP.

\begin{rhp}
  Find a $2 \times 2$ matrix-valued function $M^{(2)}(z;x,t)$ with the following properties:
\begin{enumerate} [label=(\roman*)]
  \item  $M^{(2)}(z;x,t)$ is analytic in $z \in \mathbb{C} \backslash \Sigma^{(2)}$, where $\Sigma^{(2)}=(-\infty, \lambda_{s}(\xi)] \cup L_1 \cup L_2 \cup  L_3 \cup L_4.$
  \item  $M^{(2)}(z;x,t)=I+\mathcal{O} (z^{-1}) $ as $z\to \infty$.
  \item  $M^{(2)}(z;x,t)$ achieves the CBVs  $M^{(2)}_{+}(z;x,t)$ and $M^{(2)}_{-}(z;x,t)$  on $\Sigma^{(2)}$ away from self-intersection points and
  branch points that satisfy the jump condition $M^{(2)}_{+}(z;x,t)=M^{(2)}_{-}(z;x,t)V^{(2)}(z;x,t)$, where
  \begin{equation}\label{A2V2}
    \begin{aligned}
    V^{(2)}(z;x,t)&=\left\{  \begin{array}{ll}
     (1-r(z)r^{*}(z))^{\sigma_{3}},  ~& z\in (-\infty, \lambda_r^-) \cup (\lambda_r^+, \lambda_l^-), \\ \\
     (a_+(z) a_-^*(z))^{-\sigma_{3}},  ~&z \in (\lambda_r^-, \lambda_r^+), \\ \\
     \begin{pmatrix} 0 & -r_{-}^{*}(z) \\ r_{+}(z) & 0 \end{pmatrix}, ~&z \in (\lambda_l^-, \eta(\xi)), \\ \\
     \begin{pmatrix} 0 & -r_{-}^{*}(z) \\ r_{+}(z) & e^{- 2 \mathrm{i} t g_+(z)} \end{pmatrix}, ~&z \in (\eta(\xi), \lambda_s(\xi)), \\ \\
     \begin{pmatrix} 1 & 0 \\ r(z)e^{2\mathrm{i} tg(z)} & 1 \end{pmatrix}, ~ &z \in L_1, \\ \\
     \begin{pmatrix} 1 & \frac{-r^*(z)}{1-r(z)r^*(z)}e^{-2\mathrm{i} tg(z)} \\0 & 1 \end{pmatrix},  ~&z \in L_2,\\ \\
     \begin{pmatrix} 1 & 0 \\ \frac{r(z)}{1-r(z)r^*(z)}e^{2\mathrm{i} tg(z)} & 1 \end{pmatrix},  ~&z \in L_3,\\ \\
     \begin{pmatrix} 1 & -r^*(z)e^{-2\mathrm{i} tg(z)} \\ 0  & 1 \end{pmatrix}, ~ &z \in L_4.\end{array} \right.\\
    \end{aligned}
  \end{equation}
\end{enumerate}
\end{rhp}

As before, it is necessary to determine a scalar function $\delta(z)$ that is analytic in  $z \in \mathbb{C} \backslash(-\infty, \lambda_s(\xi)]$ and satisfies the jump conditions
\begin{equation*}
  \begin{aligned}
\left\{  \begin{array}{ll}
   \delta _{+}(z)=(1-r(z)r^{*}(z))\delta _{-}(z),  ~& z\in (-\infty, \lambda_r^-) \cup (\lambda_r^+, \lambda_l^-), \\ \\
   \delta _{+}(z)=(a_+(z) a_-^*(z))^{-1}\delta _{-}(z),  ~& z\in (\lambda_r^-, \lambda_r^+), \\ \\
   \delta _{+}(z)\delta _{-}(z)=r_{+}(z),  ~& z\in (\lambda_l^-, \lambda_s(\xi)).
  \end{array} \right.\\
  \end{aligned}
\end{equation*}
By Plemelj formula, we have
\begin{equation} \label{A2delta}
  \begin{aligned}
  \delta(z) = \exp & \left\{ \frac{\mathcal{R} (z;\lambda_s(\xi), \lambda_l^-)}{2 \pi \mathrm{i}}  \left[ \left( \int_{-\infty}^{\lambda_r^-} + \int_{\lambda_r^+}^{\lambda_l^-} \right) \frac{\ln (1-|r(\zeta)|^2)}{\mathcal{R} (\zeta;\lambda_s(\xi), \lambda_l^-) } \,\frac{\dif\zeta}{\zeta - z}  \right. \right. \\
  & \left. \left.  + \int_{\lambda_r^-}^{\lambda_r^+} \frac{-\ln(a_+(\zeta)a_-^*(\zeta))}{\mathcal{R} (\zeta;\lambda_s(\xi), \lambda_l^-)} \,\frac{\dif\zeta}{\zeta - z}
  + \int_{\lambda_l^-}^{\lambda_s(\xi)} \frac{\ln (r_+(\zeta))}{\mathcal{R}_+ (\zeta;\lambda_s(\xi), \lambda_l^-)} \, \frac{\dif\zeta}{\zeta - z} \right]  \right\}
  \end{aligned}
\end{equation}
with
\begin{equation}\label{A2deltaINF}
  \begin{aligned}
  \delta_{\infty }= \exp & \left\{ \frac{ \mathrm{i}}{2\pi} \left[ \left( \int_{-\infty}^{\lambda_r^-} + \int_{\lambda_r^+}^{\lambda_l^-} \right)  \frac{\ln (1-|r(\zeta)|^2)}{\mathcal{R} (\zeta;\lambda_s(\xi), \lambda_l^-) } \dif \zeta  \right. \right. \\
  & \left. \left. +  \int_{\lambda_r^-}^{\lambda_r^+} \frac{-\ln(a_+(\zeta)a_-^*(\zeta))}{\mathcal{R} (\zeta;\lambda_s(\xi), \lambda_l^-)} \dif \zeta
  + \int_{\lambda_l^-}^{\lambda_s(\xi)} \frac{\ln (r_+(\zeta))}{\mathcal{R}_+ (\zeta;\lambda_s(\xi), \lambda_l^-)} \dif \zeta \right]  \right\}.
  \end{aligned}
\end{equation}
This results in the following transformation
\begin{equation}
	M^{(3)}(z;x,t)=\delta _{\infty}^{\sigma_{3}} M^{(2)}(z;x,t) \delta(z) ^{-\sigma_{3}}
\end{equation}
and the following RHP.
\begin{rhp} \label{lrrhp} Find a $2 \times 2$ matrix-valued function $M^{(3)}(z;x,t)$ with the following properties:
\begin{enumerate} [label=(\roman*)]
  \item $M^{(3)}(z;x,t)$ is analytic in $z \in \mathbb{C} \backslash \Sigma^{(3)}$, where $\Sigma^{(3)}= [\lambda_l^-, \lambda_s(\xi)] \cup L_1 \cup L_2 \cup  L_3 \cup L_4.$
  \item $M^{(3)}(z;x,t)=I+\mathcal{O} (z^{-1}) $ as $z\to \infty$.
  \item $M^{(3)}(z;x,t)$ achieves the CBVs $M^{(3)}_{+}(z;x,t)$ and $M^{(3)}_{-}(z;x,t)$ on $\Sigma^{(3)}$ away from self-intersection points and
branch points that satisfy the jump condition $M^{(3)}_{+}(z;x,t)=M^{(3)}_{-}(z;x,t)V^{(3)}(z;x,t)$, where
\begin{equation}\label{A2V3}
  \begin{aligned}
V^{(3)}(z;x,t)&=\left\{  \begin{array}{ll}
   \begin{pmatrix} 0 & -1 \\ 1 & 0 \end{pmatrix}, ~&z \in (\lambda_l^-, \eta(\xi)),\\ \\
   \begin{pmatrix} 0 & -1 \\ 1 & \delta _{+}(z)\delta _{-}^{-1}(z)e^{- 2 \mathrm{i} t g_+(z)} \end{pmatrix}, ~&z \in (\eta(\xi), \lambda_s(\xi)),\\ \\
   \begin{pmatrix} 1 & 0 \\ r(z) \delta ^{-2}(z) e^{2\mathrm{i} tg(z)} & 1 \end{pmatrix}, ~ &z \in L_1, \\ \\
   \begin{pmatrix} 1 & \frac{-r^*(z)}{1-r(z)r^*(z)} \delta ^2(z) e^{-2\mathrm{i}t g(z)} \\0 & 1 \end{pmatrix},  ~&z \in L_2,\\ \\
   \begin{pmatrix} 1 & 0 \\ \frac{r(z)}{1-r(z)r^*(z)} \delta ^{-2}(z) e^{2\mathrm{i} tg(z)} & 1 \end{pmatrix},  ~&z \in L_3,\\ \\
   \begin{pmatrix} 1 & -r^*(z) \delta ^2(z) e^{-2\mathrm{i} tg(z)} \\ 0  & 1 \end{pmatrix}, ~ &z \in L_4. \end{array} \right.\\
  \end{aligned}
\end{equation}
\end{enumerate}
\end{rhp}
The above jump matrices $V^{(3)}(z;x,t)$ on the steepest descent contours converge to the identity matrix, while the jump matrices on the interval $(\lambda_l^-, \lambda_s(\xi))$
converge exponentially to a constant matrix. Furthermore, the latter convergence is uniform outside a fixed neighborhood $\mathcal{U}$ of $\lambda_s(\xi)$ due to
$g(z)=\mathcal{O}((z-\lambda_s(\xi))^{3/2})$ near $\lambda_s(\xi)$. As before, consider the limiting problem
\begin{equation}
  M^{(\mathrm{mod})}_+(z;x,t)=M^{(\mathrm{mod})}_-(z;x,t) V^{(\mathrm{mod})} \qquad  \mathrm{for}  \quad z \in (\lambda_l^-, \lambda_s(\xi)),
\end{equation}
where
\begin{equation}\label{A2vmod}
  V^{(\mathrm{mod})}=\begin{pmatrix} 0 & -1 \\ 1 & 0 \end{pmatrix}.
\end{equation}
This problem is exactly solvable and the solution is
\begin{equation} \label{A2mod}
  M^{(\mathrm{mod})}(z;x,t) = \mathcal{E} (z;\lambda_s(\xi),\lambda_l^-).
\end{equation}

The construction of the local parametrix inside $\mathcal{U}$ is the same as that in Section \ref{ABrsh},  except that the conformal mapping from $\mathcal{U}$ to a neighborhood of the origin is replaced by
\begin{equation}\label{B2f}
  \frac{2}{3}(-f(z))^{3/2}=\mathrm{i\, sgn}(\mathrm{Im}\,z)\,g(z)=(\lambda_s(\xi) - z)^{3/2}(z-\lambda_l^-)^{1/2}.
\end{equation}
Then we construct the global parametrix, which leads to a similar error RHP as RHP \ref{rhp-rare} (with $\hat{\delta}$ replaced by $\delta$) and thus the same error estimate $M^{(\mathrm{err})}_1(x,t)= \mathcal{O}(t^{-1})$.
Using equations (\ref{rec}), (\ref{A2ginf}), (\ref{A2deltaINF}) and (\ref{A2mod}), it is concluded that the long-time asymptotics of the defocusing NLS equation (\ref{(NLS)})
in the rarefaction wave region is given by
\begin{equation}
  q(x,t)=-\frac{x+2\lambda_l^- t}{3t}e^{-\mathrm{i}(2{\lambda_l^-}^2t+2\lambda_l^- x-x^2/t)/3}e^{-\mathrm{i}\phi_{lr}(\xi) }+\mathcal{O}(t^{-1}),
\end{equation}
where
\begin{equation} \label{Bphilr}
  \begin{aligned}
    \phi_{lr}(\xi)=\frac{1}{\pi}&\left[\left(\int_{-\infty}^{\lambda_r^-} +\int_{\lambda_r^+}^{\lambda_l^-}\right) \right.
   \frac{\ln (1-|r(\zeta)|^2)}{\sqrt{(\zeta-\lambda_s(\xi))(\zeta-\lambda_l^-)} } \dif \zeta\\
   &+ \left. \int_{\lambda_r^-}^{\lambda_r^+} \frac{-\ln(a_+(\zeta)a_-^*(\zeta))}{\sqrt{(\zeta-\lambda_s(\xi))(\zeta-\lambda_l^-)}} \dif \zeta
   + \int_{\lambda_l^-}^{\lambda_s(\xi)} \frac{\arg (r_+(\zeta))}{\sqrt{(\lambda_s(\xi)-\zeta)(\zeta-\lambda_l^-)}} \dif \zeta\right].
  \end{aligned}
\end{equation}

\subsubsection{Vacuum region: $-2\lambda_l^-<\xi <-2\lambda_r^+$}
\
\newline
\indent

\begin{figure}[htbp]
  \centering
  \subfigure[]{\includegraphics[width=5.5cm]{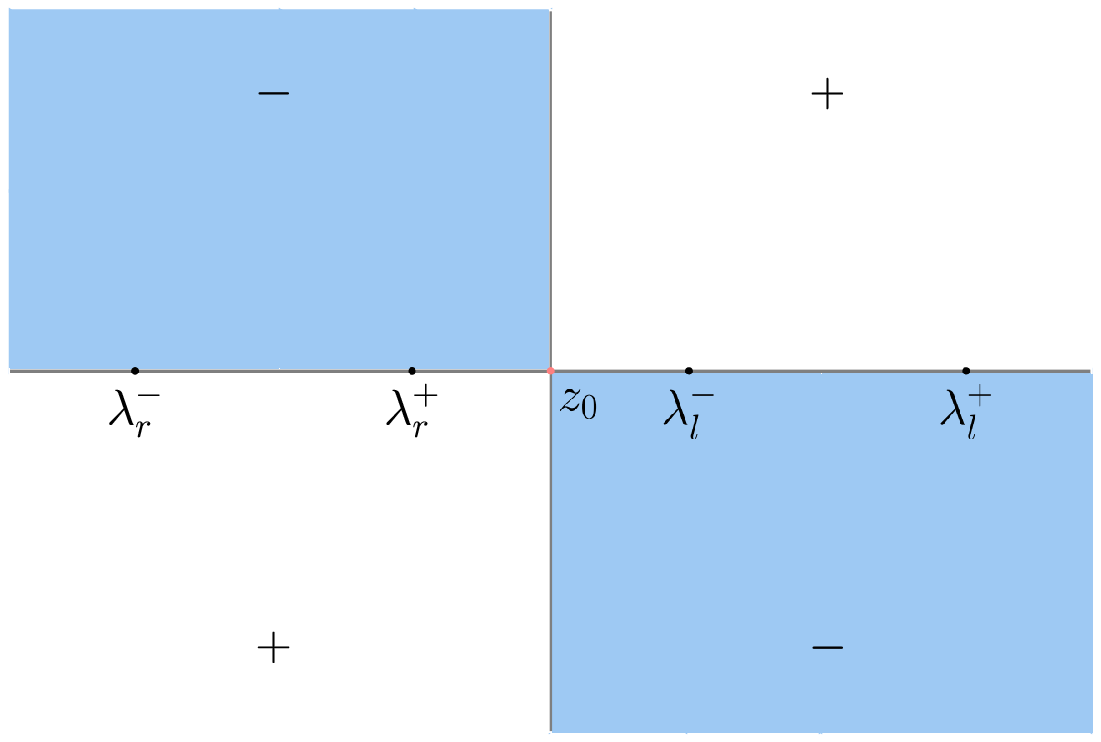}}
  \subfigure[]{\includegraphics[width=5cm]{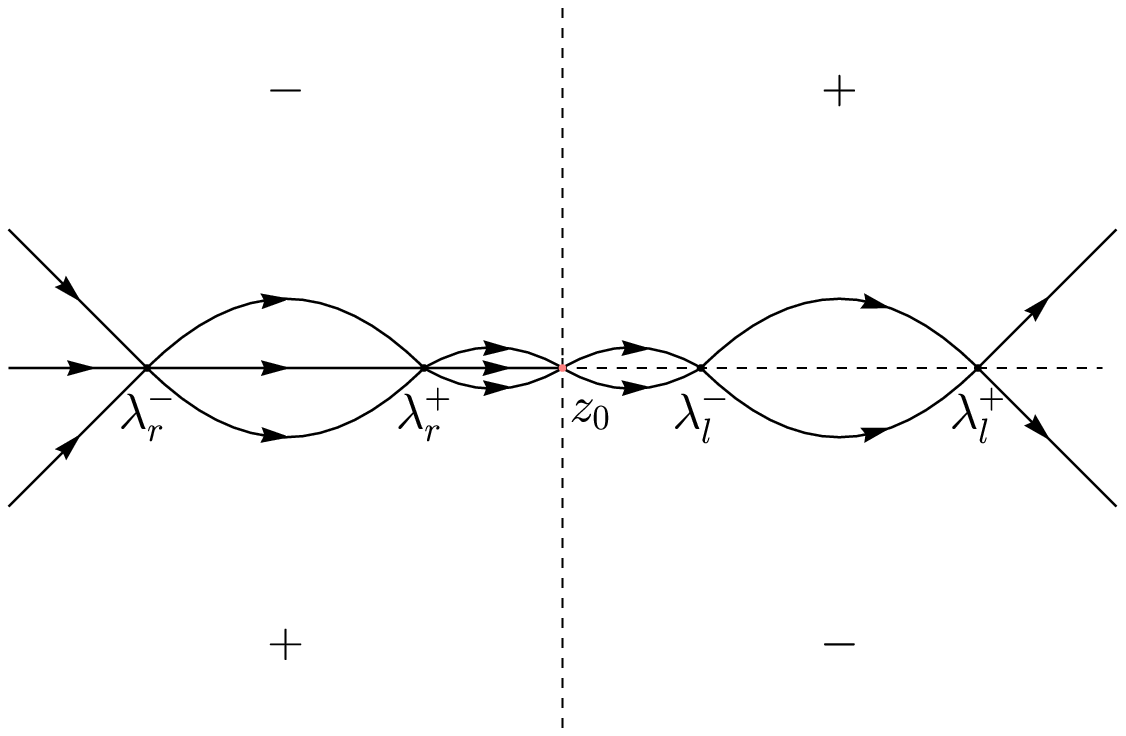}}
  \subfigure[]{\includegraphics[width=5.5cm]{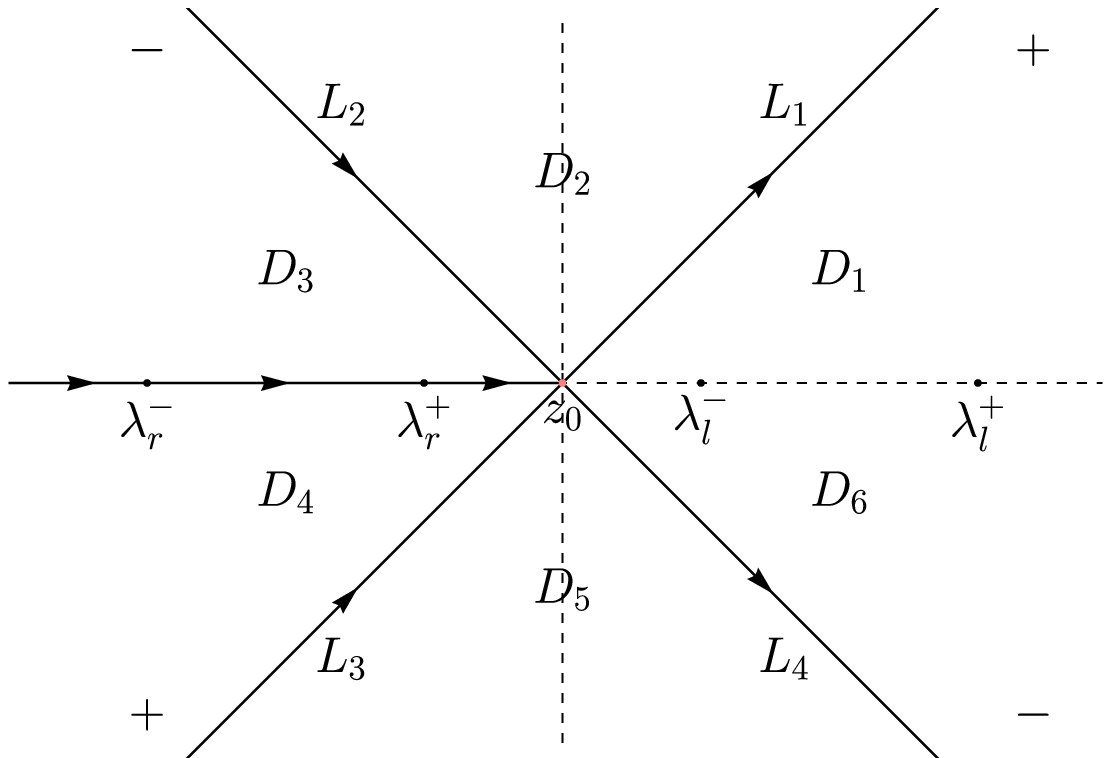}}
\caption{{\protect\small (a) The sign structure of Im$(g(z))$; (b) Opening lenses; (c) The jump contours of $M^{(2)}(z;x,t)$.}}
\label{figB3}
\end{figure}

If $\xi > -2\lambda_l^-$, the stationary phase point $\lambda_s(\xi)$ of the rarefaction $g$-function is less than $\lambda_l^-$, which leads to exponentially large  diagonal entries  $e^{\pm 2 \mathrm{i}t g_+(z)}$  of the jump matrix $V^{(2)}(z;x,t)$ on $(\lambda_s(\xi), \lambda_l^-)$ in (\ref{A2V2}).
It seems that a new $g$-function should be introduced. However, the original phase function $\theta(z)$ defined by (\ref{theta}) is sufficient since the off-diagonal entries of the jump matrices $V^{(2)}(z;x,t)$ across all steepest contours through the original stationary phase point $z_0(\xi)=-{\xi}/{2}$
decay exponentially to zero as $t \to \infty$. This implies that the asymptotic behavior in this region is essentially the same as zero background \cite{deift1994long}. Furthermore, $z_0(\xi)$ determines the boundaries of this vacuum region:
\begin{equation} \label{B3b}
  \lambda_r^+ < z_0(\xi) < \lambda_l^- \qquad \mathrm{iff}  \qquad  -2\lambda_l^- = \xi_{B2}<\xi < \xi_{B3} = -2\lambda_r^+.
\end{equation}

At this time, open lenses through the stationary phase point $z_0(\xi)$ (as shown in Figure \ref{figB3}) by the transformations $M(z;x,t) =M^{(1)}(z;x,t) \to M^{(2)}(z;x,t)$, which leads to the following RHP.

\begin{rhp}
  Find a $2 \times 2$ matrix-valued function $M^{(2)}(z;x,t)$ with the following properties:
\begin{enumerate} [label=(\roman*)]
  \item  $M^{(2)}(z;x,t)$ is analytic in $z \in \mathbb{C} \backslash \Sigma^{(2)}$, where $\Sigma^{(2)}=(-\infty, z_0(\xi)] \cup L_1 \cup L_2 \cup  L_3 \cup L_4.$
  \item  $M^{(2)}(z;x,t)=I+\mathcal{O} (z^{-1}) $ as $z\to \infty$.
  \item  $M^{(2)}(z;x,t)$ achieves the CBVs  $M^{(2)}_{+}(z;x,t)$ and $M^{(2)}_{-}(z;x,t)$  on $\Sigma^{(2)}$ away from self-intersection points
  and branch points that satisfy the jump condition $M^{(2)}_{+}(z;x,t)=M^{(2)}_{-}(z;x,t)V^{(2)}(z;x,t)$, where
  \begin{equation}\label{A3V2}
    \begin{aligned}
    V^{(2)}(z;x,t)&=\left\{  \begin{array}{ll}
     (1-r(z)r^{*}(z))^{\sigma_{3}},  ~& z\in (-\infty, \lambda_r^-) \cup (\lambda_r^+, z_0(\xi)) , \\ \\
     (a_+(z) a_-^*(z))^{-\sigma_{3}},  ~&z \in (\lambda_r^-, \lambda_r^+), \\ \\
     \begin{pmatrix} 1 & 0 \\ r(z)e^{2\mathrm{i} t\theta(z) } & 1 \end{pmatrix}, ~ &z \in L_1, \\ \\
     \begin{pmatrix} 1 & \frac{-r^*(z)}{1-r(z)r^*(z)}e^{-2\mathrm{i} t\theta(z)} \\0 & 1 \end{pmatrix},  ~&z \in L_2,\\ \\
     \begin{pmatrix} 1 & 0 \\ \frac{r(z)}{1-r(z)r^*(z)}e^{2\mathrm{i} t\theta(z)} & 1 \end{pmatrix},  ~&z \in L_3,\\ \\
     \begin{pmatrix} 1 & -r^*(z)e^{-2\mathrm{i} t \theta(z)} \\ 0  & 1 \end{pmatrix}, ~ &z \in L_4.\end{array} \right.\\
    \end{aligned}
  \end{equation}
\end{enumerate}
\end{rhp}
This RHP is roughly the same as that of zero background \cite{deift1994long}, except for the jump matrix on the interval $(\lambda_r^-, \lambda_r^+)$. Introduce a scalar function
$\delta(z)$ as
\begin{equation} \label{A3delta}
  \delta(z) = \exp \left\{   \frac{1}{2 \pi \mathrm{i}} \left[ \left(\int_{-\infty}^{\lambda_r^-} + \int_{\lambda_r^+}^{z_0(\xi)} \right) \frac{\ln (1-|r(\zeta)|^2)}{\zeta - z} \dif \zeta
 + \int_{\lambda_r^-}^{\lambda_r^+} \frac{-\ln(a_+(\zeta)a_-^*(\zeta))}{\zeta - z} \dif \zeta  \right]  \right\}
\end{equation}
which is analytic in $z \in \mathbb{C} \backslash (-\infty, z_0(\xi)]$ with asymptotic behavior $\delta(z)=1+\mathcal{O} (z^{-1}) $ as $z\to \infty$. $\delta(z)$ achieves the CBVs on $(-\infty, z_0(\xi))$ with the jump condition
\begin{equation*}
  \begin{aligned}
\left\{  \begin{array}{ll}
   \delta _{+}(z)=(1-r(z)r^{*}(z))\delta _{-}(z),  ~& z\in (-\infty, \lambda_r^-) \cup (\lambda_r^+, z_0(\xi)), \\ \\
   \delta _{+}(z)=(a_+(z) a_-^*(z))^{-1}\delta _{-}(z),  ~& z\in (\lambda_r^-, \lambda_r^+).
  \end{array} \right.\\
  \end{aligned}
\end{equation*}

This results in the following transformation (noting that $\delta_{\infty}=1$)
\begin{equation}
	M^{(3)}(z;x,t)= M^{(2)}(z;x,t) \delta(z) ^{-\sigma_{3}}
\end{equation}
and the following RHP.
\begin{rhp} Find a $2 \times 2$ matrix-valued function $M^{(3)}(z;x,t)$ with the following properties:
\begin{enumerate} [label=(\roman*)]
  \item $M^{(3)}(z;x,t)$ is analytic in $z \in \mathbb{C} \backslash \Sigma^{(3)}$, where $\Sigma^{(3)}= L_1 \cup L_2 \cup  L_3 \cup L_4.$
  \item $M^{(3)}(z;x,t)=I+\mathcal{O} (z^{-1}) $ as $z\to \infty$.
  \item $M^{(3)}(z;x,t)$ achieves CBVs $M^{(3)}_{+}(z;x,t)$ and $M^{(3)}_{-}(z;x,t)$ on $\Sigma^{(3)}$ away from the self-intersection point with jump condition $M^{(3)}_{+}(z;x,t)=M^{(3)}_{-}(z;x,t)V^{(3)}(z;x,t)$, where
\begin{equation}\label{A3V3}
  \begin{aligned}
V^{(3)}(z;x,t)&=\left\{  \begin{array}{ll}
   \begin{pmatrix} 1 & 0 \\ r(z) \delta ^{-2}(z) e^{2\mathrm{i} t \theta(z)} & 1 \end{pmatrix}, ~ &z \in L_1, \\ \\
   \begin{pmatrix} 1 & \frac{-r^*(z)}{1-r(z)r^*(z)} \delta ^2 (z) e^{-2\mathrm{i}t \theta(z)} \\0 & 1 \end{pmatrix},  ~&z \in L_2,\\ \\
   \begin{pmatrix} 1 & 0 \\ \frac{r(z)}{1-r(z)r^*(z)} \delta ^{-2} (z) e^{2\mathrm{i} t\theta(z)} & 1 \end{pmatrix},  ~&z \in L_3,\\ \\
   \begin{pmatrix} 1 & -r^*(z) \delta ^2(z) e^{-2\mathrm{i} t\theta(z)} \\ 0  & 1 \end{pmatrix}, ~ &z \in L_4. \end{array} \right.\\
  \end{aligned}
\end{equation}
\end{enumerate}
\end{rhp}
The above RHP for $M^{(3)}(z;x,t)$ is the same as in the case of zero background (see \cite{deift1994long}), except the behavior of $\delta(z)$ as $z \to z_0(\xi)=-\xi/2$:
\begin{equation}
  \delta(z)=\left(z+\frac{\xi}{2}\right)^{\mathrm{i}\nu(-\xi/2)}e^{\widetilde{\chi}(z)},
\end{equation}
where
\begin{equation} \label{B3nu}
  \begin{aligned}
    &\nu(-\xi/2) =-\frac{1}{2\pi} \ln (1-|r(-\xi/2)|^2) \in \mathbb{R},\\
    &\widetilde{\chi}(z) =-\frac{1}{2 \pi \mathrm{i}}\left[\left(\int_{-\infty}^{\lambda_r^-} + \int_{\lambda_r^+}^{-\xi/2} \right) \ln ({z - \zeta }) \,\dif \, {\ln (1-|r(\zeta)|^2)}
    + \int_{\lambda_r^-}^{\lambda_r^+} \frac{\ln(a_+(\zeta)a_-^*(\zeta))}{\zeta + \xi/2} \dif \zeta \right] \in \mathrm{i}\mathbb{R}.
  \end{aligned}
\end{equation}
\par
Thus it is concluded that the long-time asymptotic behavior of the defocusing NLS equation (\ref{(NLS)}) in the vacuum region is formulated by
\begin{equation}
   q(x,t)=\frac{\nu(-\xi)}{\sqrt{t}}e^{\mathrm{i}{t}{\xi^2}/{2}-\mathrm{i} \nu(-\xi/2) \ln t + \mathrm{i}\phi_{va} (\xi)}
   +o(t^{-\frac{1}{2}})=\mathcal{O}(t^{-\frac{1}{2}}),
\end{equation}
where
\begin{equation} \label{phiva}
   \phi_{va}(\xi)=-2 \ln2 ~ \nu(-\xi/2) + \frac{\pi}{4} + \arg \Gamma(\mathrm{i} \nu(-\xi/2))- \arg(r(-\xi/2))-2 \mathrm{i} \widetilde{\chi}(-\xi/2).
\end{equation}

\subsubsection{Rarefaction wave region: $-2\lambda_r^+ <\xi <-\frac{\lambda_r^+ + 3\lambda_r^-}{2}$} \label{Arrare}
\
\newline
\indent

\begin{figure}[htbp]
  \centering
  \subfigure[]{\includegraphics[width=5.5cm]{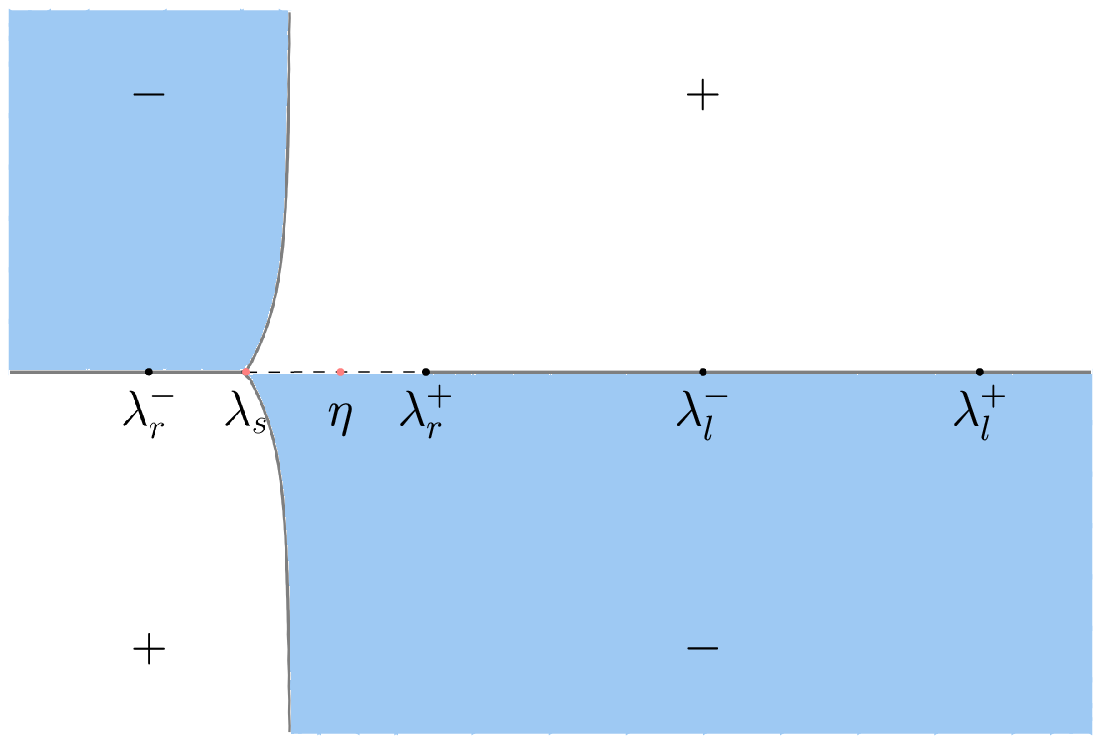}}
  \subfigure[]{\includegraphics[width=5cm]{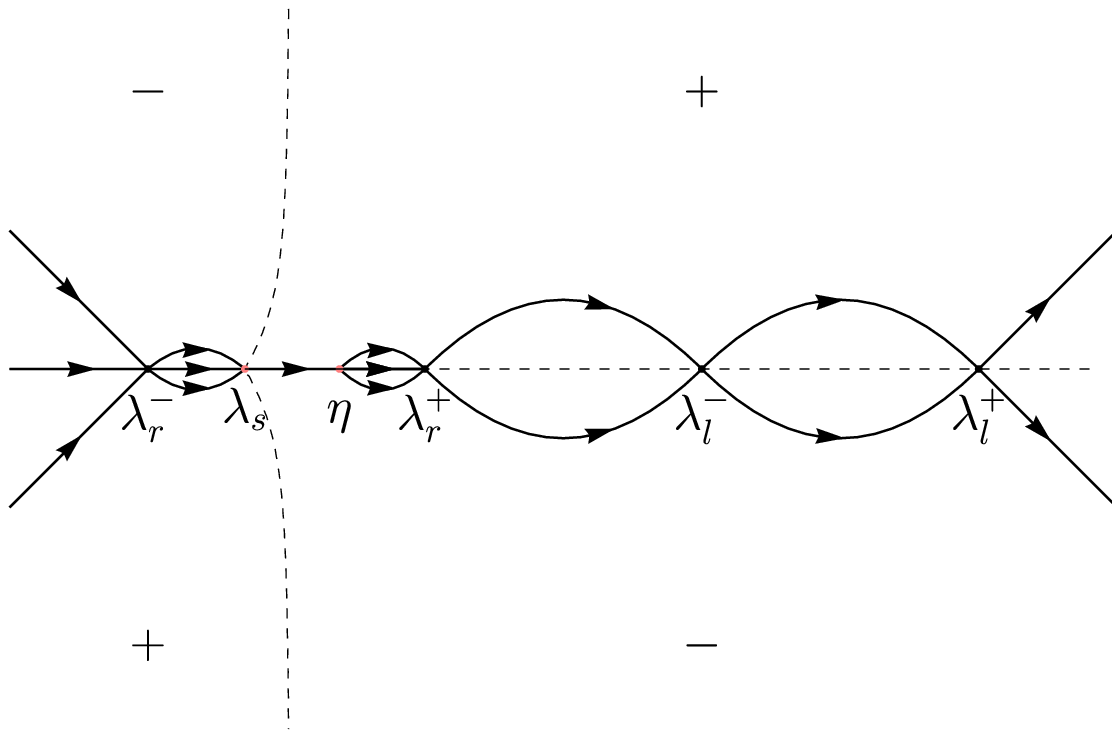}}
  \subfigure[]{\includegraphics[width=5.5cm]{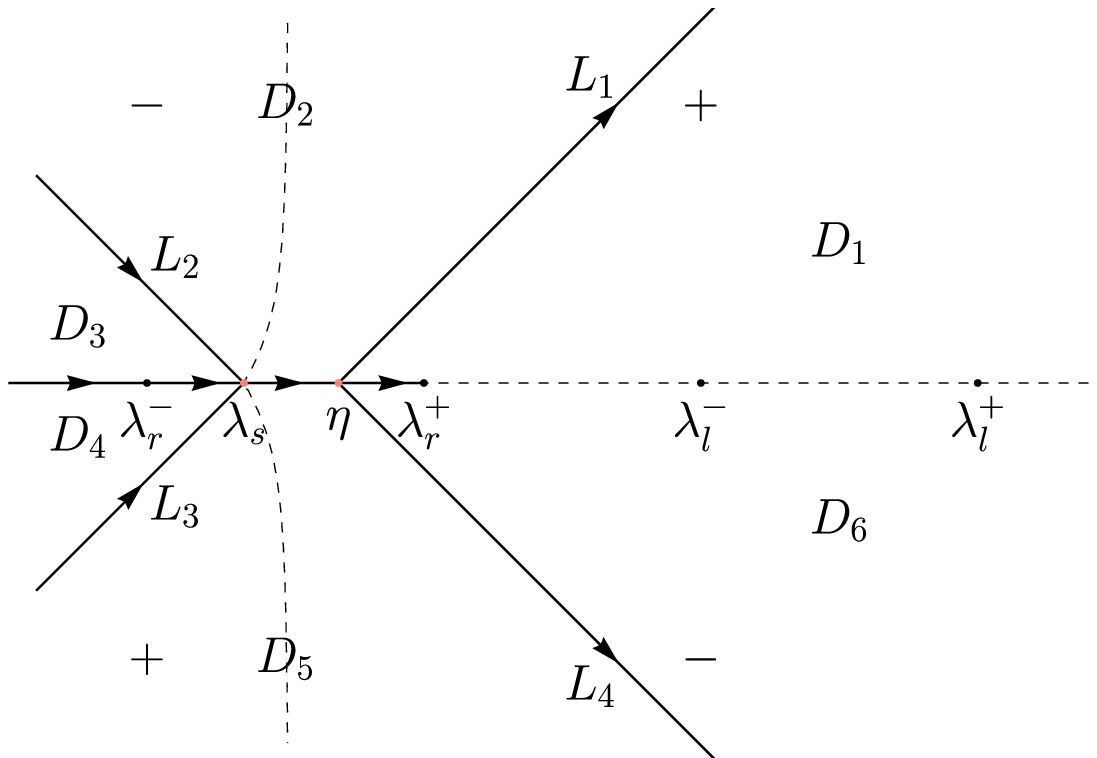}}
\caption{{\protect\small (a) The sign structure of Im$(g(z))$; (b) Opening lenses; (c) The jump contours of $M^{(2)}(z;x,t)$.}}
\label{figB4}
\end{figure}

As $\xi$ increases such that $\xi > -2\lambda_r^+ $, the stationary phase point $z_0(\xi)$ moves inside $\mathcal{I}_r$,
which leads to exponential oscillation  $e^{\pm 2 \mathrm{i} t \theta(z)}$ in the jump matrix $V^{(2)}(z;x,t)$ on the interval $(z_0(\xi), \lambda_r^+) \subset
\mathcal{I}_r$. For a new soft edge $\lambda_s(\xi)=-(\lambda_r^+ + 2 \xi)/3$, introduce the one-band $g$-function similar to (\ref{A2g}) as
\begin{equation} \label{A4g}
  g(z)= \mathcal{R} (z;\lambda_r^+,\lambda_s(\xi))(z-\lambda_s(\xi)),
\end{equation}
with the differential
\begin{equation} \label{A4dg}
  \dif g(z)= 2 \frac{(z-\lambda_s(\xi))(z-\eta(\xi))}{\mathcal{R} (z;\lambda_r^+,\lambda_s(\xi))} \,\dif z,
\end{equation}
where  the other stationary phase point $\eta(\xi)=(\lambda_s(\xi)+3\lambda_r^+)/4$ is determined by $\int_{\lambda_s(\xi)}^{\lambda_r^+} \,\dif g(z) =0$.
Furthermore, $g(z)$ is analytic in $z \in \mathbb{C} \backslash [\lambda_s(\xi), \lambda_r^+]$ with
\begin{equation}\label{A4ginf}
  g_{\infty}=-(2{\lambda_r^+}^2+2\lambda_r^+ \xi-\xi^2)/6.
\end{equation}
The boundaries of this rarefaction wave region are determined by
\begin{equation} \label{B4b}
  \lambda_r^-  < \lambda_s(\xi) < \lambda_r^+ \qquad \mathrm{iff}  \qquad -2\lambda_r^+ = \xi_{B3} <\xi < \xi_{B4}= -\frac{\lambda_r^+ + 3\lambda_r^-}{2}.
\end{equation}

Similar to the previous rarefaction wave region, open lenses from $\lambda_s(\xi)$ and $\eta(\xi)$ as shown in Figure \ref{figB4}, which leads to the following new RHP.

\begin{rhp}
  Find a $2 \times 2$ matrix-valued function $M^{(2)}(z;x,t)$ with the following properties:
\begin{enumerate} [label=(\roman*)]
  \item  $M^{(2)}(z;x,t)$ is analytic in $z \in \mathbb{C} \backslash \Sigma^{(2)}$, where $\Sigma^{(2)}=(-\infty, \lambda_r^+] \cup L_1 \cup L_2 \cup  L_3 \cup L_4.$
  \item  $M^{(2)}(z;x,t)=I+\mathcal{O} (z^{-1}) $ as $z\to \infty$.
  \item  $M^{(2)}(z;x,t)$ achieves the CBVs  $M^{(2)}_{+}(z;x,t)$ and $M^{(2)}_{-}(z;x,t)$  on $\Sigma^{(2)}$ away from self-intersection  points and
  branch points that satisfy the jump condition $M^{(2)}_{+}(z;x,t)=M^{(2)}_{-}(z;x,t)V^{(2)}(z;x,t)$, where
  \begin{equation}\label{A4V2}
    \begin{aligned}
    V^{(2)}(z;x,t)&=\left\{  \begin{array}{ll}
     (1-r(z)r^{*}(z))^{\sigma_{3}},  ~& z\in (-\infty, \lambda_r^-), \\ \\
     (a_+(z) a_-^*(z))^{-\sigma_{3}},  ~&z \in  (\lambda_r^-, \lambda_s(\xi)), \\ \\
     \begin{pmatrix} \frac{e^{2 \mathrm{i} t g_+(z)}}{a_+(z) a_-^*(z)} & -1 \\ 1 & 0 \end{pmatrix}, ~&z \in (\lambda_s(\xi),\eta(\xi)), \\ \\
     \begin{pmatrix} 0 & -1 \\ 1 & 0 \end{pmatrix}, ~&z \in (\eta(\xi),\lambda_r^+), \\ \\
     \begin{pmatrix} 1 & 0 \\ r(z)e^{2\mathrm{i} tg(z)} & 1 \end{pmatrix}, ~ &z \in L_1, \\ \\
     \begin{pmatrix} 1 & \frac{-r^*(z)}{1-r(z)r^*(z)}e^{-2\mathrm{i} tg(z)} \\0 & 1 \end{pmatrix},  ~&z \in L_2,\\ \\
     \begin{pmatrix} 1 & 0 \\ \frac{r(z)}{1-r(z)r^*(z)}e^{2\mathrm{i} tg(z)} & 1 \end{pmatrix},  ~&z \in L_3,\\ \\
     \begin{pmatrix} 1 & -r^*(z)e^{-2\mathrm{i} tg(z)} \\ 0  & 1 \end{pmatrix}, ~ &z \in L_4.\end{array} \right.\\
    \end{aligned}
  \end{equation}
\end{enumerate}
\end{rhp}

As before, introduce a scalar function
\begin{equation} \label{A4delta}
  \begin{aligned}
  \delta(z) = \exp & \left\{ \frac{\mathcal{R} (z;\lambda_r^+,\lambda_s(\xi))}{2 \pi \mathrm{i}}  \left(\int_{-\infty}^{ \lambda_r^-}  \frac{\ln (1-|r(\zeta)|^2)}{\mathcal{R} (\zeta;\lambda_r^+,\lambda_s(\xi)) } \,\frac{\dif\zeta}{\zeta - z}    + \int_{\lambda_r^-}^{\lambda_s(\xi)} \frac{-\ln(a_+(\zeta)a_-^*(\zeta))}{\mathcal{R} (\zeta;\lambda_r^+,\lambda_s(\xi))} \,\frac{\dif\zeta}{\zeta - z} \right)  \right\}
  \end{aligned}
\end{equation}
with the following properties:
\begin{enumerate} [label=(\roman*)]
  \item $\delta(z)$ is analytic in $z \in \mathbb{C} \backslash (-\infty, \lambda_r^+ ]$;
  \item $\delta(z)=\delta_\infty+\mathcal{O} (z^{-1}) $ as $z\to \infty$, where
    \begin{equation}\label{A4deltaINF}
     \begin{aligned}
     \delta_{\infty }= \exp & \left\{ \frac{ \mathrm{i}}{2\pi} \left(\int_{-\infty}^{ \lambda_r^-}  \frac{\ln (1-|r(\zeta)|^2)}{\mathcal{R} (\zeta;\lambda_r^+,\lambda_s(\xi)) } \dif \zeta
     + \int_{\lambda_r^-}^{\lambda_s(\xi)} \frac{-\ln(a_+(\zeta)a_-^*(\zeta))}{\mathcal{R} (\zeta;\lambda_r^+,\lambda_s(\xi))} \dif \zeta \right)  \right\};
     \end{aligned}
    \end{equation}
  \item $\delta(z)$ achieves CBVs on  $(-\infty, \lambda_r^+ )$ satisfying the jump conditions
  \begin{equation*}
    \begin{aligned}
  \left\{  \begin{array}{ll}
     \delta _{+}(z)=(1-r(z)r^{*}(z))\delta _{-}(z),  ~& z\in (-\infty, \lambda_r^-) , \\ \\
     \delta _{+}(z)=(a_+(z) a_-^*(z))^{-1}\delta _{-}(z),  ~& z\in (\lambda_r^-, \lambda_s(\xi)), \\ \\
     \delta _{+}(z)\delta _{-}(z)=1,  ~& z\in (\lambda_s(\xi), \lambda_r^+).
    \end{array} \right.\\
    \end{aligned}
  \end{equation*}
\end{enumerate}
This results in the following transformation
\begin{equation}
	M^{(3)}(z;x,t)=\delta _{\infty}^{\sigma_{3}} M^{(2)}(z;x,t) \delta(z) ^{-\sigma_{3}}
\end{equation}
and the RHP below.
\begin{rhp} \label{rrrhp} Find a $2 \times 2$ matrix-valued function $M^{(3)}(z;x,t)$ with the following properties:
\begin{enumerate} [label=(\roman*)]
  \item $M^{(3)}(z;x,t)$ is analytic in $z \in \mathbb{C} \backslash \Sigma^{(3)}$, where $\Sigma^{(3)}= [\lambda_s(\xi), \lambda_r^+] \cup L_1 \cup L_2 \cup  L_3 \cup L_4.$
  \item $M^{(3)}(z;x,t)=I+\mathcal{O} (z^{-1}) $ as $z\to \infty$.
  \item $M^{(3)}(z;x,t)$ achieves the CBVs $M^{(3)}_{+}(z;x,t)$ and $M^{(3)}_{-}(z;x,t)$ on $\Sigma^{(3)}$ away from self-intersection points and
branch points that satisfy the jump condition $M^{(3)}_{+}(z;x,t)=M^{(3)}_{-}(z;x,t)V^{(3)}(z;x,t)$, where
\begin{equation}\label{A4V3}
  \begin{aligned}
V^{(3)}(z;x,t)&=\left\{  \begin{array}{ll}
   \begin{pmatrix} \frac{\delta ^{-1}_{+}(z)\delta _{-}(z)}{a_+(z) a_-^*(z)}  e^{2 \mathrm{i} t g_+(z)} & -1 \\ 1 & 0 \end{pmatrix}, ~&z \in (\lambda_s(\xi),\eta(\xi)), \\ \\
   \begin{pmatrix} 0 & -1 \\ 1 & 0 \end{pmatrix}, ~&z \in (\eta(\xi),\lambda_r^+), \\ \\
   \begin{pmatrix} 1 & 0 \\ r(z) \delta ^{-2}(z) e^{2\mathrm{i} tg(z)} & 1 \end{pmatrix}, ~ &z \in L_1, \\ \\
   \begin{pmatrix} 1 & \frac{-r^*(z)}{1-r(z)r^*(z)} \delta ^2(z) e^{-2\mathrm{i}t g(z)} \\0 & 1 \end{pmatrix},  ~&z \in L_2,\\ \\
   \begin{pmatrix} 1 & 0 \\ \frac{r(z)}{1-r(z)r^*(z)} \delta ^{-2}(z) e^{2\mathrm{i} tg(z)} & 1 \end{pmatrix},  ~&z \in L_3,\\ \\
   \begin{pmatrix} 1 & -r^*(z) \delta ^2(z) e^{-2\mathrm{i} tg(z)} \\ 0  & 1 \end{pmatrix}, ~ &z \in L_4. \end{array} \right.\\
  \end{aligned}
\end{equation}
\end{enumerate}
\end{rhp}

The jump matrices $V^{(3)}(z;x,t)$ uniformly converge to the identity matrix or a constant matrix outside a fixed neighborhood $\mathcal{U}$ of $\lambda_s(\xi)$,
while the convergence is not uniform inside $\mathcal{U}$ due to the exponential phase function $g(z)=\mathcal{O}((z-\lambda_s(\xi))^{3/2})$. So it is necessary to construct an outer
parametrix
\begin{equation} \label{A4mod}
  M^{(\mathrm{mod})}(z;x,t) = \mathcal{E} (z;\lambda_r^+, \lambda_s(\xi))
\end{equation}
solving the limiting problem
\begin{equation}
  M^{(\mathrm{mod})}_+(z;x,t)=M^{(\mathrm{mod})}_-(z;x,t) V^{(\mathrm{mod})} \qquad  \mathrm{for}  \quad z \in (\lambda_s(\xi),\lambda_r^+),
\end{equation}
where
\begin{equation}\label{A4vmod}
  V^{(\mathrm{mod})}=\begin{pmatrix} 0 & -1 \\ 1 & 0 \end{pmatrix}.
\end{equation}
The construction of the local parametrix inside $\mathcal{U}$ is the same as that in Section \ref{ALSH},  except that the conformal mapping from $\mathcal{U}$ to a neighborhood of the origin is replaced by
\begin{equation}\label{B4f}
  \frac{2}{3}f^{3/2}(z)=-\mathrm{i\, sgn}(\mathrm{Im}\,z)\,g(z)=( z- \lambda_s(\xi))^{3/2}(\lambda_r^+ - z)^{1/2}.
\end{equation}
Then we construct the global parametrix, which leads to a similar error RHP as RHP \ref{a2err} (with $\hat{\delta}$ replaced by $\delta$) and thus the same error estimate $M^{(\mathrm{err})}_1(x,t)= \mathcal{O}(t^{-1})$.
Using equations (\ref{rec}), (\ref{A4ginf}), (\ref{A4deltaINF}) and (\ref{A4mod}), it is concluded that the long-time asymptotics of the defocusing NLS equation (\ref{(NLS)})
in the rarefaction wave region is formulated by
\begin{equation}
  q(x,t)=\frac{\xi+2\lambda_r^+}{3}e^{-\mathrm{i}t(2{\lambda_r^+}^2+2\lambda_r^+ \xi-\xi^2)/3}e^{-\mathrm{i}\phi_{rr}(\xi) }+\mathcal{O}(t^{-1}),
\end{equation}
where
\begin{equation} \label{phirr}
   \phi_{rr}(\xi)=\frac{1}{\pi} \left(\int_{-\infty}^{\lambda_r^-} \frac{\ln (1-|r(\zeta)|^2)}{\sqrt{(\zeta-\lambda_r^+)(\zeta-\lambda_s(\xi))} } \dif \zeta
   +\int_{\lambda_r^-}^{\lambda_s(\xi)} \frac{-\ln(a_+(\zeta)a_-^*(\zeta))}{\sqrt{(\zeta-\lambda_r^+)(\zeta-\lambda_s(\xi))}} \dif \zeta \right).
\end{equation}

\subsubsection{The right plane wave region: $\xi > -\frac{\lambda_r^+ + 3\lambda_r^-}{2}$} \label{Arpw}
\
\newline
\indent
As $\xi$ increases such that $\xi > -\frac{\lambda_r^+ + 3\lambda_r^-}{2}$, the stationary phase point $\lambda_s(\xi)$  moves outside $\mathcal{I}_r$, which leads to exponentially
large diagonal entries $e^{\pm 2 \mathrm{i}t g_+(z)}$ of the jump matrix $V^{(2)}(z;x,t)$ on $(\lambda_s(\xi), \lambda_r^-)$.
Corresponding to the explicit eigenfunction $\Psi_r^{\mathrm{p}}(z;x,t)$ defined by (\ref{psip}), the $g$-function is the same as (\ref{A5g}) for the right plane wave regions in all six cases in (\ref{Classification}). As a result, the leading-order asymptotics in the right plane wave region  are consistent with the initial condition for $t = 0, x > 0$, up to the phase shift  $e^{-\mathrm{i}\phi_{rp}(\xi)}$.
The factorizations (\ref{f1})-(\ref{f4}) imply that only jump matrices on  $\mathcal{I}_r$  (sometimes consisting of $\mathcal{I}_r \backslash \overline{\mathcal{I}_l} $  and $\mathcal{I}_l \cap \mathcal{I}_r $) contribute to the asymptotic behaviors.
The boundary of the right plane wave region is characterized by
\begin{equation} \label{B5b}
  \eta_- (\xi) < \lambda_r^- \qquad \mathrm{iff}  \qquad  \xi > \xi_{B4}=-\frac{\lambda_r^+ + 3\lambda_r^-}{2},
\end{equation}
where $\eta_- (\xi)$ is given by (\ref{rpweta}). All transformations from $M(z;x,t)$ to  $M^{(\mathrm{err})}(z;x,t)$ are the same in the six cases, since the stationary phase point $\eta_-(\xi)$ of the right plane wave $g$-function
lies on the left side of both intervals $\mathcal{I}_l $ and $\mathcal{I}_r$, and the upper/lower matrix factorizations in (\ref{f1}) and (\ref{f2}) ((\ref{f3}) and (\ref{f4})) are the same.
Thus the unified form of the long-time asymptotics of the defocusing NLS equation (\ref{(NLS)}) in the right plane wave region is given by (\ref{A5}).

\subsection{\rm Case C:  $\lambda_r^+>\lambda_l^+ >\lambda_r^->\lambda_l^-$}
\
\newline
\indent
In this case, the two intervals $\mathcal{I}_l$ and $\mathcal{I}_r$ intersect but are not completely contained within each other, so all the factorizations in
(\ref{f1})-(\ref{f4}) should be considered. In comparison with Case A, only the positions of points $ \lambda_l^+$ and $\lambda_r^-$ are swapped. This means that this case is very similar to Case A, except that the center region becomes a plane wave region, as we will see later.
As the self-similar variable $\xi=x/t$ increases, the stationary phase points of the corresponding
$g$-functions change continuously at different intervals on $\mathbb{R}$, which implies that there are five different regions: the left plane wave region, the left
dispersive shock wave region, the middle plane wave region, the right dispersive shock wave region and the right plane wave region.

\subsubsection{The left plane wave region: $\xi<-\frac{2\lambda_r^++\lambda_l^++\lambda_l^-}{2}+\frac{2(\lambda_r^+ - \lambda_l^+)
(\lambda_r^+ - \lambda_l^-)}{-2\lambda_r^+ + \lambda_l^+ + \lambda_l^-}$}\label{Clpw}
\
\newline
\indent
In Section \ref{Blpw}, we have given the uniform form of the left plane wave region and all that remains is to determine the boundary of the region:
\begin{equation} \label{C1b}
  \eta_+ (\xi) > \lambda_r^+ \qquad \mathrm{iff}  \qquad  \xi< \xi_{C1}= -\frac{2\lambda_r^++\lambda_l^++\lambda_l^-}{2}+\frac{2(\lambda_r^+ - \lambda_l^+)
  (\lambda_r^+ - \lambda_l^-)}{-2\lambda_r^+ + \lambda_l^+ + \lambda_l^-},
\end{equation}
where $\eta_+ (\xi) $ is given by (\ref{lpweta}).

\subsubsection{The left dispersive shock wave region: $-\frac{2\lambda_r^++\lambda_l^++\lambda_l^-}{2}+\frac{2(\lambda_r^+ - \lambda_l^+)
(\lambda_r^+ - \lambda_l^-)}{-2\lambda_r^+ + \lambda_l^+ + \lambda_l^-} < \xi < -\frac{\lambda_r^++2\lambda_l^++\lambda_l^-}{2}$} \label{Clsh}
\
\newline
\indent
When the stationary phase point $\eta_+ (\xi)$ moves inside $\mathcal{I}_r$, the exponential oscillation $e^{\pm 2 \mathrm{i} t \theta(z)}$ appears in the jump matrix $V^{(2)}(z;x,t)$
on the interval $(\eta_+(\xi), \lambda_r^+) \subset \mathcal{I}_r$. To remove this oscillation, the modified $g$-function defined by (\ref{A2G}) should be introduced.
The boundaries of this region are characterized by degeneration below
\begin{equation}\label{C2b}
  \begin{aligned}
  &\xi \to \xi_{C1}= v_2(\lambda_r^+, \lambda_r^+, \lambda_l^+, \lambda_l^-)= -\frac{2\lambda_r^++\lambda_l^++\lambda_l^-}{2}+\frac{2(\lambda_r^+ - \lambda_l^+)
  (\lambda_r^+ - \lambda_l^-)}{-2\lambda_r^+ + \lambda_l^+ + \lambda_l^-}, \quad &\mathrm{as} \quad \lambda_s(\xi) \to \lambda_r^+,\\
  &\xi \to  \xi_{C2}= v_2(\lambda_r^+, \lambda_l^+, \lambda_l^+, \lambda_l^-)=-\frac{\lambda_r^++2\lambda_l^++\lambda_l^-}{2}, \quad &\mathrm{as} \quad \lambda_s(\xi) \to \lambda_l^+,
  \end{aligned}
\end{equation}
where $v_2$ is the Whitham velocity given by (\ref{vv2}).

Because of the same $g$-function, we will summarize the RHP analysis of the left dispersive shock wave region.  In a similar way, open lenses from intervals to the steepest descent contours through $\lambda_s(\xi)$ and $\eta_+(\xi)$ by the transformations $M(z;x,t) \mapsto M^{(1)}(z;x,t) \mapsto M^{(2)}(z;x,t)$, which yields the following RHP.

\begin{rhp}
  Find a $2 \times 2$ matrix-valued function $M^{(2)}(z;x,t)$ with the following properties:
\begin{enumerate} [label=(\roman*)]
  \item  $M^{(2)}(z;x,t)$ is analytic in $z \in \mathbb{C} \backslash \Sigma^{(2)}$ where $\Sigma^{(2)}=(-\infty, \lambda_r^+] \cup L_1 \cup L_2 \cup  L_3 \cup L_4.$
  \item  $M^{(2)}(z;x,t)=I+\mathcal{O} (z^{-1}) $ as $z\to \infty$.
  \item  $M^{(2)}(z;x,t)$ achieves the CBVs  $M^{(2)}_{+}(z;x,t)$ and $M^{(2)}_{-}(z;x,t)$  on $\Sigma^{(2)}$ away from self-intersection points and
  branch points that satisfy the jump condition $M^{(2)}_{+}(z;x,t)=M^{(2)}_{-}(z;x,t)V^{(2)}(z;x,t)$, where
  \begin{equation}\label{C2V2}
    \begin{aligned}
    V^{(2)}(z;x,t)&=\left\{  \begin{array}{ll}
     (1-r(z)r^{*}(z))^{\sigma_{3}},  ~& z\in (-\infty, \lambda_s(\xi))\backslash (\overline{\mathcal{I}_l \cup \mathcal{I}_r} ), \\ \\
     (a_+(z) a_-^*(z))^{-\sigma_{3}},  ~&z \in (-\infty, \lambda_s(\xi)) \cap (\mathcal{I}_r \backslash \overline{\mathcal{I}_l} ), \\ \\
     \begin{pmatrix} 0 & -r_{-}^{*}(z) e^{-\mathrm{i} \gamma}  \\ r_{+}(z) e^{\mathrm{i} \gamma}  & 0 \end{pmatrix}, ~&z \in \mathcal{I}_l \backslash \overline{\mathcal{I}_r}, \\ \\
     \begin{pmatrix} \frac{e^{ 2 \mathrm{i} t g_+(z)}}{a_+(z) a_-^*(z)} & -1 \\1 & 0 \end{pmatrix}, ~&z \in (\lambda_s(\xi), \eta_+(\xi)), \\ \\
     \begin{pmatrix} 0 & -1 \\1 & 0 \end{pmatrix}, ~&z \in (\eta_+(\xi), \lambda_r^+) \cup (\mathcal{I}_l \cap \mathcal{I}_r) , \\ \\
     \begin{pmatrix} 1 & 0 \\ r(z)e^{2\mathrm{i} tg(z)} & 1 \end{pmatrix}, ~ &z \in L_1, \\ \\
     \begin{pmatrix} 1 & \frac{-r^*(z)}{1-r(z)r^*(z)}e^{-2\mathrm{i} tg(z)} \\0 & 1 \end{pmatrix},  ~&z \in L_2,\\ \\
     \begin{pmatrix} 1 & 0 \\ \frac{r(z)}{1-r(z)r^*(z)}e^{2\mathrm{i} tg(z)} & 1 \end{pmatrix},  ~&z \in L_3,\\ \\
     \begin{pmatrix} 1 & -r^*(z)e^{-2\mathrm{i} tg(z)} \\ 0  & 1 \end{pmatrix}, ~ &z \in L_4.\end{array} \right.\\
    \end{aligned}
  \end{equation}
\end{enumerate}
\end{rhp}
Hence, the general $\delta(z)$ is now given by
\begin{equation} \label{C2delta}
  \begin{aligned}
     \delta(z)= \exp & \left\{  \frac{\mathcal{R} (z;\boldsymbol{\lambda})}{2 \pi \mathrm{i}}  \left( \int_{(-\infty, \lambda_s(\xi))\backslash (\overline{\mathcal{I}_l \cup \mathcal{I}_r} )}  \frac{\ln (1-|r(\zeta)|^2)}{\mathcal{R} (\zeta;\boldsymbol{\lambda}) } \,\frac{\dif\zeta}{\zeta - z} \right. \right. \\
      & \left. \left.+ \int_{(-\infty, \lambda_s(\xi)) \cap (\mathcal{I}_r \backslash \overline{\mathcal{I}_l} )} \frac{-\ln(a_+(\zeta)a_-^*(\zeta))}{\mathcal{R} (\zeta;\boldsymbol{\lambda})} \,\frac{\dif\zeta}{\zeta - z} + \int_{\mathcal{I}_l \backslash \overline{\mathcal{I}_r}} \frac{\ln (r_+(\zeta))}{\mathcal{R}_+ (\zeta;\boldsymbol{\lambda})} \, \frac{\dif\zeta}{\zeta - z} \right)  \right\},
    \end{aligned}
\end{equation}
where $\boldsymbol{\lambda}=(\lambda_r^+, \lambda_s(\xi), \lambda_l^+, \lambda_l^-)$. Then the corresponding coefficients $\varphi^{(1)}$ and $\varphi^{(0)}$ are as follows:
\begin{equation}\label{C2phi1}
  \begin{aligned}
   \varphi^{(1)}=\frac{1}{2 \pi} & \left( \int_{(-\infty, \lambda_s(\xi))\backslash (\overline{\mathcal{I}_l \cup \mathcal{I}_r} )} \frac{\ln (1-|r(\zeta)|^2) }{\mathcal{R} (\zeta;\boldsymbol{\lambda}) } \dif \zeta \right. \\
    &\left. +\int_{(-\infty, \lambda_s(\xi)) \cap (\mathcal{I}_r \backslash \overline{\mathcal{I}_l} )} \frac{-\ln(a_+(\zeta)a_-^*(\zeta))}{\mathcal{R} (\zeta;\boldsymbol{\lambda})} \dif \zeta + \int_{\mathcal{I}_l \backslash \overline{\mathcal{I}_r}} \frac{\ln (r_+(\zeta))}{\mathcal{R}_+ (\zeta;\boldsymbol{\lambda})} \,{\dif\zeta} \right),\\
  \end{aligned}
\end{equation}
\begin{equation}\label{C2phi0}
  \begin{aligned}
    \varphi^{(0)}=\frac{1}{2 \pi} & \left(\int_{(-\infty, \lambda_s(\xi))\backslash (\overline{\mathcal{I}_l \cup \mathcal{I}_r} )} \frac{\ln (1-|r(\zeta)|^2)  (\zeta+ V) }{\mathcal{R} (\zeta;\boldsymbol{\lambda}) } \dif \zeta \right. \\
    &+ \left. \int_{(-\infty, \lambda_s(\xi)) \cap (\mathcal{I}_r \backslash \overline{\mathcal{I}_l} )} \frac{-\ln(a_+(\zeta)a_-^*(\zeta))(\zeta+ V )}{\mathcal{R} (\zeta;\boldsymbol{\lambda})} \dif \zeta
    + \int_{\mathcal{I}_l \backslash \overline{\mathcal{I}_r}} \frac{\ln (r_+(\zeta)) (\zeta+ V)}{\mathcal{R}_+ (\zeta;\boldsymbol{\lambda})} \,{\dif\zeta} \right),
    \end{aligned}
\end{equation}
where $V=-\frac{1}{2}\sigma_1(\boldsymbol{\lambda})=-\frac{1}{2}(\lambda_r^+ + \lambda_s(\xi)+ \lambda_l^++\lambda_l^-)$.

Using the modified function $\hat{\delta}(z)$ defined by (\ref{AB2mdelta}) (with ${\delta}(z)$ replaced by (\ref{C2delta})) and the modified transformation
\begin{equation}
	M^{(3)}(z;x,t)=\hat{\delta} _{\infty}^{\sigma_{3}} M^{(2)}(z;x,t) \hat{\delta}(z) ^{-\sigma_{3}},
\end{equation}
 the RHP \ref{lsrhp} will be obtained again and the following steps are the same as Section \ref{ALSH}. Finally, it is concluded that the long-time asymptotics of the defocusing NLS equation (\ref{(NLS)})
in the left dispersive shock wave region is given by:
\begin{equation} \label{lswq}
  q(x,t)=\frac{\lambda_r^+ - \lambda_s(\xi) + \lambda_l^+ -\lambda_l^-}{2}\frac{\Theta (0)~\Theta (2\mathcal{A}(\infty)+(\gamma + \hat{\gamma})/{2 \pi})}{\Theta ((\gamma + \hat{\gamma})/{2 \pi})~\Theta (2\mathcal{A}(\infty))} e^{2\mathrm{i} (t g_\infty+\hat{g}_\infty-\varphi^{(0)})}+\mathcal{O}(t^{-1}),
\end{equation}
where $g_\infty$, $\hat{g}_\infty$ and $\varphi^{(0)}$ are given by (\ref{A2GINF}), (\ref{AB2ghatINF}) (with $\varphi^{(1)}$ replaced by (\ref{C2phi1})) and (\ref{C2phi0}), respectively.

\subsubsection{The middle plane wave region: $ -\frac{\lambda_r^++2\lambda_l^++\lambda_l^-}{2} < \xi< -\frac{\lambda_r^++2\lambda_r^- +\lambda_l^-}{2}$}\label{Cmpw}
\
\newline
\indent
\begin{figure}[htbp]
  \centering
  \subfigure[]{\includegraphics[width=5.5cm]{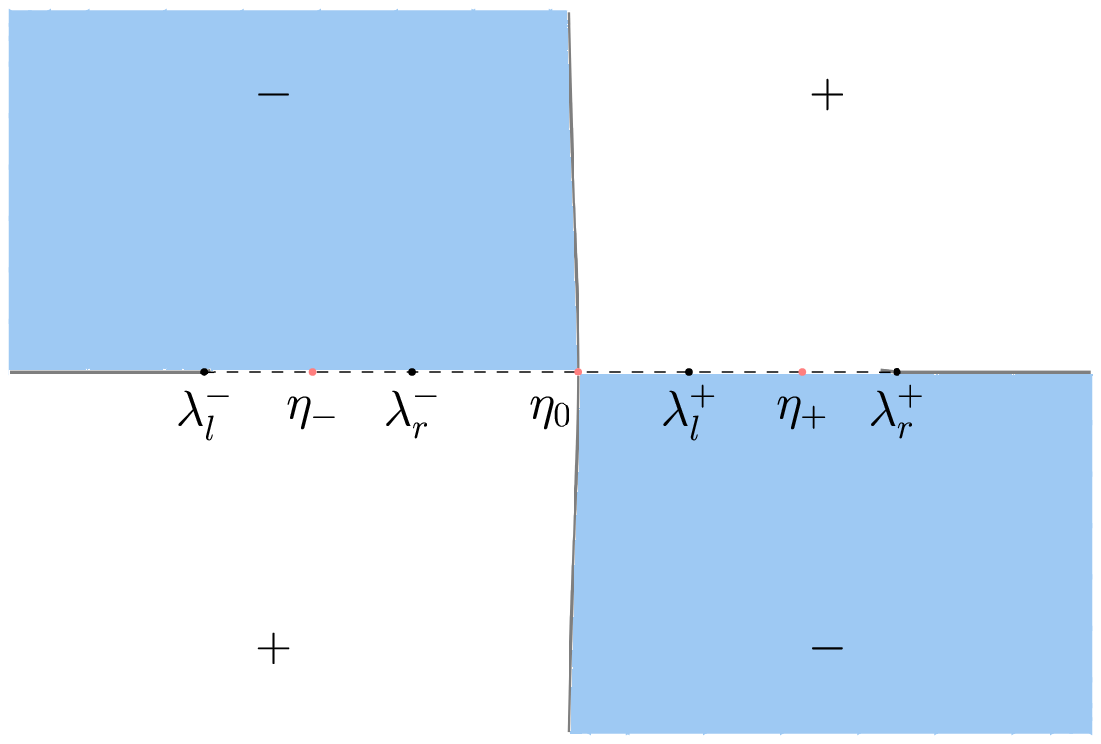}}
  \subfigure[]{\includegraphics[width=5cm]{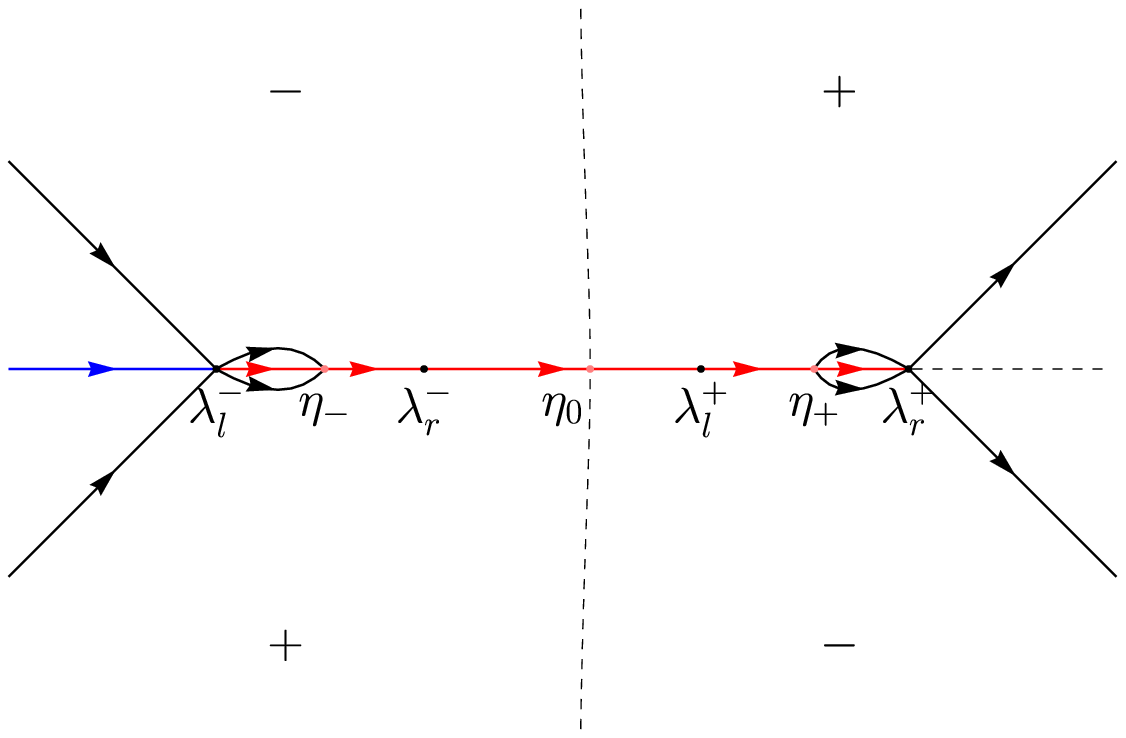}}
  \subfigure[]{\includegraphics[width=5.5cm]{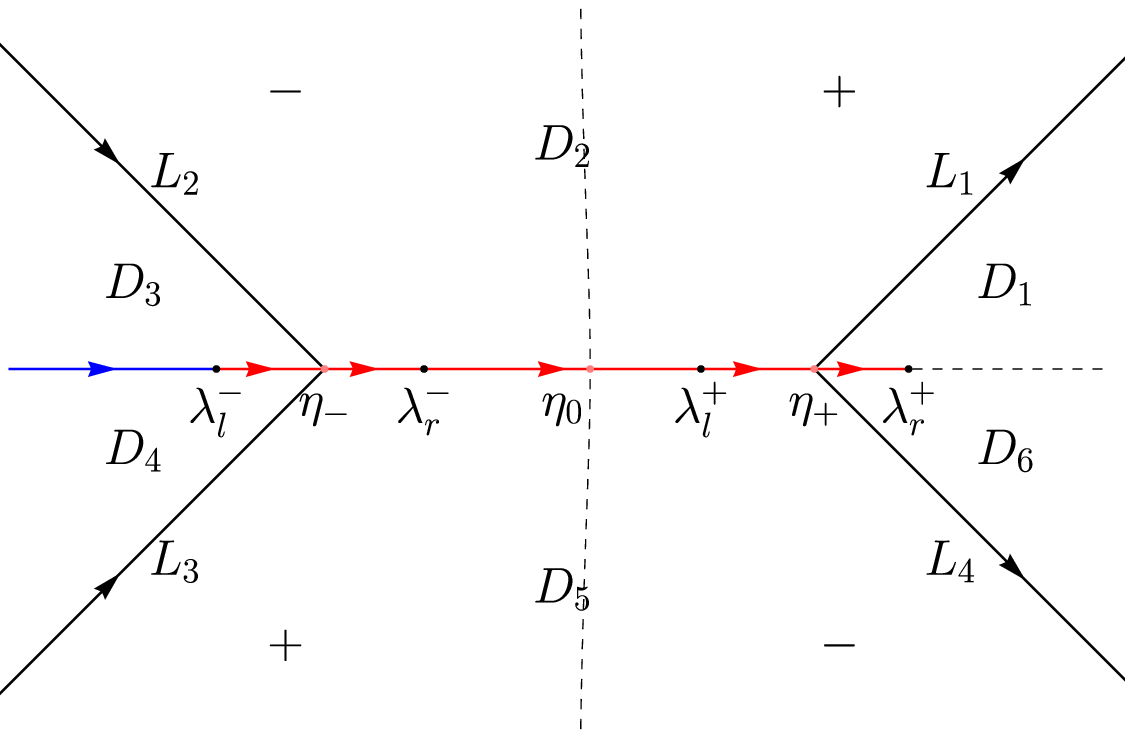}}
\caption{{\protect\small (a) The sign structure of Im$(g(z))$; (b) Opening lenses; (c) The jump contours of $M^{(2)}(z;x,t)$.}}
\label{figC3}
\end{figure}

When $\xi > -\frac{\lambda_r^++2\lambda_l^++\lambda_l^-}{2}$, the stationary phase point $\lambda_s(\xi)$ of the shock $g$-function is less than $\lambda_l^+$, which leads to exponentially large  diagonal entries $e^{\pm 2 \mathrm{i}t g_+(z)}$ of the jump matrix $V^{(2)}(z;x,t)$ on $(\lambda_s(\xi), \lambda_l^+)$ in (\ref{C2V2}).
To avoid this situation,  another new one-band $g$-function should be introduced, which can be regarded as the degeneration of the shock (two-band) $g$-function:
\begin{equation} \label{C3g}
  g(z)= \mathcal{R} (z;\lambda_r^+,\lambda_l^-)(z+\frac{1}{2}(\lambda_r^++\lambda_l^-)+\xi),
\end{equation}
which is consistent with the genus-zero case of the two hard edges ($\lambda_r^+$ and $\lambda_l^-$) in the Whitham modulation theory.
The function $g(z)$ is analytic in $z \in \mathbb{C} \backslash [\lambda_l^-, \lambda_r^+]$ with
\begin{equation}\label{C3ginf}
   g_{\infty}=- \xi(\lambda_r^+ + \lambda_l^- )/2 - (\lambda_r^+ + \lambda_l^- )^2/2- (\lambda_r^+ - \lambda_l^- )^2/8,
\end{equation}
and the differential of $g(z)$ is
\begin{equation} \label{C3dg}
  \dif g(z)= 2 \frac{(z-\eta_+(\xi))(z-\eta_-(\xi))}{\mathcal{R} (z;\lambda_l^+,\lambda_l^-)} \,\dif z,
\end{equation}
where $\eta_{\pm}(\xi)$ are the stationary phase points given by
\begin{equation}\label{cpweta}
  \eta_{\pm}(\xi)=\frac{\lambda_r^++\lambda_l^--\xi \pm \sqrt{(\lambda_r^++\lambda_l^-+\xi)^2+2(\lambda_r^+-\lambda_l^-)^2}}{4}.
\end{equation}
Denote $\eta_0(\xi)=-\frac{1}{2}(\lambda_r^++\lambda_l^- +2\xi)$, and the boundary of the middle plane wave region is characterized by:
\begin{equation}\label{C3b}
  \lambda_r^- < \eta_0(\xi) < \lambda_l^+    \qquad \mathrm{iff}  \qquad -\frac{\lambda_r^++2\lambda_l^++\lambda_l^-}{2} = \xi_{C2} < \xi< \xi_{C3}= -\frac{\lambda_r^++2\lambda_r^- +\lambda_l^-}{2} .
\end{equation}

Then open lenses from the real axis to the steepest descent contours through $\eta_+(\xi)$ and $\eta_-(\xi)$ (as shown in Figure \ref{figC3}) by the transformations (\ref{g-fun}) and (\ref{open}), i.e.,
$M(z;x,t) \mapsto M^{(1)}(z;x,t) \mapsto M^{(2)}(z;x,t)$, which yields the following new RHP.

\begin{rhp}
  Find a $2 \times 2$ matrix-valued function $M^{(2)}(z;x,t)$ with the following properties:
\begin{enumerate} [label=(\roman*)]
  \item  $M^{(2)}(z;x,t)$ is analytic in $z \in \mathbb{C} \backslash \Sigma^{(2)}$, where $\Sigma^{(2)}=(-\infty, \lambda_r^+] \cup L_1 \cup L_2 \cup  L_3 \cup L_4.$
  \item  $M^{(2)}(z;x,t)=I+\mathcal{O} (z^{-1}) $ as $z\to \infty$.
  \item  $M^{(2)}(z;x,t)$ achieves the CBVs  $M^{(2)}_{+}(z;x,t)$ and $M^{(2)}_{-}(z;x,t)$  on $\Sigma^{(2)}$ away from self-intersection points and
  branch points that satisfy the jump condition $M^{(2)}_{+}(z;x,t)=M^{(2)}_{-}(z;x,t)V^{(2)}(z;x,t)$, where
  \begin{equation}\label{C3V2}
    \begin{aligned}
    V^{(2)}(z;x,t)&=\left\{  \begin{array}{ll}
     (1-r(z)r^{*}(z))^{\sigma_{3}},  ~& z\in (-\infty, \lambda_l^- ), \\ \\
     \begin{pmatrix} 0 & -r_{-}^{*}(z) \\ r_{+}(z) & e^{-2\mathrm{i}tg_+(z)} \chi_{(\eta_-(\xi), \lambda_r^-)}  \end{pmatrix}, ~&z \in (\lambda_l^-, \lambda_r^-), \\ \\
     \begin{pmatrix} 0 & -1 \\ 1 & 0 \end{pmatrix}, ~&z \in (\lambda_r^-, \lambda_l^+), \\ \\
     \begin{pmatrix} \frac{e^{ 2 \mathrm{i} t g_+(z)}}{a_+(z) a_-^*(z)} \chi_{(\lambda_l^+, \eta_+(\xi))} & -1 \\ 1 & 0 \end{pmatrix}, ~&z \in (\lambda_l^+, \lambda_r^+), \\ \\
     \begin{pmatrix} 1 & 0 \\ r(z)e^{2\mathrm{i} tg(z)} & 1 \end{pmatrix}, ~ &z \in L_1, \\ \\
     \begin{pmatrix} 1 & \frac{-r^*(z)}{1-r(z)r^*(z)}e^{-2\mathrm{i} tg(z)} \\0 & 1 \end{pmatrix},  ~&z \in L_2,\\ \\
     \begin{pmatrix} 1 & 0 \\ \frac{r(z)}{1-r(z)r^*(z)}e^{2\mathrm{i} tg(z)} & 1 \end{pmatrix},  ~&z \in L_3,\\ \\
     \begin{pmatrix} 1 & -r^*(z)e^{-2\mathrm{i} tg(z)} \\ 0  & 1 \end{pmatrix}, ~ &z \in L_4.\end{array} \right.\\
    \end{aligned}
  \end{equation}
  Here, $\chi_{(a, b)}$ is the characteristic function of the set $(a, b)$ and  $\chi_{(a, b)} \equiv 0$ when $a>b$.
\end{enumerate}
\end{rhp}

Define
\begin{equation} \label{C3delta}
  \begin{aligned}
  \delta(z) = \exp & \left\{ \frac{\mathcal{R} (z;\lambda_r^+, \lambda_l^-)}{2 \pi \mathrm{i}}  \left(\int_{-\infty}^{\lambda_l^-}   \frac{\ln (1-|r(\zeta)|^2)}{\mathcal{R} (\zeta;\lambda_r^+, \lambda_l^-) } \,\frac{\dif\zeta}{\zeta - z}
    + \int_{\lambda_l^-}^{\lambda_r^-} \frac{\ln (r_+(\zeta))}{\mathcal{R}_+ (\zeta;\lambda_r^+, \lambda_l^-)} \, \frac{\dif\zeta}{\zeta - z} \right)  \right\}
  \end{aligned}
\end{equation}
whose large $z$ asymptotic behavior is $\delta(z)=\delta_{\infty}+\mathcal{O} (z^{-1}) $ as $z \to \infty$, where
\begin{equation}\label{C3deltaINF}
  \begin{aligned}
  \delta_{\infty }= \exp & \left\{ \frac{ \mathrm{i}}{2\pi} \left(\int_{-\infty}^{\lambda_l^-}  \frac{\ln (1-|r(\zeta)|^2)}{\mathcal{R} (\zeta;\lambda_r^+, \lambda_l^-) } \dif \zeta
  + \int_{\lambda_l^-}^{\lambda_r^-} \frac{\ln (r_+(\zeta))}{\mathcal{R}_+ (\zeta;\lambda_r^+, \lambda_l^-)} \dif \zeta \right)  \right\}.
  \end{aligned}
\end{equation}
Using the function ${\delta}(z)$  and the transformation
\begin{equation}
	M^{(3)}(z;x,t)=\delta_{\infty}^{\sigma_{3}} M^{(2)}(z;x,t) \delta(z) ^{-\sigma_{3}},
\end{equation}
 the following RHP is obtained immediately.

\begin{rhp} \label{rhpmid} Find a $2 \times 2$ matrix-valued function $M^{(3)}(z;x,t)$ with the following properties:
  \begin{enumerate} [label=(\roman*)]
    \item $M^{(3)}(z;x,t)$ is analytic in $z \in \mathbb{C} \backslash \Sigma^{(3)}$, where $\Sigma^{(3)}=  [\lambda_l^-,\lambda_r^+]\cup L_1 \cup L_2 \cup  L_3 \cup L_4.$
    \item $M^{(3)}(z;x,t)=I+\mathcal{O} (z^{-1}) $ as $z\to \infty$.
    \item $M^{(3)}(z;x,t)$ achieves the CBVs $M^{(3)}_{+}(z;x,t)$ and $M^{(3)}_{-}(z;x,t)$ on $\Sigma^{(3)}$ away from the self-intersection point and
  branch points that satisfy the jump condition $M^{(3)}_{+}(z;x,t)=M^{(3)}_{-}(z;x,t)V^{(3)}(z;x,t)$, where
  \begin{equation}\label{C3V3}
    \begin{aligned}
  V^{(3)}(z;x,t)&=\left\{  \begin{array}{ll}
     \begin{pmatrix} \frac{\delta ^{-1}_{+}(z)\delta _{-}(z)}{a_+(z) a_-^*(z)}{e^{ 2 \mathrm{i} t g_+(z)}} \chi_{(\lambda_l^+, \eta_+(\xi))} & -1 \\1 & {\delta _{+}(z)\delta^{-1} _{-}(z)}e^{-2\mathrm{i}tg_+(z)} \chi_{(\eta_-(\xi), \lambda_r^-)} \end{pmatrix}, ~&z \in (\lambda_l^-, \lambda_r^+), \\ \\
     \begin{pmatrix} 1 & 0 \\ r(z) {\delta} ^{-2}(z) e^{2\mathrm{i} tg(z)} & 1 \end{pmatrix}, ~ &z \in L_1, \\ \\
     \begin{pmatrix} 1 & \frac{-r^*(z)}{1-r(z)r^*(z)} {\delta} ^2(z) e^{-2\mathrm{i}t g(z)} \\0 & 1 \end{pmatrix},  ~&z \in L_2,\\ \\
     \begin{pmatrix} 1 & 0 \\ \frac{r(z)}{1-r(z)r^*(z)} {\delta} ^{-2}(z) e^{2\mathrm{i} tg(z)} & 1 \end{pmatrix},  ~&z \in L_3,\\ \\
     \begin{pmatrix} 1 & -r^*(z) {\delta} ^2(z) e^{-2\mathrm{i} tg(z)} \\ 0  & 1 \end{pmatrix}, ~ &z \in L_4. \end{array} \right.\\
    \end{aligned}
  \end{equation}
  \end{enumerate}
\end{rhp}

The above jump matrices subsequently uniformly converge to constant matrices, and stationary phase points both lie on the cut $[\lambda_l^-,\lambda_r^+]$ of the phase function $g(z)$ where $\mathrm{Im}~g_{\pm}(z) \neq 0$.
Therefore, there exists a constant $c>0$ such that the diagonal entries of the above jump matrix on $(\lambda_l^-, \lambda_r^+)$ and the off-diagonal entries of the above jump matrices on $L_i$ $(i=1, \cdots, 4),$ are
$\mathcal{O}(e^{-ct})$, and no local parametrix is needed. Then a global parametrix satisfying the uniformly limiting problem is found as follows
\begin{equation}
  M^{(\mathrm{par})}_+(z;x,t)=M^{(\mathrm{par})}_-(z;x,t) V^{(\mathrm{par})} \qquad  \mathrm{for}  \quad z \in (\lambda_l^-,\lambda_r^+),
\end{equation}
where
\begin{equation}\label{C3vpar}
  V^{(\mathrm{par})}=\begin{pmatrix} 0 & -1 \\ 1 & 0 \end{pmatrix}.
\end{equation}
The solution is exactly given by
\begin{equation} \label{C3par}
  M^{(\mathrm{par})}(z;x,t) = \mathcal{E} (z;\lambda_r^+, \lambda_l^-).
\end{equation}
Define the error matrix as
$M^{(\mathrm{err})}(z;x,t)= M^{(3)}(z;x,t) {M^{(\mathrm{par})}(z;x,t)}^{-1}$,
which results in the following error RHP.

\begin{rhp} Find a $2 \times 2$ matrix-valued function $M^{(\mathrm{err})}(z;x,t)$ with the following properties:
  \begin{enumerate} [label=(\roman*)]
    \item $M^{(\mathrm{err})}(z;x,t)$ is analytic in $z \in \mathbb{C} \backslash \Sigma^{(\mathrm{err})}$, where $\Sigma^{(\mathrm{err})}=L_1 \cup L_2 \cup  L_3 \cup L_4 $.
    \item $M^{(\mathrm{err})}(z;x,t)=I+\mathcal{O} (z^{-1}) $ as $z\to \infty$.
    \item $M^{(\mathrm{err})}(z;x,t)$ achieves the CBVs $M^{(\mathrm{err})}_{+}(z;x,t)$ and $M^{(\mathrm{err})}_{-}(z;x,t)$ on $\Sigma^{(\mathrm{err})}$, which satisfy the jump condition $M^{(\mathrm{err})}_{+}(z;x,t)=M^{(\mathrm{err})}_{-}(z;x,t)V^{(\mathrm{err})}(z;x,t)$, where
  \begin{equation}\label{C3Verr}
    \begin{aligned}
  V^{(\mathrm{err})}(z;x,t)&=\left\{  \begin{array}{ll}
    M^{(\mathrm{par})}(z;x,t) \begin{pmatrix} 1 & 0 \\ r(z) \delta ^{-2}(z) e^{2\mathrm{i} tg(z)} & 1 \end{pmatrix} {M^{(\mathrm{par})}(z;x,t)}^{-1}, ~ &z \in L_1, \\ \\
    M^{(\mathrm{par})}(z;x,t) \begin{pmatrix} 1 & \frac{-r^*(z)}{1-r(z)r^*(z)} \delta ^2(z) e^{-2\mathrm{i}t g(z)} \\0 & 1 \end{pmatrix} {M^{(\mathrm{par})}(z;x,t)}^{-1},  ~&z \in L_2,\\ \\
    M^{(\mathrm{par})}(z;x,t) \begin{pmatrix} 1 & 0 \\ \frac{r(z)}{1-r(z)r^*(z)} \delta ^{-2}(z) e^{2\mathrm{i} tg(z)} & 1 \end{pmatrix} {M^{(\mathrm{par})}(z;x,t)}^{-1},  ~&z \in L_3,\\ \\
    M^{(\mathrm{par})}(z;x,t) \begin{pmatrix} 1 & -r^*(z) \delta ^2(z) e^{-2\mathrm{i} tg(z)} \\ 0  & 1 \end{pmatrix} {M^{(\mathrm{par})}(z;x,t)}^{-1}, ~ &z \in L_4. \end{array} \right.\\
    \end{aligned}
  \end{equation}
  \end{enumerate}
\end{rhp}
Note that no local parametrix is needed, so the error estimate is exactly $M^{(\mathrm{err})}_1(x,t)= \mathcal{O}(e^{-ct})$ for some constant $c>0$.
Using equations (\ref{rec}), (\ref{C3ginf}), (\ref{C3deltaINF}) and (\ref{C3par}),
it is concluded that the long-time asymptotics of the defocusing NLS equation (\ref{(NLS)}) in the middle plane wave region is a plane wave of the form
\begin{equation}
    q(x,t)=\frac{\lambda_r^+-\lambda_l^-}{2} e^{\mathrm{i}(kx-\omega t)}e^{-\mathrm{i}\phi_{mp} }+\mathcal{O}(e^{-ct}),
\end{equation}
where
\begin{equation}
  k=-(\lambda_r^++\lambda_l^-),~~~~ \omega=\frac{(\lambda_r^++\lambda_l^-)^2}{2}+\frac{(\lambda_r^+ - \lambda_l^-)^2}{4},
\end{equation}
and
\begin{equation} \label{phimp}
  \phi_{mp}=\frac{1}{\pi} \left( \int_{-\infty}^{\lambda_l^-} \frac{\ln (1-|r(\zeta)|^2)}{\sqrt{(\zeta-\lambda_r^+)(\zeta-\lambda_l^-)} } \dif \zeta
  +\int_{\lambda_l^-}^{\lambda_r^-} \frac{\arg (r_+(\zeta))}{\sqrt{(\lambda_r^+-\zeta)(\zeta-\lambda_l^-)}} \dif \zeta \right). \\
\end{equation}

\subsubsection{The right dispersive shock wave region: $ -\frac{\lambda_r^++2\lambda_r^- +\lambda_l^-}{2} < \xi < -\frac{\lambda_r^+ + \lambda_r^- + 2\lambda_l^-}{2}+\frac{2(\lambda_r^+ - \lambda_l^-)
(\lambda_r^- - \lambda_l^-)}{\lambda_r^+ + \lambda_r^- -2\lambda_l^-} $} \label{Crsh}
\
\newline
\indent
If $\xi > -\frac{\lambda_r^++2\lambda_r^- +\lambda_l^-}{2}$, the point $\eta_0(\xi)$ of the one-band $g$-function (\ref{C3g}) is less than $\lambda_r^-$, which leads to the emergence of the interval $(\eta_0(\xi), \lambda_r^-)$, on which the jump matrix (\ref{C3V2}) grows exponentially with respect to a sufficiently large $t$.
To remove the exponentially large entries, it is necessary to modify the $g$-function to include a gap $(\lambda_s(\xi), \lambda_r^-)$, where $\lambda_s(\xi)$ is the soft edge determined by (\ref{vv3}).
The $g$-function is exactly given by (\ref{AB4g}).
The boundaries of this region are characterized by degeneration below
\begin{equation}\label{C4b}
  \begin{aligned}
  &\xi \to \xi_{C3}= v_3(\lambda_r^+, \lambda_r^-, \lambda_r^-, \lambda_l^-)= -\frac{\lambda_r^++2\lambda_r^- +\lambda_l^-}{2}, \quad &\mathrm{as} \quad \lambda_s(\xi) \to \lambda_r^-,\\
  &\xi \to \xi_{C4}=  v_3(\lambda_r^+, \lambda_r^-, \lambda_l^-, \lambda_l^-)= -\frac{\lambda_r^+ + \lambda_r^- + 2\lambda_l^-}{2}+\frac{2(\lambda_r^+ - \lambda_l^-)
  (\lambda_r^- - \lambda_l^-)}{\lambda_r^+ + \lambda_r^- -2\lambda_l^-}, \quad &\mathrm{as} \quad \lambda_s(\xi) \to \lambda_l^-.
  \end{aligned}
\end{equation}
where $v_3$ is the Whitham velocity given by (\ref{vv3}).

Due to the same $g$-function, we will summarize the RHP analysis of the right dispersive shock wave region.  Similarly, open lenses from intervals to the steepest descent contours through $\lambda_s(\xi)$ and $\eta_-(\xi)$ by the transformations
$M(z;x,t) \mapsto M^{(1)}(z;x,t) \mapsto M^{(2)}(z;x,t)$, which exactly results in RHP \ref{rsrhp}, since the upper/lower matrix factorizations in (\ref{f1}) and (\ref{f2}) ((\ref{f3}) and (\ref{f4})) are the same, respectively.
The following steps are exactly the same as in Section {\ref{ABrsh}} and  the asymptotic behaviors with error estimates in the right dispersive shock wave region are given by (\ref{rsab}).

\subsubsection{The right plane wave region: $\xi > -\frac{\lambda_r^+ + \lambda_r^- + 2\lambda_l^-}{2}+\frac{2(\lambda_r^+ - \lambda_l^-)
(\lambda_r^- - \lambda_l^-)}{\lambda_r^+ + \lambda_r^- -2\lambda_l^-}$} \label{Crpw}
\
\newline
\indent
When the stationary phase $\lambda_s(\xi)$  coincides with $\lambda_l^-$, the two-band $g$-function degenerates to the one-band one. This is the same as that
in Section \ref{ABrpw}, where we have given the uniform form of the right plane wave region. So we only present the boundary of the right plane wave region
\begin{equation} \label{C5b}
  \eta_- (\xi) < \lambda_l^- \qquad \mathrm{iff}  \qquad  \xi > \xi_{C4} = -\frac{\lambda_r^+ + \lambda_r^- + 2\lambda_l^-}{2}+\frac{2(\lambda_r^+ - \lambda_l^-)
  (\lambda_r^- - \lambda_l^-)}{\lambda_r^+ + \lambda_r^- -2\lambda_l^-},
\end{equation}
where $\eta_- (\xi)$ is given by (\ref{rpweta}). The boundary is consistent with (\ref{C4b}).

\subsection{\rm Case D:  $\lambda_l^+>\lambda_r^+>\lambda_l^->\lambda_r^-$}
\
\newline
\indent
In this case, the two intervals $\mathcal{I}_l$ and $\mathcal{I}_r$ intersect but are not completely contained within each other, thus all the factorizations in
(\ref{f1})-(\ref{f4}) should be considered. In comparison with Case B, only the positions of points $ \lambda_l^+$ and $\lambda_r^-$ are swapped, which means that this case is very similar to Case B, except that the center region becomes a plane wave region, as we will see later.
As the self-similar variable $\xi=x/t$ increases, the stationary phase points of the corresponding
$g$-functions change continuously at different intervals on $\mathbb{R}$, which implies that there are five different regions: the left plane wave region,
the left rarefaction wave region, the middle plane wave region, the right rarefaction wave region and the right plane wave region.

\subsubsection{The left plane wave region: $ \xi<-\frac{3\lambda_l^++\lambda_l^-}{2}  $}\label{Dlpw}
\
\newline
\indent
The uniform form of the left plane wave region is similar to the case in Section \ref{Blpw}, and only the boundary of the region needs to be determined, which is easy to formulate as
\begin{equation} \label{D1b}
  \eta_+(\xi) > \lambda_l^+    \qquad \mathrm{iff}  \qquad \xi< \xi_{D1}=-\frac{3\lambda_l^++\lambda_l^-}{2},
\end{equation}
where $\eta_+ (\xi) $ is given by (\ref{lpweta}).

\subsubsection{The left rarefaction wave region: $-\frac{3\lambda_l^++\lambda_l^-}{2} <\xi < -\frac{3\lambda_r^++\lambda_l^-}{2} $} \label{Dlrare}
\
\newline
\indent
When $\xi > -\frac{3\lambda_l^++\lambda_l^-}{2} $, the stationary phase point $\eta_+(\xi)$ of the $g$-function (\ref{A1g}) moves inside $\mathcal{I}_l$, which leads to exponentially large diagonal entries $e^{\pm 2 \mathrm{i} t g_+(z)}$ of the jump matrix $V^{(2)}(z;x,t)$ on the interval $(\eta_+(\xi), \lambda_l^+) \subset
\mathcal{I}_l$. So it is necessary to introduce the rarefaction $g$-function defined by (\ref{A2g}). The boundaries of this rarefaction wave region are characterized by:
\begin{equation} \label{D2b}
  \lambda_r^+ < \lambda_s(\xi) < \lambda_l^+ \qquad \mathrm{iff}  \qquad -\frac{3\lambda_l^++\lambda_l^-}{2}  = \xi_{D1} <\xi < \xi_{D2}= -\frac{3\lambda_r^++\lambda_l^-}{2} ,
\end{equation}
where $\lambda_s(\xi)=-(\lambda_l^- + 2 \xi)/3$  is the soft edge defined in Whitham modulation theory.
\par
Due to the same rarefaction $g$-function and the same  upper/lower matrix factorizations in (\ref{f3}) and (\ref{f4}), we will summarize the RHP analysis of the left rarefaction wave region. In the similar way, open lenses from intervals to the steepest descent contours through $\lambda_s(\xi)$ and $\eta(\xi)$ by the transformations
$M(z;x,t) \mapsto M^{(1)}(z;x,t) \mapsto M^{(2)}(z;x,t)$, which results in the following RHP.

\begin{rhp}
  Find a $2 \times 2$ matrix-valued function $M^{(2)}(z;x,t)$ with the following properties:
\begin{enumerate} [label=(\roman*)]
  \item  $M^{(2)}(z;x,t)$ is analytic in $z \in \mathbb{C} \backslash \Sigma^{(2)}$, where $\Sigma^{(2)}=(-\infty, \lambda_{s}(\xi)] \cup L_1 \cup L_2 \cup  L_3 \cup L_4.$
  \item  $M^{(2)}(z;x,t)=I+\mathcal{O} (z^{-1}) $ as $z\to \infty$.
  \item  $M^{(2)}(z;x,t)$ achieves the CBVs  $M^{(2)}_{+}(z;x,t)$ and $M^{(2)}_{-}(z;x,t)$  on $\Sigma^{(2)}$ away from self-intersection points and
  branch points that satisfy the jump condition $M^{(2)}_{+}(z;x,t)=M^{(2)}_{-}(z;x,t)V^{(2)}(z;x,t)$, where
  \begin{equation}\label{D2V2}
    \begin{aligned}
    V^{(2)}(z;x,t)&=\left\{  \begin{array}{ll}
     (1-r(z)r^{*}(z))^{\sigma_{3}},  ~& z\in (-\infty, \lambda_{s}(\xi))\backslash (\overline{\mathcal{I}_l \cup \mathcal{I}_r} ), \\ \\
     (a_+(z) a_-^*(z))^{-\sigma_{3}},  ~&z \in  \mathcal{I}_r \backslash \overline{\mathcal{I}_l} , \\ \\
     \begin{pmatrix} 0 & -r_{-}^{*}(z) \\ r_{+}(z) & e^{- 2 \mathrm{i} t g_+(z)} \chi_{(\eta(\xi), \lambda_s(\xi))} \end{pmatrix}, ~&z \in (-\infty, \lambda_{s}(\xi))\cap (\mathcal{I}_l \backslash \overline{\mathcal{I}_r} ), \\ \\
     \begin{pmatrix} 0 & -1 \\ 1 & 0 \end{pmatrix}, ~&z \in \mathcal{I}_l \cap \mathcal{I}_r \\ \\
     \begin{pmatrix} 1 & 0 \\ r(z)e^{2\mathrm{i} tg(z)} & 1 \end{pmatrix}, ~ &z \in L_1, \\ \\
     \begin{pmatrix} 1 & \frac{-r^*(z)}{1-r(z)r^*(z)}e^{-2\mathrm{i} tg(z)} \\0 & 1 \end{pmatrix},  ~&z \in L_2,\\ \\
     \begin{pmatrix} 1 & 0 \\ \frac{r(z)}{1-r(z)r^*(z)}e^{2\mathrm{i} tg(z)} & 1 \end{pmatrix},  ~&z \in L_3,\\ \\
     \begin{pmatrix} 1 & -r^*(z)e^{-2\mathrm{i} tg(z)} \\ 0  & 1 \end{pmatrix}, ~ &z \in L_4.\end{array} \right.\\
    \end{aligned}
  \end{equation}
  Here, $\chi_{(a, b)}$ is the characteristic function of the set $(a, b)$.
\end{enumerate}
\end{rhp}

Now, the general function $\delta(z)$ is given by
\begin{equation} \label{D2delta}
  \begin{aligned}
  \delta(z) = \exp & \left\{ \frac{\mathcal{R} (z;\lambda_s(\xi), \lambda_l^-)}{2 \pi \mathrm{i}}  \left(\int_{(-\infty, \lambda_{s}(\xi))\backslash (\overline{\mathcal{I}_l \cup \mathcal{I}_r} )}  \frac{\ln (1-|r(\zeta)|^2)}{\mathcal{R} (\zeta;\lambda_s(\xi), \lambda_l^-) } \,\frac{\dif\zeta}{\zeta - z} \right. \right.  \\
  & \left. \left.  + \int_{\mathcal{I}_r \backslash \overline{\mathcal{I}_l}} \frac{-\ln(a_+(\zeta)a_-^*(\zeta))}{\mathcal{R} (\zeta;\lambda_s(\xi), \lambda_l^-)} \,\frac{\dif\zeta}{\zeta - z}
  + \int_{(-\infty, \lambda_{s}(\xi))\cap (\mathcal{I}_l \backslash \overline{\mathcal{I}_r} )} \frac{\ln (r_+(\zeta))}{\mathcal{R}_+ (\zeta;\lambda_s(\xi), \lambda_l^-)} \, \frac{\dif\zeta}{\zeta - z} \right)  \right\}
  \end{aligned}
\end{equation}
whose large $z$ asymptotic behavior is $\delta(z)=\delta_{\infty}+\mathcal{O} (z^{-1}) $ as $z \to \infty$, where
\begin{equation}\label{D2deltaINF}
  \begin{aligned}
  \delta_{\infty }= \exp & \left\{ \frac{ \mathrm{i}}{2\pi} \left(\int_{(-\infty, \lambda_{s}(\xi))\backslash (\overline{\mathcal{I}_l \cup \mathcal{I}_r} )}  \frac{\ln (1-|r(\zeta)|^2)}{\mathcal{R} (\zeta;\lambda_s(\xi), \lambda_l^-) } \dif \zeta \right. \right.  \\
  &\left. \left.  + \int_{\mathcal{I}_r \backslash \overline{\mathcal{I}_l}}  \frac{-\ln(a_+(\zeta)a_-^*(\zeta))}{\mathcal{R} (\zeta;\lambda_s(\xi), \lambda_l^-)} \dif \zeta
  + \int_{(-\infty, \lambda_{s}(\xi))\cap (\mathcal{I}_l \backslash \overline{\mathcal{I}_r} )} \frac{\ln (r_+(\zeta))}{\mathcal{R}_+ (\zeta;\lambda_s(\xi), \lambda_l^-)} \dif \zeta \right)  \right\}.
  \end{aligned}
\end{equation}
Using the function ${\delta}(z)$  and the transformation
\begin{equation}
	M^{(3)}(z;x,t)={\delta} _{\infty}^{\sigma_{3}} M^{(2)}(z;x,t) {\delta}(z) ^{-\sigma_{3}},
\end{equation}
the RHP \ref{lrrhp} is obtained again and the following steps are exactly the same as Section \ref{Alrare}.
Using equations (\ref{rec}), (\ref{A2ginf}),  (\ref{A2mod}) and (\ref{D2deltaINF}), the long-time asymptotics of the defocusing NLS equation (\ref{(NLS)}) with error estimates in the left rarefaction wave regions is derived as
\begin{equation}\label{lrwq}
  \begin{aligned}
   q(x,t)=-\frac{\xi+2\lambda_l^-}{3}e^{-\mathrm{i}t(2{\lambda_l^-}^2+2\lambda_l^- \xi-\xi^2)/3}e^{-\mathrm{i}\phi_{lr}(\xi) }+\mathcal{O}(t^{-1}),
  \end{aligned}
\end{equation}
where
\begin{equation}\label{philrlr}
  \begin{aligned}
    \phi_{lr}(\xi)=&\frac{1}{\pi} \left(\int_{(-\infty, \lambda_{s}(\xi))\backslash (\overline{\mathcal{I}_l \cup \mathcal{I}_r} )}  \frac{\ln (1-|r(\zeta)|^2)}{\sqrt{(\zeta-\lambda_s(\xi))(\zeta-\lambda_l^-)} } \dif \zeta \right.  \\
    & \left. + \int_{\mathcal{I}_r \backslash \overline{\mathcal{I}_l}}  \frac{-\ln(a_+(\zeta)a_-^*(\zeta))}{\sqrt{(\zeta-\lambda_s(\xi))(\zeta-\lambda_l^-)}} \dif \zeta
    + \int_{(-\infty, \lambda_{s}(\xi))\cap (\mathcal{I}_l \backslash \overline{\mathcal{I}_r} )} \frac{\arg(r_+(\zeta))}{\sqrt{(\lambda_s(\xi)-\zeta)(\zeta-\lambda_l^-)}} \dif \zeta \right)  .
  \end{aligned}
\end{equation}

\subsubsection{The middle plane wave region: $ -\frac{3\lambda_r^+ + \lambda_l^-}{2} < \xi< -\frac{\lambda_r^+ + 3\lambda_l^-}{2}$}\label{Dmpw}
\
\newline
\indent
When $\xi > -\frac{3\lambda_r^+ + \lambda_l^-}{2} $, the stationary phase point $\lambda_s(\xi)$ of the rarefaction $g$-function is less than $\lambda_r^+$, which leads to the exponentially large diagonal entries $e^{\pm 2 \mathrm{i}t g_+(z)}$ of the jump matrix $V^{(2)}(z;x,t)$ on $(\lambda_s(\xi), \lambda_r^+)$ in (\ref{D2V2}).
To avoid this case, the one-band $g$-function defined by (\ref{C3g}) should be introduced.
The boundary of the middle plane wave region is characterized by
\begin{equation} \label{D3b}
  \lambda_l^- < \eta_0(\xi) < \lambda_r^+    \qquad \mathrm{iff}  \qquad -\frac{3\lambda_r^+ + \lambda_l^-}{2} = \xi_{D2} < \xi_{D3}= \xi< -\frac{\lambda_r^+ + 3\lambda_l^-}{2},
\end{equation}
where $\eta_0(\xi)=-\frac{1}{2}(\lambda_r^++\lambda_l^- +2\xi)$.

Due to the same one-band $g$-function and the same  upper/lower matrix factorizations in (\ref{f3}) and (\ref{f4}), we will summarize the RHP analysis of the middle plane wave region. As before, open lenses from intervals to the steepest descent contours through  $\eta_+(\xi)$ and $\eta_-(\xi)$ by the transformations
$M(z;x,t) \mapsto M^{(1)}(z;x,t) \mapsto M^{(2)}(z;x,t)$, which results in the following RHP.

\begin{rhp}
  Find a $2 \times 2$ matrix-valued function $M^{(2)}(z;x,t)$ with the following properties:
\begin{enumerate} [label=(\roman*)]
  \item  $M^{(2)}(z;x,t)$ is analytic in $z \in \mathbb{C} \backslash \Sigma^{(2)}$, where $\Sigma^{(2)}=(-\infty, \lambda_r^+] \cup L_1 \cup L_2 \cup  L_3 \cup L_4.$
  \item  $M^{(2)}(z;x,t)=I+\mathcal{O} (z^{-1}) $ as $z\to \infty$.
  \item  $M^{(2)}(z;x,t)$ achieves the CBVs  $M^{(2)}_{+}(z;x,t)$ and $M^{(2)}_{-}(z;x,t)$  on $\Sigma^{(2)}$ away from self-intersection points and
  branch points that satisfy the jump condition $M^{(2)}_{+}(z;x,t)=M^{(2)}_{-}(z;x,t)V^{(2)}(z;x,t)$, where
  \begin{equation}\label{D3V2}
    \begin{aligned}
    V^{(2)}(z;x,t)&=\left\{  \begin{array}{ll}
     (1-r(z)r^{*}(z))^{\sigma_{3}},  ~& z\in (-\infty, \eta_0(\xi))\backslash (\overline{\mathcal{I}_l \cup \mathcal{I}_r} ), \\ \\
     (a_+(z) a_-^*(z))^{-\sigma_{3}},  ~&z \in  (-\infty, \eta_0(\xi) )\cap(\mathcal{I}_r \backslash \overline{\mathcal{I}_l}) , \\ \\
     \begin{pmatrix} 0 & -r_{-}^{*}(z) \\ r_{+}(z) & e^{-2\mathrm{i}tg_+} \chi_{(\eta_-(\xi), \lambda_r^-)}  \end{pmatrix}, ~&z \in (-\infty, \eta_0(\xi) )\cap(\mathcal{I}_l \backslash \overline{\mathcal{I}_r}), \\ \\
     \begin{pmatrix} 0 & -1 \\ 1 & 0 \end{pmatrix}, ~&z \in \mathcal{I}_l \cap \mathcal{I}_r, \\ \\
     \begin{pmatrix} \frac{e^{ 2 \mathrm{i} t g_+(z)}}{a_+(z) a_-^*(z)} \chi_{(\lambda_l^+, \eta_+(\xi))} & -1 \\ 1 & 0 \end{pmatrix}, ~&z \in (\eta_0(\xi), +\infty)\cap(\mathcal{I}_r \backslash \overline{\mathcal{I}_l}), \\ \\
     \begin{pmatrix} 1 & 0 \\ r(z)e^{2\mathrm{i} tg(z)} & 1 \end{pmatrix}, ~ &z \in L_1, \\ \\
     \begin{pmatrix} 1 & \frac{-r^*(z)}{1-r(z)r^*(z)}e^{-2\mathrm{i} tg(z)} \\0 & 1 \end{pmatrix},  ~&z \in L_2,\\ \\
     \begin{pmatrix} 1 & 0 \\ \frac{r(z)}{1-r(z)r^*(z)}e^{2\mathrm{i} tg(z)} & 1 \end{pmatrix},  ~&z \in L_3,\\ \\
     \begin{pmatrix} 1 & -r^*(z)e^{-2\mathrm{i} tg(z)} \\ 0  & 1 \end{pmatrix}, ~ &z \in L_4.\end{array} \right.\\
    \end{aligned}
  \end{equation}
  Here, $\chi_{(a, b)}$ is the characteristic function of the set $(a, b)$ and  $\chi_{(a, b)} \equiv 0$ when $a>b$.
\end{enumerate}
\end{rhp}

Define
\begin{equation} \label{D3delta}
  \begin{aligned}
    \delta(z) &= \exp \left\{    \frac{\mathcal{R} (z;\lambda_r^+, \lambda_l^-)}{2 \pi \mathrm{i}}  \left(\int_{(-\infty, \eta_0(\xi))\backslash (\overline{\mathcal{I}_l \cup \mathcal{I}_r} )}  \frac{\ln (1-|r(\zeta)|^2)}{\mathcal{R} (\zeta;\lambda_r^+, \lambda_l^-) } \,\frac{\dif\zeta}{\zeta - z}  \right. \right. \\
    & \left. \left.  + \int_{(-\infty, \eta_0(\xi) )\cap(\mathcal{I}_r \backslash \overline{\mathcal{I}_l})} \frac{-\ln(a_+(\zeta)a_-^*(\zeta))}{\mathcal{R} (\zeta;\lambda_r^+, \lambda_l^-)} \,\frac{\dif\zeta}{\zeta - z}
    +  \int_{(-\infty, \eta_0(\xi) )\cap(\mathcal{I}_l \backslash \overline{\mathcal{I}_r})} \frac{\ln (r_+(\zeta))}{\mathcal{R}_+ (\zeta;\lambda_r^+, \lambda_l^-)} \, \frac{\dif\zeta}{\zeta - z} \right)  \right\}
    \end{aligned}
\end{equation}
whose large $z$ asymptotic behavior is $\delta(z)=\delta_{\infty}+\mathcal{O} (z^{-1}) $ as $z \to \infty$, where
\begin{equation}\label{D3deltaINF}
  \begin{aligned}
    \delta_{\infty }= \exp & \left\{ \frac{ \mathrm{i}}{2\pi} \left(\int_{(-\infty, \eta_0(\xi))\backslash (\overline{\mathcal{I}_l \cup \mathcal{I}_r} )}  \frac{\ln (1-|r(\zeta)|^2)}{\mathcal{R} (\zeta;\lambda_r^+, \lambda_l^-) } \dif \zeta \right. \right.  \\
    & \left. \left. + \int_{(-\infty, \eta_0(\xi) )\cap(\mathcal{I}_r \backslash \overline{\mathcal{I}_l})}  \frac{-\ln(a_+(\zeta)a_-^*(\zeta))}{\mathcal{R} (\zeta;\lambda_r^+, \lambda_l^-)} \dif \zeta
    + \int_{(-\infty, \eta_0(\xi) )\cap(\mathcal{I}_l \backslash \overline{\mathcal{I}_r})} \frac{\ln (r_+(\zeta))}{\mathcal{R}_+ (\zeta;\lambda_r^+, \lambda_l^-)} \dif \zeta \right)  \right\}.
    \end{aligned}
\end{equation}
Using the function ${\delta}(z)$  and the transformation
\begin{equation}
	M^{(3)}(z;x,t)=\delta_{\infty}^{\sigma_{3}} M^{(2)}(z;x,t) \delta(z) ^{-\sigma_{3}},
\end{equation}
the RHP \ref{rhpmid} is obtained again and the following steps are exactly the same as Section \ref{Cmpw}.
Using equations (\ref{rec}), (\ref{C3ginf}),  (\ref{C3par}) and (\ref{D3deltaINF}), the long-time asymptotic behaviors of the defocusing NLS equation (\ref{(NLS)}) with error estimates in middle plane wave regions are obtained as follows
 \begin{equation} \label{mplq}
  q(x,t)=\frac{\lambda_r^+-\lambda_l^-}{2} e^{\mathrm{i}(kx-\omega t)}e^{-\mathrm{i}\phi_{mp} }+\mathcal{O}(e^{-ct}),
\end{equation}
where
\begin{equation} \label{mplqk}
k=-(\lambda_r^++\lambda_l^-),~~~~ \omega=\frac{(\lambda_r^++\lambda_l^-)^2}{2}+\frac{(\lambda_r^+ - \lambda_l^-)^2}{4},
\end{equation}
and
\begin{equation} \label{phimpmp}
  \begin{aligned}
    \phi_{mp}=    &\frac{ 1}{\pi} \left(\int_{(-\infty, \eta_0(\xi))\backslash (\overline{\mathcal{I}_l \cup \mathcal{I}_r} )}  \frac{\ln (1-|r(\zeta)|^2)}{\sqrt{(\zeta-\lambda_r^+)(\zeta-\lambda_l^-)} } \dif \zeta  \right.  \\
    & \left. + \int_{(-\infty, \eta_0(\xi) )\cap(\mathcal{I}_r \backslash \overline{\mathcal{I}_l})}  \frac{-\ln(a_+(\zeta)a_-^*(\zeta))}{\sqrt{(\zeta-\lambda_r^+)(\zeta-\lambda_l^-)} } \dif \zeta
    + \int_{(-\infty, \eta_0(\xi) )\cap(\mathcal{I}_l \backslash \overline{\mathcal{I}_r})} \frac{\arg (r_+(\zeta))}{\sqrt{(\lambda_r^+-\zeta)(\zeta-\lambda_l^-)}} \dif \zeta \right) .
    \end{aligned}
\end{equation}

\subsubsection{The right rarefaction wave region: $ -\frac{\lambda_r^+ + 3\lambda_l^-}{2} <\xi <-\frac{\lambda_r^+ + 3\lambda_r^-}{2}$} \label{Drrare}
\
\newline
\indent
When $\xi >-\frac{\lambda_r^+ + 3\lambda_l^-}{2}$, the point $\eta_0(\xi)$ of the one-band $g$-function (\ref{C3g}) is less than $\lambda_l^-$, which leads to the emergence of the interval $(\eta_0(\xi), \lambda_l^-)$, on which the jump matrix (\ref{D3V2}) grows exponentially with respect to  $t$ large enough.
In order to remove the exponentially large entries, introduce the rarefaction $g$-function defined by (\ref{A4g}). The boundaries of this rarefaction wave region are characterized by
\begin{equation} \label{D4b}
  \lambda_r^- < \lambda_s(\xi) < \lambda_l^- \qquad \mathrm{iff}  \qquad -\frac{\lambda_r^+ + 3\lambda_l^-}{2} =\xi_{D3} <\xi < \xi_{D4}= -\frac{\lambda_r^+ + 3\lambda_r^-}{2} ,
\end{equation}
where $\lambda_s(\xi)=-(\lambda_r^+ + 2 \xi)/3$  is the soft edge defined in Whitham modulation theory.

Due to the same rarefaction $g$-function and the same  upper/lower matrix factorizations in (\ref{f1}) and (\ref{f2}), we will summarize the RHP analysis of the right rarefaction wave region. Similarly, open lenses from intervals to the steepest descent contours through $\lambda_s(\xi)$ and $\eta(\xi)$ by the transformations
$M(z;x,t) \mapsto M^{(1)}(z;x,t) \mapsto M^{(2)}(z;x,t)$, which results in the following RHP.

\begin{rhp}
  Find a $2 \times 2$ matrix-valued function $M^{(2)}(z;x,t)$ with the following properties:
\begin{enumerate} [label=(\roman*)]
  \item  $M^{(2)}(z;x,t)$ is analytic in $z \in \mathbb{C} \backslash \Sigma^{(2)}$ where $\Sigma^{(2)}=(-\infty, \lambda_r^+] \cup L_1 \cup L_2 \cup  L_3 \cup L_4.$
  \item  $M^{(2)}(z;x,t)=I+\mathcal{O} (z^{-1}) $ as $z\to \infty$.
  \item  $M^{(2)}(z;x,t)$ achieves the CBVs  $M^{(2)}_{+}(z;x,t)$ and $M^{(2)}_{-}(z;x,t)$  on $\Sigma^{(2)}$ away from the self-intersection  points and
  branch points that satisfy the jump condition $M^{(2)}_{+}(z;x,t)=M^{(2)}_{-}(z;x,t)V^{(2)}(z;x,t)$, where
  \begin{equation}\label{D4V2}
    \begin{aligned}
    V^{(2)}(z;x,t)&=\left\{  \begin{array}{ll}
     (1-r(z)r^{*}(z))^{\sigma_{3}},  ~& z\in (-\infty, \lambda_r^-), \\ \\
     (a_+(z) a_-^*(z))^{-\sigma_{3}},  ~&z \in  (\lambda_r^-, \lambda_s(\xi)), \\ \\
     \begin{pmatrix} \frac{e^{ 2 \mathrm{i} t g_+(z)}}{a_+(z) a_-^*(z)} \chi_{(\lambda_s(\xi), \eta(\xi))} & -1 \\ 1 & 0 \end{pmatrix}, ~&z \in (\lambda_s(\xi), +\infty)\cap(\mathcal{I}_r \backslash \overline{\mathcal{I}_l}), \\ \\
     \begin{pmatrix} 0 & -1 \\ 1 & 0 \end{pmatrix}, ~&z \in \mathcal{I}_l \cap \mathcal{I}_r, \\ \\
     \begin{pmatrix} 1 & 0 \\ r(z)e^{2\mathrm{i} tg(z)} & 1 \end{pmatrix}, ~ &z \in L_1, \\ \\
     \begin{pmatrix} 1 & \frac{-r^*(z)}{1-r(z)r^*(z)}e^{-2\mathrm{i} tg(z)} \\0 & 1 \end{pmatrix},  ~&z \in L_2,\\ \\
     \begin{pmatrix} 1 & 0 \\ \frac{r(z)}{1-r(z)r^*(z)}e^{2\mathrm{i} tg(z)} & 1 \end{pmatrix},  ~&z \in L_3,\\ \\
     \begin{pmatrix} 1 & -r^*(z)e^{-2\mathrm{i} tg(z)} \\ 0  & 1 \end{pmatrix}, ~ &z \in L_4.\end{array} \right.\\
    \end{aligned}
  \end{equation}
  Here, $\chi_{(a, b)}$ is the characteristic function of the set $(a, b)$.
\end{enumerate}
\end{rhp}

Using the function ${\delta}(z)$ defined by (\ref{A2delta})  and the transformation
\begin{equation}
	M^{(3)}(z;x,t)=\delta_{\infty}^{\sigma_{3}} M^{(2)}(z;x,t) \delta(z) ^{-\sigma_{3}},
\end{equation}
we will obtain the RHP \ref{rrrhp} again and the following steps are exactly the same as Section \ref{Arrare}. Hence, the asymptotic behaviors with error estimates in the right rarefaction wave region are given by (\ref{B4}).

\subsubsection{The right plane wave region: $\xi > -\frac{\lambda_r^+ + 3\lambda_r^-}{2}$}
\
\newline
\indent
When $\xi > -\frac{\lambda_r^+ + 3\lambda_r^-}{2}$, the stationary phase point $\lambda_s(\xi)$  moves outside $\mathcal{I}_r$, which leads to the exponentially
large diagonal entries $e^{\pm 2 \mathrm{i}t g_+(z)}$ of the jump matrix (\ref{D4V2}) on $(\lambda_s(\xi), \lambda_r^-)$. This is the same as Section \ref{ABrpw}, thus the only thing is to determine the boundary of the right plane wave region, which is
\begin{equation} \label{D5b}
  \eta_- (\xi) < \lambda_r^- \qquad \mathrm{iff}  \qquad  \xi > \xi_{D4}= -\frac{\lambda_r^+ + 3\lambda_r^-}{2}
\end{equation}
where $\eta_{-}(\xi)$ is given by (\ref{rpweta}).

\subsection{\rm Case E:   $\lambda_r^+>\lambda_l^+> \lambda_l^->\lambda_r^-
$}
\
\newline
\indent
In the previous sections, we have given all the possible regions and the corresponding long-time asymptotic behaviors with error estimates. These results can be regarded as stitching some specific regions together, while some regions are fixed, such as the left and right plane wave regions, which are determined by the initial data. So in the remaining two cases, i.e., Case E and Case F, we just need to determine the boundaries of each region.

In this case, the  interval $\mathcal{I}_l$ is completely contained within  $\mathcal{I}_r$, so we only need to consider the factorizations in (\ref{f1}), (\ref{f3}) and (\ref{f4}).
As the self-similar variable $\xi=x/t$ increases, the stationary phase points of the corresponding
$g$-functions change continuously at different intervals on $\mathbb{R}$, which implies that there are five different regions: the left plane wave region,
dispersive shock wave region, the middle plane wave region, rarefaction wave region and the right plane wave region.

\subsubsection{The left plane wave region: $  \xi<-\frac{2\lambda_r^++\lambda_l^++\lambda_l^-}{2}+\frac{2(\lambda_r^+ - \lambda_l^+)
(\lambda_r^+ - \lambda_l^-)}{-2\lambda_r^+ + \lambda_l^+ + \lambda_l^-} $}\label{Elpw}
\
\newline
\indent

The asymptotic form of the left plane wave region has given in Section \ref{Blpw} and the remaining thing is to determine the boundary of the region, which is
\begin{equation} \label{E1b}
  \eta_+ (\xi) > \lambda_r^+ \qquad \mathrm{iff}  \qquad  \xi< \xi_{E1}= -\frac{2\lambda_r^++\lambda_l^++\lambda_l^-}{2}+\frac{2(\lambda_r^+ - \lambda_l^+)
  (\lambda_r^+ - \lambda_l^-)}{-2\lambda_r^+ + \lambda_l^+ + \lambda_l^-},
\end{equation}
where $\eta_+ (\xi)$ is given by (\ref{lpweta}).

\subsubsection{Dispersive shock wave region: $-\frac{2\lambda_r^++\lambda_l^++\lambda_l^-}{2}+\frac{2(\lambda_r^+ - \lambda_l^+)
(\lambda_r^+ - \lambda_l^-)}{-2\lambda_r^+ + \lambda_l^+ + \lambda_l^-} < \xi < -\frac{\lambda_r^++2\lambda_l^++\lambda_l^-}{2}$} \label{Elsh}
\
\newline
\indent
When the stationary phase point $\eta_+ (\xi)$ moves inside $\mathcal{I}_r$, the exponential oscillation $e^{\pm 2 \mathrm{i} t \theta(z)}$ appears in the jump matrix $V^{(2)}(z;x,t)$
on the interval $(\eta_+(\xi), \lambda_r^+) $. This region is the same as Section \ref{Clsh} and
the boundaries of this region are also characterized by the degeneration below
\begin{equation}\label{E2b}
  \begin{aligned}
  &\xi \to \xi_{E1}= v_2(\lambda_r^+, \lambda_r^+, \lambda_l^+, \lambda_l^-)= -\frac{2\lambda_r^++\lambda_l^++\lambda_l^-}{2}+\frac{2(\lambda_r^+ - \lambda_l^+)
  (\lambda_r^+ - \lambda_l^-)}{-2\lambda_r^+ + \lambda_l^+ + \lambda_l^-}, \quad &\mathrm{as} \quad \lambda_s(\xi) \to \lambda_r^+,\\
  &\xi \to \xi_{E2}= v_2(\lambda_r^+, \lambda_l^+, \lambda_l^+, \lambda_l^-)=-\frac{\lambda_r^++2\lambda_l^++\lambda_l^-}{2}, \quad &\mathrm{as} \quad \lambda_s(\xi) \to \lambda_l^+,
  \end{aligned}
\end{equation}
where $v_2$ is the Whitham velocity given by (\ref{vv2}). The long-time asymptotics of the defocusing NLS equation (\ref{(NLS)}) is given by (\ref{lswq}).

\subsubsection{The middle plane wave region: $ -\frac{\lambda_r^++2\lambda_l^++\lambda_l^-}{2} < \xi<  -\frac{\lambda_r^+ + 3\lambda_l^-}{2}$}\label{Empw}
\
\newline
\indent
When $\xi > -\dfrac{\lambda_r^++2\lambda_l^++\lambda_l^-}{2}$, the stationary phase point $\lambda_s(\xi)$ of the shock $g$-function is less than $\lambda_l^+$, which leads to the exponentially large diagonal entries
$e^{\pm 2 \mathrm{i}t g_+(z)}$ of the jump matrix $V^{(2)}(z;x,t)$ on $(\lambda_s(\xi), \lambda_l^+)$. This region is similar to the case in Section \ref{Cmpw} and
the boundaries of this region are  characterized by
\begin{equation} \label{E3b}
  \lambda_l^- < \eta_0(\xi) < \lambda_l^+    \qquad \mathrm{iff}  \qquad -\frac{\lambda_r^++2\lambda_l^++\lambda_l^-}{2} =\xi_{E2} < \xi< \xi_{E3} = -\frac{\lambda_r^+ + 3\lambda_l^-}{2},
\end{equation}
where $\eta_0(\xi)=-\frac{1}{2}(\lambda_r^++\lambda_l^- +2\xi)$. The long-time asymptotics of the defocusing NLS equation (\ref{(NLS)}) is given by (\ref{mplq})-(\ref{phimpmp}).

\subsubsection{Rarefaction wave region: $ -\frac{\lambda_r^+ + 3\lambda_l^-}{2} <\xi <-\frac{\lambda_r^+ + 3\lambda_r^-}{2}$} \label{Errare}
\
\newline
\indent
When $\xi >-\frac{\lambda_r^+ + 3\lambda_l^-}{2}$, the point $\eta_0(\xi)$ of the one-band $g$-function (\ref{C3g}) is less than $\lambda_l^-$, which leads to the emergence of the interval $(\eta_0(\xi), \lambda_l^-)$, on which the jump matrix (\ref{D3V2}) grows exponentially with respect to $t$ large enough.
This region is the same as Section \ref{Drrare} and the boundaries of this rarefaction wave region are also characterized by
\begin{equation} \label{E4b}
  \lambda_r^- < \lambda_s(\xi) < \lambda_l^- \qquad \mathrm{iff}  \qquad -\frac{\lambda_r^+ + 3\lambda_l^-}{2} =\xi_{E3} <\xi < \xi_{E4}=-\frac{\lambda_r^+ + 3\lambda_r^-}{2} ,
\end{equation}
where $\lambda_s(\xi)=-(\lambda_r^+ + 2 \xi)/3$. The long-time asymptotics of the defocusing NLS equation (\ref{(NLS)}) is given by (\ref{B4}).

\subsubsection{The right plane wave region: $\xi > -\frac{\lambda_r^+ + 3\lambda_r^-}{2}$}
\
\newline
\indent
When $\xi > -\frac{\lambda_r^+ + 3\lambda_r^-}{2}$, the stationary phase point $\lambda_s(\xi)$  moves outside $\mathcal{I}_r$, which leads to the exponentially
large diagonal entries $e^{\pm 2 \mathrm{i}t g_+(z)}$ of the jump matrix (\ref{D4V2}) on $(\lambda_s(\xi), \lambda_r^-)$. This is the same as Section \ref{ABrpw}, so the only thing is to determine the boundary of the right plane wave region, which is
\begin{equation} \label{E5b}
  \eta_- (\xi) < \lambda_r^- \qquad \mathrm{iff}  \qquad  \xi > \xi_{E4} = -\frac{\lambda_r^+ + 3\lambda_r^-}{2},
\end{equation}
where $\eta_{-}(\xi)$ is given by (\ref{rpweta}).

\subsection{\rm Case F:   $\lambda_l^+>\lambda_r^+>\lambda_r^->\lambda_l^-$}
\
\newline
\indent
In this case, the interval $\mathcal{I}_r$ is completely contained within  $\mathcal{I}_l$, so only the   factorizations in (\ref{f1}), (\ref{f2}) and (\ref{f4}) should be considered. In comparison with Case E,  the positions of intervals $\mathcal{I}_l$ and $\mathcal{I}_r$ are swapped. This implies that this case is very similar to Case E, and the regions appear in exactly the opposite order, as we will see later.
As the self-similar variable $\xi=x/t$ increases, the stationary phase points of the corresponding
$g$-functions change continuously at different intervals on $\mathbb{R}$, which implies that there are five different regions: the left plane wave region, rarefaction wave region, the middle plane wave region, dispersive shock wave region and the right plane wave region. Note that the special initial data included in this case for the defocusing NLS equation (\ref{(NLS)}) have been fully studied by Jenkins \cite{Jenkins2015}.

\subsubsection{The left plane wave region: $ \xi<-\frac{3\lambda_l^++\lambda_l^-}{2}  $}\label{Flpw}
\
\newline
\indent
In Section \ref{Blpw}, we have given the uniform form of the left plane wave region, and here we only need to determine the boundary of the region, which is
\begin{equation} \label{F1b}
  \eta_+(\xi) > \lambda_l^+    \qquad \mathrm{iff}  \qquad \xi< \xi_{F1} =-\frac{3\lambda_l^++\lambda_l^-}{2} .
\end{equation}
where $\eta_+ (\xi) $ is exactly given in (\ref{lpweta}).

\subsubsection{Rarefaction wave region: $-\frac{3\lambda_l^++\lambda_l^-}{2} <\xi < -\frac{3\lambda_r^++\lambda_l^-}{2} $} \label{Flrare}
\
\newline
\indent
When $\xi > -\frac{3\lambda_l^++\lambda_l^-}{2} $, the stationary phase point $\eta_+(\xi)$ of the $g$-function (\ref{A1g}) moves inside $\mathcal{I}_l$. This contributes exponentially large diagonal entries $e^{\pm 2 \mathrm{i} t g_+(z)}$ of the jump matrix $V^{(2)}(z;x,t)$ on the interval $(\eta_+(\xi), \lambda_l^+) \subset
\mathcal{I}_l$. This is the same as Section \ref{Dlrare} and thus the boundaries of this rarefaction wave region are also characterized by:
\begin{equation} \label{F2b}
  \lambda_r^+ < \lambda_s(\xi) < \lambda_l^+ \qquad \mathrm{iff}  \qquad -\frac{3\lambda_l^++\lambda_l^-}{2} = \xi_{F1} <\xi < \xi_{F2} = -\frac{3\lambda_r^++\lambda_l^-}{2} ,
\end{equation}
where $\lambda_s(\xi)=-(\lambda_l^- + 2 \xi)/3$.   The long-time asymptotic behaviors are given by (\ref{lrwq}) and (\ref{philrlr}).

\subsubsection{Middle plane wave region: $ -\frac{3\lambda_r^+ + \lambda_l^-}{2} < \xi<  -\frac{\lambda_r^++2\lambda_r^- +\lambda_l^-}{2}$}\label{Fmpw}
\
\newline
\indent
When $\xi > -\frac{3\lambda_r^+ + \lambda_l^-}{2} $, the stationary phase point $\lambda_s(\xi)$ of the rarefaction $g$-function is less than $\lambda_r^+$. This leads to exponentially large  diagonal entries $e^{\pm 2 \mathrm{i}t g_+(z)}$ of the jump matrix $V^{(2)}(z;x,t)$ on $(\lambda_s(\xi), \lambda_r^+)$ in (\ref{D2V2}).
Such a region is similar to the case in Section \ref{Cmpw} and the boundary of the middle plane wave region is characterized by
\begin{equation} \label{F3b}
  \lambda_r^- < \eta_0(\xi) < \lambda_r^+    \qquad \mathrm{iff}  \qquad -\frac{3\lambda_r^+ + \lambda_l^-}{2}  = \xi_{F2} < \xi< \xi_{F3} = -\frac{\lambda_r^++2\lambda_r^- +\lambda_l^-}{2},
\end{equation}
where $\eta_0(\xi)=-\frac{1}{2}(\lambda_r^++\lambda_l^- +2\xi)$. The long-time asymptotics of the defocusing NLS equation (\ref{(NLS)}) is given by (\ref{mplq})-(\ref{phimpmp}).

\subsubsection{Dispersive shock wave region: $ -\frac{\lambda_r^++2\lambda_r^- +\lambda_l^-}{2} < \xi < -\frac{\lambda_r^+ + \lambda_r^- + 2\lambda_l^-}{2}+\frac{2(\lambda_r^+ - \lambda_l^-)
(\lambda_r^- - \lambda_l^-)}{\lambda_r^+ + \lambda_r^- -2\lambda_l^-} $} \label{Frsh}
\
\newline
\indent
If $\xi > -\dfrac{\lambda_r^++2\lambda_r^- +\lambda_l^-}{2}$, the point $\eta_0(\xi)$ of the one-band $g$-function (\ref{C3g}) is less than $\lambda_r^-$, which leads to the emergence of the interval $(\eta_0(\xi), \lambda_r^-)$, on which the jump matrix (\ref{C3V2}) grows exponentially with respect to  $t$ large enough.
This is the same as Section \ref{Crsh} and thus the boundaries of this region are also characterized by the degeneration below
\begin{equation} \label{F4b}
  \begin{aligned}
  &\xi \to \xi_{F3}= v_3(\lambda_r^+, \lambda_r^-, \lambda_r^-, \lambda_l^-)= -\frac{\lambda_r^++2\lambda_r^- +\lambda_l^-}{2}, \quad &\mathrm{as} \quad \lambda_s(\xi) \to \lambda_r^-,\\
  &\xi \to \xi_{F4}= v_3(\lambda_r^+, \lambda_r^-, \lambda_l^-, \lambda_l^-)= -\frac{\lambda_r^+ + \lambda_r^- + 2\lambda_l^-}{2}+\frac{2(\lambda_r^+ - \lambda_l^-)
  (\lambda_r^- - \lambda_l^-)}{\lambda_r^+ + \lambda_r^- -2\lambda_l^-}, \quad &\mathrm{as} \quad \lambda_s(\xi) \to \lambda_l^-.
  \end{aligned}
\end{equation}
where $v_3$ is the Whitham velocity given by (\ref{vv3}). The long-time asymptotics of the defocusing NLS equation (\ref{(NLS)}) is given by (\ref{rsab}).

\subsubsection{The right plane wave region: $\xi > -\frac{\lambda_r^+ + \lambda_r^- + 2\lambda_l^-}{2}+\frac{2(\lambda_r^+ - \lambda_l^-)
(\lambda_r^- - \lambda_l^-)}{\lambda_r^+ + \lambda_r^- -2\lambda_l^-}$} \label{Frpw}
\
\newline
\indent
When the stationary phase $\lambda_s(\xi)$  coincides with $\lambda_l^-$, the two-band $g$-function degenerates to the one-band one. This is the same as that
in Section \ref{ABrpw}, where we have given the uniform form of the right plane wave region. The only thing is to determine the boundary of the right plane wave region, which is
\begin{equation} \label{F5b}
  \eta_- (\xi) < \lambda_l^- \qquad \mathrm{iff}  \qquad  \xi > \xi_{F4}  = -\frac{\lambda_r^+ + \lambda_r^- + 2\lambda_l^-}{2}+\frac{2(\lambda_r^+ - \lambda_l^-)
  (\lambda_r^- - \lambda_l^-)}{\lambda_r^+ + \lambda_r^- -2\lambda_l^-},
\end{equation}
where $\eta_- (\xi)$ is given by (\ref{rpweta}).

\section{Conclusions}

In conclusion, it is significant and challenging to study the Riemann problem of the defocusing nonlinear Schr\"{o}dinger hydrodynamics from the aspects of both physics and mathematics.
In recent years, little work has been done to investigate the long-time asymptotics of the defocusing nonlinear Schr\"{o}dinger equation with step-like initial data, and the full asymptotic analysis of the general step-like initial data like (\ref{(initial)}) is very interesting and remains lacking. Thus in this work, we have carried out rigorous asymptotic analysis for the Riemann problem of the defocusing nonlinear Schr\"{o}dinger hydrodynamics based on the Whitham modulation theory and Riemann-Hilbert formulation. First, the complete classification, including six cases of the asymptotic solutions to the defocusing nonlinear Schr\"{o}dinger equation with initial data (\ref{(initial)}), is given according to the orders of the Riemann invariants by Whitham modulation theory. Then, the long-time asymptotic behaviors with the leading-order terms and error estimates for each region of the six cases are formulated by the Deift-Zhou nonlinear steepest descent method for oscillatory Riemann-Hilbert problems. Finally, the long-time asymptotic solutions in each case are displayed, and it is shown that the theoretical results are in excellent agreement with the results from Whitham modulation theory and the numerical simulations.\\
\par
\par
{\bf Acknowledgements}
\par
This work was supported by the National Natural Science Foundation of China through grant 11971067.

\appendix

\section{Riemann-Hilbert Formalism for the Inverse Scattering Problem} \label{AAC}

The Lax pair of the defocusing nonlinear Schr\"{o}dinger equation (\ref{(NLS)}) is firstly given by Zakharov and Shabat \cite{zakharov1971exact}, that is,
\begin{subequations}\label{lax}
  \begin{align}
  &\Psi_x + \mathrm{i} z\sigma_3\Psi=Q(x,t)\Psi,\label{xpair}\\
  &\Psi_t + \mathrm{i} z^2\sigma_3\Psi=\widetilde{Q} (x,t)\Psi,\label{tpair}
  \end{align}
\end{subequations}
where $\Psi(x,t)$ is a $2\times 2$ matrix-valued function, $z \in \mathbb{C}$ is the spectral parameter, and $Q(x,t)$ and $\widetilde{Q} (x,t)$ are expressed in terms of the potential function $q(x,t)$ as
\begin{subequations}\label{QQ}
  \begin{align}
  & Q(x,t)= \begin{pmatrix}
    0 & q(x,t) \\ \overline{q(x,t)} & 0
  \end{pmatrix},\\
  & \widetilde{Q} (x,t)=z Q(x,t) - \frac{\mathrm{i}}{2}({Q^2(x,t)}+Q_x(x,t)) \sigma_3.
  \end{align}
\end{subequations}
The Lax pair is compatible (i.e., $\Psi_{xt}=\Psi_{tx}$) if and only if the potential function $q(x,t)$ solves the defocusing NLS equation (\ref{(NLS)}).

The defocusing NLS function (\ref{(NLS)}) has plane wave solutions:
\begin{equation}
  q^{\mathrm{p}}_j (x,t)=A_j e^{-2\mathrm{i} \mu_j x-\mathrm{i}(A_j^2+2\mu_j^2) t}, \qquad j\in\{l, r\},
\end{equation}
which are consistent with the initial data (\ref{(initial)}). The nature extensions of (\ref{(initial)}) to $t>0$  given by
\begin{equation}
  \int_{0}^{\pm \infty} |q(x,t) - q^{\mathrm{p}}_j (x,t) | \,dx < \infty, \qquad \forall t \geqslant 0, \qquad j\in\{l, r\},
\end{equation}
make sure that the Riemann-Hilbert formulation in this work is complete.
\par
Now, introduce a pair of Jost functions by the asymptotics conditions
\begin{equation}
  \Psi_l(z;x,t) \sim \Psi^{\mathrm{p}} _l(z;x,t) \quad \mathrm{as} \quad x \to -\infty, \qquad \Psi_r(z;x,t) \sim \Psi^{\mathrm{p}}_r(z;x,t) \quad \mathrm{as} \quad x \to \infty,
\end{equation}
where $\Psi^{\mathrm{p}}_j(z;x,t)$, $j \in \{l ,r\},$ are the solutions of the Lax pair (\ref{lax}) with $q(x,t)$ replaced by $q^{\mathrm{p}}_j (x,t)$
\begin{equation} \label{psip}
  \Psi^{\mathrm{p}}_j(z;x,t)=e^{-{\mathrm{i}}(\mu_j x + (A_j^2 /2 +\mu_j^2 ) t)\sigma_3} \mathcal{E}_j (z) e^{-{\mathrm{i}} \mathcal{R}(z;\lambda_j^+,  \lambda_j^-)(x+(z+\mu_j)t)\sigma_3},
\end{equation}
where
\begin{equation}\label{EPsilon}
  \mathcal{E}_j (z):= \mathcal{E} (z;\lambda_j^+,  \lambda_j^-), \quad \beta_j(z):=\beta(z;\lambda_j^+,  \lambda_j^-), \qquad j \in \{l ,r\}.
\end{equation}
\par
The Jost functions can be expressed as the solutions of the Volterra integral equations as
\begin{equation}
  \Psi_j(z;x,t) = \Psi^{\mathrm{p}}_j(z;x,t) + \int_{\mp \infty }^{x} \Psi^{\mathrm{p}}_j(z;x,t) ({\Psi^{\mathrm{p}}_j})^{-1}(z;s,t)((Q-Q^{\mathrm{p}}_j)(s,t))\Psi_j(z;s,t) \,ds,
\end{equation}
where
\begin{equation}
  Q^{\mathrm{p}}_j(x,t)= \begin{pmatrix}
    0 & q^{\mathrm{p}}_j(x,t) \\ \overline{q^{\mathrm{p}}_j(x,t)} & 0
  \end{pmatrix},  \qquad j \in \{l ,r\}.
\end{equation}
As a consequence, the columns $\Psi_l^{(1)}(z;x,t)$ and $\Psi_r^{(2)}(z;x,t)$ are analytic in $z \in \mathbb{C}^+$, while the columns $\Psi_l^{(2)}(z;x,t)$ and $\Psi_r^{(1)}(z;x,t)$ are analytic in $z \in \mathbb{C}^-$.
Furthermore, for the exact initial data (\ref{(initial)}) we can calculate the Jost functions explicitly and thus it is proved that $\Psi_j(z;x,t)$ can be analytically continued to $\mathbb{C} \backslash \overline{\mathcal{I}_j} $, $ j\in \{l, r\}$.
Define the scattering matrix
\begin{equation}
  S(z):= {\Psi_r}^{-1}(z;x,t) \Psi_l(z;x,t) = \mathcal{E}_r^{-1} (z) \mathcal{E}_l (z)=\begin{pmatrix} a(z) & b^*(z) \\ b(z) & a^*(z) \end{pmatrix},
\end{equation}
with the scattering functions
\begin{equation}
  \begin{aligned}
  a(z)&=a^*(z)=\frac{\beta_l(z)\beta^{-1}_r(z)+\beta^{-1}_l(z)\beta_r(z)}{2},\\
  b(z)&=-b^*(z)=\frac{\beta_l(z)\beta^{-1}_r(z)-\beta^{-1}_l(z)\beta_r(z)}{2\mathrm{i}},\\
\end{aligned}
\end{equation}
and the reflection coefficient
\begin{equation}
  r(z):=\frac{b(z)}{a(z)}=-\mathrm{i} \frac{\beta_l^2(z)-\beta_r^2(z)}{\beta_l^2(z)+\beta_r^2(z)}.
\end{equation}
Then, $a(z)$, $b(z)$ and $r(z)$ are all analytic in $z\in \mathbb{C} \backslash ({(\overline{\mathcal{I}_l \cup \mathcal{I}_r} ) \backslash (\mathcal{I}_l \cap \mathcal{I}_r)})$ and satisfy the jump conditions below
\begin{equation}\label{rraa1}
  \begin{aligned}
  &a_+(z)=a_+^*(z)=-b_-(z)=b^*_-(z), \quad b_+(z)=-b^*_+(z)=a_-(z)=a_-^*(z), \qquad z\in \mathcal{I}_l \backslash \overline{\mathcal{I}_r},\\
  &a_+(z)=a_+^*(z)=b_-(z)=-b^*_-(z), \quad b_+(z)=-b^*_+(z)=-a_-(z)=-a_-^*(z), \qquad z\in \mathcal{I}_r \backslash \overline{\mathcal{I}_l},\\
  \end{aligned}
\end{equation}
and thus it follows that
\begin{equation}\label{rraa2}
  \begin{aligned}
  &r_+(z){r^*_-(z)}=1, \qquad z\in ({\mathcal{I}_l \cup \mathcal{I}_r} ) \backslash (\overline{\mathcal{I}_l \cap \mathcal{I}_r}),\\
  &r_+(z)+r^*_-(z)=-(a_+(z)a_-^*(z))^{-1}, \qquad z\in \mathcal{I}_r \backslash \overline{\mathcal{I}_l}.
\end{aligned}
\end{equation}
Furthermore, we have
\begin{equation}\label{absr}
    |r(z)|<1, \quad  z\in  \mathbb{C} \backslash ({(\overline{\mathcal{I}_l \cup \mathcal{I}_r} ) \backslash (\mathcal{I}_l \cap \mathcal{I}_r)}), \qquad  |r_\pm(z)|=1, \quad  z\in  (\overline{\mathcal{I}_l \cup \mathcal{I}_r} ) \backslash ({\mathcal{I}_l \cap \mathcal{I}_r}).
\end{equation}
By the definition of $a(z)$, one has $a(z)\neq0$, which means that there are no solitons for the initial data (\ref{(initial)}).
\par
The analytic and asymptotic behaviors of the Jost functions imply that one can construct the  $2 \times 2$ matrix-valued piecewise analytic function of the form
\begin{equation}\label{Mrhp}
  M(z;x,t):=\left\{  \begin{array}{ll}
  \left[ \dfrac{\Psi^{(1)}_l(z;x,t)}{a(z)} \quad \Psi^{(2)}_r(z;x,t) \right]e^{\mathrm{i}t\theta(z)\sigma_3},z \in \mathbb{C}^+, \\
  \\
  \left[ \Psi^{(1)}_r(z;x,t) \quad \dfrac{\Psi^{(2)}_l(z;x,t)}{a^*(z)} \right]e^{\mathrm{i}t\theta(z)\sigma_3}, z \in \mathbb{C}^-,
  \end{array} \right.
\end{equation}
where
\begin{equation} \label{theta}
  \theta(z)=\theta(z; \xi)=z^2+ \xi z \qquad  \textrm{with} \qquad \xi =\frac{x}{t}.
\end{equation}
The matrix-valued function $M(z;x,t)$ solves the basic RHP:
\begin{rhp} \label{RHP-rhp-0}
Find a $2 \times 2$ matrix-valued function $M(z;x,t)$ with the following properties:
  \begin{enumerate} [label=(\roman*)]
    \item $M(z;x,t)$ is analytic in $z \in \mathbb{C} \backslash \mathbb{R}$.
    \item $M(z;x,t)=I + \mathcal{O} (z^{-1}) $ as $z \to \infty$.
    \item $M(z;x,t)$ achieves the CBVs $M_+(z;x,t)$ and $M_-(z;x,t)$ on $\mathbb{R}$ away from branch points  which satisfy the jump condition $M_+(z;x,t)=M_-(z;x,t)V(z;x,t)$ where
  \begin{equation}\label{MV}
    \begin{aligned}
    V(z;x,t):=\left\{  \begin{array}{ll}
      \begin{pmatrix} 1-r(z)r^*(z) & -r^*(z)e^{-2\mathrm{i}t\theta(z)} \\ r(z)e^{2\mathrm{i}t\theta(z)} & 1 \end{pmatrix}, &z \in \mathbb{R} \backslash ({\overline{\mathcal{I}_l \cup \mathcal{I}_r} }), \\ \\
      \begin{pmatrix} 0 & -r_-^*(z)e^{-2\mathrm{i}t\theta(z)} \\ r_+(z)e^{2\mathrm{i}t\theta(z)} & 1 \end{pmatrix}, &z \in \mathcal{I}_l \backslash \overline{\mathcal{I}_r}, \\ \\
      \begin{pmatrix}(a_+(z)a_-^*(z))^{-1} & -e^{-2\mathrm{i}t\theta(z)} \\ e^{2\mathrm{i}t\theta(z)} & 0 \end{pmatrix}, &z \in \mathcal{I}_r \backslash \overline{\mathcal{I}_l}, \\ \\
      \begin{pmatrix} 0 & -e^{-2\mathrm{i}t\theta(z)} \\ e^{2\mathrm{i}t\theta(z)} & 0 \end{pmatrix}, &z \in \mathcal{I}_l \cap \mathcal{I}_r. \end{array} \right.
    \end{aligned}
  \end{equation}
  \end{enumerate}
\end{rhp}
If $M(z;x,t)$ is the solution of the basic Riemann-Hilbert Problem above, then the solution $q(x,t)$ of the Cauchy problem (\ref{(NLS)})-(\ref{(initial)}) is given by the reconstruction formula
\begin{equation} \label{recf}
  q(x,t)= 2 \mathrm{i} \lim_{z \to \infty} (z M(z;x,t))_{12}.
\end{equation}

\section{The local model problems and error estimates} \label{ABC}

\begin{figure}[htbp]
  \centering
  \subfigure[]{\includegraphics[width=6cm]{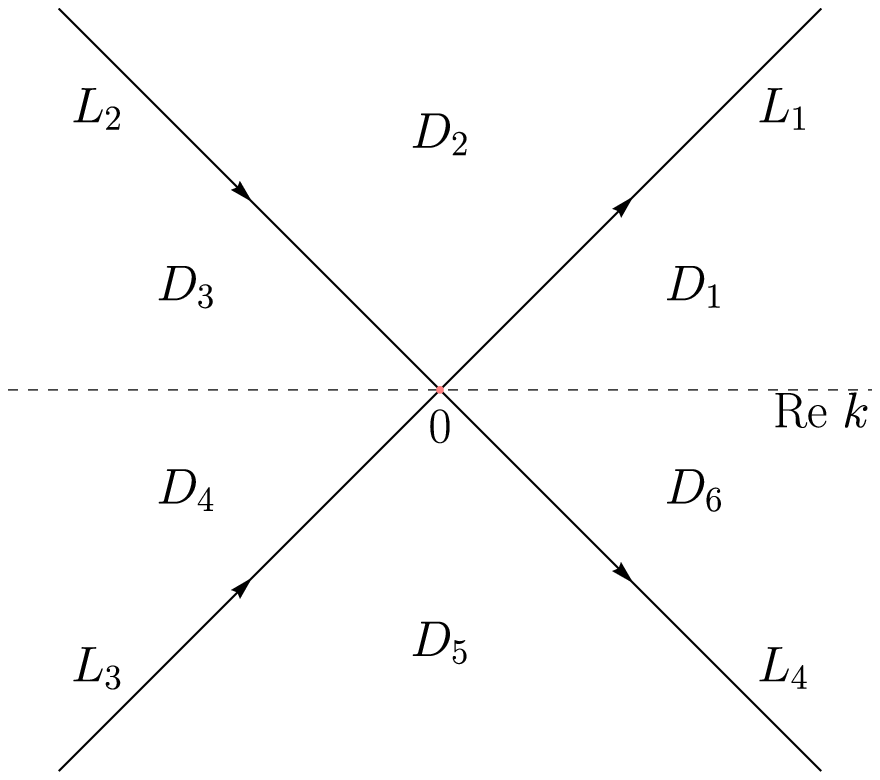}} \qquad \qquad
  \subfigure[]{\includegraphics[width=6cm]{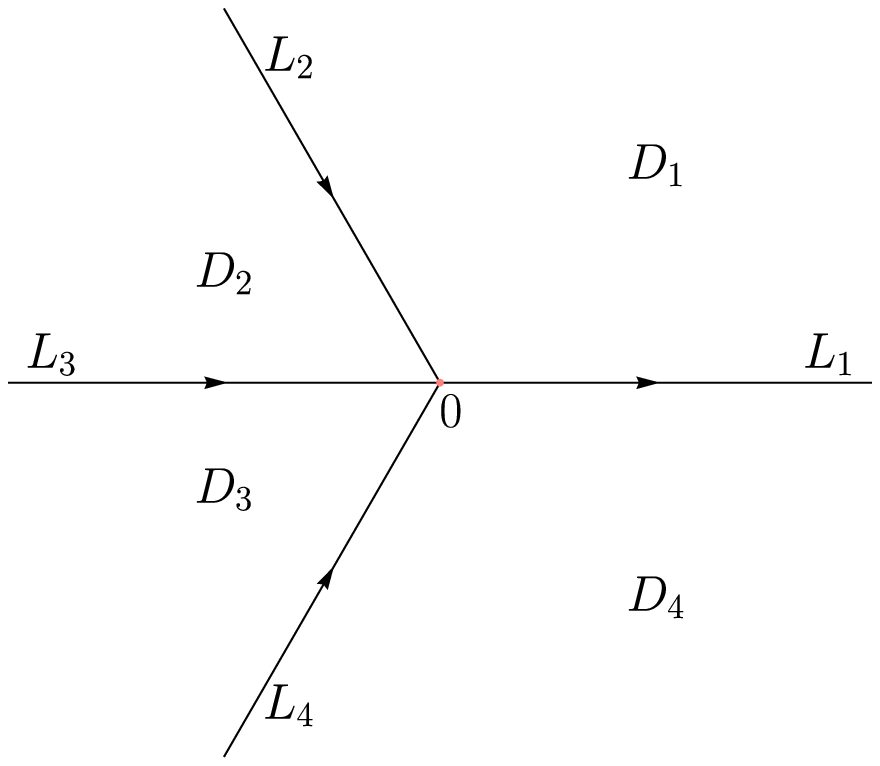}}
\caption{{\protect\small (a) The jump contours of parabolic cylinder model problem; (b) The jump contours of Airy model problem.}}
\label{figAPPEN}
\end{figure}

This appendix provides the local model RHPs (including the parabolic cylinder model and the Airy model)  for the construction of the local parametrices and the corresponding error estimates.
These models are standard and we only summarize them for use without rigorous proofs. Please see \cite{deift1994long,deift1999uniform} for details.

\subsection{The parabolic cylinder model problem}\label{pc}
\
\newline
\indent
Consider the oriented contours $L_i$ $(i=1, 2, 3, 4)$, which divide the complex $k$-plane into six sectors $D_i$ $(i=1,2, \cdots, 6)$, with the real axis, as shown in Figure \ref{figAPPEN}(a). Fix $r \in \mathbb{C}$ and let $\nu=\nu(r):= -\frac{1}{2\pi}\ln(1-rr^*)$. Then we consider the following RHP.

\begin{rhp}[The parabolic cylinder model problem] Find a $2 \times 2$ matrix-valued function $M^{(\mathrm{PC})}(k; r)$ with the following properties:
  \begin{enumerate} [label=(\roman*)]
    \item $M^{(\mathrm{PC})}(k;r)$ is analytic in $k \in \mathbb{C} \backslash \Sigma^{(\mathrm{PC})}$, where $\Sigma^{(\mathrm{PC})}=  L_1 \cup L_2 \cup  L_3 \cup L_4.$
    \item $M^{(\mathrm{PC})}(k;r)=I+ M_1^{(\mathrm{PC})}(r)~k^{-1} + \mathcal{O} (k^{-2}) $ uniformly as $k\to \infty$.
    \item $M^{(\mathrm{PC})}(k;r)$ achieves the CBVs $M^{(\mathrm{PC})}_{+}(k;r)$ and $M^{(\mathrm{PC})}_{-}(k;r)$ on $\Sigma^{(\mathrm{PC})}$  which satisfy the jump condition $M^{(\mathrm{PC})}_{+}(k;r)=M^{(\mathrm{PC})}_{-}(k;r)V^{(\mathrm{PC})}(k;r)$ where
  \begin{equation}\label{PCV}
    \begin{aligned}
  V^{(\mathrm{PC})}(k;r)&=\left\{  \begin{array}{ll}
     \begin{pmatrix} 1 & 0 \\ r k^{-2\mathrm{i}\nu} e^{\mathrm{i} k^2/2} & 1 \end{pmatrix}, ~ &k \in L_1, \\ \\
     \begin{pmatrix} 1 & \frac{-r^*}{1-rr^*}k^{2\mathrm{i}\nu} e^{-\mathrm{i} k^2/2}  \\0 & 1 \end{pmatrix},  ~&k \in L_2,\\ \\
     \begin{pmatrix} 1 & 0 \\ \frac{r}{1-rr^*}k^{-2\mathrm{i}\nu} e^{\mathrm{i} k^2/2} & 1 \end{pmatrix},  ~&k \in L_3,\\ \\
     \begin{pmatrix} 1 & -r^*k^{2\mathrm{i}\nu} e^{-\mathrm{i} k^2/2}  \\ 0  & 1 \end{pmatrix}, ~ &k \in L_4. \end{array} \right.\\
    \end{aligned}
  \end{equation}
  \end{enumerate}
\end{rhp}
There exists exactly one solution: when $r=0$, the solution is trivial; when $r \neq 0$, the solution is given in \cite{deift1994long} as follows:
\begin{equation}
  M^{(\mathrm{PC})}(k;r)=\Xi (k;r) P(k;r) e^{\mathrm{i} k^2 \sigma_3 /4}k^{-\mathrm{i}\nu \sigma_3},
\end{equation}
where
\begin{equation}
  \begin{aligned}
    \Xi (k;r)&=\left\{  \begin{array}{ll}
   \begin{pmatrix} e^{-3\pi \nu/4} D_{\mathrm{i} \nu}( e^{-3\pi \mathrm{i} /4} k) &  -\mathrm{i} \beta  e^{\pi(\nu - \mathrm{i}) /4} D_{-\mathrm{i} \nu-1}( e^{-\pi \mathrm{i} /4} k) \\ \mathrm{i} \nu e^{-3\pi(\nu + \mathrm{i}) /4} D_{\mathrm{i} \nu-1}( e^{-3\pi \mathrm{i} /4} k)/\beta   & e^{\pi \nu/4} D_{-\mathrm{i} \nu}( e^{-\pi \mathrm{i} /4} k)\end{pmatrix}, ~ &k \in \mathbb{C}^+, \\ \\
   \begin{pmatrix} e^{pi \nu/4} D_{\mathrm{i} \nu}( e^{\pi \mathrm{i} /4} k) &  -\mathrm{i} \beta  e^{-3\pi(\nu - \mathrm{i}) /4} D_{-\mathrm{i} \nu-1}( e^{3\pi \mathrm{i} /4} k) \\ \mathrm{i} \nu e^{\pi(\nu + \mathrm{i}) /4} D_{\mathrm{i} \nu-1}( e^{\pi \mathrm{i} /4} k)/\beta   & e^{-3\pi \nu/4} D_{-\mathrm{i} \nu}( e^{3\pi \mathrm{i} /4} k)\end{pmatrix}, ~ &k \in \mathbb{C}^-,
   \end{array} \right.\\
   P(k;r)&=\left\{  \begin{array}{ll}
    \begin{pmatrix} 1 & 0 \\ -r  & 1 \end{pmatrix}, ~ &k \in D_1, \\ \\
    \begin{pmatrix} 1 & \frac{r^*}{1-rr^*}  \\0 & 1 \end{pmatrix},  ~&k \in D_3,\\ \\
    \begin{pmatrix} 1 & 0 \\ \frac{r}{1-rr^*} & 1 \end{pmatrix},  ~&k \in D_4,\\ \\
    \begin{pmatrix} 1 & -r^*  \\ 0  & 1 \end{pmatrix}, ~ &k \in D_6, \\ \\
    I, ~ &k \in D_2 \cup D_5,\end{array} \right.\\
  \end{aligned}
\end{equation}
and $\beta$ is the complex constant
\begin{equation}
  \beta=\beta(r)=\frac{\sqrt{2\pi}e^{\mathrm{i}\pi/4}e^{-\pi\nu/2}}{r \Gamma(-\mathrm{i}\nu)},
\end{equation}
and $D_a(\cdot )$ denotes the standard (entire) parabolic cylinder function.

As $k \to \infty$, the asymptotic behavior of the solution is
\begin{equation}\label{pcinf}
  M^{(\mathrm{PC})}(k;r)=I+ \frac{1}{k} \begin{pmatrix} 0 & -\mathrm{i} \beta(r) \\ \mathrm{i} \nu(r) /\beta(r) & 0 \end{pmatrix}  + \mathcal{O} (k^{-2}).
\end{equation}

\subsection{The Airy model problem} \label{airy}
\
\newline
\indent
Consider the oriented contours $L_i$ $(i=1, 2, 3, 4)$, which divide the complex $k$-plane into four sectors $D_i$ $(i=1, 2, 3, 4)$, as shown in Figure \ref{figAPPEN}(b). Then we consider the following RHP.

\begin{rhp}[Airy model problem] Find a $2 \times 2$ matrix-valued function $M^{(\mathrm{Ai})}(k)$ with the following properties:
  \begin{enumerate} [label=(\roman*)]
    \item $M^{(\mathrm{Ai})}(k)$ is analytic in $k \in \mathbb{C} \backslash \Sigma^{(\mathrm{Ai})}$, where $\Sigma^{(\mathrm{Ai})}=  L_1 \cup L_2 \cup  L_3 \cup L_4.$
    \item $M^{(\mathrm{Ai})}(k)=I+ \mathcal{O} (k^{-1}) $ as $k\to \infty$.
    \item $M^{(\mathrm{Ai})}(k)$ achieves the CBVs $M^{(\mathrm{Ai})}_{+}(k)$ and $M^{(\mathrm{Ai})}_{-}(k)$ on $\Sigma^{(\mathrm{Ai})}$  which satisfy the jump condition $M^{(\mathrm{Ai})}_{+}(k)=M^{(\mathrm{Ai})}_{-}(k)V^{(\mathrm{Ai})}(k)$, where
  \begin{equation}\label{AiV}
    \begin{aligned}
  V^{(\mathrm{Ai})}(k)&=\left\{  \begin{array}{ll}
     \begin{pmatrix} 1 &  e^{-\frac{4}{3}k^{3/2}}  \\0 & 1 \end{pmatrix},  ~&k \in L_1,\\ \\
     \begin{pmatrix} 1 & 0 \\  e^{\frac{4}{3}k^{3/2}} & 1 \end{pmatrix},  ~&k \in L_2 \cup L_4,\\ \\
     \begin{pmatrix} 0 & 1  \\ -1  & 0 \end{pmatrix}, ~ &k \in L_3. \end{array} \right.\\
    \end{aligned}
  \end{equation}
  \end{enumerate}
\end{rhp}

The solution of the Airy model problem is given in \cite{deift1999uniform} as follows
\begin{equation}
  \begin{aligned}
    M^{(\mathrm{Ai})}(k)&=\left\{  \begin{array}{ll}
    \begin{pmatrix} \mathrm{Ai}(k) & \mathrm{Ai}(e^{4\pi \mathrm{i}/3} k) \\ \mathrm{Ai}^{'}(k) & e^{4\pi \mathrm{i}/3} \mathrm{Ai}^{'}(e^{4\pi \mathrm{i}/3} k) \end{pmatrix} e^{-\pi \mathrm{i} \sigma_3 /6}   e^{\frac{2}{3} k^{3/2} \sigma_3  }  , ~ &k \in D_1, \\ \\
    \begin{pmatrix} \mathrm{Ai}(k) & \mathrm{Ai}(e^{4\pi \mathrm{i}/3} k) \\ \mathrm{Ai}^{'}(k) & e^{4\pi \mathrm{i}/3} \mathrm{Ai}^{'}(e^{4\pi \mathrm{i}/3} k) \end{pmatrix} e^{-\pi \mathrm{i} \sigma_3 /6} \begin{pmatrix} 1 & 0 \\ -1 & 0 \end{pmatrix}    e^{\frac{2}{3} k^{3/2} \sigma_3  }  , ~ &k \in D_2, \\ \\
    \begin{pmatrix} \mathrm{Ai}(k) & -e^{4\pi \mathrm{i}/3} \mathrm{Ai}(e^{2\pi \mathrm{i}/3} k) \\ \mathrm{Ai}^{'}(k) & -\mathrm{Ai}^{'}(e^{2\pi \mathrm{i}/3} k) \end{pmatrix} e^{-\pi \mathrm{i} \sigma_3 /6} \begin{pmatrix} 1 & 0 \\ 1 & 0 \end{pmatrix}    e^{\frac{2}{3} k^{3/2} \sigma_3  }  , ~ &k \in D_3, \\ \\
    \begin{pmatrix} \mathrm{Ai}(k) & -e^{4\pi \mathrm{i}/3} \mathrm{Ai}(e^{2\pi \mathrm{i}/3} k) \\ \mathrm{Ai}^{'}(k) & -\mathrm{Ai}^{'}(e^{2\pi \mathrm{i}/3} k) \end{pmatrix} e^{-\pi \mathrm{i} \sigma_3 /6}   e^{\frac{2}{3} k^{3/2} \sigma_3  }  , ~ &k \in D_4,
    \end{array} \right.\\
  \end{aligned}
\end{equation}
and as $k \to \infty$, the full asymptotic expansion is given by
\begin{equation} \label{aiinf}
  \mathrm{Ai}(k)=\frac{e^{\pi \mathrm{i}/12}}{2 \sqrt{\pi}} k^{-\sigma_3/4} \left( \begin{pmatrix} 1 & 1 \\ -1 & 1 \end{pmatrix} e^{-\pi \mathrm{i} \sigma_3/4} + \sum_{n = 1}^{\infty} \begin{pmatrix} (-1)^n s_n  & s_n \\ (-1)^{n+1} t_n & t_n \end{pmatrix} e^{-\pi \mathrm{i} \sigma_3/4} (\frac{2}{3} k ^{3/2})^{-n}\right),
\end{equation}
where
\begin{equation}
  s_n= \frac{\Gamma(3n+1/2)}{54^n n! \Gamma(n+1/2)}, \qquad t_n=-\frac{6n+1}{6n-1} s_n.
\end{equation}

\subsection{Error estimates} \label{appb}
\
\newline
\indent
This appendix gives the error estimates by using the Cauchy singular integrals. Define the Cauchy projection operator on the error contour $\Sigma^{(\mathrm{err})}$
\begin{equation}
  (\mathcal{C} (f) )(z) = \frac{1}{2\pi \mathrm{i}} \int_{\Sigma^{(\mathrm{err})}} \frac{f(\zeta)}{\zeta-z}\dif \zeta,
\end{equation}
and define $\mathcal{C}_- (f) $ to be the nontangential limit of $\mathcal{C} (f) $ from the right side of $\Sigma^{(\mathrm{err})}$. Given the error matrix $V^{(\mathrm{err})}$, let
\begin{equation}
  \mathcal{C}_{V^{(\mathrm{err})}}f  = \mathcal{C}_- (f(V^{(\mathrm{err})}-I)).
\end{equation}
By Beals and Coifman theory \cite{beals1984scattering}, the solution of the error RHP is given by
\begin{equation}
  M^{(\mathrm{err})}(z;x,t)=I + \frac{1}{2\pi \mathrm{i}} \int_{\Sigma^{(\mathrm{err})}} \frac{\mu (\zeta;x,t)(V^{(\mathrm{err})}(\zeta;x,t)-I)}{\zeta-z}\dif \zeta,
\end{equation}
where $\mu $ is the unique solution of
\begin{equation}
  \mu = I + \mathcal{C}_{V^{(\mathrm{err})}}\mu,
\end{equation}
that is,
\begin{equation}
  \mu-I=(I-\mathcal{C}_{V^{(\mathrm{err})}} )^{-1} \mathcal{C}_{V^{(\mathrm{err})}}I,
\end{equation}
so
\begin{equation}\label{errmu}
  \begin{aligned}
   {\| \mu(\cdot;x,t)-I \|}_{L^2(\Sigma^{(\mathrm{err})})}  &\leqslant  {\|( I-\mathcal{C}_{V^{(\mathrm{err})} } ) ^{-1}\|}_{\mathcal{L} (L^2(\Sigma^{(\mathrm{err})}))} {\| \mathcal{C}_{V^{(\mathrm{err})}}I \|}_{L^2(\Sigma^{(\mathrm{err})})} \\
   & \leqslant C \| V^{(\mathrm{err})}(\cdot;x,t)-I \| _{L^2(\Sigma^{(\mathrm{err})})}
  \end{aligned}
\end{equation}
for some constant $C$, since $ V^{(\mathrm{err})}(z;x,t)-I$ decays to 0 uniformly as $t \to \infty$. Consider the expansion of $M^{(\mathrm{err})}(z;x,t)$ for large $z$
\begin{equation}
  M^{(\mathrm{err})}(z;x,t)=I +  \frac{M^{(\mathrm{err})}_1(x,t)}{z}+ \mathcal{O}(z^{-2}),
\end{equation}
and then
\begin{equation}
  M^{(\mathrm{err})}_1(x,t)= -\frac{1}{2\pi \mathrm{i}} \int_{\Sigma^{(\mathrm{err})}} {\mu (\zeta;x,t)(V^{(\mathrm{err})}(\zeta;x,t)-I)}\dif \zeta.
\end{equation}
By Cauchy-Schwarz inequality, we have
\begin{equation}\label{errm}
  \begin{aligned}
    |M^{(\mathrm{err})}_1(x,t)| &\leqslant \frac{1}{2\pi} \left(\| \mu(\cdot;x,t)-I  \|_{L^2(\Sigma^{(\mathrm{err})})}  \| V^{(\mathrm{err})}(\cdot;x,t)-I \|_{L^2(\Sigma^{(\mathrm{err})})}+ \| V^{(\mathrm{err})}(\cdot;x,t)-I \|_{L^1(\Sigma^{(\mathrm{err})})} \right).
  \end{aligned}
\end{equation}
Outside $\mathcal{U} $, $V^{(\mathrm{err})}(z;x,t)$ decays to $I$ exponentially for large $t$, so
\begin{equation}\label{err1}
  \| V^{(\mathrm{err})}(\cdot;x,t)-I \|_{L^p(\Sigma^{(\mathrm{err})}\backslash \partial\mathcal{U})} = \mathcal{O}(e^{-ct}), \qquad p=1,2,
\end{equation}
for some constant $c>0$. Inside $\mathcal{U}$, the parabolic cylinder model problem or the Airy model problem has been used to construct the local parametrix.
In what follows, the error estimate is considered in details. From (\ref{A1Verr}) and (\ref{pcinf}), it follows
\begin{equation}\label{err21}
  \| V^{(\mathrm{err})}(\cdot;x,t)-I \|_{L^p( \partial\mathcal{U})} = \mathcal{O}(t^{-1/2}), \qquad p=1,2.
\end{equation}
Combining equations (\ref{err1}) and (\ref{err21}), we have
\begin{equation}
  \| V^{(\mathrm{err})}(\cdot;x,t)-I \|_{L^p( \Sigma^{(\mathrm{err})})} = \mathcal{O}(t^{-1/2}), \qquad p=1,2,
\end{equation}
and thus $|M^{(\mathrm{err})}_1(x,t)|=\mathcal{O}(t^{-1/2})$ together from (\ref{errmu}) and (\ref{errm}).
As for the error estimate for the latter scenario, note that from (\ref{Ab2Verr}) and (\ref{aiinf}),
\begin{equation}\label{err22}
  \| V^{(\mathrm{err})}(\cdot;x,t)-I \|_{L^p( \partial\mathcal{U})} = \mathcal{O}(t^{-1}), \qquad p=1,2,
\end{equation}
so we have
\begin{equation}
  \| V^{(\mathrm{err})}(\cdot;x,t)-I \|_{L^p( \Sigma^{(\mathrm{err})})} = \mathcal{O}(t^{-1}), \qquad p=1,2,
\end{equation}
and thus $|M^{(\mathrm{err})}_1(x,t)|=\mathcal{O}(t^{-1})$.

\section{Representation in terms of Jacobi elliptic functions} \label{CC}

This appendix expresses the square modulus of the leading-order term of the asymptotic solution (\ref{onephase}) in terms of the Jacobi elliptic function.

Introduce the Jacobi theta functions with the nome $q$ as follows
\begin{equation}
  \begin{aligned}
      &\vartheta_1(z;q)=-\mathrm{i} \sum_{n\in \mathbb{Z}} (-1)^n q^{(n+1/2)^2} e^{\mathrm{i}(2n+1)z}, \quad
      \vartheta_2(z;q)=\sum_{n\in \mathbb{Z}} q^{(n+1/2)^2} e^{\mathrm{i}(2n+1)z}, \\
      &\vartheta_3(z;q)=\sum_{n\in \mathbb{Z}} q^{n^2} e^{2\mathrm{i}nz}, \quad
      \vartheta_4(z;q)=\sum_{n\in \mathbb{Z}} (-1)^n q^{n^2} e^{2\mathrm{i}nz}.
  \end{aligned}
\end{equation}
Let the nome $q$ be
\begin{equation}
  q=e^{\mathrm{i}\pi \tau},
\end{equation}
where $\tau$ is the Riemann period defined by (\ref{tau}). Note that the function $\Theta$ defined by (\ref{Theta}) can be expressed in terms of $\vartheta_3$ function as
\begin{equation}
  \Theta (z) = \vartheta_3(\pi z).
\end{equation}
Hereafter, the nome $q$ will be suppressed from the arguments of the theta functions for brevity.
The Jacobi elliptic functions $\mathrm{cn}$ and $\mathrm{sd}$ are defined as ratios of the above Jacobi theta functions
\begin{equation} \label{sd}
  \mathrm{cn}(z,m) = \frac{\vartheta_4(0)}{\vartheta_2(0)} \frac{\vartheta_2(z/{\vartheta_3^2(0)})}{\vartheta_4(z/{\vartheta_3^2(0)})}, \quad
  \mathrm{sd}(z,m) = \frac{\mathrm{sn}(z,m)}{\mathrm{dn}(z,m)} = \frac{{\vartheta_3^2(0)}}{\vartheta_2(0)\vartheta_4(0)} \frac{\vartheta_1(z/{\vartheta_3^2(0)})}{\vartheta_3(z/{\vartheta_3^2(0)})}
\end{equation}
with the elliptic modulus $m=\vartheta_2^2(0)/\vartheta_3^2(0)$.

Introduce the real quantities
\begin{equation}
  a= \frac{\gamma+ \hat{\gamma}}{2}, \qquad b= -2 \mathrm{i} \mathcal{A}(\infty),
\end{equation}
where the real quantities $\gamma$, $\hat{\gamma}$ and the pure imaginary quantity $\mathcal{A}(\infty)$ are given by (\ref{gamma}), (\ref{gammah}) and (\ref{Abelinf}), respectively.
Using the relation between $\Theta$ and $\vartheta_3$ and the addition formula
\begin{equation}
  \vartheta_3(a+\mathrm{i} b) \vartheta_3(a-\mathrm{i} b) \vartheta_3^2(0) = \vartheta_1^2(a) \vartheta_1^2(\mathrm{i} b) + \vartheta_3^2(a) \vartheta_3^2(\mathrm{i} b),
\end{equation}
we can write the square modulus of the leading-order term of the solution (\ref{onephase}) in the form
\begin{equation}
  |q_{\mathrm{as}}(x,t)|^2 =\frac{( \lambda_r^+ - \lambda_s(\xi) + \lambda_l^+ -\lambda_l^- )^2}{4} \left(1+ \frac{{\vartheta_1^2(a)}\vartheta_1^2(\mathrm{i} b)}{\vartheta_3^2(a) \vartheta_3^2(\mathrm{i} b)}\right).
\end{equation}
Then, by the definition (\ref{sd}) and the relation between $\mathrm{cn}$ and $\mathrm{sd}$
\begin{equation}
  \mathrm{cn}(z+K(m),m)= - \sqrt{1-m^2} \mathrm{sd}(z,m),
\end{equation}
one has
\begin{equation} \label{qqq}
  |q_{\mathrm{as}}(x,t)|^2 =\frac{( \lambda_r^+ - \lambda_s(\xi) + \lambda_l^+ -\lambda_l^- )^2}{4} \left(1+ \frac{m}{\sqrt{1-m^2}} \frac{\vartheta_1^2(\mathrm{i} b)}{\vartheta_3^2(\mathrm{i} b)}
  \mathrm{cn}^2(a\vartheta_3^2(0)+ K(m),m) \right).
\end{equation}
It remains to compute the ratio of the theta functions in the above equation as well as the constant $\vartheta_3^2(0)$. From \cite{doob1971handbook}, we have
\begin{equation} \label{t0}
  \vartheta_2^2(0)= \frac{2 m K(m)}{\pi}, \quad \vartheta_3^2(0)= \frac{2 K(m)}{\pi}, \quad \vartheta_4^2(0)= \frac{2 \sqrt{1-m^2} K(m)}{\pi}.
\end{equation}
Writing $y=\pi\mathcal{A}(\infty)$ and using the duplication fomulae
\begin{equation}
  \begin{aligned}
   &\vartheta_1(2y)\vartheta_2(0)\vartheta_3(0)\vartheta_4(0)=2\vartheta_1(y)\vartheta_2(y)\vartheta_3(y)\vartheta_4(y),\\
   &\vartheta_3(2y)\vartheta_4^2(0)\vartheta_3(0)=\vartheta_3^2(y)\vartheta_4^2(y)-\vartheta_1^2(y)\vartheta_2^2(y),
  \end{aligned}
\end{equation}
one has
\begin{equation} \label{ins1}
  \begin{aligned}
  \frac{\vartheta_1^2(\mathrm{i} b)}{\vartheta_3^2(\mathrm{i} b)} = \frac{\vartheta_1^2(2y)}{\vartheta_3^2(2y)}=4\frac{\vartheta_4^2(0)}{\vartheta_3^2(0)}
  \left( \frac{\vartheta_1(y)\vartheta_2(y)\vartheta_3(y)\vartheta_4(y)}{\vartheta_3^2(y)\vartheta_4^2(y)-\vartheta_1^2(y)\vartheta_2^2(y)} \right)^2
  = 4 \frac{\sqrt{1-m^2}}{m} \frac{\mathcal{G}^2 }{(\mathcal{G}^2-1)^2},
  \end{aligned}
\end{equation}
where $\mathcal{G} = (\vartheta_3(y)\vartheta_4(y))/(\vartheta_1(y)\vartheta_2(y))$.

To compute $\mathcal{G} $, the idea is to express it in terms of a meromorphic function. Firstly, recall the identities between the theta functions
\begin{equation}\label{ccc}
  \begin{aligned}
 &\vartheta_1(z)= - \mathrm{i} e^{\mathrm{i}\pi \tau /4} e^{\mathrm{i}z} \vartheta_3(z+\pi/2+\pi\tau/2), \\
  &\vartheta_2(z)= e^{\mathrm{i}\pi \tau /4} e^{\mathrm{i}z} \vartheta_3(z+\pi\tau/2),  \qquad \vartheta_4(z)=  e^{\mathrm{i}z} \vartheta_3(z+\pi/2).
\end{aligned}
\end{equation}
Consider the meromorphic function on the Riemann surface $X$ as
\begin{equation}
  G(P)=\mathrm{i} e^{-\mathrm{i}\pi \tau /2-2\pi\mathrm{i}\int_{\lambda_r^+}^{P} \omega_1}
  \frac{\Theta (\int_{\lambda_l^+}^{P} \omega_1 + \frac{1}{2}+\frac{\tau}{2})\Theta (\int_{\lambda_l^-}^{P} \omega_1 + \frac{1}{2}+\frac{\tau}{2})}
  {\Theta (\int_{\lambda_r^+}^{P} \omega_1 + \frac{1}{2}+\frac{\tau}{2})\Theta (\int_{\lambda_s(\xi)}^{P} \omega_1 + \frac{1}{2}+\frac{\tau}{2})}.
\end{equation}
Then considering the restriction on the upper sheet $X_+$, it follows
\begin{equation}
  \int_{\lambda_l^+}^{z} \omega_1=\int_{\lambda_r^+}^{z} \omega_1 +  \frac{1}{2}+\frac{\tau}{2}+ \mathbb{Z}, \quad
  \int_{\lambda_l^-}^{z} \omega_1=\int_{\lambda_r^+}^{z} \omega_1 +\frac{\tau}{2}+ \mathbb{Z}, \quad
  \int_{\lambda_s(\xi)}^{z} \omega_1= \int_{\lambda_r^+}^{z} \omega_1 +  \frac{1}{2} + \mathbb{Z},
\end{equation}
and using the periodicity properties of $\Theta $ and (\ref{ccc}), one has
\begin{equation}
  G(z)=\frac{\vartheta_3(\pi \int_{\lambda_r^+}^{z} \omega_1)\vartheta_4(\pi \int_{\lambda_r^+}^{z} \omega_1)}{\vartheta_1(\pi \int_{\lambda_r^+}^{z} \omega_1)\vartheta_2(\pi \int_{\lambda_r^+}^{z} \omega_1)}
\end{equation}
and $\mathcal{G}=G(\infty_+)$ is what we want to compute.

By the construction, $G(z)$ is meromorphic on $X_+$, with simple poles at $\lambda_r^+$ and $\lambda_s(\xi)$ and simple zeros at $\lambda_l^+$ and $\lambda_l^-$.
Then $G(z)$ can be expressed in the form
\begin{equation}
  G(z)=\mathcal{G}\frac{(z-\lambda_l^+)^{1/2}(z-\lambda_l^-)^{1/2}}{(z-\lambda_r^+)^{1/2}(z-\lambda_s(\xi))^{1/2}},
\end{equation}
by Liouville's theorem. To determine $\mathcal{G}$, we evaluate the residue of $G(z)$ at $\lambda_r^+$ in two different ways.
On the one hand, using the above expression, yields
\begin{equation} \label{res1}
  \mathrm{Res}(G(z),\lambda_r^+ )=\mathcal{G} \frac{(\lambda_r^+-\lambda_l^+)^{1/2}(\lambda_r^+-\lambda_l^-)^{1/2}}{(\lambda_r^+-\lambda_s(\xi))^{1/2}}.
\end{equation}
On the other hand, using the local coordinate near $\lambda_r^+$, yields
\begin{equation} \label{res2}
  \mathrm{Res}(G(z),\lambda_r^+ )=\frac{\vartheta_3(0)\vartheta_4(0)}{2\pi d \vartheta_1^{'}(0) \vartheta_2(0)},
\end{equation}
where $d$ is given by (\ref{d}). Matching the expressions (\ref{res1}) and (\ref{res2}), we have
\begin{equation}
  \mathcal{G}=-\mathrm{i}\frac{(\lambda_r^+ - \lambda_s(\xi))^{1/2}}{(\lambda_l^+-\lambda_l^-)^{1/2}},
\end{equation}
where the identity $\vartheta_1^{'}(0)=\vartheta_2(0)\vartheta_3(0)\vartheta_4(0)$ has been used. Inserting the expression of $\mathcal{G}$ into equation (\ref{ins1}), yields
\begin{equation}
  \frac{\vartheta_1^2(\mathrm{i} b)}{\vartheta_3^2(\mathrm{i} b)} =- 4 \frac{\sqrt{1-m^2}}{m} \frac{(\lambda_r^+ - \lambda_s(\xi))(\lambda_l^+-\lambda_l^-)}{((\lambda_r^+ - \lambda_s(\xi))+(\lambda_l^+-\lambda_l^-))^2}.
\end{equation}
Inserting this expression into equation (\ref{qqq}) and simplifying the expression by identifying $a$ with the help of the Riemann bilinear relations
\begin{equation}
  a=2 \pi \mathrm{i} d (x- Vt + \varphi^{(1)}),
\end{equation}
it finally follows
\begin{equation}
  |q_{\mathrm{as}}(x,t)|^2= \rho_2 - (\rho_2- \rho_3) \mathrm{cn}^2 \left(\sqrt{\rho_1-\rho_3}\left(x- Vt + \varphi^{(1)}\left(\frac{x}{t}\right)\right) -K(m),m \right),
\end{equation}
where
\begin{equation}\label{rho}
  \rho_1= \frac{(\lambda_r^+ + \lambda_s(\xi) -\lambda_l^+ - \lambda_l^-)^2}{4},
  \rho_2= \frac{(\lambda_r^+ - \lambda_s(\xi) +\lambda_l^+ - \lambda_l^-)^2}{4},
  \rho_3= \frac{(\lambda_r^+ - \lambda_s(\xi) -\lambda_l^+ + \lambda_l^-)^2}{4},
\end{equation}
and $\varphi^{(1)}$ and $V$  are given by (\ref{phi111}) and (\ref{VVV}), respectively.

\bibliographystyle{plain}

\end{document}